\documentclass{compositio}
\usepackage[T1]{fontenc}
\usepackage{amsmath}
\usepackage{amssymb}
\usepackage{amsfonts}
\usepackage{eucal}
\usepackage[all]{xy}
\usepackage[frenchb]{babel}
\usepackage{pslatex}

\newtheorem{prop}{Proposition}[subsection]
\newtheorem{theo}[prop]{Théor\`eme}
\newtheorem*{theo**}{Théorème}
\newtheorem{coro}[prop]{Corollaire}
\newtheorem{lemm}[prop]{Lemme}

\theoremstyle{definition}

\newtheorem{vide}[prop]{}
\newtheorem{defi}[prop]{Définition}

\theoremstyle{remark}
\newtheorem{rema}[prop]{Remarques}

\newtheorem{nota}[prop]{Notations}

\numberwithin{equation}{prop}

\newcommand{\riso}{ \overset{\sim}{\longrightarrow}\, }
\newcommand{\liso}{ \overset{\sim}{\longleftarrow}\, }

\newcommand{\Spf}{\mathrm{Spf}\,}

\renewcommand{\sp}{\mathrm{sp}}

\renewcommand{\AA}{{\mathcal{A}}}

\newcommand{\FF}{{\mathcal{F}}}
\newcommand{\B}{{\mathcal{B}}}

\newcommand{\E}{{\mathcal{E}}}
\newcommand{\G}{{\mathcal{G}}}
\renewcommand{\H}{{\mathcal{H}}}
\newcommand{\M}{{\mathcal{M}}}

\newcommand{\D}{{\mathcal{D}}}
\newcommand{\I}{{\mathcal{I}}}
\newcommand{\PP}{{\mathcal{P}}}

\renewcommand{\O}{{\mathcal{O}}}

\newcommand{\V}{\mathcal{V}}
\renewcommand{\S}{\mathcal{S}}
\newcommand{\T}{{\mathcal{T}}}
\newcommand{\Y}{\mathcal{Y}}
\newcommand{\ZZ}{\mathcal{Z}}
\newcommand{\X}{\mathfrak{X}}

\newcommand{\U}{\mathfrak{U}}

\newcommand{\A}{\mathbb{A}}

\newcommand{\DD}{\mathbb{D}}
\renewcommand{\L}{\mathbb{L}}
\newcommand{\R}{\mathbb{R}}
\newcommand{\Q}{\mathbb{Q}}
\newcommand{\Z}{\mathbb{Z}}
\newcommand{\N}{\mathbb{N}}

\newcommand{\hdag}{  \phantom{}{^{\dag} }    }

\begin{document}
\title{$\mathcal{D}$-modules arithmétiques associés aux isocristaux surconvergents. Cas lisse}
\author{Daniel Caro}
\email{daniel.caro@durham.ac.uk}
\address{Department of Mathematical Sciences,
Durham University,
Science Laboratories,
South Rd \newline
DURHAM DH1 3LE,
UNITED KINGDOM}
\classification{14F10, 14F30}
\keywords{arithmetical $\D$-modules, Frobenius, dual functor, direct image}
\thanks{L'auteur a bénéficié du soutien du réseau européen TMR \textit{Arithmetic Algebraic Geometry}
(contrat numéro UE MRTN-CT-2003-504917).}
\begin{abstract}
Let $\mathcal{V}$ be a mixed characteristic complete discrete valuation ring, $\mathcal{P}$ a separated smooth formal
scheme over $\mathcal{V}$, $P$ its special fiber, $X$ a smooth closed subscheme of $P$, $T$ a divisor in $P$ such that
$T _X = T \cap X$ is a divisor in $X$ and $\smash{\D}^\dag _{\mathcal{P}}(\hdag T)$ the weak completion of the sheaf of
differential operators on $\mathcal{P}$ with overconvergent singularities along $T$. We construct a fully faithful
functor denoted by $ \sp _{X \hookrightarrow \mathcal{P},T,+}$ from the category of isocrystal on $X \setminus T _X$
overconvergent along $T _X$ into the category of coherent $\smash{\D}^\dag _{\mathcal{P}}(\hdag T) \otimes _\mathbb{Z}
\mathbb{Q} $-modules with support in $X$. Next, we prove the commutation of $ \sp _{X \hookrightarrow \mathcal{P},T,+}$
with (extraordinary) inverse images and dual functors. These properties are compatible with Frobenius.
\end{abstract}

\maketitle


\tableofcontents

\section*{Introduction}

Soient $\V$ un anneau de valuation discrète complet d'inégales caractéristiques $(0,\,p)$
et de corps résiduel $k$,
$\PP$ un $\V$-schéma formel séparé et lisse,
$P$ sa fibre spéciale,
$T$ un diviseur de $P$, $U := P \setminus T$
et $F$ la puissance $s$-ième du Frobenius absolu de $P$, avec $s$ un entier fixé.
On désigne par
$\smash{\D} ^\dag _{\PP} (\hdag T ) _{\Q}$, l'{\it anneau des opérateurs différentiels sur $\PP$
de niveau fini à singularités surconvergentes le long de $T $} (voir \cite[4.2.5]{Be1}).
Lorsque le diviseur $T$ est vide, on ne l'indique pas.
On obtient la notion de
$\smash{\D}$-modules arithmétiques sur $U$ à singularités surconvergentes le long de $T$
(lorsque $\PP$ est propre, il n'est pas nécessaire de préciser $T$),
i.e., de $\smash{\D} ^\dag _{\PP} (\hdag T ) _{\Q}$-modules (toujours à gauche par défaut).
Pour illustrer ce fait, rappelons que Berthelot a construit
un foncteur canonique, noté $\sp _*$, pleinement fidèle de
la catégorie des isocristaux sur $U$ surconvergent le long de $T$ dans celle
des $\smash{\D} ^\dag _{\PP} (\hdag T) _{\Q}$-modules cohérents
(voir \cite[4.4.5]{Be1}).
Il a aussi défini la catégorie
des $F$-$\smash{\D} ^\dag _{\PP} (\hdag T ) _{\Q}$-modules (cohérents) de la façon suivante :
les objets sont les couples $(\M,\phi)$, où $\M$ est un $\smash{\D} ^\dag _{\PP} (\hdag T ) _{\Q}$-module (cohérent)
et $\phi$ un isomorphisme
$\smash{\D} ^\dag _{\PP} (\hdag T ) _{\Q}$-linéaire $F ^* \M \riso \M$.
Les flèches $(\M,\phi) \rightarrow (\M',\phi')$ sont les morphismes
$\smash{\D} ^\dag _{\PP} (\hdag T ) _{\Q}$-linéaires $\M \rightarrow \M '$ commutant à Frobenius.
En construisant de même les $F$-isocristaux surconvergents, il
a établi dans \cite[4.6]{Be2} que le foncteur $\sp _*$ commute aux actions de Frobenius.
La conjecture de Berthelot \cite[5.3.6.D]{Beintro2} implique
que l'image essentielle par $\sp _*$ des $F$-isocristaux sur
$U$ surconvergents le long de $T$
est incluse dans la catégorie des
$F \text{-}\smash{\D} ^\dag _{\PP,\Q}$-modules holonomes
(cela a un sens grâce à l'homomorphisme canonique
$\smash{\D} ^\dag _{\PP,\Q}\rightarrow \smash{\D} ^\dag _{\PP} (\hdag T) _{\Q}$).
Afin de valider ce genre de propriétés de finitude, le théorème de désingularisation
de de Jong (\cite{dejong}) fournit un outil très puissant en permettant de se ramener
au cas où le diviseur $T$ est lisse
(voir par exemple la preuve par Berthelot de la cohérence différentielle de l'algèbre des fonctions à singularités
 surconvergentes \cite{Becohdiff}).
La difficulté technique de cet outil est qu'en altérant $P$, on obtienne des variétés lisses non nécessairement relevables
mais se plongeant comme sous-schéma fermé dans un $\V$-schéma formel lisse.
D'où l'intérêt concernant une extension à cette situation géométrique que l'on appellera
{\og cas lisse non relevable \fg} ou simplement {\og cas lisse \fg}.
Cette généralisation est une première étape dans le problème de la construction
de $\smash{\D}$-modules arithmétiques associés aux isocristaux surconvergents.
Dans \cite{caro_devissge_surcoh}, nous donnerons une seconde construction
dans une situation géométrique différente et complémentaire de la présente.

Détaillons à présent les résultats de cet article.

Dans une première partie, nous établissons la transitivité, pour la composition de morphismes propres, des morphismes
d'adjonction entre images directes et images inverses extraordinaires.
Nous en déduisons en particulier
que l'isomorphisme {\it canonique} (i.e., celui construit directement à la main : \ref{f+og+comm})
de composition des images directes utilisé ici est
  le même (dans le cas de morphismes propres)
  que celui construit dans \cite[1.2.15]{caro_courbe-nouveau},
  qui utilise le délicat théorème de dualité relative.
  D'où une unification de ces deux constructions (voir \ref{remacomp+}).

Soit $X $ un sous-schéma fermé {\it $k$-lisse} de $P$ tel que $T \cap X$ soit un diviseur de $X$.
Dans la deuxième partie de cet article,
nous construisons un foncteur pleinement fidèle, noté $\sp _{X \hookrightarrow \PP,T,+}$, de la
catégorie des isocristaux sur $X \setminus (T\cap X)$ surconvergents le long de $T \cap X$
dans celle des
$\smash{\D} ^\dag _{\PP} (\hdag T) _{\Q}$-modules cohérents à support dans $X$.
Expliquons comment est construit $\sp _{X \hookrightarrow \PP,T,+}$.
Lorsque $X \hookrightarrow P $ se relève en
un morphisme de $\V$-schémas formels lisses $u$ : $\X \hookrightarrow \PP$,
on le définit en posant $\sp _{X \hookrightarrow \PP,T,+}:= u _{T,+} \circ \sp _*$,
où $\sp $ : $ \PP _K \rightarrow \PP$ le morphisme de spécialisation de la fibre générique de $\PP$
(en tant qu'espace analytique rigide) dans $\PP$
et où $u _{T,+}$ est l'image directe par $u$ à singularités surconvergentes le long de $T$.
Or, comme $X$ est lisse, $X \hookrightarrow P $ se relève localement.
Il s'agit alors de procéder par recollement.

Nous avions construit un isomorphisme de commutation du foncteur pleinement fidèle
$\sp _*$
aux foncteurs duaux respectifs (voir \cite{caro_comparaison}). Pour les {\it recoller}, il s'agit
alors d'établir la compatibilité de cet isomorphisme de commutation aux foncteurs images inverses
(extraordinaires) par une immersion ouverte.  Cela fait l'objet du troisième chapitre.

Dans une quatrième partie, nous prouvons la commutation de
$\sp _{X \hookrightarrow \PP,T,+}$ aux foncteurs restrictions, images inverses extraordinaires
et foncteur dual.
La vérification la plus technique est celle concernant le foncteur dual.
Nous construisons pour cela un isomorphisme de commutation
des foncteurs duaux aux images inverses extraordinaires par une immersion.
Pour valider sa transitivité, nous utilisons celle prouvée dans le premier chapitre.
Dans sa procédure de recollement,
nous nous servons aussi du troisième chapitre.

Le théorème de pleine fidélité de Kedlaya du foncteur associant
à des $F$-isocristaux surconvergents
les $F$-isocristaux convergents correspondants
(voir \cite{kedlaya_full_faithfull})
et celui de Tsuzuki sur les restrictions (voir \cite[4.1.1]{tsumono})
simplifie, lorsque l'on se restreint aux $F$-isocristaux surconvergents,
la vérification de la commutation
de $\sp _{X \hookrightarrow \PP,T,+}$ aux foncteurs cohomologiques ci-dessus.
Ces théorèmes nous ramènent sans difficulté au cas où $X \hookrightarrow P$ se relève, ce qui évite
les tracas engendrés par les recollements (transitivité etc.).
Ainsi, {\it lorsque l'on dispose de structures de Frobenius},
le premier et le troisième chapitre deviennent superflus.

Pour illustrer la maniabilité de la construction de $ \sp _{X \hookrightarrow \PP,T,+}$
(et son utilité), nous obtenons la cohérence, en tant que $\smash{\D} ^\dag _{\PP,\Q}$-module,
de l'image par $ \sp _{X \hookrightarrow \PP,T,+}$ de la catégorie
$F$-isocristaux {\it unités} (i.e. dont les pentes de Frobenius sont nulles)
sur $X \setminus (T\cap X)$ surconvergents le long de $T\cap X$.
La preuve utilise le théorème de désingularisation
de de Jong (\cite{dejong}) et le théorème de monodromie valable
pour les $F$-isocristaux surconvergents unités (voir \cite{tsumono}).
Cette technique est prometteuse et n'attend que des théorèmes de monodromie plus généraux.

\section*{Remerciements}
Je remercie P. Berthelot pour une erreur qu'il a décelée dans la précédente
procédure de recollement permettant de construire $ \sp _{X \hookrightarrow \PP,T,+}$.

\section*{Notations}
Tout au long de cet article, nous garderons les notations
suivantes :
les schémas formels seront notés par des lettres calligraphiques ou
gothiques et leur fibre spéciale par les lettres romanes
correspondantes.
De plus, la lettre $\V$ désignera un anneau de valuation discrète complet,
de corps résiduel $k$ de caractéristique $p>0$, de corps de
fractions $K$ de caractéristique $0$, d'idéal maximal $\mathfrak{m}$ et $\pi$ une uniformisante.
On fixe un entier naturel $s\geq 1$ et on désigne par $F$ la puissance
$s$-ième de l'endomorphisme de Frobenius.
Les modules sont par défaut des modules à gauche.
Si $X$ est un schéma ou schéma formel, on note
$d _X$ la dimension de Krull de $X$ et $\omega _{X}= \Omega _X ^{d _X}$
le faisceau des différentielles de degré maximum.
S'il existe un système de coordonnées locales $t _1,\dots,t_d$ (sur un schéma ou schéma formel),
on notera, pour tous $i =1,\dots ,d$, $\tau _i = 1 \otimes t _i - t_i \otimes 1$ et
  $\partial _i$ les dérivations correspondantes. Pour tout $\underline{k} \in \N ^d$,
  on pose $\underline{\tau} ^{\{ \underline{k}\}}:=\tau _1 ^{\{ k_1\}} \cdots \tau _d ^{\{ k_d\}}$,
  $\underline{\partial} ^{<\underline{k}>}:= \partial _1 ^{<k_1>}\cdots \partial _d ^{<k_d>}$.
Si $\E$ est un faisceau abélien, on écrit $\E _\Q:=\E \otimes _\Z \Q$.

$\bullet$ Les indices $\mathrm{qc}$, $\mathrm{tdf}$, $\mathrm{coh}$ et $\mathrm{parf}$
signifient respectivement \textit{quasi-cohérent}, \textit{de Tor-dimension finie}, {\it cohérent}
et {\it parfait}
tandis que $D ^\mathrm{b}$, $D ^+$ et $D ^-$ désignent respectivement les catégories dérivées
des complexes à cohomologie bornée,
bornée inférieurement et bornée supérieurement. Enfin, si $\mathcal{A}$ est un faisceau d'anneaux,
les symboles $\overset{ ^g}{} \mathcal{A}$, $\overset{^d}{} \mathcal{A}$ puis $\overset{^*}{}\mathcal{A}$
se traduisent respectivement par
$\mathcal{A}$-module {\og à gauche\fg}, {\og à droite\fg} puis {\og à droite ou à gauche\fg}
(par exemple,
$D (\overset{ ^g}{} \mathcal{A})$ indique la catégorie dérivée des complexes de $\mathcal{A}$-modules à gauche).
Si $\AA$ et $\B$ sont deux faisceaux d'anneaux (sur un même espace topologique),
on notera $D _{(.,\,\mathrm{qc})} ( \overset{^*}{}\AA ,\overset{^*}{} \B)$,
la sous catégorie pleine de
$D   ( \overset{^*}{} \AA  ,\overset{^*}{} \B)$ formée des complexes parfaits à droite (resp. de Tor-dimension finie à
droite).
De même en mettant le point à droite et en remplaçant {\og droite\fg} par {\og gauche\fg} ou
{\og $\mathrm{tdf}$\fg} par {\og $\mathrm{qc}$\fg} etc.
Comme les catégories de complexes sont par défaut les catégories dérivées,
on écrit
$\mathrm{Hom} _{\AA}(-,-)$ pour $\mathrm{Hom} _{D (\AA)}(-,-)$.

Soit $f$ : $\PP ' \rightarrow \PP$ un morphisme de $\V$-schémas formels lisses,
$T$ un diviseur de $P$ tel que $ f ^{-1} (T)$ soit un diviseur de $P'$.
On notera
$f _i$ : $P ' _i \rightarrow P _i$, la réduction de $f $ modulo $\mathfrak{m} ^{i+1}$.
En particulier, on obtient $f _0$ : $P' \rightarrow P$.
Pour tout entier $m\geq 0$, nous reprenons les constructions de Noot-Huyghe (voir \cite[2.1]{huyghe2})
concernant les faisceaux de la forme $\B _{P _i} ^{(m)} ( T)$ ou $\B _{\PP} ^{(m)} ( T)$.
En notant
$\smash{\D} _{\PP} ^{(m)}(T):=\B _{\PP} ^{(n _m)} ( T) \otimes _{\O _{\PP}} \smash{\D} ^{(m)} _{\PP}$,
avec $(n _m ) _{m\in \N}$ une suite d'entiers telle que $n _m \geq m$,
on construit la catégorie des complexes quasi-cohérents
$\smash{\underset{^{\longrightarrow}}{LD}} ^{\mathrm{b}} _{\Q
,\mathrm{qc}} ( \smash{\widehat{\D}} _{\PP} ^{(\bullet)}(T))$
(voir \cite[1.1]{caro_courbe-nouveau}
et \cite[4.2]{Beintro2}).
On prendra par défaut la suite $n _m =m$.
On écrira $\O _{\PP } (\hdag T) _\Q$ pour le
{\it faisceau des fonctions sur $\PP$ à singularités surconvergentes le long de $T$} (voir \cite[4.2.4.2]{Be1}).
Les produits tensoriels internes pour ces complexes quasi-cohérents seront notés
$-\smash{\overset{\L}{\otimes}}^{\dag}_{\O _{\PP } (\hdag T) _{\Q}}-$ (voir \ref{notaotimesdag}).
On note
$\underset{\longrightarrow}{\lim}$ : $\smash{\underset{^{\longrightarrow}}{LD}} ^{\mathrm{b}} _{\Q
,\mathrm{qc}} ( \smash{\widehat{\D}} _{\PP} ^{(\bullet)}(T))
\rightarrow D (\smash{\D} ^\dag _{\PP,\,\Q} (\hdag T))$
le foncteur
canonique (voir \cite[4.2.2]{Beintro2} lorsque
le diviseur est vide, mais la construction est similaire).
De plus, nous désignerons respectivement par $f ^! _T$, $f _{T+}$
les foncteurs image inverse extraordinaire, image directe par $f$ à singularités surconvergentes le long de $T$
(voir \cite[3.4, 3.5, 4.3]{Beintro2} et \cite[1.1.5]{caro_courbe-nouveau}).
Pour ces deux opérations cohomologiques,
lorsque l'on aura affaire à des complexes de bimodules, pour préciser quelle structure nous choisissons
dans leur calcul, on mettra un indice {\og $g$\fg} (resp. {\og $d$\fg})
{\it en bas} pour indiquer la structure {\it gauche} (resp. {\it droite}).
En outre, si $Z$ est un sous-schéma fermé de $P$,
$\R \underline{\Gamma} ^\dag _Z $ sera le foncteur cohomologique local
à support strict dans $Z$ (au sens de \cite[2.2.6]{caro_surcoherent})
et $(\hdag Z)$ le foncteur restriction (\cite[2.2.6]{caro_surcoherent}).
Si $T' \subset T$ sont deux diviseurs de $P$, on notera abusivement
$(\hdag T)$ à la place de $(\hdag T ,\,T')$.
Pour tout diviseur $T$ de $P$, nous désignerons par $\DD ^\dag _{\PP, T}$, $\DD _{\PP, T}$ ou simplement par $\DD _{T}$,
le foncteur dual $\smash{\D} ^\dag _{\PP} (\hdag T) _{\Q}$-linéaire (voir \ref{nota-morptraceOdagTpre}).
Lorsque $T$ est l'ensemble vide, nous omettrons de l'indiquer dans les opérations cohomologiques ci-dessus.

\section*{Convention}
Soit
$\AA$ un faisceau d'anneaux et $\E$ un $\AA$-bimodule à gauche.
Pour calculer les termes de la forme $- \otimes _{\AA} \E$ ou
$\mathcal{H}om _{\AA} (-,\E)$, on prendra par défaut la structure gauche de $\E$.
De plus, si $\M$ est un $\AA$-bimodule à droite,
pour calculer $\M \otimes _{\AA} -$ (resp. $\mathcal{H}om _{\AA} (-,\M)$),
nous choisissons par défaut la structure droite (resp. gauche) de $\M$. De même pour les complexes.

\section{Transitivité du morphisme d'adjonction entre image directe et image inverse extraordinaire}

Nous vérifions ici que la transitivité pour la composition de morphismes propres
de l'isomorphisme de dualité relative implique celle
des morphismes d'adjonction entre image directe et image inverse extraordinaire.
Cette transitivité nous permettra dans un prochain chapitre d'établir la transitivité
de l'isomorphisme de commutation des foncteurs duaux aux images inverses extraordinaires
par une immersion (voir \ref{trans-theta-f2}).

\subsection{Cas des schémas}
Soient $S$ un schéma noethérien de dimension de Krull finie,
   $f$ : $Y \rightarrow X$ et $g$ : $Z \rightarrow Y$ deux $S$-morphismes propres de $S$-schémas lisses.
   On rappelle que $ \smash{\D} ^{(m)}   _{Y \rightarrow X}$ est le
   $(\smash{\D} ^{(m)}   _{Y},f ^{-1}\smash{\D} ^{(m)}   _{X} )$-bimodule $f ^* \smash{\D} ^{(m)}   _{X} $,
   où $f^*$ désigne l'image inverse en tant que $\O _{X}$-module. Pour la définition
   des images inverses extraordinaires et images directes, nous renvoyons à \cite[2.2 et 2.3]{Beintro2}.
On pourra identifier $f ^! \smash{\D} ^{(m)}   _{X} $ avec $ \smash{\D} ^{(m)}   _{Y \rightarrow X}[d _{Y/X}]$.
\begin{vide}\label{nota-thetaf}
Soient $\AA _X$ un faisceau d'anneaux sur $X$.
   Comme $Y$ est noethérien, le foncteur $f_*$ est de dimension cohomologique finie et induit
   le foncteur $\R f _*$ :
   $D ( f ^{-1} \AA _X  \overset{^\mathrm{d}}{})\rightarrow D(\AA    _{X} \overset{^\mathrm{d}}{})$.
Pour tous $\E \in D  _{\mathrm{tdf},\mathrm{qc}} (\overset{^\mathrm{g}}{} \AA   _{X})$ et
  $\FF \in D ( f ^{-1} \AA _X \overset{^\mathrm{d}}{})$, on notera
  $$\mathrm{proj} _f \ : \ \R f _* \FF \otimes _{\AA _X} ^\L \E \riso
  \R f _* (\FF \otimes ^\L _{f ^{-1} \AA _X} f ^{-1} \E )$$
  l'isomorphisme de projection ou ceux qui s'en déduisent par fonctorialité.
  De même, en remplaçant $f$ par $g$ etc. Il pourra aussi se noter $\mathrm{proj}$.

  De plus, nous désignerons par $\otimes$ les isomorphismes de la forme
  \cite[2.1.12.(i)]{caro_comparaison}
  ou ceux qui s'en déduisent par fonctorialité.

\end{vide}

\begin{vide}[Isomorphismes canoniques de composition]\label{isocanocomps}
Les images directes par $f$ de complexes de $\smash{\D} _Y$-modules à droite (resp. à gauche) seront notés
$f _+ ^{\mathrm{d}}$ (resp. $f _+ ^{\mathrm{g}}$). Si aucune ambiguïté n'est à craindre, on les notera
simplement $f _+$.

Pour tout $\G \in D ^\mathrm{b} _\mathrm{coh} (\smash{\D} ^{(m)}   _{Z} \overset{^\mathrm{d}}{} )$, on dispose d'un
isomorphisme canonique $f ^{\mathrm{d}} _+ \circ g ^{\mathrm{d}} _+ (\G) \riso f \circ g  ^{\mathrm{d}} _+ (\G)$ fonctoriel
en $\G$.
Celui-ci est le composé suivant :
\begin{gather}\notag
  f ^{\mathrm{d}} _+ \circ g ^{\mathrm{d}} _+ (\G) = \R f _* ( \R g _* (\G \otimes ^\L _{\smash{\D} ^{(m)}   _{Z}}
\smash{\D} ^{(m)}   _{Z \rightarrow Y} ) \otimes ^\L _{\smash{\D} ^{(m)}   _{Y}} \smash{\D} ^{(m)}   _{Y \rightarrow X})
\\
\notag
\tilde{\underset{\mathrm{proj} _g}{\longrightarrow}}
\R f _*  \R g _* (\G \otimes ^\L _{\smash{\D} ^{(m)}   _{Z}}
\smash{\D} ^{(m)}   _{Z \rightarrow Y} \otimes ^\L _{g ^{-1}\smash{\D} ^{(m)}   _{Y}} g ^{-1} \smash{\D} ^{(m)}   _{Y \rightarrow X})
\riso
\R f _*  \R g _* (\G \otimes ^\L _{\smash{\D} ^{(m)}   _{Z}}\smash{\D} ^{(m)}   _{Z \rightarrow X})
\\
\label{f+og+comm}
\riso
f \circ g ^{\mathrm{d}} _+ (\G),
\end{gather}
où l'avant dernier isomorphisme découle
de l'isomorphisme canonique
$\smash{\D} ^{(m)}   _{Z \rightarrow Y} \otimes ^\L _{g ^{-1}
\smash{\D} ^{(m)}   _{Y }} g ^{-1} \smash{\D} ^{(m)}   _{Y \rightarrow X}
\riso
\smash{\D} ^{(m)}   _{Z \rightarrow X}$.
En tordant les structures de droite à gauche,
il en découle, pour tout
$\G \in D ^\mathrm{b} _\mathrm{coh} (\overset{^\mathrm{g}}{} \smash{\D} ^{(m)}   _{Z}  )$,
l'isomorphisme canonique
$f ^{\mathrm{g}} _+ \circ g ^{\mathrm{g}} _+ (\G) \riso f \circ g ^{\mathrm{g}} _+ (\G)$.

De plus, si $\E \in D ^\mathrm{b} _\mathrm{coh} (\overset{^\mathrm{g}}{} \smash{\D} ^{(m)}   _{X}  )$,
on dispose de $g ^! \circ f ^!  (\E) \riso (f \circ g) ^! (\E)$ via les isomorphismes :
\begin{gather}
\notag
  g ^! f ^! (\E ) =
g ^! \smash{\D} ^{(m)}   _{Y}
\otimes ^\L _{g ^{-1}\smash{\D} ^{(m)}   _{Y}}
g ^{-1} (f ^! \smash{\D} ^{(m)}   _{X }
\otimes ^\L _{f^{-1}\smash{\D} ^{(m)}   _{X }} f^{-1} \E)
\\ \notag
\riso   g ^! \smash{\D} ^{(m)}   _{Y}
\otimes ^\L _{g ^{-1}\smash{\D} ^{(m)}   _{Y}}
g ^{-1} f ^! \smash{\D} ^{(m)}   _{X }
\otimes ^\L _{g ^{-1} f^{-1}\smash{\D} ^{(m)}   _{X }}
g ^{-1} f^{-1} \E
=
g ^! f ^! \smash{\D} ^{(m)}   _{X }
\otimes ^\L _{g ^{-1} f^{-1}\smash{\D} ^{(m)}   _{X }}
g ^{-1} f^{-1} \E
\\
\notag
\riso
(f \circ g) ^! \smash{\D} ^{(m)}   _{X } \otimes ^\L _{g ^{-1} f^{-1}\smash{\D} ^{(m)}   _{X }}
g ^{-1} f^{-1} \E
\riso
(f \circ g) ^! (\E).
\end{gather}
\end{vide}

\begin{vide}\label{thetaf}
  Soient
  $\E \in D ^\mathrm{b} _{\mathrm{qc},\mathrm{tdf}} (\overset{^\mathrm{g}}{} \smash{\D} ^{(m)}   _{X})$ et
  $\FF \in D ^\mathrm{b} _\mathrm{coh} (\overset{^\mathrm{g}}{} \smash{\D} ^{(m)}   _{Y})$.
  D'après Virrion (voir \cite{Vir04}),
  on dispose de l'isomorphisme
\begin{equation}\label{defchi'f}
  f ^{\mathrm{d}} _+ \R \mathcal{H}om _{\smash{\D} ^{(m)}   _{Y}} (\FF, \smash{\D} ^{(m)}   _{Y}) [d _Y]
  \riso
  \R \mathcal{H}om _{\smash{\D} ^{(m)}   _{X}} (f ^{\mathrm{g}} _+ (\FF), \smash{\D} ^{(m)}   _{X})[d _X],
\end{equation}
  que l'on notera
  ici $\chi ' _f$.
  Elle en déduit le composé suivant :
  \begin{gather}\notag
    \R \mathcal{H}om _{\smash{\D} ^{(m)}   _{X}} (f _+ (\FF), \E)
    \\ \notag
    \tilde{\underset{\otimes}{\longleftarrow}}
    \R \mathcal{H}om _{\smash{\D} ^{(m)}   _{X}} (f _+ (\FF), \smash{\D} ^{(m)}   _{X})
    \otimes ^\L _{\smash{\D} ^{(m)}   _{X}} \E
    \tilde{\underset{\chi '_f}{\longleftarrow}}
    f _+ \R \mathcal{H}om _{\smash{\D} ^{(m)}   _{Y}} (\FF, \smash{\D} ^{(m)}   _{Y})
    \otimes ^\L _{\smash{\D} ^{(m)}   _{X}} \E [d_{Y/X}]
    \\
    \notag
    =\R f _* (
    \R \mathcal{H} om _{\smash{\D} ^{(m)}   _{Y}} (\FF, \smash{\D} ^{(m)}   _{Y}) \otimes^\L _{\smash{\D} ^{(m)}   _{Y}} f ^! \smash{\D} ^{(m)}   _{X}
    )
    \otimes ^\L _{\smash{\D} ^{(m)}   _{X}} \E
    \tilde{\underset{\otimes}{\longrightarrow}}
    \R f _* \R \mathcal{H} om _{\smash{\D} ^{(m)}   _{Y}} (\FF, f ^! \smash{\D} ^{(m)}   _{X})
    \otimes ^\L _{\smash{\D} ^{(m)}   _{X}} \E
    \\
    \notag
    \tilde{\underset{\mathrm{proj} _f}{\longrightarrow}}
    \R f _* (\R \mathcal{H} om _{\smash{\D} ^{(m)}   _{Y}} (\FF, f ^! \smash{\D} ^{(m)}   _{X})
    \otimes ^\L _{f ^{\text{-}1}\smash{\D} ^{(m)}   _{X}} f ^{\text{-}1} \E )
    \\ \label{consdeltaf}
    \tilde{\underset{\otimes}{\longrightarrow}}
   \R f _* \R \mathcal{H} om _{\smash{\D} ^{(m)}   _{Y}} (\FF, f ^! \E ),
  \end{gather}
que l'on notera
\begin{equation}
\label{defdelta}
  \delta _f \ :\
  \R \mathcal{H}om _{\smash{\D} ^{(m)}   _{X}} (f _+ (\FF), \E)
  \riso
  \R f _* \R \mathcal{H} om _{\smash{\D} ^{(m)}   _{Y}} (\FF, f ^! \E ).
\end{equation}

\end{vide}

\begin{lemm}
\label{lemmrema-thetaf}
Soient $\E \in D ^\mathrm{b} _{\mathrm{qc},\mathrm{tdf}} (\overset{^\mathrm{g}}{} \smash{\D} ^{(m)}   _{X})$ et
  $\FF \in D ^\mathrm{b} _\mathrm{coh} (\overset{^\mathrm{g}}{} \smash{\D} ^{(m)}   _{Y})$.
Le diagramme suivant :
\begin{equation}
  \label{rema-thetaf}
\xymatrix  @R=0,3cm   @C=2cm {
{\R \mathcal{H}om _{\smash{\D} ^{(m)}   _{X}} (f _+ (\FF), \E)}
\ar[d] ^-{\delta _f}
&
{ \R \mathcal{H}om _{\smash{\D} ^{(m)}   _{X}} (f _+ (\FF), \smash{\D} ^{(m)}   _{X})
    \otimes ^\L _{\smash{\D} ^{(m)}   _{X}} \E}
    \ar[l] ^-{\otimes}
    \ar[d] ^-{\delta _f}
    \\
{\R f _* \R \mathcal{H} om _{\smash{\D} ^{(m)}   _{Y}} (\FF, f ^! \E )}
&
{\R f _* \R \mathcal{H} om _{\smash{\D} ^{(m)}   _{Y}} (\FF, f ^! \smash{\D} ^{(m)}   _{X})
    \otimes ^\L _{\smash{\D} ^{(m)}   _{X}} \E}
    \ar[l] ^-{\otimes \circ \mathrm{proj} _f}
}
\end{equation}
est commutatif.
\end{lemm}
\begin{proof}
  Par construction de $\delta _f$ (\ref{consdeltaf}),
cela découle du fait que, pour tout $\G \in D ( f ^{-1}\smash{\D} ^{(m)}   _{X} \overset{^\mathrm{d}}{})$,
le morphisme composé canonique :
  $\R f _* \G \otimes _{\smash{\D} ^{(m)}   _{X}} ^\L \smash{\D} ^{(m)}   _{X} \overset{\mathrm{proj} _f}{\longrightarrow}
  \R f _* (\G \otimes ^\L _{f ^{-1}\smash{\D} ^{(m)}   _{X}} f ^{-1} \smash{\D} ^{(m)}   _{X} )\riso
  \R f _* \G $ est l'isomorphisme canonique
  $\R f _* \G \otimes _{\smash{\D} ^{(m)}   _{X}} ^\L \smash{\D} ^{(m)}   _{X} \riso
  \R f _* \G $.
\end{proof}

\begin{prop}\label{pre-tr-transif}
 Soient
  $\E \in D ^\mathrm{b} _{\mathrm{qc},\mathrm{tdf}} (\overset{^\mathrm{g}}{} \smash{\D} ^{(m)}   _{X})$ et
  $\G \in D ^\mathrm{b} _\mathrm{coh} (\overset{^\mathrm{g}}{} \smash{\D} ^{(m)}   _{Z})$.
  Le diagramme suivant
\begin{equation}
  \label{pre-tr-transifdiag0}
  \xymatrix  @R=0,3cm    {
  {\R \mathcal{H}om _{\smash{\D} ^{(m)}   _{X}} (f _+ \circ g _+ (\G), \E)}
  \ar[d] ^-{\delta _g \circ \delta _f} _-\sim
  &
  {\R \mathcal{H}om _{\smash{\D} ^{(m)}   _{X}} (f \circ g _+ (\G), \E)}
  \ar[l] _-\sim
  \ar[d] ^-{\delta _{f \circ g}} _-\sim
  \\
 {\R f _* \R g _* \R \mathcal{H} om _{\smash{\D} ^{(m)}   _{Z}} (\G, g ^! f ^! \E )}
 &
 {\R f \circ g _* \R \mathcal{H} om _{\smash{\D} ^{(m)}   _{Z}} (\G, f \circ g ^! \E ),}
 \ar[l] _-\sim
  }
\end{equation}
  où les isomorphismes horizontaux dérivent par fonctorialité des isomorphismes canoniques de composition
  de \ref{isocanocomps},
  est commutatif.
\end{prop}
\begin{proof}
Considérons le diagramme
\begin{equation}
  \label{pre-tr-transif-diag1}
  \xymatrix  @R=0,3cm   @C=-2cm {
  &
  {\R \mathcal{H}om _{\smash{\D} ^{(m)}   _{X}} (f \circ g _+ (\G), \E)}
  \ar[ld]
  \ar'[d][dd] ^-{\delta _{f \circ g}}
  &&
  {\R \mathcal{H}om _{\smash{\D} ^{(m)}   _{X}} (f \circ g _+ (\G), \smash{\D} ^{(m)}   _{X}) \otimes ^\L _{\smash{\D} ^{(m)}   _{X}} \E}
  \ar[ld]
  \ar[ll] ^-{\otimes}
  \ar[dd] ^-{\delta _{f \circ g}}
  \\
  {\R \mathcal{H}om _{\smash{\D} ^{(m)}   _{X}} (f _+ \circ g _+ (\G), \E)}
  \ar[dd] _-{\delta _g \circ \delta _f}
  &&
  {\R \mathcal{H}om _{\smash{\D} ^{(m)} _{X}} (f _+ \circ g _+ (\G),
  \smash{\D} ^{(m)}   _{X}) \otimes ^\L _{\smash{\D} ^{(m)}   _{X}} \E }
  \ar[ll] _(0.4){\otimes}
  \ar[dd] ^(0.25){\delta _g \circ \delta _f}
  &
  \\
  &
  {\R f \circ g _* \R \mathcal{H} om _{\smash{\D} ^{(m)} _{Z}}
  (    \G, f \circ g ^!  \E  )                               }
  &&
  {\R f \circ g _* \R \mathcal{H} om _{\smash{\D} ^{(m)} _{Z}}
  (\G, f \circ g ^!  \smash{\D} ^{(m)}   _{X} )
  \otimes ^\L _{\smash{\D} ^{(m)}   _{X}}  \E }
  \ar'[l][ll]^(0.1){\otimes \circ \mathrm{proj} _{f\circ g}  }
  \\
  {\R f _* \R g _* \R \mathcal{H} om _{\smash{\D} ^{(m)}   _{Z}} (\G, g ^! f ^! \E )}
  \ar[ur]
  &&
  {\R f _* \R g _* \R \mathcal{H} om _{\smash{\D} ^{(m)} _{Z}}
  (\G, g ^! f ^! \smash{\D} ^{(m)}  _{X} )
  \otimes ^\L _{\smash{\D} ^{(m)}  _{X}} \E ,}
  \ar[ur]
}
\end{equation}
où le carré de gauche est \ref{pre-tr-transifdiag0}
tandis celui de droite se déduit de \ref{pre-tr-transifdiag0} appliqué à
$\E = \smash{\D} ^{(m)} _{X}$.
Comme les flèches sont des isomorphismes,
il suffit donc de vérifier que les carrés de droite, du haut, du fond et le composé de celui
de devant avec celui du bas sont commutatifs.

D'après \ref{rema-thetaf}, le carré du fond est commutatif. Par fonctorialité, celui du haut l'est aussi.
Afin de vérifier celle du carré de droite, nous aurons besoin du lemme suivant.
\begin{lemm}\label{pre-tr-transif-lemm}
Pour tout $\G \in D ^\mathrm{b} _\mathrm{coh} (\overset{^\mathrm{g}}{} \smash{\D} ^{(m)}   _{Z})$,
le diagramme qui suit
  $$\xymatrix  @R=0,3cm   @C = -1,8cm {
  &
  {\R \mathcal{H}om _{\smash{\D} ^{(m)}   _{X}} (f _+ \circ g _+ (\G),\smash{\D} ^{(m)}   _{X}) }
  &&
  {\R \mathcal{H}om _{\smash{\D} ^{(m)}   _{X}} (f \circ g _+ (\G), \smash{\D} ^{(m)}   _{X}) }
  \ar[ll]
  \\
  &
  {f _+ \R \mathcal{H}om _{\smash{\D} ^{(m)}   _{Y}} ( g _+ (\G),\smash{\D} ^{(m)}   _{Y}) [d_{Y/X}]}
  \ar[u] ^-{\chi' _f}
  \ar@<9,5ex>[ddl] ^-{\otimes \circ \mathrm{proj} _g \circ \otimes}
  \\
  &
  {f _+  g _+ \R \mathcal{H}om _{\smash{\D} ^{(m)}   _{Z}} (  \G,\smash{\D} ^{(m)}   _{Z}) [d_{Z/X}]}
  \ar[u] ^-{\chi' _g}
  \ar @<-10ex> [ddl] _-{\otimes} 
  &&
  {f \circ g _+ \R \mathcal{H}om _{\smash{\D} ^{(m)}   _{Z}} (  \G,\smash{\D} ^{(m)}   _{Z}) [d_{Z/X}]}
  \ar[ll]
  \ar[uu] ^-{\chi' _{f \circ g}}
  \ar[ddl] ^-{\otimes}
  \\
  {\R f _* \R \mathcal{H} om _{\smash{\D} ^{(m)}   _{Y}} (g _+ (\G), f ^! \smash{\D} ^{(m)}   _{X})}
  \ar[d] ^-{\delta_g}
  \\
  {\R f _* \R g _* \R \mathcal{H} om _{\smash{\D} ^{(m)}   _{Z}} ( \G, g ^! f ^! \smash{\D} ^{(m)}   _{X})}
  &&
  {\R f\circ g  _* \R \mathcal{H} om _{\smash{\D} ^{(m)}   _{Z}} ( \G, f\circ g ^! \smash{\D} ^{(m)}   _{X})}
  \ar[ll]
}$$
est commutatif.
\end{lemm}
\begin{proof}
  Le rectangle du fond est commutatif par transitivité en $f$ des isomorphismes $\chi' _f$
  (voir \cite[2.5.9]{Beintro2}).
  Par construction de $\delta _g$ (voir \ref{thetaf}), le diagramme
\begin{equation}
  \label{pre-tr-transif-lemm-diag1}
\xymatrix  @R=0,3cm   @C=2,5cm  {
{\R \mathcal{H} om _{\smash{\D} ^{(m)}   _{Y}} (g _+ (\G), f ^! \smash{\D} ^{(m)}   _{X})}
\ar[d] ^-{\delta_g}
&
{\R \mathcal{H} om _{\smash{\D} ^{(m)}   _{Y}} (g _+ (\G), \smash{\D} ^{(m)}   _{Y})
\otimes ^\L _{\smash{\D} ^{(m)}   _{Y}} f ^! \smash{\D} ^{(m)}   _{X} }
\ar[l] ^-{\otimes}
\\
{ \R g _* \R \mathcal{H} om _{\smash{\D} ^{(m)}   _{Z}} ( \G, g ^! f ^! \smash{\D} ^{(m)}   _{X})}
&
{g _+ \R \mathcal{H} om _{\smash{\D} ^{(m)}   _{Y}} (\G, \smash{\D} ^{(m)}   _{Y})
\otimes ^\L _{\smash{\D} ^{(m)}   _{Y}} f ^! \smash{\D} ^{(m)}   _{X} [d_{Z/Y}].}
\ar[u] ^-{\chi' _g}
\ar[l]^(0.55){\otimes \circ \mathrm{proj} _g \circ \otimes}
}
\end{equation}
  est commutatif. En appliquant $\R f _*$ à \ref{pre-tr-transif-lemm-diag1}, on obtient
  le carré de gauche de \ref{pre-tr-transif-lemm}, qui est donc commutatif.
  Le carré du bas de \ref{pre-tr-transif-lemm} correspond (via \ref{f+og+comm}) au contour de

  \begin{equation}
    \label{pre-tr-transif-lemm-diag2}
    \xymatrix  @R=0,3cm    {
    {f _+  g _+ \R \mathcal{H}om _{\smash{\D} ^{(m)}   _{Z}} (  \G,\smash{\D} ^{(m)}   _{Z}) [d_{Z/X}]}
    \ar@{=}[r]
    \ar[d] ^-{\otimes \circ \mathrm{proj} _g \circ \otimes}
    &
    {\R f _* (\R g _* (\R \mathcal{H}om _{\smash{\D} ^{(m)}   _{Z}} (  \G,\smash{\D} ^{(m)}   _{Z})
    \otimes ^\L _{\smash{\D} ^{(m)}   _{Z}} g ^! \smash{\D} ^{(m)}   _{Y})
    \otimes ^\L _{\smash{\D} ^{(m)}   _{Y}} f ^! \smash{\D} ^{(m)}   _{X}) }
    \ar[d] ^-{\mathrm{proj} _g}
    \\
    {\R f _* \R g _* \R \mathcal{H} om _{\smash{\D} ^{(m)}   _{Z}} (\G, g ^! f ^! \smash{\D} ^{(m)}   _{X} )}
    \ar[d]
    &
    {\R f _* \R g _* (\R \mathcal{H} om _{\smash{\D} ^{(m)}   _{Z}} (\G, \smash{\D} ^{(m)}   _{Z} )
    \otimes ^\L _{\smash{\D} ^{(m)}   _{Z}} g ^! f ^! \smash{\D} ^{(m)}   _{X} )}
    \ar[l] ^-{\otimes}
    \ar[d]
    \\
    {\R f \circ g _* \R \mathcal{H} om _{\smash{\D} ^{(m)}   _{Z}} (\G, f\circ g ^! \smash{\D} ^{(m)}   _{X} )}
    &
    {\R f\circ g _* (\R \mathcal{H} om _{\smash{\D} ^{(m)}   _{Z}} (\G, \smash{\D} ^{(m)}   _{Z} )
    \otimes ^\L _{\smash{\D} ^{(m)}   _{Z}} f\circ g ^! \smash{\D} ^{(m)}   _{X} )}
    \ar[l] ^-{\otimes}
    }
  \end{equation}
  Le carré du bas de \ref{pre-tr-transif-lemm-diag2} est commutatif par fonctorialité.
  On vérifie par fonctorialité ou transitivité la commutativité du diagramme suivant
$$\xymatrix  @R=0,3cm    @C=0,2cm {
{\R g _* \R \mathcal{H}om _{\smash{\D} ^{(m)}   _{Z}} (  \G, g ^! \smash{\D} ^{(m)}   _{Y} )
    \otimes ^\L _{\smash{\D} ^{(m)}   _{Y}} f ^! \smash{\D} ^{(m)}   _{X}}
    \ar[d] ^-{\mathrm{proj} _g}
    &
    *[F]{\R g _* (\R \mathcal{H}om _{\smash{\D} ^{(m)}   _{Z}} (  \G,\smash{\D} ^{(m)}   _{Z})
    \otimes ^\L _{\smash{\D} ^{(m)}   _{Z}} g ^! \smash{\D} ^{(m)}   _{Y})
    \otimes ^\L _{\smash{\D} ^{(m)}   _{Y}} f ^! \smash{\D} ^{(m)}   _{X}}
    \ar[d] ^-{\mathrm{proj} _g}
    \ar[l] ^-{\otimes}
    \\
    {\R g _* (\R \mathcal{H}om _{\smash{\D} ^{(m)}   _{Z}} (  \G,g ^! \smash{\D} ^{(m)}   _{Y})
    \otimes ^\L _{g ^{-1}  \smash{\D} ^{(m)}   _{Y}} g ^{-1} f ^! \smash{\D} ^{(m)}   _{X} )}
    \ar[d] ^-{\otimes}
    &
    *[F]{\R g _* (\R \mathcal{H}om _{\smash{\D} ^{(m)}   _{Z}} (  \G,\smash{\D} ^{(m)}   _{Z})
    \otimes ^\L _{\smash{\D} ^{(m)}   _{Z}} g ^!
    f ^! \smash{\D} ^{(m)}   _{X} )}
    \ar[l] ^-{\otimes}
    \ar[d] ^-{\otimes}
    \\
    {\R g _* \R \mathcal{H}om _{\smash{\D} ^{(m)}   _{Z}} (  \G,g ^! \smash{\D} ^{(m)}   _{Y}
    \otimes ^\L _{g ^{-1}  \smash{\D} ^{(m)}   _{Y}} g ^{-1} f ^! \smash{\D} ^{(m)}   _{X} )}
    &
    *[F]{\R g _* (\R \mathcal{H}om _{\smash{\D} ^{(m)}   _{Z}} (  \G,g ^! f ^!\smash{\D} ^{(m)}   _{X} )),}
    \ar@{=}[l]
    }$$
    dont les termes encadrés correspondent aux sommets du carré du haut de \ref{pre-tr-transif-lemm-diag2}.
  En lui appliquant $\R f _*$, on obtient le carré du haut de \ref{pre-tr-transif-lemm-diag2}.
\end{proof}
Le foncteur $-\otimes ^\L _{\smash{\D} ^{(m)}   _{X}} \E $ appliqué au contour du diagramme de \ref{pre-tr-transif-lemm}
donne le carré de droite de \ref{pre-tr-transif-diag1}, qui est donc commutatif.

Il reste à prouver que le carré de devant composé avec celui du bas de \ref{pre-tr-transif-diag1} est commutatif.
Pour cela considérons le diagramme
\begin{equation}
\label{pre-tr-transif-lemm-diag3}
 \xymatrix  @R=0,3cm    @C=0,1cm  {
  {\R \mathcal{H}om _{\smash{\D} ^{(m)}   _{X}} (f _+ \circ g _+ (\G), \E)}
  \ar[d] ^-{\delta _f}
  &
  {\R \mathcal{H}om _{\smash{\D} ^{(m)}   _{X}}
  (f _+ \circ g _+ (\G),\smash{\D} ^{(m)}   _{X}) \otimes ^\L _{\smash{\D} ^{(m)}   _{X}} \E}
  \ar[l] ^-{\otimes}
  \ar[d] ^-{\delta _f}
  \\
  {\R f _* \R \mathcal{H}om _{\smash{\D} ^{(m)}   _{Y}} (g _+ (\G), f ^! \E)}
  &
  {\R f _* \R \mathcal{H}om _{\smash{\D} ^{(m)}   _{Y}} ( g _+ (\G),f ^! \smash{\D} ^{(m)}   _{X}) \otimes ^\L _{\smash{\D} ^{(m)}   _{X}} \E}
  \ar[l] ^-{\otimes\circ \mathrm{proj} _f}
  \ar@{=}[d]
  \\
  {\R f _* (\R \mathcal{H}om _{\smash{\D} ^{(m)}   _{Y}} (g _+ (\G), f ^! \smash{\D} ^{(m)}   _{X})
  \otimes ^\L _{f ^{-1} \smash{\D} ^{(m)}   _{X}} f ^{-1} \E)}
  \ar[u] ^-{\otimes}
  &
  {\R f _* \R \mathcal{H}om _{\smash{\D} ^{(m)}   _{Y}}
  ( g _+ (\G),f ^! \smash{\D} ^{(m)}   _{X}) \otimes ^\L _{\smash{\D} ^{(m)}   _{X}} \E}
  \ar[l] ^-{\mathrm{proj} _f}
  \\
  {\R f _* (\R \mathcal{H}om _{\smash{\D} ^{(m)}   _{Y}} (g _+ (\G), \smash{\D} ^{(m)}   _{Y} )
  \otimes ^\L _{\smash{\D} ^{(m)}   _{Y}} f ^!  \E)}
  \ar[u] ^-{\otimes}
  \ar[d] ^-{\delta _g}
  &
  {\R f _* (\R \mathcal{H}om _{\smash{\D} ^{(m)}   _{Y}} ( g _+ (\G),\smash{\D} ^{(m)}   _{Y} )
  \otimes ^\L _{\smash{\D} ^{(m)}   _{Y}} f ^! \smash{\D} ^{(m)}   _{X}) \otimes ^\L _{\smash{\D} ^{(m)}   _{X}} \E}
  \ar[l] ^-{\mathrm{proj} _f}
  \ar[u] ^-{\otimes}
  \ar[d] ^-{\delta _g}
  \\
  {\R f _* (\R g _* \R \mathcal{H}om _{\smash{\D} ^{(m)}   _{Z }} (\G, g ^! \smash{\D} ^{(m)}   _{Y} )
  \otimes ^\L _{\smash{\D} ^{(m)}   _{Y}} f ^!  \E)}
  \ar[d] ^-{\otimes \circ \mathrm{proj} _g}
  &
  {\R f _* (\R g _* \R \mathcal{H}om _{\smash{\D} ^{(m)}   _{Z}} ( \G , g ^! \smash{\D} ^{(m)}   _{Y} )
  \otimes ^\L _{\smash{\D} ^{(m)}   _{Y}} f ^! \smash{\D} ^{(m)}   _{X}) \otimes ^\L _{\smash{\D} ^{(m)}   _{X}} \E}
  \ar[l] ^-{\mathrm{proj} _f}
  \ar[d] ^-{\otimes \circ \mathrm{proj} _g}
  \\
  {\R f _* \R g _* \R \mathcal{H}om _{\smash{\D} ^{(m)}   _{Z}} (\G, g ^! f ^! \E )}
  &
  {\R f _* \R g _* \R \mathcal{H}om _{\smash{\D} ^{(m)}   _{Z}} ( \G , g ^! f ^! \smash{\D} ^{(m)}   _{X} )
  \otimes ^\L _{\smash{\D} ^{(m)}   _{X}} \E}
  \\
  {\R f \circ g _* \R \mathcal{H}om _{\smash{\D} ^{(m)}   _{Z}} (\G, f\circ g  ^! \E )}
  \ar[u] ^-{\mathrm{can}}
  &
  {\R f \circ g _* \R \mathcal{H}om _{\smash{\D} ^{(m)}   _{Z}} ( \G , f\circ g ^! \smash{\D} ^{(m)}   _{X} )
  \otimes ^\L _{\smash{\D} ^{(m)}   _{X}} \E.}
  \ar[u] ^-{\mathrm{can}}
  \ar[l] ^-{\otimes \circ \mathrm{proj} _{f\circ g}} }
  \end{equation}
Grâce à \ref{lemmrema-thetaf},
la flèche composée de gauche de \ref{pre-tr-transif-lemm-diag3} :
$\R f _* \R \mathcal{H}om _{\smash{\D} ^{(m)}   _{Y}} (g _+ (\G), f ^! \E)
\rightarrow
\R f _* \R g _* \R \mathcal{H}om _{\smash{\D} ^{(m)}   _{Z}} (\G, g ^! f ^! \E )$
est celle induite par $\delta _g$ (plus précisément $\R f _* \delta _g$).
De manière analogue, le composé de droite \ref{pre-tr-transif-lemm-diag3} :
$\R f _* \R \mathcal{H}om _{\smash{\D} ^{(m)}   _{Y}}
  ( g _+ (\G),f ^! \smash{\D} ^{(m)}   _{X}) \otimes ^\L _{\smash{\D} ^{(m)}   _{X}} \E
  \rightarrow
\R f _* \R g _* \R \mathcal{H}om _{\smash{\D} ^{(m)}   _{Z}} ( \G , g ^! f ^! \smash{\D} ^{(m)}   _{X} )
  \otimes ^\L _{\smash{\D} ^{(m)}   _{X}} \E$
  se déduit fonctoriellement de $\delta _g$ (est égale à $\R f _* \delta _g \otimes \mathrm{Id}$.
  Il en résulte que la commutativité du carré de devant composé avec celui du bas de \ref{pre-tr-transif-diag1}
  est équivalente à celle de
  \ref{pre-tr-transif-lemm-diag3}. \'Etablissons alors cette dernière.

D'après \ref{rema-thetaf}, le carré du haut de \ref{pre-tr-transif-lemm-diag3} est commutatif.
Par fonctorialité ou définition, il en est de même des trois autres.
Il reste à vérifier la commutativité du rectangle. Pour cela considérons le diagramme suivant.
\begin{equation}
\label{pre-tr-transif-lemm-diag4}
  \xymatrix  @R=0,3cm   @C=0,3cm  {
  {\R f _* (\R g _* \R \mathcal{H}om _{\smash{\D} ^{(m)}   _{Z}} (\G, g ^! \smash{\D} ^{(m)}   _{Y} )
  \underset{\smash{\D} ^{(m)}   _{Y}}{\otimes ^\L}  f ^!  \E)}
  \ar@{=}[r]
  \ar[d] ^-{\mathrm{proj} _g}
  &
  {\R f _* (\R g _* \R \mathcal{H}om _{\smash{\D} ^{(m)}   _{Z}} (\G, g ^! \smash{\D} ^{(m)}   _{Y} )
  \underset{\smash{\D} ^{(m)}   _{Y}}{\otimes ^\L}  f ^! \smash{\D} ^{(m)}   _{X} \underset{f ^{\text{-}1}\smash{\D} ^{(m)}   _{X}}{\otimes ^\L} f ^{\text{-}1} \E)}
  \ar[d] ^-{\mathrm{proj} _g}
  \\
  {\R f _* \R g _* (\R \underset{\smash{\D} ^{(m)}   _{Z}}{\mathcal{H}om} (\G, g ^! \smash{\D} ^{(m)}   _{Y})
  \underset{g ^{\text{-}1}\smash{\D} ^{(m)}   _{Y}}{\otimes ^\L}  g ^{\text{-}1} f ^!  \E)}
  \ar[d] ^-{\otimes}
  &
  {\R f _* (\R g _* (\R \underset{\smash{\D} ^{(m)}   _{Z}}{\mathcal{H}om}  (\G, g ^! \smash{\D} ^{(m)}   _{Y} )
  \underset{g ^{\text{-}1}\smash{\D} ^{(m)}   _{Y}}{\otimes ^\L}  g ^{\text{-}1} f ^! \smash{\D} ^{(m)}   _{X})
  \underset{f ^{\text{-}1}\smash{\D} ^{(m)}   _{X}}{\otimes ^\L} f ^{\text{-}1} \E)}
  \ar[d] ^-{\otimes}
  \\
  {\R f _* \R g _* \R \mathcal{H}om _{\smash{\D} ^{(m)}   _{Z}} (\G, g ^! f ^!  \E)}
  \ar[d]
  &
  {\R f _* (\R g _* \R \mathcal{H}om _{\smash{\D} ^{(m)}   _{Z}} (\G, g ^!   f ^! \smash{\D} ^{(m)}   _{X})
  \underset{f ^{\text{-}1}\smash{\D} ^{(m)}   _{X}}{\otimes ^\L} f ^{\text{-}1} \E)}
  \ar[d] ^-{\mathrm{proj} _g}
  \\
  {\R f _* \R g _* \R \underset{\smash{\D} ^{(m)}   _{Z}}{\mathcal{H}om}
  (\G, g ^! f ^! \smash{\D} ^{(m)}   _{X}
  \underset{g ^{\text{-}1}f ^{\text{-}1}\smash{\D} ^{(m)}   _{X}}{\otimes ^\L}
  g ^{\text{-}1} f ^{\text{-}1} \E )}
  \ar[d]
  &
  {\R f _* \R g _* (  \R \underset{\smash{\D} ^{(m)}   _{Z}}{\mathcal{H}om}
  (\G, g ^!   f ^! \smash{\D} ^{(m)}   _{X})
  \underset{g ^{\text{-}1}f ^{\text{-}1}\smash{\D} ^{(m)}   _{X}}{\otimes ^\L}
  g ^{\text{-}1} f ^{\text{-}1} \E)}
  \ar[l]^{\otimes}
  \ar[d]
  \\
  {\R f \circ g _* \R \underset{\smash{\D} ^{(m)}   _{Z}}{\mathcal{H}om}
  (\G, g ^! f ^! \smash{\D} ^{(m)}   _{X}
  \underset{f\circ g ^{\text{-}1}\smash{\D} ^{(m)}   _{X}}{\otimes ^\L}
  f\circ g ^{\text{-}1} \E )}
  \ar[d]
  &
  {\R f \circ g _* (  \R \underset{\smash{\D} ^{(m)}   _{Z}}{\mathcal{H}om}
  (\G, g ^!   f ^! \smash{\D} ^{(m)}   _{X})
  \underset{f\circ g ^{\text{-}1}\smash{\D} ^{(m)}   _{X}}{\otimes ^\L}
  f\circ g ^{\text{-}1} \E)}
  \ar[l]^{\otimes}
  \ar[d]
  \\
  {\R f \circ g _* \R \mathcal{H}om _{\smash{\D} ^{(m)}   _{Z}}
  (\G, f \circ g  ^! \E )}
  &
  {\R f \circ g _* (  \R \underset{\smash{\D} ^{(m)}   _{Z}}{\mathcal{H}om}
  (\G, f\circ g ^! \smash{\D} ^{(m)}   _{X})
  \underset{f\circ g ^{\text{-}1}\smash{\D} ^{(m)}   _{X}}{\otimes ^\L}
  f\circ g ^{\text{-}1} \E)}
  \ar[l]^{\otimes}
  }
\end{equation}
Les deux carrés de \ref{pre-tr-transif-lemm-diag4} sont commutatifs par fonctorialité.
Nous aurons besoin des deux lemmes ci-après, dont les preuves sont aisées,
qui établissent deux conditions de transitivité validées par
les morphismes de projection.
\begin{lemm}
  \label{trans-proj-e}
  Soient $\H \in D ( g ^{-1}\smash{\D} ^{(m)}   _{Y} \overset{^\mathrm{d}}{})$,
  $\E _1 \in D  _{(\mathrm{tdf},\mathrm{qc},\cdot)} (\overset{^\mathrm{g}}{} \smash{\D} ^{(m)}   _{Y},
  f ^{-1} \smash{\D} ^{(m)}   _{X} \overset{^\mathrm{d}}{})$
  et $\E _2 \in D  _{\mathrm{tdf},\mathrm{qc}} (f ^{-1} \smash{\D} ^{(m)}   _{X} \overset{^\mathrm{d}}{})$.
  Le diagramme
$$\xymatrix  @R=0,3cm    @C=2cm {
{\R g _* \H \otimes _{\smash{\D} ^{(m)}   _{Y}} ^\L \E _1\otimes ^\L _{f ^{-1} \smash{\D} ^{(m)}   _{X}}\E _2}
\ar[r] ^-{\mathrm{proj} _g} _-\sim
\ar[d] ^-{\mathrm{proj} _g} _-\sim
&
{\R g _* (\H \otimes _{ g ^{-1} \smash{\D} ^{(m)}   _{Y}} ^\L g ^{-1} \E _1) \otimes ^\L _{f ^{-1} \smash{\D} ^{(m)}   _{X}}\E _2}
\ar[d] ^-{\mathrm{proj} _g} _-\sim
\\
{\R g _* \H \otimes _{ g ^{-1} \smash{\D} ^{(m)}   _{Y}} ^\L g ^{-1} (\E _1\otimes ^\L _{f ^{-1} \smash{\D} ^{(m)}   _{X}}\E _2)}
\ar[r] _-\sim
&
{\R g _* (\H \otimes _{ g ^{-1} \smash{\D} ^{(m)}   _{Y}} ^\L g ^{-1} \E _1
\otimes ^\L _{g ^{-1}f ^{-1} \smash{\D} ^{(m)}   _{X}} g ^{-1}\E _2)}
}
$$
  est commutatif.
\end{lemm}

\begin{lemm}  \label{trans-proj-e2}
  Soient $\H \in D ( g ^{-1} f ^{-1} \smash{\D} ^{(m)}   _{X} \overset{^\mathrm{d}}{})$
  et $\E \in D ^\mathrm{b} _{\mathrm{qc},\mathrm{tdf}} (\overset{^\mathrm{g}}{} \smash{\D} ^{(m)}   _{X})$.
  Le diagramme
\begin{equation}
  \notag
  \xymatrix  @R=0,3cm    {
  { \R f _* (\R g _* \H \otimes ^\L _{f ^{-1} \smash{\D} ^{(m)}   _{X}} f ^{-1} \E)}
 \ar[d] ^-{\mathrm{proj} _g}
  &
  { \R f _* \R g _* \H \otimes ^\L _{\smash{\D} ^{(m)}   _{X}} \E }
  \ar[l] ^-{\mathrm{proj} _f}  \ar[d]
  \\
  { \R f _* \R g _* (\H \otimes ^\L _{g ^{-1} f ^{-1} \smash{\D} ^{(m)}   _{X}} g ^{-1}  f ^{-1} \E)}
  &
  { \R (f \circ g ) _* (\H \otimes ^\L _{(f \circ g ) ^{-1} \smash{\D} ^{(m)}   _{X}}   (f \circ g ) ^{-1} \E)}
    \ar[l] ^-{\mathrm{proj} _{f \circ g}}
}
\end{equation}
  est commutatif.
\end{lemm}

En utilisant le lemme \ref{trans-proj-e}, la transitivité des isomorphismes de la forme $\otimes$,
on vérifie la commutativité du rectangle (en haut) de \ref{pre-tr-transif-lemm-diag4}
où l'on a omis $\R f _*$.
Le diagramme \ref{pre-tr-transif-lemm-diag4} est donc commutatif.

On parvient par fonctorialité de \ref{trans-proj-e2} à la commutativité de rectangle du milieu de :
\begin{equation}
\label{pre-tr-transif-lemm-diag5}
  \xymatrix  @R=0,3cm   @C=0,3cm  {
  {\R f _* (\R g _* \R \underset{\smash{\D} ^{(m)}   _{Z}}{\mathcal{H}om}
  (\G, g ^! \smash{\D} ^{(m)}   _{Y} )
  \underset{\smash{\D} ^{(m)}   _{Y}}{\otimes ^\L}  f ^! \smash{\D} ^{(m)}   _{X}
  \underset{f ^{\text{-}1}\smash{\D} ^{(m)}   _{X}}{\otimes ^\L} f ^{\text{-}1} \E)}
  \ar[d] ^-{\mathrm{proj} _g}
  &
  {\R f _* (\R g _* \R \underset{\smash{\D} ^{(m)}   _{Z}}{\mathcal{H}om}
  (\G, g ^! \smash{\D} ^{(m)}   _{Y} )
 \!\!\!\!  \underset{\smash{\D} ^{(m)}   _{Y}}{\otimes ^\L}
 \!\!\!\!  f ^! \smash{\D} ^{(m)}   _{X})
\!\!\!\!  \underset{\smash{\D} ^{(m)}   _{X}}{\otimes ^\L}
\!\!\!\! \E}
  \ar[d] ^-{\mathrm{proj} _g}
  \ar[l] ^-{\mathrm{proj} _f}
  \\
  {\R f _* (\R g _* (\R \underset{\smash{\D} ^{(m)}   _{Z}}{\mathcal{H}om}  (\G, g ^! \smash{\D} ^{(m)}   _{Y} )
  \!\!\!\!\!\!\underset{g ^{\text{-}1}\smash{\D} ^{(m)}   _{Y}}{\otimes ^\L}
  \!\!\!\! g ^{\text{-}1} f ^! \smash{\D} ^{(m)}   _{X})
  \!\!\!\!\!\!  \underset{f ^{\text{-}1}\smash{\D} ^{(m)}   _{X}}{\otimes ^\L}
  \!\!\!\! f ^{\text{-}1} \E)}
  \ar[d] ^-{\otimes}
  &
  {\R f _* \R g _* (\R \underset{\smash{\D} ^{(m)}   _{Z}}{\mathcal{H}om}  (\G, g ^! \smash{\D} ^{(m)}   _{Y} )
  \!\!\!\!\!\!\underset{g ^{\text{-}1}\smash{\D} ^{(m)}   _{Y}}{\otimes ^\L}
  \!\!\!\! g ^{\text{-}1} f ^! \smash{\D} ^{(m)}   _{X})
  \!\!\!\! \underset{\smash{\D} ^{(m)}   _{X}}{\otimes ^\L}
  \!\!\!\! \E}
  \ar[d] ^-{\otimes}
  \ar[l] ^-{\mathrm{proj} _f}
  \\
  {\R f _* (\R g _* \R \mathcal{H}om _{\smash{\D} ^{(m)}   _{Z}} (\G, g ^!   f ^! \smash{\D} ^{(m)}   _{X})
  \underset{f ^{\text{-}1}\smash{\D} ^{(m)}   _{X}}{\otimes ^\L} f ^{\text{-}1} \E)}
  \ar[d] ^-{\mathrm{proj} _g}
  &
  {\R f _* \R g _* \R \mathcal{H}om _{\smash{\D} ^{(m)}   _{Z}} (\G, g ^!   f ^! \smash{\D} ^{(m)}   _{X})
  \underset{\smash{\D} ^{(m)}   _{X}}{\otimes ^\L} \E}
  \ar[dd]
  \ar[l] ^-{\mathrm{proj} _f}
  \\
  {\R f _* \R g _* (  \R \underset{\smash{\D} ^{(m)}   _{Z}}{\mathcal{H}om}
  (\G, g ^!   f ^! \smash{\D} ^{(m)}   _{X})
  \underset{g ^{\text{-}1}f ^{\text{-}1}\smash{\D} ^{(m)}   _{X}}{\otimes ^\L}
  g ^{\text{-}1} f ^{\text{-}1} \E)}
  \ar[d]
  &
  {}
  \\
  {\R f \circ g _* (  \R \underset{\smash{\D} ^{(m)}   _{Z}}{\mathcal{H}om}
  (\G, g ^!   f ^! \smash{\D} ^{(m)}   _{X})
  \underset{f\circ g ^{\text{-}1}\smash{\D} ^{(m)}   _{X}}{\otimes ^\L}
  f\circ g ^{\text{-}1} \E)}
  \ar[d]
  &
  {\R f \circ g _*   \R \underset{\smash{\D} ^{(m)}   _{Z}}{\mathcal{H}om}
  (\G, g ^!   f ^! \smash{\D} ^{(m)}   _{X})
  \underset{\smash{\D} ^{(m)}   _{X}}{\otimes ^\L} \E}
  \ar[d]
  \ar[l] ^-{\mathrm{proj} _{f\circ g}}
  \\
  {\R f \circ g _* (  \R \underset{\smash{\D} ^{(m)}   _{Z}}{\mathcal{H}om}
  (\G, f\circ g ^! \smash{\D} ^{(m)}   _{X})
  \underset{f\circ g ^{\text{-}1}\smash{\D} ^{(m)}   _{X}}{\otimes ^\L}
  f\circ g ^{\text{-}1} \E)}
  &
  {\R f \circ g _*   \R \underset{\smash{\D} ^{(m)}   _{Z}}{\mathcal{H}om}
  (\G, f\circ g ^! \smash{\D} ^{(m)}   _{X})
  \underset{\smash{\D} ^{(m)}   _{X}}{\otimes ^\L}  \E)}
  \ar[l] ^-{\mathrm{proj} _{f\circ g}}
  }
\end{equation}
On obtient des carrés par fonctorialité. Ainsi, \ref{pre-tr-transif-lemm-diag5} est commutatif.

Le composé de \ref{pre-tr-transif-lemm-diag4} et \ref{pre-tr-transif-lemm-diag5} correspond au
rectangle de \ref{pre-tr-transif-lemm-diag3}.
Celui-ci est donc commutatif.
On a donc vérifié la commutativité du diagramme \ref{pre-tr-transif-lemm-diag3}, i.e.,
celle du carré de devant composé avec celui du bas de \ref{pre-tr-transif-diag1}.
D'où le résultat.
\end{proof}

\begin{vide}\label{defadjH}
Soient $\E \in D ^\mathrm{b} _{\mathrm{qc},\mathrm{tdf}} (\overset{^\mathrm{g}}{} \smash{\D} ^{(m)}   _{X})$ et
  $\FF \in D ^\mathrm{b} _\mathrm{coh} (\overset{^\mathrm{g}}{} \smash{\D} ^{(m)}   _{Y})$.
 En appliquant $\R \Gamma (X, -)$ et $\mathrm{H} ^0 \circ \R \Gamma (X, -)$ à $\delta _f$, on obtient
 les bijections :
 \begin{gather}\label{defadjH1}
 \mathrm{adj} _{f,\FF,\E}\ : \
\R  \mathrm{Hom} _{\smash{\D} ^{(m)}   _{X}} (f _+ (\FF), \E) \riso
  \R\mathrm{Hom} _{\smash{\D} ^{(m)}   _{Y}} (\FF, f ^! \E ), \\
  \label{defadjH2}
   \mathrm{adj} _{f,\FF,\E}\ : \
  \mathrm{Hom} _{\smash{\D} ^{(m)}   _{X}} (f _+ (\FF), \E) \riso
  \mathrm{Hom} _{\smash{\D} ^{(m)}   _{Y}} (\FF, f ^! \E ).
 \end{gather}
 Si aucune confusion n'est à craindre, on les notera $\mathrm{adj} _f$.
\end{vide}

Il découle immédiatement de \ref{pre-tr-transif} le corollaire qui suit.
\begin{coro}
  \label{coropre-tr-transif}
 Soient
  $\E \in D ^\mathrm{b} _{\mathrm{qc},\mathrm{tdf}} (\overset{^\mathrm{g}}{} \smash{\D} ^{(m)}   _{X})$ et
  $\G \in D ^\mathrm{b} _\mathrm{coh} (\overset{^\mathrm{g}}{} \smash{\D} ^{(m)}   _{Z})$.
Le diagramme
\begin{equation}
  \label{coropre-tr-transifdiag0}
  \xymatrix  @R=0,3cm    {
  {\R \mathrm{Hom} _{\smash{\D} ^{(m)}   _{X}} (f _+ \circ g _+ (\G), \E)}
  \ar[d] ^-{\mathrm{adj} _g \circ \mathrm{adj} _f} _-\sim
  &
  {\R \mathrm{Hom} _{\smash{\D} ^{(m)}   _{X}} (f \circ g _+ (\G), \E)}
  \ar[l] _-\sim
  \ar[d] ^-{\mathrm{adj} _{f \circ g}} _-\sim
  \\
 {\R \mathrm{Hom} _{\smash{\D} ^{(m)}   _{Z}} (\G, g ^! f ^! \E )}
 &
 {\R \mathrm{Hom} _{\smash{\D} ^{(m)}   _{Z}} (\G, f \circ g ^! \E )}
 \ar[l] _-\sim
  }
\end{equation}
est commutatif. De même en enlevant {\og $\R$\fg}.
\end{coro}

\begin{vide}
Soit $\FF \in D ^\mathrm{b} _\mathrm{coh} (\overset{^\mathrm{g}}{} \smash{\D} ^{(m)}   _{Y})$.
  L'image de l'identité de $f _+ (\FF)$
  par la bijection $\mathrm{adj} _{f, \FF, f _+ (\FF)}$ de \ref{defadjH}
  donne le morphisme $\FF \rightarrow f ^! f _+ (\FF)$ noté $\mathrm{adj} _{f, \FF}$
  ou simplement $\mathrm{adj} _{\FF}$ ou $\mathrm{adj} _{f}$.

Soit $\E \in D ^\mathrm{b} _{\mathrm{coh}} (\overset{^\mathrm{g}}{} \smash{\D} ^{(m)}   _{X})$
tel que
$f ^! \E \in D ^\mathrm{b} _{\mathrm{coh}} (\overset{^\mathrm{g}}{} \smash{\D} ^{(m)}   _{Y})$.
On notera $\mathrm{adj} _{f,\E}$ : $f _+ \circ f ^! (\E) \rightarrow \E$,
  l'image inverse de l'identité de $f ^! \E $ par
  $\mathrm{adj} _{f, f ^! \E , \E  }$ ou $\mathrm{adj} _{\E}$ ou $\mathrm{adj} _{f}$.
\end{vide}

\begin{prop}
  \label{tr-transif}
  Pour tout $\G \in D ^\mathrm{b} _\mathrm{coh} (\overset{^\mathrm{g}}{} \smash{\D} ^{(m)}   _{Z}  )$ et
  pour tout
  $\E \in D ^\mathrm{b} _\mathrm{coh} (\overset{^\mathrm{g}}{} \smash{\D} ^{(m)}   _{X} )$ tel
  que $f ^! (\E) \in D ^\mathrm{b} _\mathrm{coh} (\overset{^\mathrm{g}}{} \smash{\D} ^{(m)}   _{Y})$
  et $ f \circ g ^!  (\E)
  \in D ^\mathrm{b} _\mathrm{coh} (\overset{^\mathrm{g}}{} \smash{\D} ^{(m)}   _{Z})$,
  les diagrammes
$$\xymatrix  @R=0,3cm   @C=2cm  {
{f _+  g _+  g ^!  f ^! (\E)}
\ar[r] ^-{f _+ \mathrm{adj} _{g, f ^! (\E)}}
\ar[d] _-\sim
&
{f _+  f ^! (\E)}
\ar[d] ^-{\mathrm{adj} _{f,\E}}
\\
{f \circ g _+  f \circ g ^! (\E)}
\ar[r] ^-{\mathrm{adj} _{f \circ g,\E}}
&
{\E}
}
\text{ et }
\hfill
\xymatrix  @R=0,3cm    @C=2cm {
{g ^! g _+ (\G)} \ar[r] ^-{g ^! \mathrm{adj} _{f,g _+ \G}}
&
{g ^! f ^! f _+ g _+ (\G)}
\ar[d] _-\sim
\\
{\G}
\ar[r] ^-{\mathrm{adj} _{f\circ g,\G}}
\ar[u] ^-{\mathrm{adj} _{g,\G}}
&
{ f \circ g ^! f \circ g _+  (\G)}
}
$$
sont commutatifs.
\end{prop}
\begin{proof}
Notons {\og $\mathrm{can}$\fg} les isomorphismes (ou induit fonctoriellement par ceux)
de la forme
$g ^! f ^! \riso (f \circ g) ^!$ ou
$f _+ g _+ \riso (f\circ g) _+$.
La démonstration de la commutativité du
diagramme de droite étant analogue, contentons-nous de vérifier celle de gauche.
À cette fin, considérons le diagramme ci-après :
\begin{equation}
   \label{tr-transif-diag1}
   \xymatrix  @R=0,3cm     {
   &
   {\mathrm{Hom} _{\smash{\D} ^{(m)}   _{Y}} (f^! \E , f ^! \E )}
   \ar[ld] _-{g ^!}
   \ar[d] ^-{\mathrm{adj} _{g, f ^! \E}} _-\sim
   &
   {\mathrm{Hom} _{\smash{\D} ^{(m)}   _{X}} (f _+ f^! \E ,  \E )}
   \ar[l] ^-{\mathrm{adj} _f} _-\sim
   \ar[d] ^-{f _+ \mathrm{adj} _{g, f ^! \E}} _-\sim
   \\
   {\mathrm{Hom} _{\smash{\D} ^{(m)}   _{Z}} (g ^! f^! \E ,g ^! f ^! \E )}
   &
   {\mathrm{Hom} _{\smash{\D} ^{(m)}   _{Y}} (g _+ g ^! f^! \E , f ^! \E )}
   \ar[l] ^-{\mathrm{adj} _g} _-\sim
   &
   {\mathrm{Hom} _{\smash{\D} ^{(m)}   _{X}} (f _+ g _+ g ^! f^! \E , \E )}
   \ar[l] ^-{\mathrm{adj} _f} _-\sim
   \\
   {\mathrm{Hom} _{\smash{\D} ^{(m)}   _{Z}} (f \circ g ^! \E ,g ^! f ^! \E )}
   \ar[u] ^-{\mathrm{can}} _-\sim
   &
   {\mathrm{Hom} _{\smash{\D} ^{(m)}   _{Y}} (g _+ f \circ g ^! \E , f ^! \E )}
   \ar[l] ^-{\mathrm{adj} _g} _-\sim
   \ar[u] ^-{\mathrm{can}} _-\sim
   &
   {\mathrm{Hom} _{\smash{\D} ^{(m)}   _{X}} (f _+ g _+ f \circ g ^! \E , \E )}
   \ar[l] ^-{\mathrm{adj} _f} _-\sim
   \ar[u] ^-{\mathrm{can}} _-\sim
   \\
   {\mathrm{Hom} _{\smash{\D} ^{(m)}   _{Z}} (f \circ g ^! \E ,f \circ g ^! \E )}
   \ar[u] ^-{\mathrm{can}} _-\sim
   &
   &
   {\mathrm{Hom} _{\smash{\D} ^{(m)}   _{X}} (f \circ g _+ f \circ g ^! \E , \E ),}
   \ar[ll] ^-{\mathrm{adj} _{f\circ g}} _-\sim
   \ar[u] ^-{\mathrm{can}} _-\sim
}
\end{equation}
dont les bijections horizontales sont les isomorphismes d'adjonction de \ref{defadjH2}.
Il dérive de \ref{coropre-tr-transif} la commutativité du rectangle du bas de \ref{tr-transif-diag1}.
Par fonctorialité en le terme de droite $f ^! \E$, il suffit de vérifier celle
du triangle pour l'identité de $f ^! \E$, ce qui est immédiat.
Enfin, celle des carrés s'établissent par fonctorialité.
Le diagramme \ref{tr-transif-diag1} est donc commutatif.

Or, la flèche $\mathrm{Hom} _{\smash{\D} ^{(m)}   _{X}} (f _+ g _+ g ^! f^! \E , \E )
\rightarrow
\mathrm{Hom} _{\smash{\D} ^{(m)}   _{Z}} (f \circ g ^! \E ,f \circ g ^! \E )$
de \ref{tr-transif-diag1} passant par le droite puis par le bas
(resp. par le haut puis par la gauche) , envoie
$\mathrm{adj} _{f \circ g, \E} \circ \mathrm{can} $
(resp. $\mathrm{adj} _{f,\E} \circ (f _+ \mathrm{adj} _{g, f ^! (\E)})$) sur l'identité.
\end{proof}

\subsection{Cas des schémas formels}
\begin{vide}\label{defisochgtbasepre}
Soient
  $\E , \FF \in D ^\mathrm{b} _{\mathrm{qc}} (\overset{^\mathrm{g}}{} \smash{\widehat{\D}} ^{(m)} _\X)$.
  Pour tout entier $i\geq 0$, on pose $\E _i :=
  \O _{S _i} \otimes ^\L _{\O _{\S}} \E  \riso
\smash{\D} ^{(m)} _{X _i} \otimes ^\L _{\smash{\widehat{\D}} ^{(m)} _\X} \E   $ (de même pour $\FF$).
  On dispose de l'{\it isomorphisme de changement de base} suivant :
\begin{equation}
  \label{defchgtbaseD}
  \O _{S _i} \otimes ^\L _{\O _{\S}}  \R \mathcal{H} om _{\smash{\widehat{\D}} ^{(m)} _\X} ( \E ,\FF)
  \riso
\R \mathcal{H} om _{\smash{\widehat{\D}} ^{(m)} _\X} ( \E ,\FF _i)
\liso
\R \mathcal{H} om _{\smash{\D} ^{(m)} _{X _i}} ( \E _i,\FF _i).
\end{equation}
Celui-ci est {\it transitif en $i$}, i.e., pour tout $1\leq i ' \leq i$, on construit de manière analogue
des isomorphismes
$\O _{S _{i'}} \otimes ^\L _{\O _{S _i}}
\R \mathcal{H} om _{\smash{\D} ^{(m)} _{X _i}} ( \E _i,\FF _i)
\riso
\R \mathcal{H} om _{\smash{\D} ^{(m)} _{X _{i'}}} ( \E _{i'},\FF _{i'})$,
ceux-ci étant compatibles avec \ref{defchgtbaseD}.

On remarque enfin que le composé
\begin{equation}
  \label{defchgtbaseD2}
\R \mathcal{H} om _{\smash{\widehat{\D}} ^{(m)} _\X} ( \E ,\FF)
\rightarrow   \O _{S _i} \otimes ^\L _{\O _{\S}}  \R \mathcal{H} om _{\smash{\widehat{\D}} ^{(m)} _\X} ( \E ,\FF)
  \riso
\R \mathcal{H} om _{\smash{\D} ^{(m)} _{X _i}} ( \E _i,\FF _i),
\end{equation}
dont l'isomorphisme est \ref{defchgtbaseD}, est le morphisme canonique
$\O _{S _i} \otimes ^\L _{\O _{\S}}-$ :
$\R \mathcal{H} om _{\smash{\widehat{\D}} ^{(m)} _\X} ( \E ,\FF)
\rightarrow
\R \mathcal{H} om _{\smash{\D} ^{(m)} _{X _i}} ( \E _i,\FF _i)$.
\end{vide}

\begin{vide}\label{defisochgtbase}
Soient $f$ : $ \Y \rightarrow \X$ un morphisme propre de $\V$-schémas formels lisses,
  $\E \in D ^\mathrm{b} _{\mathrm{qc}} (\overset{^\mathrm{g}}{} \smash{\widehat{\D}} ^{(m)} _\X)$.
Par construction de l'image inverse extraordinaire définie pour les complexes quasi-cohérents
(voir \cite[3.4.2.1]{Beintro2} et aussi \cite[3.2.2]{Beintro2}),
on dispose de l'{\it isomorphisme de changement de base de l'image inverse extraordinaire}
$\O _{S _i} \otimes ^\L _{\O _{\S}}   f ^! ( \E)
\riso
f _i ^! ( \E _i)$.
Lorsque $\E \in D ^\mathrm{b} _{\mathrm{coh}} (\overset{^\mathrm{g}}{} \smash{\widehat{\D}} ^{(m)} _\X)$,
on identifiera
$f ^! ( \E)$ avec
$\smash{\widehat{\D}} ^{(m)} _{\Y \rightarrow \X}
\otimes ^\L _{f ^{-1}  \smash{\widehat{\D}} ^{(m)} _{\X}} f ^{-1} \E [d _{Y/X}] $.
On bénéficie alors, pour tout $\E \in D ^\mathrm{b} _{\mathrm{coh}} (\overset{^\mathrm{g}}{} \smash{\widehat{\D}} ^{(m)} _\X)$,
 de l'{\it isomorphisme de changement de base} via les isomorphismes :
\begin{equation}
  \label{defchgtbasef!}
  \xymatrix  @R=0,3cm   {
  {\O _{S _i} \otimes ^\L _{\O _{\S}}   f ^! ( \E)}
  \ar@{=}[d]
  \ar@{.>}[rr]
  &&
{f _i ^! ( \E _i)}
  \ar@{=}[d]
  \\
  {\O _{S _i} \otimes ^\L _{\O _{\S}}   (\smash{\widehat{\D}} ^{(m)} _{\Y \rightarrow \X} \otimes ^\L
  _{f ^{-1}  \smash{\widehat{\D}} ^{(m)} _{\X}} f ^{-1} \E) [d _{Y/X}]}
  \ar[r] _-\sim
  &
{ \smash{\D} ^{(m)} _{Y _i  \rightarrow X _i } \otimes ^\L   _{f ^{-1}
  \smash{\widehat{\D}} ^{(m)} _{\X}} f ^{-1} \E [d _{Y/X}]}
  \ar[r] _-\sim
  &
 { \smash{\D} ^{(m)} _{Y _i  \rightarrow X _i } \otimes ^\L   _{f ^{-1}
  \smash{\D} ^{(m)} _{X _i}} f ^{-1} \E  _i [d _{Y/X}] .}
  }
\end{equation}
Comme pour \ref{defchgtbaseD}, celui-ci est transitif en $i$, i.e., on en tire
un isomorphisme entre les systèmes projectifs
$(\O _{S _i} \otimes ^\L _{\O _{\S}}   f ^! ( \E) ) _{i\in \N}$
et $(f _i ^! ( \E _i) )_{i\in \N}$.

De plus, pour tout
$\M \in D ^\mathrm{b} _{\mathrm{qc}} ( \smash{\widehat{\D}} ^{(m)} _\Y \overset{^\mathrm{d}}{})$,
on dispose tautologiquement de l'{\it isomorphisme de changement de base de l'image directe}
(voir \cite[3.5.1.1]{Beintro2}) :
$ f _+ ^\mathrm{d} ( \M ) \otimes ^\L _{\O _{\S}} \O _{S _i}
\riso
f _{i+}  ^\mathrm{d} (\M _i)$.
De même, lorsque
$\M \in D ^\mathrm{b} _{\mathrm{coh}} ( \smash{\widehat{\D}} ^{(m)} _\Y \overset{^\mathrm{d}}{})$, celui-ci
est défini via le diagramme commutatif :
\begin{equation}
  \label{defchgtbasef+}
  \xymatrix  @R=0,3cm    {
{ f _+ ^\mathrm{d} ( \M ) \otimes ^\L _{\O _{\S}} \O _{S _i} }
\ar@{=}[r]
\ar@{.>}[d] _-\sim
&
{(\R f _* ( \M \otimes ^\L _{\smash{\widehat{\D}} ^{(m)} _\Y}
\smash{\widehat{\D}} ^{(m)} _{\Y \rightarrow \X} ) )  \otimes ^\L _{\O _{\S}} \O _{S _i} }
\ar[r] ^-{\mathrm{proj}} _-\sim
&
{\R f _* ( \M \otimes ^\L _{\smash{\widehat{\D}} ^{(m)} _\Y}
\smash{\widehat{\D}} ^{(m)} _{\Y \rightarrow \X} \otimes ^\L _{\O _{\S}} \O _{S _i})}
\ar[d] _-\sim
\\
{f _{i+}  ^\mathrm{d} (\M _i)}
&
{\R f _* ( \M  _i \otimes ^\L _{\smash{\D} ^{(m)} _{Y _i}}
\smash{\D} ^{(m)} _{Y _i \rightarrow X _i } )}
\ar@{=}[l]
&
{\R f _* ( \M \otimes ^\L _{\smash{\widehat{\D}} ^{(m)} _\Y}
\smash{\D} ^{(m)} _{Y _i \rightarrow X _i } ).}
\ar[l] _-\sim
}
\end{equation}
Comme pour \ref{defchgtbaseD}, l'isomorphisme de gauche de \ref{defchgtbasef+} est transitif en $i$, i.e.,
induit un isomorphisme entre les systèmes projectifs
$( f _+ ^\mathrm{d} ( \M ) \otimes ^\L _{\O _{\S}} \O _{S _i} ) _{i\in \N}$
et
$(f _{i+}  ^\mathrm{d} (\M _i) ) _{i\in \N}$.
De même, en remplaçant dans \ref{defchgtbasef+} le symbole {\og d \fg} par {\og g \fg}.

\end{vide}

\begin{vide}
Soient $f$ : $ \Y \rightarrow \X$ un morphisme propre de $\V$-schémas formels lisses,
  $\E \in D ^\mathrm{b} _{\mathrm{coh}} (\smash{\widehat{\D}} ^{(m)} _\X)$
  et
  $\FF \in D ^\mathrm{b} _\mathrm{coh} (\smash{\widehat{\D}} ^{(m)} _\Y)$.
Comme pour \ref{defchi'f}, on dispose d'un isomorphisme canonique
\begin{equation}
  \label{defchi'fhat}
    \chi ' _f \ :\ f ^{\mathrm{d}} _+ \R \mathcal{H}om _{\smash{\widehat{\D}} ^{(m)}   _{\Y}}
(\FF, \smash{\widehat{\D}} ^{(m)}   _{\Y}) [d _Y]
  \riso
  \R \mathcal{H}om _{\smash{\widehat{\D}} ^{(m)}   _{\X}}
  (f ^{\mathrm{g}} _+ (\FF), \smash{\widehat{\D}} ^{(m)}   _{\X})[d _X].
\end{equation}
  Ce dernier s'inscrit
  dans le diagramme commutatif ci-après :
  \begin{equation}
    \label{chi'fchgtbase}
    \xymatrix  @R=0,3cm   {
    {\O _{S _i} \otimes ^\L _{\O _{\S}}  (
    f ^{\mathrm{d}} _+ \R \mathcal{H}om _{\smash{\widehat{\D}} ^{(m)}   _{\Y}}
(\FF, \smash{\widehat{\D}} ^{(m)}   _{\Y}) [d _Y])}
\ar[r] _-\sim ^{\chi _{f} '}
\ar[d] _-\sim
&
{\O _{S _i} \otimes ^\L _{\O _{\S}}  (
\R \mathcal{H}om _{\smash{\widehat{\D}} ^{(m)}   _{\X}}
  (f ^{\mathrm{g}} _+ (\FF), \smash{\widehat{\D}} ^{(m)}   _{\X})[d _X])}
  \ar[d] _-\sim
  \\
{  f _{i +} ^{\mathrm{d}}  \R \mathcal{H}om _{\smash{\D} ^{(m)}   _{Y _i }}
(\FF _i , \smash{\D} ^{(m)}   _{Y _i}) [d _Y]}
\ar[r] _-\sim ^{\chi ' _f}
&
{\R \mathcal{H}om _{\smash{\D} ^{(m)}   _{X _i }}
(f _{i +} ^{\mathrm{g}}  (\FF _i ), \smash{\D} ^{(m)}   _{X _i})[d _X],}
}
  \end{equation}
dont les isomorphisme verticaux sont induits par les isomorphismes de changement de base de \ref{defisochgtbasepre}
et \ref{defisochgtbase}.

De façon analogue à \ref{defdelta} (pour le morphisme de projection, nous utilisons la cohérence de $\E$),
on construit grâce à $\chi ' _f$ l'isomorphisme :
\begin{equation}
\label{defdelta2}
  \delta _f \ :\
  \R \mathcal{H}om _{\smash{\widehat{\D}} ^{(m)}   _{X}} (f _+ (\FF), \E)
  \riso
  \R f _* \R \mathcal{H} om _{\smash{\widehat{\D}} ^{(m)}   _{Y}} (\FF, f ^! \E ).
\end{equation}
Il dérive de \ref{chi'fchgtbase}, le diagramme commutatif suivant :
\begin{equation}
  \label{chgtbasedelta}
  \xymatrix  @R=0,3cm   {
  {\O _{S _i} \otimes ^\L _{\O _{\S}}  (
  \R \mathcal{H}om _{\smash{\widehat{\D}} ^{(m)}   _{\X}} (f _+ (\FF), \E))}
  \ar[r] _-\sim ^{\delta _f}
  \ar[d] _-\sim
  &
  {\O _{S _i} \otimes ^\L _{\O _{\S}}  (
  \R f _* \R \mathcal{H} om _{\smash{\widehat{\D}} ^{(m)}   _{\Y}} (\FF, f ^! \E ))}
    \ar[d] _-\sim
  \\
  {\R \mathcal{H}om _{\smash{\D} ^{(m)}   _{X _i }} (f _{i+} (\FF _i ), \E _i ) }
\ar[r] _-\sim ^{\delta _{f_i}}
&
{  \R f _* \R \mathcal{H} om _{\smash{\D} ^{(m)}   _{Y _i }} (\FF _i , f _i  ^! \E _i ) ,}
}
\end{equation}
où les isomorphisme verticaux se déduisent des isomorphismes de changement de base de \ref{defisochgtbasepre}
et \ref{defisochgtbase}.

Avec la remarque de \ref{defchgtbaseD2}, on déduit de \ref{chgtbasedelta} par application
du foncteur $H ^0 \circ \Gamma (\X, -)$, le diagramme commutatif suivant :
\begin{equation}
  \label{chgtbaseadj}
  \xymatrix  @R=0,3cm   {
  {\mathrm{Hom} _{\smash{\widehat{\D}} ^{(m)}   _{\X}} (f _+ (\FF), \E))}
  \ar[r] _-\sim ^{\mathrm{adj} _f}
  \ar[d]
  &
  {\mathrm{Hom} _{\smash{\widehat{\D}} ^{(m)}   _{\Y}} (\FF, f ^! \E ))}
    \ar[d]
  \\
  {\mathrm{Hom} _{\smash{\D} ^{(m)}   _{X _i }} (f _{i+} (\FF _i ), \E _i ) }
\ar[r] _-\sim ^{\mathrm{adj} _{f_i}}
&
{ \mathrm{Hom} _{\smash{\D} ^{(m)}   _{Y _i }} (\FF _i , f _i  ^! \E _i ) ,}
}
\end{equation}
dont les morphismes verticaux sont ceux induits fonctoriellement par
$\O _{S _i} \otimes ^\L _{\O _{\S}} -$.

Lorsque l'on suppose seulement
$\E \in D ^\mathrm{b} _{\mathrm{qc}} (\smash{\widehat{\D}} ^{(m)} _\X)$,
on construit l'isomorphisme d'adjonction
$\mathrm{adj} _f$ :
$\mathrm{Hom} _{\smash{\widehat{\D}} ^{(m)}   _{\X}} (f _+ (\FF), \E))
\riso
\mathrm{Hom} _{\smash{\widehat{\D}} ^{(m)}   _{\Y}} (\FF, f ^! \E ))$
via la commutativité du diagramme ci-dessous :
\begin{equation}
  \label{chgtbaseadjqc}
  \xymatrix  @R=0,3cm    @C=2cm {
  {\mathrm{Hom} _{\smash{\widehat{\D}} ^{(m)}   _{\X}} (f _+ (\FF), \E))}
  \ar@{.>}[r]  ^{\mathrm{adj} _f}
  \ar[d] _-\sim
  &
  {\mathrm{Hom} _{\smash{\widehat{\D}} ^{(m)}   _{\Y}} (\FF, f ^! \E ))}
    \ar[d] _-\sim
  \\
  {\underset{\underset{i}{\longleftarrow}}{\lim} \,
  \mathrm{Hom} _{\smash{\D} ^{(m)}   _{X _i }} (f _{i+} (\FF _i ), \E _i ) }
\ar[r] _-\sim ^{\underset{\underset{i}{\longleftarrow}}{\lim} \, \mathrm{adj} _{f_i}}
&
{\underset{\underset{i}{\longleftarrow}}{\lim} \, \mathrm{Hom} _{\smash{\D} ^{(m)}   _{Y _i }} (\FF _i , f _i  ^! \E _i ),}
}
\end{equation}
dont les isomorphismes verticaux sont dus à l'équivalence de catégorie (\cite[3.2.3]{Beintro2}).
Grâce à la commutativité de \ref{chgtbaseadj}, les deux constructions de l'isomorphisme d'adjonction
$\mathrm{adj} _f$ sont compatibles.

Pour tous morphismes propres
$g$ : $\ZZ \rightarrow \Y$ et $f $ : $\Y \rightarrow \X$ de $\V$-schémas formels lisses,
pour tout $\E \in D ^\mathrm{b} _{\mathrm{qc}} (\smash{\widehat{\D}} ^{(m)} _\X)$ et
tout $\G \in D ^\mathrm{b} _\mathrm{coh} (\smash{\widehat{\D}} ^{(m)} _\Y)$,
on obtient par complétion (voir \ref{chgtbaseadjqc})
de \ref{pre-tr-transifdiag0} le diagramme commutatif suivant
\begin{equation}
  \label{pre-tr-transifdiag0c}
  \xymatrix  @R=0,3cm    {
  { \mathrm{Hom} _{\smash{\widehat{\D}} ^{(m)}   _{\X}} (f _+ \circ g _+ (\G), \E)}
  \ar[d] ^-{\mathrm{adj} _g \circ \mathrm{adj} _f} _-\sim
  &
  { \mathrm{Hom} _{\smash{\widehat{\D}} ^{(m)}   _{\X}} (f \circ g _+ (\G), \E)}
  \ar[l] _-\sim
  \ar[d] ^-{\mathrm{adj} _{f \circ g}} _-\sim
  \\
 { \mathrm{Hom} _{\smash{\widehat{\D}} ^{(m)}   _{\ZZ}} (\G, g ^! f ^! \E )}
 &
 { \mathrm{Hom} _{\smash{\widehat{\D}} ^{(m)}   _{\ZZ}} (\G, f \circ g ^! \E ).}
 \ar[l] _-\sim
  }
\end{equation}

\end{vide}

\begin{vide}
  \label{chgtniv}
  Soient $f$ : $ \Y \rightarrow \X$ un morphisme propre de $\V$-schémas formels lisses,
  $m' \geq m$ deux entiers,
  $\E \in D ^\mathrm{b} _{\mathrm{qc}} (\smash{\widehat{\D}} ^{(m)} _\X)$,
  $\E '\in D ^\mathrm{b} _{\mathrm{qc}} (\smash{\widehat{\D}} ^{(m')} _\X)$,
  $\E \rightarrow \E '$ un morphisme $\smash{\widehat{\D}} ^{(m)} _\X$-linéaire,
  $\FF \in D ^\mathrm{b} _\mathrm{coh} (\smash{\widehat{\D}} ^{(m)} _\Y)$
  et
  $\FF ' := \smash{\widehat{\D}} ^{(m')} _\Y \otimes ^\L _{\smash{\widehat{\D}} ^{(m)} _\Y} \FF$.
  On bénéficie, via le diagramme commutatif ci-dessous,
  du {\it morphisme de changement de niveau}, noté $\mathrm{niv}$ :
  \begin{equation}
    \label{chgtnivfor1}
    \xymatrix  @R=0,3cm   {
    {\R \mathcal{H}om _{\smash{\widehat{\D}} ^{(m)}   _{\X}} (f _+ ^{(m)} (\FF), \E)}
    \ar[r]
    \ar@{.>}[d] ^-{\mathrm{niv}}
    &
    {\R \mathcal{H}om _{\smash{\widehat{\D}} ^{(m)}   _{\X}} (f _+ ^{(m)} (\FF), \E ')}
\\
{\R \mathcal{H}om _{\smash{\widehat{\D}} ^{(m')}   _{\X}} (  f _+ ^{(m ')} (\FF '), \E ')}
    \ar[r] _-\sim
  &
{ \R \mathcal{H}om _{\smash{\widehat{\D}} ^{(m')}   _{\X}} (
  \smash{\widehat{\D}} ^{(m')}   _{\X} \otimes _{\smash{\widehat{\D}} ^{(m)}   _{\X}} ^\L
  f _+ ^{(m)} (\FF), \E '),}
  \ar[u] _-\sim
 }
  \end{equation}
  dont l'isomorphisme du bas découle de \cite[3.5.3.(ii)]{Beintro2}.
On vérifie ensuite la commutativité du diagramme
\begin{equation}
  \label{chgtnivfor11}
  \xymatrix  @R=0,3cm   {
{\O _{S _i} \otimes ^\L _{\O _{\S}}  (
\R \mathcal{H}om _{\smash{\widehat{\D}} ^{(m)}   _{\X}} (f _+ ^{(m)} (\FF), \E)   )}
\ar[r] ^-{\mathrm{niv}} \ar[d] _-\sim
&
{\O _{S _i} \otimes ^\L _{\O _{\S}}  (
\R \mathcal{H}om _{\smash{\widehat{\D}} ^{(m')}   _{\X}} (  f _+ ^{(m ')} (\FF '), \E ')  )}
\ar[d] _-\sim
\\
{\R \mathcal{H}om _{\smash{\D} ^{(m)}   _{X _i }} (f _{i+} ^{(m)} (\FF _i ), \E _i)  }
\ar[r] ^-{\mathrm{niv}}
&
{\R \mathcal{H}om _{\smash{\D} ^{(m')}   _{X _i }} (f _{i+} ^{(m')} (\FF '_i ), \E ' _i)  ,}
}
\end{equation}
dont le morphisme de changement de niveau du bas se construit de manière analogue à \ref{chgtnivfor1}
et dont les isomorphisme verticaux se déduisent des isomorphismes de changement de base de \ref{defisochgtbasepre}
et \ref{defisochgtbase}.
Il en résulte celle du suivant :
\begin{equation}
  \label{chgtnivfor12}
  \xymatrix  @R=0,3cm   {
{\mathrm{Hom} _{\smash{\widehat{\D}} ^{(m)}   _{\X}} (f _+ ^{(m)} (\FF), \E)   }
\ar[r] ^-{\mathrm{niv}} \ar[d]
&
{\mathrm{Hom} _{\smash{\widehat{\D}} ^{(m')}   _{\X}} (  f _+ ^{(m ')} (\FF '), \E ')  }
\ar[d]
\\
{\mathrm{Hom} _{\smash{\D} ^{(m)}   _{X _i }} (f _{i+} ^{(m)} (\FF _i ), \E _i)  }
\ar[r] ^-{\mathrm{niv}}
&
{\mathrm{Hom} _{\smash{\D} ^{(m')}   _{X _i }} (f _{i+} ^{(m')} (\FF '_i ), \E ' _i)  .}
}
\end{equation}

On a aussi le morphisme de changement de niveau suivant :
  \begin{equation}
      \label{chgtnivfor2}
\mathrm{niv}\ :\  \R f _* \R \mathcal{H} om _{\smash{\widehat{\D}} ^{(m)}   _{\Y}} (\FF, f ^{! ^{(m)}} \E )
 \rightarrow
 \R f _* \R \mathcal{H} om _{\smash{\widehat{\D}} ^{(m)}   _{\Y}} (\FF, f ^{! ^{(m')}} \E ' )
 \liso
\R f _* \R \mathcal{H} om _{\smash{\widehat{\D}} ^{(m')}   _{\Y}} (\FF ', f ^{! ^{(m')}} \E ' ).
  \end{equation}
  De plus, le diagramme
  \begin{equation}
      \label{chgtnivfor21}
      \xymatrix  @R=0,3cm   {
{\O _{S _i} \otimes ^\L _{\O _{\S}}  (
 \R f _* \R \mathcal{H} om _{\smash{\widehat{\D}} ^{(m)}   _{\Y}} (\FF, f ^{! ^{(m)}} \E ) )}
\ar[r] ^-{\mathrm{niv}} \ar[d] _-\sim
&
{\O _{S _i} \otimes ^\L _{\O _{\S}}  (
\R f _* \R \mathcal{H} om _{\smash{\widehat{\D}} ^{(m')}   _{\Y}} (\FF ', f ^{! ^{(m')}} \E ' ) )}
\ar[d] _-\sim
\\
{ \R f _* \R \mathcal{H} om _{\smash{\D} ^{(m)}   _{Y _i}} (\FF _i , f _i ^{! ^{(m)}} \E _i )}
\ar[r] ^-{\mathrm{niv}}
&
{ \R f _* \R \mathcal{H} om _{\smash{\D} ^{(m')}   _{Y '_i}} (\FF '_i , f _i ^{! ^{(m')}} \E '_i ),}
}
  \end{equation}
dont le morphisme du bas se définit comme dans \ref{chgtnivfor2}
et dont les isomorphisme verticaux se déduisent des isomorphismes de changement de base de \ref{defisochgtbasepre}
et \ref{defisochgtbase},
est commutatif.
On obtient ainsi la commutativité du suivant :
  \begin{equation}
      \label{chgtnivfor22}
      \xymatrix  @R=0,3cm   {
{\mathrm{Hom} _{\smash{\widehat{\D}} ^{(m)}   _{\Y}} (\FF, f ^{! ^{(m)}} \E ) }
\ar[r] ^-{\mathrm{niv}} \ar[d]
&
{\mathrm{Hom} _{\smash{\widehat{\D}} ^{(m')}   _{\Y}} (\FF ', f ^{! ^{(m')}} \E ' ) }
\ar[d]
\\
{\mathrm{Hom} _{\smash{\D} ^{(m)}   _{Y _i}} (\FF _i , f _i ^{! ^{(m)}} \E _i )}
\ar[r] ^-{\mathrm{niv}}
&
{\mathrm{Hom} _{\smash{\D} ^{(m')}   _{Y '_i}} (\FF '_i , f _i ^{! ^{(m')}} \E '_i ).}
}
\end{equation}
Considérons enfin le diagramme
\begin{equation}
 \label{chghtnivcube}
 \xymatrix  @R=0,3cm   @C=0,3cm  {
 &
 {\mathrm{Hom} _{\smash{\widehat{\D}} ^{(m')}   _{\X}} (  f _+ ^{(m ')} (\FF '), \E ')  }
 \ar[rr] _-\sim ^-{\mathrm{adj} _f} \ar'[d][dd]
 &&
 {\mathrm{Hom} _{\smash{\widehat{\D}} ^{(m')}   _{\Y}} (\FF ', f ^{! ^{(m')}} \E ' ) }
 \ar[dd]
 \\
 {\mathrm{Hom} _{\smash{\widehat{\D}} ^{(m)}   _{\X}} (f _+ ^{(m)} (\FF), \E)   }
 \ar[rr] _-(0.3)\sim ^-(0.3){\mathrm{adj} _f} \ar[dd] \ar[ru] ^-{\mathrm{niv}}
 &&
 {\mathrm{Hom} _{\smash{\widehat{\D}} ^{(m)}   _{\Y}} (\FF, f ^{! ^{(m)}} \E )}
 \ar[dd] \ar[ru] ^-{\mathrm{niv}}
 \\
&
{\mathrm{Hom} _{\smash{\D} ^{(m)}   _{X _i }} (f _{i+} ^{(m')} (\FF '_i ), \E ' _i)  }
\ar'[r] _-\sim ^-{\mathrm{adj} _{f _i}} [rr]
&&
{\mathrm{Hom} _{\smash{\D} ^{(m')}   _{Y '_i}} (\FF '_i , f _i ^{! ^{(m')}} \E '_i )}
\\
{\mathrm{Hom} _{\smash{\D} ^{(m)}   _{X _i }} (f _{i+} ^{(m)} (\FF _i ), \E _i)  }
\ar[rr]_-\sim ^-{\mathrm{adj} _{f _i}} \ar[ru] ^-{\mathrm{niv}}
&&
{\mathrm{Hom} _{\smash{\D} ^{(m)}   _{Y _i}} (\FF _i , f _i ^{! ^{(m)}} \E _i ).}
\ar[ru] ^-{\mathrm{niv}}
}
\end{equation}
Par \ref{chgtnivfor12} (resp. \ref{chgtnivfor22}, resp. \ref{chgtbaseadj}),
le carré de gauche (resp. de droite, resp. de devant et de derrière) sont commutatifs.
La commutativité de l'isomorphisme de dualité relative de $f _i$ au changement de niveau implique
celle de $\delta _{f _i}$ puis celle de $\mathrm{adj} _{f_i}$. Le carré du bas est ainsi commutatif.
Or, en passant à la limite projective sur $i$, les flèches verticales de
\ref{chghtnivcube} deviennent des isomorphismes. Il en résulte la commutativité du carré supérieur
de \ref{chghtnivcube}.

\end{vide}

\begin{vide}
Soit $f$ : $\Y \rightarrow \X$ un morphisme propre de $\V$-schémas formels lisses.
Pour tous $\E ^{(\bullet)} \in \smash{\underset{^{\longrightarrow}}{LD}} ^{\mathrm{b}} _{\Q,\mathrm{qc}}
( \smash{\widehat{\D}} _{\X} ^{(\bullet)})$,
$\FF ^{(\bullet)} \in \smash{\underset{^{\longrightarrow}}{LD}} ^{\mathrm{b}} _{\Q,\mathrm{coh}}
( \smash{\widehat{\D}} _{\Y} ^{(\bullet)})$, par construction des
catégories de la forme
$\smash{\underset{^{\longrightarrow}}{LD}} ^{\mathrm{b}} _{\Q,\mathrm{qc}}
( \smash{\widehat{\D}} _{\X} ^{(\bullet)})$,
de la commutativité au changement de niveau des isomorphismes d'adjonction de
\ref{chgtbaseadjqc} (i.e. commutativité du carré supérieur de \ref{chghtnivcube}),
l'isomorphisme canonique d'adjonction suivant :
\begin{equation}\label{defadjdag}
  \mathrm{adj} _f \ : \
  \mathrm{Hom} _{\smash{\underset{^{\longrightarrow}}{LD}} ^{\mathrm{b}} _{\Q,\mathrm{qc}}
( \smash{\widehat{\D}} _{\X} ^{(\bullet)})} (f ^{(\bullet)} _+ (\FF ^{(\bullet)}), \E ^{(\bullet)})
\riso
  \mathrm{Hom} _{\smash{\underset{^{\longrightarrow}}{LD}} ^{\mathrm{b}} _{\Q,\mathrm{qc}}
( \smash{\widehat{\D}} _{\Y} ^{(\bullet)})} (\FF ^{(\bullet)}, f ^{! (\bullet) } \E ^{(\bullet)}).
\end{equation}

Pour tout morphisme propre
$g$ : $\ZZ \rightarrow \Y$ de $\V$-schémas formels lisses, pour tout
$\E ^{(\bullet)} \in \smash{\underset{^{\longrightarrow}}{LD}} ^{\mathrm{b}} _{\Q,\mathrm{qc}}
( \smash{\widehat{\D}} _{\X} ^{(\bullet)})$ et
 tout $\G ^{(\bullet)} \in \smash{\underset{^{\longrightarrow}}{LD}} ^{\mathrm{b}} _{\Q,\mathrm{coh}}
( \smash{\widehat{\D}} _{\ZZ} ^{(\bullet)})$,
il dérive de \ref{pre-tr-transifdiag0c} le diagramme commutatif :
\begin{equation}
  \label{pre-tr-transifdiag0dag}
  \xymatrix  @R=0,3cm    {
  { \mathrm{Hom} _{\smash{\underset{^{\longrightarrow}}{LD}} ^{\mathrm{b}} _{\Q,\mathrm{qc}}
( \smash{\widehat{\D}} _{\X} ^{(\bullet)})} (f _+ \circ g _+ (\G^{(\bullet)}), \E^{(\bullet)})}
  \ar[d] ^-{\mathrm{adj} _g \circ \mathrm{adj} _f} _-\sim
  &
  { \mathrm{Hom} _{\smash{\underset{^{\longrightarrow}}{LD}} ^{\mathrm{b}} _{\Q,\mathrm{qc}}
( \smash{\widehat{\D}} _{\X} ^{(\bullet)})} (f \circ g _+ (\G ^{(\bullet)}), \E ^{(\bullet)})}
  \ar[l] _-\sim
  \ar[d] ^-{\mathrm{adj} _{f \circ g}} _-\sim
  \\
 { \mathrm{Hom} _{\smash{\underset{^{\longrightarrow}}{LD}} ^{\mathrm{b}} _{\Q,\mathrm{qc}}
( \smash{\widehat{\D}} _{\ZZ} ^{(\bullet)})} (\G^{(\bullet)}, g ^! f ^! \E ^{(\bullet)})}
 &
 { \mathrm{Hom} _{\smash{\underset{^{\longrightarrow}}{LD}} ^{\mathrm{b}} _{\Q,\mathrm{qc}}
( \smash{\widehat{\D}} _{\ZZ} ^{(\bullet)})} (\G ^{(\bullet)}, f \circ g ^! \E ^{(\bullet)})}
 \ar[l] _-\sim
  }
\end{equation}

\end{vide}

\begin{vide}\label{videadjdag=}
Soient $f$ : $\Y \rightarrow \X$ un morphisme propre de $\V$-schémas formels lisses,
$\E ^{(\bullet)} \in \smash{\underset{^{\longrightarrow}}{LD}} ^{\mathrm{b}} _{\Q,\mathrm{coh}}
( \smash{\widehat{\D}} _{\X} ^{(\bullet)})$,
$\FF ^{(\bullet)} \in \smash{\underset{^{\longrightarrow}}{LD}} ^{\mathrm{b}} _{\Q,\mathrm{coh}}
( \smash{\widehat{\D}} _{\Y} ^{(\bullet)})$,
$\E := \underset{\longrightarrow}{\lim} \E ^{(\bullet)}$
et
$\FF := \underset{\longrightarrow}{\lim} \FF ^{(\bullet)}$ (le foncteur
$\underset{\longrightarrow}{\lim}$ est celui défini dans \cite[4.2.2]{Beintro2}).
On bénéficie à nouveau de l'isomorphisme :
$\chi '_f$ : $  f ^{\mathrm{d}} _+
\R \mathcal{H}om _{\smash{\D} ^{\dag}   _{\Y, \Q}} (\FF, \smash{\D} ^{\dag}   _{\Y, \Q}) [d _Y]
 \riso
 \R \mathcal{H}om _{\smash{\D} ^{\dag}   _{\X, \Q}} (f ^{\mathrm{g}} _+ (\FF), \smash{\D} ^{\dag}   _{\X, \Q})[d _X]$
 (voir \cite{Vir04}).
 On construit ensuite de manière analogue à \ref{defdelta}, l'isomorphisme
 $\delta _f \ :\
  \R \mathcal{H}om _{\smash{\D} ^{\dag}   _{\X,\Q}} (f _+ (\FF), \E)
  \riso
  \R f _* \R \mathcal{H} om _{\smash{\D} ^{\dag}   _{\Y,\Q}} (\FF, f ^! \E )$,
  puis l'isomorphisme d'adjonction :
  $\mathrm{adj} _f \ : \
  \mathrm{Hom} _{\smash{\D} ^{\dag}   _{\X,\Q}} (f _+ (\FF), \E)
  \riso
  \mathrm{Hom} _{\smash{\D} ^{\dag}   _{\Y,\Q}} (\FF, f ^! \E ).$

  Le diagramme suivant
\begin{equation}
  \label{adjdag=}
  \xymatrix  @R=0,3cm   {
  {\mathrm{Hom} _{\smash{\underset{^{\longrightarrow}}{LD}} ^{\mathrm{b}} _{\Q,\mathrm{qc}}
( \smash{\widehat{\D}} _{\X} ^{(\bullet)})} (f ^{(\bullet)} _+ (\FF ^{(\bullet)}), \E ^{(\bullet)})}
\ar[r] _-\sim ^-{\mathrm{adj}}
\ar[d] _-\sim ^{\underset{\longrightarrow}{\lim}}
&
{  \mathrm{Hom} _{\smash{\underset{^{\longrightarrow}}{LD}} ^{\mathrm{b}} _{\Q,\mathrm{qc}}
( \smash{\widehat{\D}} _{\Y} ^{(\bullet)})} (\FF ^{(\bullet)}, f ^{! (\bullet) } \E ^{(\bullet)} )}
\ar[d]  ^{\underset{\longrightarrow}{\lim}}
\\
{ \mathrm{Hom} _{\smash{\D} ^{\dag}   _{\X,\Q}} (f _+ (\FF), \E)}
\ar[r] _-\sim ^-{\mathrm{adj}}
&
{  \mathrm{Hom} _{\smash{\D} ^{\dag}   _{\Y,\Q}} (\FF, f ^! \E ),}
}
\end{equation}
où l'isomorphisme du haut est \ref{defadjdag},
est commutatif.

En effet, quitte à compliquer les notations, supposons que
$\E ^{(m)} \in D ^\mathrm{b} _{\mathrm{coh}} (\smash{\widehat{\D}} ^{(m)} _\X)$,
le morphisme
$\smash{\widehat{\D}} ^{(m')} _\X \otimes ^\L _{\smash{\widehat{\D}} ^{(m)} _\X}\E ^{(m)}
 \rightarrow \E ^{(m')}$
est un isomorphisme et de même en remplaçant $\E$ par $\FF$.
On dispose du diagramme commutatif :
\begin{equation}
  \label{chgtnivmdag}
  \xymatrix  @R=0,3cm   {
{ f ^{(m)\mathrm{d}} _+
\R \mathcal{H}om _{\smash{\D} ^{(m)}   _{\Y}} (\FF^{(m)} , \smash{\D} ^{(m)}   _{\Y}) [d _Y]
\otimes _{\smash{\D} ^{(m)}   _{\X }}  \smash{\D} ^{\dag}   _{\X,\Q }}
\ar[r] _-\sim ^{\chi '_f}
\ar[d] _-\sim
&
{ \R \mathcal{H}om _{\smash{\D} ^{(m)}   _{\X}} (f ^{(m)\mathrm{g}} _+ (\FF ^{(m)}), \smash{\D} ^{(m)}   _{\X })[d _X]
\otimes _{\smash{\D} ^{(m)}   _{\X }}  \smash{\D} ^{\dag}   _{\X,\Q }}
\ar[d] _-\sim
\\
{ f ^{\mathrm{d}} _+
\R \mathcal{H}om _{\smash{\D} ^{\dag}   _{\Y, \Q}} (\FF, \smash{\D} ^{\dag}   _{\Y, \Q}) [d _Y]}
\ar[r] _-\sim ^{\chi '_f}
&
{ \R \mathcal{H}om _{\smash{\D} ^{\dag}   _{\X, \Q}} (f ^{\mathrm{g}} _+ (\FF), \smash{\D} ^{\dag}   _{\X, \Q})[d _X],}
}
\end{equation}
dont les isomorphismes verticaux se déduisent de \cite[3.5.3.(ii)]{Beintro2} et de la
commutation des foncteurs duaux à l'extension des scalaires.
Il en dérive le carré commutatif suivant :
\begin{equation}
  \label{adjdag=2}
  \xymatrix  @R=0,3cm   {
  { \mathrm{Hom} _{\smash{\D} ^{(m)}   _{\X}} (f _+ ^{(m)}(\FF ^{(m)}), \E ^{(m)})}
\ar[r] _-\sim ^-{\mathrm{adj}_f}
\ar[d] _-\sim
&
{  \mathrm{Hom} _{\smash{\D} ^{(m)}   _{\Y}} (\FF ^{(m)}, f ^{!(m)} \E ^{(m)} )}
\ar[d] _-\sim
  \\
{ \mathrm{Hom} _{\smash{\D} ^{\dag}   _{\X,\Q}} (f _+ (\FF), \E)}
\ar[r] _-\sim ^-{\mathrm{adj}_f}
&
{  \mathrm{Hom} _{\smash{\D} ^{\dag}   _{\Y,\Q}} (\FF, f ^! \E ).}
}
\end{equation}
Avec les notations de \cite[4.2.2]{Beintro2}, il en découle, pour tous $\lambda \in L$, $\chi \in M$,
le diagramme commutatif suivant :
\begin{equation}
  \label{adjdag=3}
  \xymatrix  @R=0,3cm   @C=0,3cm {
  { \mathrm{Hom} _{D ( \smash{\D} ^{(\bullet)}   _{\X})}
  (f _+ ^{(\bullet)}(\FF ^{(\bullet)}), \lambda ^* \chi ^* \E ^{(\bullet)})}
\ar[r] _-\sim ^-{\mathrm{adj}_f}
\ar[d] ^-{\underset{\longrightarrow}{\lim}}
&
{  \mathrm{Hom} _{D(\smash{\D} ^{(\bullet)}   _{\Y})}
(\FF ^{(\bullet)}, f ^{!(\bullet)} \lambda ^* \chi ^* \E ^{(\bullet)} )}
\ar[d] ^-{\underset{\longrightarrow}{\lim}}
\ar[r] _-\sim
&
{  \mathrm{Hom} _{D(\smash{\D} ^{(\bullet)}   _{\Y})}
(\FF ^{(\bullet)}, \lambda ^* \chi ^* f ^{!(\bullet)} \E ^{(\bullet)} )}
\ar[d] ^-{\underset{\longrightarrow}{\lim}}
  \\
 { \mathrm{Hom} _{\smash{\D} ^{\dag}   _{\X,\Q}}
  (\underset{\longrightarrow}{\lim} f _+ ^{(\bullet)}(\FF ^{(\bullet)}),
  \underset{\longrightarrow}{\lim} \lambda ^* \chi ^* \E ^{(\bullet)})}
\ar[r] _-\sim ^-{\mathrm{adj}_f}
\ar[d] _-\sim
&
{  \mathrm{Hom} _{\smash{\D} ^{\dag}   _{\Y,\Q}}
(\FF,
\underset{\longrightarrow}{\lim} f ^{!(\bullet)} \lambda ^* \chi ^* \E ^{(\bullet)} )}
\ar[d] _-\sim
\ar[r] _-\sim
&
{  \mathrm{Hom} _{\smash{\D} ^{\dag}   _{\Y,\Q}}
(\FF,
\underset{\longrightarrow}{\lim} \lambda ^* \chi ^*  f ^{!(\bullet)} \E ^{(\bullet)} )}
\ar[dl] _-\sim
\\
{ \mathrm{Hom} _{\smash{\D} ^{\dag}   _{\X,\Q}} (f _+ (\FF), \E)}
\ar[r] _-\sim ^-{\mathrm{adj}_f}
&
{  \mathrm{Hom} _{\smash{\D} ^{\dag}   _{\Y,\Q}} (\FF, f ^! \E ),}
}
\end{equation}
où $\underset{\longrightarrow}{\lim}$ est le foncteur de \cite[4.2.2]{Beintro2} (ou \cite[4.2.4]{Beintro2}).

En passant à la limite inductive sur $\lambda \in L$, $\chi \in M$, on déduit de \ref{adjdag=3}
le diagramme \ref{adjdag=}.

\end{vide}

\begin{vide}
Soient $f$ : $\Y  \rightarrow \X $ un morphisme propre de $\V$-schémas formels lisses,
$T$ un diviseur de $X $ tel que $T ':= f  ^{-1}  (T)$ soit un diviseur de $Y$.

D'après \cite[1.2.10]{caro_courbe-nouveau}, pour tous
$\E  \in D ^{\mathrm{b}} _\mathrm{coh}(\overset{^{\mathrm{g}}}{} \smash{\D} ^\dag _{\X } (\hdag T ) _{\Q} )$
et
$\FF  \in D ^{\mathrm{b}} _\mathrm{coh}(\smash{\D} ^\dag _{\Y } (\hdag T ' ) _{\Q} \overset{^{\mathrm{d}}}{})$,
on dispose de l'isomorphisme canonique d'adjonction :
\begin{equation}
  \label{adjnonreldiv}
  \mathrm{adj} _{f,T} \  :\ \mathrm{Hom} _{\smash{\D} ^\dag _{\X } (\hdag T ) _{\Q}} ( f _{T+}  (\FF ), \E )
\riso
{\mathrm{Hom} _{\smash{\D} ^\dag _{\Y } (\hdag T ') _{\Q}} ( \FF , f _{T} ^{!} \E )}.
\end{equation}

Il en découle comme d'habitude le morphisme d'adjonction :
$\mathrm{adj} _{f,T,\FF}$ : $\FF \rightarrow f _{T} ^{!} f _{T+}  (\FF )$.
De même, si
$f _{T} ^{!} \E \in D ^{\mathrm{b}} _\mathrm{coh}(\smash{\D} ^\dag _{\Y } (\hdag T ' ) _{\Q} \overset{^{\mathrm{d}}}{})$,
il en dérive le morphisme d'adjonction :
$\mathrm{adj} _{f,T,\E}$ : $f _{T+} f _{T} ^{!}  (\E ) \rightarrow \E$.
Lorsqu'aucune confusion n'est à craindre,
on les notera simplement $\mathrm{adj} _{f,T}$, $\mathrm{adj} _f$ etc.
Lorsque le diviseur est vide, on ne l'indiquera pas. Le morphisme d'adjonction $\mathrm{adj} _f$
de \ref{videadjdag=} est le même que \ref{adjnonreldiv}.
\end{vide}

\begin{prop}\label{proptr-transifdiag0dag}
  Soient $g$ : $ \ZZ \rightarrow \Y$ et $f$ : $\Y \rightarrow \X$ deux morphismes propres de
  $\V$-schémas formels lisses, $T$ un diviseur de $X$ tel que $T':=f ^{-1} (T)$ (resp. $T'':=g ^{-1} (T')$)
  soit un diviseur de $Y$ (resp. $Z$).

  Pour tout $\G \in D ^\mathrm{b} _\mathrm{coh} (\overset{^\mathrm{g}}{} \smash{\D} ^{\dag}   _{\ZZ} (T'')_{\Q} )$ et
  pour tout
  $\E \in D ^\mathrm{b} _\mathrm{coh} (\overset{^\mathrm{g}}{} \smash{\D} ^{\dag}   _{\X} (T)_{\Q} )$ tel
  que $f ^! _T (\E) \in D ^\mathrm{b} _\mathrm{coh} (\overset{^\mathrm{g}}{} \smash{\D} ^{\dag}_{\Y} (T')_{\Q})$
    et $ f \circ g ^! _T (\E)
  \in D ^\mathrm{b} _\mathrm{coh} (\overset{^\mathrm{g}}{} \smash{\D} ^{\dag} _{\ZZ} (T'')_{\Q})$,
  le diagramme
\begin{equation}
  \label{tr-transifdiag0dag}
  \xymatrix  @R=0,3cm    {
  { \mathrm{Hom} _{\smash{\D} ^{\dag}   _{\X} (T)_{\Q} }
  (f _{T+} \circ g _{T'+} (\G), \E)}
  \ar[d] ^-{\mathrm{adj} _{g,T'} \circ \mathrm{adj} _{f,T}} _-\sim
  &
  { \mathrm{Hom} _{\smash{\D} ^{\dag}   _{\X} (T)_{\Q}} (f \circ g _+ (\G), \E)}
  \ar[l] _-\sim
  \ar[d] ^-{\mathrm{adj} _{f \circ g,T}} _-\sim
  \\
 { \mathrm{Hom} _{\smash{\D} ^{\dag}   _{\ZZ} (T'')_{\Q} } (\G, g ^! f ^! \E )}
 &
 { \mathrm{Hom} _{\smash{\D} ^{\dag}   _{\ZZ} (T'')_{\Q} } (\G, f \circ g ^! \E ),}
 \ar[l] _-\sim
  }
\end{equation}
  où les isomorphismes horizontaux sont induits par les isomorphismes canoniques de composition (\ref{isocanocomps}),
  est commutatif.
\end{prop}
\begin{proof}
Par \cite[4.3.12]{Beintro2}, on peut supposer le diviseur $T$ vide.
Cela découle alors de \ref{pre-tr-transifdiag0dag} et \ref{adjdag=}.
\end{proof}

\begin{rema}
\label{remacomp+}
  Avec les notations et hypothèses de \ref{proptr-transifdiag0dag},
l'hypothèse {\og $f ^! _T (\E) \in D ^\mathrm{b} _\mathrm{coh} (\overset{^\mathrm{g}}{} \smash{\D} ^{\dag}_{\Y} (T')_{\Q})$
    et $ f \circ g ^! _T (\E)
  \in D ^\mathrm{b} _\mathrm{coh} (\overset{^\mathrm{g}}{} \smash{\D} ^{\dag} _{\ZZ} (T'')_{\Q})$ \fg}
  est toujours satisfaite lorsque $f$ et $g $ sont, en plus d'être propres, lisses.
  De plus, ces hypothèses sont aussi validées lorsque l'on travaille avec des
  complexes surcohérents (voir \cite{caro_surcoherent}) etc.
  Dans ces conditions, il découle de \ref{proptr-transifdiag0dag}
  que l'isomorphisme canonique de composition des images directes utilisé ici
  (voir \ref{f+og+comm}) est
  le même que celui construit dans \cite[1.2.15]{caro_courbe-nouveau}. On a ainsi obtenu
  une unification.
\end{rema}

\begin{lemm}\label{f!f+coh}
  Soient $f$ : $\Y \rightarrow \X$ un morphisme propre de $\V$-schémas formels séparés lisses,
  $T$ un diviseur de $X$ tel que $T':=f ^{-1} (T)$ soit un diviseur de $Y$,
  et $\FF \in D ^\mathrm{b} _\mathrm{coh} (\overset{^\mathrm{g}}{} \smash{\D} ^{\dag}   _{\Y} (T')_{\Q})$.

  Alors, $f _T ^! \circ f _{T,+} (\FF) \in
  D ^\mathrm{b} _\mathrm{coh} (\overset{^\mathrm{g}}{} \smash{\D} ^{\dag}   _{\Y} (T')_{\Q})$.
\end{lemm}
\begin{proof}
Quitte à compliquer les notations, supposons le diviseur $T$ vide.
  On note $\delta _{f}=(\mathrm{id}, f)$ : $ \Y \hookrightarrow \Y \times \X$,
  $\E := \delta _{f,+}  (\FF)$,
  $p$ : $\Y \times \X \rightarrow \X$ la deuxième projection,
  $g= f \times \mathrm{id}\ : \ \Y \times \X \rightarrow \X \times \X$
  et $\delta  \ : \ \X \hookrightarrow \X \times \X$ l'immersion fermée diagonale.
D'après \cite[5.3.7]{EGAI}, le carré suivant
\begin{equation}
  \label{f!f+cohdiag1}
  \xymatrix  @R=0,3cm  {
  {\Y \times \X } \ar[r] ^-g \ar[dr] ^-{p}
  &
  {\X \times \X}
  \\
  {\Y} \ar[r] ^-f \ar[u] ^-{\delta _f}
  &
  {\X} \ar[u] ^-{\delta}
  }
\end{equation}
  est cartésien.
Comme $g _+  (\E) = g _+ \circ \delta _{f,+} (\FF) \riso \delta _+ \circ f_+ (\FF)$,
$g _+ (\E)$ est à support dans $X$ (via $\delta$). Comme le carré de \ref{f!f+cohdiag1}
est cartésien, $g ^! \circ g _+ (\E)$ est à support dans $Y$ (via $\delta _f$).
Or, $p ^! \circ p _+ (\E) \riso p ^! \circ \delta ^! \circ \delta _+ \circ p _+ (\E)
\riso g ^! \circ g _+ (\E)$.
Ainsi, $p ^! \circ p _+ (\E) $ est à support dans $Y$.
Or, $p _+ (\E) \riso f _+ (\FF)
\in D ^\mathrm{b} _\mathrm{coh} (\overset{^\mathrm{g}}{} \smash{\D} ^{\dag}   _{\X,\Q})$
(car comme $f$ est propre, la cohérence est préservée par image directe).
Comme $p$ est lisse, il en découle que
$p ^! \circ p _+ (\E) \in
D ^\mathrm{b} _\mathrm{coh} (\overset{^\mathrm{g}}{} \smash{\D} ^{\dag} _{\Y \times \X,\Q})$.
D'après l'analogue $p$-adique de Berthelot du théorème de Kashiwara, il en résulte que
$\delta _f   ^! (p ^! \circ p _+ (\E))
\in D ^\mathrm{b} _\mathrm{coh} (\overset{^\mathrm{g}}{} \smash{\D} ^{\dag} _{\Y,\Q})  $.
On conclut en remarquant :
$\delta _f   ^! (p ^! \circ p _+ (\E))
=\delta _f   ^! \circ p ^! \circ p _+ \circ \delta _{f,+} (\FF)
\riso
f ^! \circ f _+ (\FF)$.

\end{proof}

\begin{prop}
  \label{tr-transifdag}
Soient $g$ : $ \ZZ \rightarrow \Y$ et $f$ : $\Y \rightarrow \X$ deux morphismes propres de
  $\V$-schémas formels séparés et lisses,
  $T$ un diviseur de $X$ tel que $T':=f ^{-1} (T)$ (resp. $T'':=g ^{-1} (T')$)
  soit un diviseur de $Y$ (resp. $Z$).

  Pour tout $\G \in D ^\mathrm{b} _\mathrm{coh} (\overset{^\mathrm{g}}{} \smash{\D} ^{\dag}   _{\ZZ} (T'')_{\Q} )$ et
  pour tout
  $\E \in D ^\mathrm{b} _\mathrm{coh} (\overset{^\mathrm{g}}{} \smash{\D} ^{\dag}   _{\X} (T)_{\Q})$ tel
  que $f ^! _T (\E) \in D ^\mathrm{b} _\mathrm{coh} (\overset{^\mathrm{g}}{} \smash{\D} ^{\dag}_{\Y} (T')_{\Q})$
  et $ f \circ g ^! _T (\E)
  \in D ^\mathrm{b} _\mathrm{coh} (\overset{^\mathrm{g}}{} \smash{\D} ^{\dag} _{\ZZ} (T'')_{\Q})$,
  on dispose des diagrammes commutatifs :
\begin{equation}
  \label{tr-transifdagdiag}
\xymatrix  @R=0,3cm    {
{f _{T,+}  g _{T',+}  g _{T'} ^!  f _T ^! (\E)}
\ar[r] ^-{\mathrm{adj} _{g,T'}}
\ar[d] _-\sim
&
{f _{T,+}  f ^! _T (\E)}
\ar[d] ^-{\mathrm{adj} _{f,T}}
\\
{f \circ g _{T,+}  f \circ g ^! _T (\E)}
\ar[r] ^-{\mathrm{adj} _{f \circ g, T}}
&
{\E}
}
\text{ et }
\hfill
\xymatrix  @R=0,3cm   {
{g ^! _{T'} g _{T',+} (\G)} \ar[r] ^-{\mathrm{adj} _{f,T}}
&
{g ^! _{T'} f ^! _T f _{T,+} g _{T',+} (\G)}
\ar[d] _-\sim
\\
{\G}
\ar[r] ^-{\mathrm{adj} _{f\circ g ,T}}
\ar[u] ^-{\mathrm{adj} _{g,T'}}
&
{ f \circ g _T ^! f \circ g _{T,+}  (\G).}
}
\end{equation}
\end{prop}
\begin{proof}
    Par \cite[4.3.12]{Be1} et via \ref{f!f+coh} (pour le diagramme de droite),
  on se ramène au cas où le diviseur $T$ est vide.
  En reprenant la preuve de \ref{tr-transif} (qui utilise \ref{coropre-tr-transif}),
  cela dérive alors de \ref{proptr-transifdiag0dag}.
\end{proof}

\begin{rema}\label{tr-transifdagdiagrema}
  L'hypothèse que $X$ soit séparé est inutile pour valider la commutativité du diagramme de gauche de
\ref{tr-transifdagdiag} car celle-ci n'utilise pas le lemme \ref{f!f+coh}.
Il en est de même pour le diagramme de droite lorsque $f$ est une immersion fermée
ou lorsque les termes de ce diagramme sont $\smash{\D} ^{\dag}   _{\ZZ} (T'')_{\Q}$-cohérents.
\end{rema}

\section{Foncteur $\sp _+$ associant un $\smash{\D}$-module arithmétique à un isocristal surconvergent. Cas lisse}

Nous adopterons les notations suivantes :
si $f$ : $\PP ' \rightarrow \PP $ est un morphisme de $\V$-schémas
formels lisses et $T $ est un diviseur de $P $ tel que $f ^{-1}(T )$ soit
un diviseur de $P '$, on notera
$ f_K$ : $ \smash{\PP '} _K \rightarrow \smash{\PP } _K$
le morphisme d'espaces analytiques rigides associé à $f$, tandis que
$\sp $ : $ \PP _{K} \rightarrow \PP $ (ou $\sp $ : $ \PP '_{K} \rightarrow \PP '$)
sera le morphisme de spécialisation (voir \cite{Berig}).

Si $j  $ : $Y  \hookrightarrow X $ est une immersion ouverte de $k$-schémas
de type fini, telle qu'il existe une immersion fermée $X \hookrightarrow \PP $, avec
$\PP $ un $\V$-schéma formel lisse sur un voisinage de $X$,
alors $j ^\dag$ signifiera le foncteur \textit{faisceau des germes de sections
surconvergentes le long de $X  \setminus Y $} (\cite[2.1.1]{Berig}).

\subsection{Isomorphismes de recollement : cas formel}

\begin{vide}\label{relev-comp-adj-immf0}
Considérons le diagramme commutatif de morphismes de $\Z _{(p)}$-schémas
  \begin{equation}
    \notag
    \xymatrix  @R=0,3cm {
     { X}
     \ar[r] ^f
     \ar[d]
     &
     {Y}
     \ar[r]^g
     \ar[d]
     &
     {Z}
     \ar[d]
     \\
     {S}
     \ar[r]
     &
     {T}
     \ar[r]
     &
     {U}
     }
  \end{equation}
  dans lequel $X$, $Y$ et $Z$ sont respectivement lisses sur $S$, $T$ et $U$.
De plus, on suppose que les schémas $S$, $T$ et $U$ sont munis respectivement pour un entier $m$ donné
de $m$-PD-idéaux quasi-cohérents $m$-PD-nilpotents
$(\mathfrak{a} _S,\,\mathfrak{b} _S,\, \alpha _S) \subset \O _S$,
$(\mathfrak{a} _T,\,\mathfrak{b} _T,\, \alpha _T) \subset \O _T$
et
$(\mathfrak{a} _U,\,\mathfrak{b} _U,\, \alpha _U) \subset \O _U$ de telle manière que
les morphismes $S \rightarrow T$ et $T\rightarrow U$ soient des $m$-PD-morphismes.
   On note $X _0$, $Y _0$ et $Z _0$ les réductions respectives de $X$, $Y$ et $Z$ modulo
   $\mathfrak{a} _S$, $\mathfrak{a} _T$ et $\mathfrak{a} _U$. On désigne ensuite par
   $f _0$ : $X _0 \rightarrow Y _0$ et $g _0$ : $ Y _0 \rightarrow Z _0$ les morphismes
   induits par $f$ et $g$.
  Soient respectivement $H _X$, $H _Y$ et $H _Z$
  des diviseurs de $X _0$, $Y _0$ et $Z _0$ induisant les factorisations
  $f _0 ^{-1} ( H _Y) \hookrightarrow  H _X$ et $g _0 ^{-1} (H _Z) \hookrightarrow H _Y$.
  En désignant par $n _m \geq m$ un entier,
  avec les notations de \cite[2.1]{huyghe2}, on dispose alors de morphismes
  $f  ^{*} \B _Y ^{(n_m)} ( H _Y)\rightarrow \B _X ^{(n_m)}( H _X)$
  et
  $g  ^{*} \B _Z ^{(n_m)} ( H _Z)\rightarrow \B _Y ^{(n_m)}( H _Y)$
  (\cite[2.1.2]{huyghe2}).

En notant $\widetilde{X}$, $\widetilde{Y}$ et $\widetilde{Z}$ les espaces annelés
$(X,\, \B _X^{(n_m)}( H _X))$, $(Y,\, \B _Y ^{(n_m)}( H _Y))$ et $(Z,\, \B _Z^{(n_m)}( H _Z))$,
on obtient des morphismes d'espaces annelés
$\tilde{f} $ : $ \widetilde{X} \rightarrow \widetilde{Y}$
et
$\tilde{g} $ : $ \widetilde{Y} \rightarrow \widetilde{Z}$.
Leur réduction modulo $\mathfrak{a}$ donne des morphismes que l'on notera
$\tilde{f} _0 $ : $ \widetilde{X} _0\rightarrow \widetilde{Y} _0$
et
$\tilde{g}_0 $ : $ \widetilde{Y} _0\rightarrow \widetilde{Z}_0$.

On pose
$\smash{\widetilde{\PP}} ^n _{X/ S(m)}:= \B _X ^{(n_m)}( H _X) \otimes _{\O _X } \PP ^n _{X/ S(m)} $,
$\smash{\widetilde{\PP}} ^n _{Y/ T(m)}:= \B _Y ^{(n_m)}( H _Y) \otimes _{\O _Y } \PP ^n _{Y/ T(m)} $ et
$\smash{\widetilde{\PP}} ^n _{Z/ U(m)}:= \B _Z ^{(n_m)}( H _Z) \otimes _{\O _Z } \PP ^n _{Z/ U(m)} $.

Le faisceau $\smash{\widetilde{\PP}} ^n _{X/ S(m)}$
est muni de deux structures canoniques (à droite et à gauche) de $\B _X ^{(n_m)}( H _X)$-algèbres
(de même pour $\smash{\widetilde{\PP}} ^n _{Y/ T(m)}$ et $\smash{\widetilde{\PP}} ^n _{Z/ U(m)}$).
Lorsqu'aucune confusion n'est à craindre, on omet d'indiquer la base ($S$, $T$ ou $U$).
On notera $\widetilde{\D} ^{(m)} _{X/S} = \B _X ^{(n_m)}( H _X) \otimes _{\O _X } \smash{\D} ^{(m)} _{X/S}$,
  $\widetilde{\D} ^{(m)} _{Y/T} = \B _Y ^{(n_m)}( H _Y) \otimes _{\O _Y } \smash{\D} ^{(m)} _{Y/T}$ et
$\tilde{f} ^!$ : $D ^- (\widetilde{\D} ^{(m)} _{Y/T})
\rightarrow D ^-(\widetilde{\D} ^{(m)} _{X/S})$, le foncteur image inverse
extraordinaire.

D'après \cite[2.1.5]{Be1} et avec les notations de \cite[2.1.2]{Be1},
pour tout $S$-morphisme $f'$ : $ X \rightarrow Y$
ayant la même restriction $X _0 \rightarrow Y$ que $f$,
on dispose, pour $n$ assez grand,
de la factorisation $\delta ^n  _{f ,f'}$ : $X \rightarrow \Delta ^n _{Y/T(m)} $ (grâce à la $m$-PD-nilpotence de
$(\mathfrak{a} _S,\,\mathfrak{b} _S,\, \alpha _S)$) rendant commutatif le diagramme
\begin{equation}
  \label{ber1.2.1.5}
\xymatrix  @R=0,3cm {
    {X _0}
    \ar@{^{(}->}[r]
    \ar[d]
    &
    {X}
    \ar@{.>}[d] _{\delta ^n  _{f ,f'}}
    \ar[dr]^{(f',f)}
    \\
    {Y}
    \ar@{^{(}->}[r]
    &
    {\Delta ^n _{Y/T}}
    \ar[r] ^-{s _n}
    &
    {Y \times _T Y}
    \ar@<0,5ex>[r]^(0.6){p _1}
    \ar@<-0,5ex>[r]_(0.6){p _2}
    &
    {Y ,}
    }
\end{equation}
où $p _1$ et $p _2$ désignent respectivement les projections à gauche et à droite et $s _n$ est le morphisme canonique.
La structure de $\widetilde{\D} ^{(m)} _{Y/T}$-module de $\B _Y ^{(n_m)}( H _Y)$
compatible à sa structure d'algèbre donne
l'isomorphisme de $\O _{\Delta ^n _{Y/T}}$-algèbres
$(p _2 \circ s _n)^* (\B _Y ^{(n_m)}( H _Y)) \riso (p _1 \circ s _n) ^* (\B _Y ^{(n_m)}( H _Y))$.
Via $\smash{\delta ^n  _{f ,f'}} ^*$, on obtient
$f ^{\prime *} (\B _Y ^{(n_m)}( H _Y)) \riso f ^*(\B _Y ^{(n_m)}( H _Y))$
s'inscrivant dans le diagramme commutatif
\begin{equation}
\label{ff'b}
  \xymatrix  @R=0,3cm {
{f ^{\prime *} (\B _Y ^{(n_m)}( H _Y))}
\ar[rr] _{\sim}
\ar[rd]
&&
{f ^{*} (\B _Y ^{(n_m)}( H _Y))}
\ar[dl]
\\
&
{\B _X ^{(n_m)}( H _X),}
}
\end{equation}
dont les flèches obliques sont les morphismes canoniques (\cite[2.1.2]{huyghe2}).

Soit
$\widetilde{\Delta} ^n _{X/S(m)}$ l'espace annelé
$ (\Delta ^n _{X/S(m)}, \O _{\widetilde{\Delta} ^n _{X/S(m)}})$,
avec $\O _{\widetilde{\Delta} ^n _{X/S(m)}} := (p _1 \circ s _n)^* (\B _X ^{(n_m)}( H _X))$.
On note
$\tilde{p} _1$ : $\widetilde{\Delta} ^n _{X/S(m)} \rightarrow \widetilde{X}$
le morphisme induit par l'application continue $p _1 \circ s _n$ et par le morphisme d'anneaux canonique
$(p _1 \circ s _n) ^{-1}(\B _X ^{(n_m)}( H _X)) \rightarrow \O _{\widetilde{\Delta} ^n _{X/S(m)}}$
et $\tilde{p} _2$ : $\widetilde{\Delta} ^n _{X/S(m)} \rightarrow \widetilde{X}$
dont le morphisme d'espaces topologiques est $p _2 \circ s _n$ et dont le morphisme d'anneaux est
$(p _2 \circ s _n)^{-1} (\B _X ^{(n_m)}( H _X))
\rightarrow (p _2 \circ s _n)^* (\B _X ^{(n_m)}( H _X))\riso  (p _1 \circ s _n)^* (\B _X ^{(n_m)}( H _X))
= \O _{\widetilde{\Delta} ^n _{X/S(m)}}$.
De même, lorsque l'on remplace respectivement $S$ et $X$ par $T$ et $Y$ ou $U$ et $Z$.

On désigne par $\delta ^n  _{\smash{\tilde{f}} ,\tilde{f}'}$ :
$\widetilde{X} \rightarrow \widetilde{\Delta} ^n _{Y/T(m)}$
le morphisme induit par $\smash{\delta ^n  _{f ,f'}}$ et
par $\smash{\delta ^n  _{f ,f'}} ^{-1} (\O _{\widetilde{\Delta} ^n _{Y/T(m)}})
\rightarrow \smash{\delta ^n  _{f ,f'}} ^* (\O _{\widetilde{\Delta} ^n _{Y/T(m)}})
=\smash{\delta ^n  _{f ,f'}} ^* (p _1 \circ s _n)^* (\B _Y ^{(n_m)}( H _Y))
 \riso f ^* (\B _Y ^{(n_m)}( H _Y))
 \rightarrow \B _X ^{(n_m)}( H _X)$.
Via \ref{ff'b}, on vérifie $ \tilde{f} = \smash{\tilde{p}} _1 \delta ^n  _{\tilde{f},\smash{\tilde{f}} '} $
et $ \smash{\tilde{f}} '= \smash{\tilde{p}} _2 \delta ^n  _{\tilde{f},\smash{\tilde{f}} '} $.

Pour tout complexe $\FF$ de $\B _Y ^{(n_m)}( H _Y) \otimes _{\O _Y } \smash{\D} ^{(m)} _{Y/T}$-modules borné supérieurement,
via le diagramme commutatif suivant valable pour $n$ assez grand
\begin{equation}
  \label{deltantilde}
  \xymatrix  @R=0,3cm {
    { \widetilde{X} _0}
    \ar[d]
    \ar[r]
    &
    {\widetilde{X}}
    \ar@<0,5ex>[dr]^{\widetilde{f}}
    \ar@<-0,5ex>[dr]_{\widetilde{f}'}
    \ar[d] _{\delta ^n  _{\tilde{f},\smash{\tilde{f}} '}}
    \\
    {\widetilde{Y} }
    \ar[r]
    &
    {\widetilde{\Delta} ^n _{Y/T(m)}}
    \ar@<0,5ex>[r]^(0.6){\tilde{p} _1}
    \ar@<-0,5ex>[r]_(0.6){\tilde{p} _2}
    &
    {\widetilde{Y} ,}
    }
\end{equation}
  l'isomorphisme $\smash{\widetilde{\PP}} ^n _{Y/ T(m)}$-linéaire
  $\smash{\widetilde{\PP}} ^n _{Y/ T(m)} \otimes _{\B _Y ^{(n_m)}( H _Y)} \FF
  \riso
  \FF \otimes _{\B _Y ^{(n_m)}( H _Y) } \smash{\widetilde{\PP}} ^n _{Y/ T(m)}$
  induit un isomorphisme
  $\widetilde{\D} ^{(m)} _{X/S}$-linéaire canonique fonctoriel en $\FF$
  $$\tau _{\tilde{f},\tilde{f}'} ^{\FF} \ :\ \tilde{f} ^{'!} (\FF) \riso \tilde{f} ^! (\FF),$$
  tel que $\tau ^{\FF} _{\tilde{f},\tilde{f}} = \mathrm{Id}$, et que, si
  $\tilde{f} ''$ : $X \rightarrow Y $ est un morphisme
  dont la restriction $X _0 \rightarrow Y$ coïncide avec celle de $f$, on ait
  la formule de transitivité
  $\tau ^{\FF} _{\tilde{f}, \tilde{f}''} = \tau ^{\FF} _{\tilde{f},\tilde{f}'}\circ
  \tau ^{\FF} _{\tilde{f}',\tilde{f}''}$.
  Si aucune confusion n'est à craindre, nous
  noterons simplement
  $\tau _{\tilde{f},\tilde{f}'} $
  ou $\tau ^{\FF} $ voire $\tau$.

\begin{vide}
[Changement de niveau]
  Soient $n _{m'} \geq m '\geq m$ deux entiers,
  $\widetilde{\D} ^{(m')} _{Y/T} = \B _Y ^{(n _{m'})}( H _Y) \otimes _{\O _Y } \smash{\D} ^{(m')} _{Y/T}$ et
  $\FF$ un complexe de $\widetilde{\D} ^{(m')} _{Y/T} $-modules borné supérieurement.
Lorsque $Y$ a des coordonnées locales $t _1, \dots t _d$,
avec les notations de \cite[2.1.2]{Be1},
le morphisme canonique $\PP ^n _{X/ S(m')} \rightarrow \PP ^n _{X/ S(m)}$
envoie $\underline{\tau} ^{\{\underline{k}\} _{(m')}}$ sur
$\underline{q} !/\underline{q} ' ! \underline{\tau} ^{\{\underline{k}\} _{(m)}}$
tandis que le morphisme $\smash{\D} ^{(m)} _{Y/T}\rightarrow \smash{\D} ^{(m')} _{Y/T}$
envoie
$\underline{\partial} ^{<\underline{k}> _{(m)}}$
sur
$\underline{q} !/\underline{q} ' ! \underline{\partial} ^{<\underline{k}> _{(m')}}$.
Grâce à la formule \cite[1.1.16.1]{caro_comparaison}, il en découle,
en notant $oub _{m ,m'}$ : $ D (\widetilde{\D} ^{(m')} _{Y/T} )
\rightarrow D (\widetilde{\D} ^{(m)} _{Y/T} )$
 le foncteur oubli, la formule :
$$\tau ^{oub _{m ,m'} (\FF )} = oub _{m ,m'} \tau ^{\FF}.$$

\end{vide}

\begin{vide}
[Compatibilité aux images inverses extraordinaires]
Soit $g '$ : $Y \rightarrow Z$, un second $U$-morphisme
  ayant la même restriction $Y _0 \rightarrow Z$ que $g$.
Pour $n$ assez grand,
il découle respectivement du diagramme commutatif de droite et de gauche suivant
  \begin{equation}
\label{transit-crist1}
    \xymatrix  @R=0,3cm {
    { \widetilde{X} _0}
    \ar[d] ^-{\tilde{f} _0}
    \ar[r]
    &
    {\widetilde{X}}
    \ar[d]^{\widetilde{f}}
    \\
    {\widetilde{Y} _0}
    \ar[r]
    \ar[d]
    &
    {\widetilde{Y}}
    \ar[d] _{\delta ^n  _{\tilde{g},\smash{\tilde{g}} '}}
    \ar@<0,5ex>[dr]^{\widetilde{g}}
    \ar@<-0,5ex>[dr]_{\widetilde{g}'}
    \\
    {\widetilde{Z}}
    \ar[r]
    &
    {\widetilde{\Delta} ^n _{Z/U(m)}}
    \ar@<0,5ex>[r]^(0.6){\tilde{p} _1}
    \ar@<-0,5ex>[r]_(0.6){\tilde{p} _2}
    &
    {\widetilde{Z} ,}
    }
    \hspace{2cm}
    \xymatrix  @R=0,3cm {
    { \widetilde{X} _0}
    \ar[d]
    \ar[r]
    &
    {\widetilde{X}}
    \ar@<0,5ex>[dr]^{\widetilde{f}}
    \ar@<-0,5ex>[dr]_{\widetilde{f}'}
    \ar[d] _{\delta ^n  _{\tilde{f},\smash{\tilde{f}} '}}
    \\
    {\widetilde{Y} }
    \ar[r]
    \ar[d] ^-{\widetilde{g}}
    &
    {\widetilde{\Delta} ^n _{Y/T(m)}}
    \ar@{.>}[d]
    \ar@<0,5ex>[r]^(0.6){\tilde{p} _1}
    \ar@<-0,5ex>[r]_(0.6){\tilde{p} _2}
    &
    {\widetilde{Y}       }
    \ar[d]^{\widetilde{g}}
    \\
    {\widetilde{Z}}
    \ar[r]
    &
    {\widetilde{\Delta} ^n _{Z/U(m)},}
    \ar@<0,5ex>[r]^(0.6){\tilde{p} _1}
    \ar@<-0,5ex>[r]_(0.6){\tilde{p} _2}
    &
    {\widetilde{Z} }
    }
  \end{equation}
que l'on dispose de l'égalité
$ \tau _{\tilde{g} \circ \tilde{f} , \tilde{g} \circ \tilde{f}'} =
\tau _{\tilde{f},\tilde{f}'} \circ \tilde{g} ^!$
et $ \tilde{f} ^! \circ \tau _{\tilde{g} ,\tilde{g}'} =
\tau _{\tilde{g} \circ \tilde{f}, \tilde{g}'\circ \tilde{f}}$.

Ainsi, les isomorphismes $\tau $ sont compatibles avec la composition des images inverses.

\end{vide}

\begin{vide}
[Changement de base]
Considérons le diagramme commutatif
\begin{equation}
  \notag
  \xymatrix  @R=0,3cm {
  {S'} \ar[r] \ar[d] & {T'} \ar[d]\\
  {S} \ar[r] & {T}}
\end{equation}
de $m$-PD-morphismes. On note $X ' = X \times _S S'$,
$Y '= Y \times _T T'$, $g$ et $g'$ : $X' \rightarrow Y'$ les morphismes induits respectivement par $f$ et $f'$.
Il découle de la propriété universelle des enveloppes à puissances divisées partielles
(\cite[1.4.1]{Be1}) que l'on dispose d'un morphisme
$\Delta ^n _{Y'/T'(m)} \rightarrow \Delta ^n _{Y/T(m)} $ rendant commutatif le diagramme
\begin{equation}
  \label{tauchgtbase1}
  \xymatrix  @R=0,3cm {
  {X'}
  \ar[r] ^-{\delta ^n  _{g ,g'}} \ar[d]
  &
  {\Delta ^n _{Y'/T'}}
  \ar[r] ^-{s _n}
  \ar@{.>}[d]
  &
  {Y' \times _{T'} Y'}
  \ar[d]
    \ar@<0,5ex>[r]^(0.6){p _1}
    \ar@<-0,5ex>[r]_(0.6){p _2}
    &
    {Y '}
    \ar[d]
    \\
    {X} \ar[r] ^-{\delta ^n  _{f ,f'}}
  &
  {\Delta ^n _{Y/T}}
  \ar[r] ^-{s _n}
  &
  {Y \times _{T} Y}
    \ar@<0,5ex>[r]^(0.6){p _1}
    \ar@<-0,5ex>[r]_(0.6){p _2}
    &
    {Y }
    }
\end{equation}
En notant respectivement $H _{X'}$ et $H _{Y'}$ les images inverses de $H _X$ et $H _Y$,
$\smash{\widetilde{X}}' = (X', \B _X ^{(n_m)}( H _{X'}))$
et $\smash{\widetilde{Y}}' = (Y', \B _Y ^{(n_m)}( H _{Y'}))$,
les morphismes $g$ et $g'$ induisent canoniquement des morphismes
$\smash{\widetilde{X}}' \rightarrow \smash{\widetilde{Y}}'$ notés $\tilde{g}$ et
$\smash{\tilde{g}} '$.
On déduit ensuite de la commutativité du diagramme \ref{tauchgtbase1} celle de
\begin{equation}
    \label{tauchgtbase2}
  \xymatrix  @R=0,3cm {
  {\smash{\widetilde{X}}'}
  \ar[r] ^-{\delta ^n  _{\tilde{g},\smash{\widetilde{g}} '}}
  \ar[d]
  &
  {\smash{\widetilde{\Delta}} ^n _{Y'/T'}}
  \ar@{.>}[d]
    \ar@<0,5ex>[r]^(0.6){\smash{\widetilde{p}} _1}
    \ar@<-0,5ex>[r]_(0.6){\smash{\widetilde{p}} _2}
    &
    {\smash{\widetilde{Y}} '}
    \ar[d]
    \\
    {\widetilde{X}} \ar[r] ^-{\delta ^n  _{\tilde{f},\smash{\tilde{f}} '}}
  &
  {\smash{\widetilde{\Delta}}  ^n _{Y/T}}
    \ar@<0,5ex>[r]^(0.6){\smash{\widetilde{p}} _1}
    \ar@<-0,5ex>[r]_(0.6){\smash{\widetilde{p}} _2}
    &
    {\widetilde{Y} .}
    }
\end{equation}
Il en résulte, notant $\theta $ : $\smash{\widetilde{X}} '
\rightarrow \smash{\widetilde{X}} $ le morphisme canonique,
$\tau _{\tilde{f},\tilde{f}'} \riso \tau _{\tilde{g},\tilde{g}'}$.
Autrement dit, les isomorphismes de la forme $\tau$ sont compatibles aux changements de base.

\end{vide}

  \end{vide}

\begin{rema}\label{tauquasinilp}
  Avec les notations de \ref{relev-comp-adj-immf0},
  soit $s$ le morphisme canonique $\Delta  _{Y/T} \rightarrow Y \times _T Y$.
On obtient
$\widetilde{\Delta} _{Y/T(m)}$ l'espace annelé
$ (\Delta  _{Y/T(m)}, \O _{\widetilde{\Delta}  _{Y/T(m)}})$ en posant
$\O _{\widetilde{\Delta}  _{Y/T(m)}} := (p _1 \circ s )^* (\B _Y ^{(n_m)}( H _Y))$.
De même, on construit $\delta   _{\smash{\tilde{f}} ,\tilde{f}'}$ :
$\widetilde{X} \rightarrow \widetilde{\Delta}  _{Y/T(m)}$, celui-ci s'insérant dans
le diagramme commutatif \ref{deltantilde} où on a enlevé les indices $n$.
Si $p$ est localement nilpotent sur $U$ (et donc sur $T$ et $S$),
pour tout complexe quasi-nilpotent (voir \ref{lemmquasi-nilp} et \ref{espilonhat})
$\FF$ de $\B _Y ^{(n_m)}( H _Y) \otimes _{\O _Y } \smash{\D} ^{(m)} _{Y/T}$-modules borné supérieurement,
on en déduit, de manière analogue à ceux de \ref{relev-comp-adj-immf0},
des isomorphismes $\tau ^{\FF} _{\tilde{f},\tilde{f}'}$, ceux-ci vérifiant les conditions de transitivité,
de changement de base et de commutation à la composition des images inverses.
\end{rema}

  \begin{vide}
  Avec les notations et hypothèses de \ref{relev-comp-adj-immf0} (nous omettons d'indiquer la base),
  on pose
$$\widetilde{\D} ^{(m)} _{\smash{X\overset{{}_f}{\rightarrow}Y}}=
\tilde{f} ^*  ( \widetilde{\D} ^{(m)} _Y )\text{ et }
\widetilde{\D} ^{(m)} _{Y\overset{{}_f}{\leftarrow}X}=
\tilde{f} ^* _d ( \widetilde{\D} ^{(m)} _Y \otimes _{\O _Y} \omega _Y ^{-1}) \otimes _{\O _X} \omega _X
\riso
\omega _X \otimes _{\O _X}  \tilde{f} ^* _g ( \widetilde{\D} ^{(m)} _Y \otimes _{\O _Y} \omega _Y ^{-1}).$$
Supposons $S=T$, $S$ noethérien et de dimension de Krull finie, $f$ et $f'$
quasi-compacts et quasi-séparés. Les isomorphismes
$\tau _{\tilde{f},\tilde{f}'} $
: $ \widetilde{\D} ^{(m)} _{X\overset{{}_{f'}}{\rightarrow}Y}
\riso
\widetilde{\D} ^{(m)} _{X\overset{{}_f}{\rightarrow}Y}$
et
$\tau _{\tilde{f},\tilde{f}'} $
: $ \widetilde{\D} ^{(m)} _{Y\overset{{}_{f'}}{\leftarrow}X}
\riso
\widetilde{\D} ^{(m)} _{Y\overset{{}_f}{\leftarrow}X}$ fournissent,
pour tout complexe de $\widetilde{\D} ^{(m)} _X$-modules à gauche (resp. à droite) borné inférieurement $\E$
(resp. $\M$),
les isomorphismes
$$ \sigma _{\tilde{f},\tilde{f}'} \ : \  f '_{H _Y,\,H _X +} (\E ) \riso f _{H _Y,\,H _X +} (\E )
\text{ et }
\sigma _{\tilde{f},\tilde{f}'} \ : \  f '_{H _Y,\,H _X +} (\M ) \riso f _{H _Y,\,H _X +} (\M )).$$
De plus, pour tout morphisme $\tilde{f} ''$ : $X \rightarrow Y $
  dont la restriction $X _0 \rightarrow Y$ coïncide avec celle de $f$, on a
  la formule de transitivité
  $\sigma _{\tilde{f}, \tilde{f}''} = \sigma _{\tilde{f},\tilde{f}'}\circ \sigma _{\tilde{f}',\tilde{f}''}$.
\end{vide}

\begin{rema}\label{rematautau}
  Soient $\FF \in D ^- (\widetilde{\D} ^{(m)} _Y)$,
  $\tilde{p} _{1 *} \smash{\widetilde{\PP}} ^n _{Y(m)}$ (resp.
  $\tilde{p} _{2 *} \smash{\widetilde{\PP}} ^n _{Y(m)}$) le faisceau
  $\smash{\widetilde{\PP}} ^n _{Y(m)}$ vu comme
  $\B _Y ^{(m)}( H _Y)$-algèbre pour la structure gauche (resp. droite).
  Considérons les diagrammes
\begin{equation}
  \notag
  \xymatrix  @R=0,3cm {
  {(\tilde{p} _{2 *} \smash{\widetilde{\PP}} ^n _{Y(m)} \otimes _{\B _Y ^{(m)}( H _Y)} \widetilde{\D} ^{(m)} _Y)
  \otimes _{\widetilde{\D} ^{(m)} _Y} \FF }
  \ar[r] _(0.59){\sim}
  \ar[d] ^-{\smash{\tilde{\epsilon}} _n ^{\widetilde{\D} ^{(m)} _Y} \otimes id}_{\sim}
  &
  {\tilde{p} _{2 *} \smash{\widetilde{\PP}} ^n _{Y(m)} \otimes _{\B _Y ^{(m)}( H _Y)} \FF }
  \ar[d] ^-{\smash{\tilde{\epsilon}} _n ^{\FF}} _{\sim}
  \\
  {(\tilde{p} _{1 *} \smash{\widetilde{\PP}} ^n _{Y(m)} \otimes _{\B _Y ^{(m)}( H _Y)} \widetilde{\D} ^{(m)} _Y)
  \otimes _{\widetilde{\D} ^{(m)} _Y} \FF }
  \ar[r]_(0.59){\sim}
  &
  {\tilde{p} _{1 *} \smash{\widetilde{\PP}} ^n _{Y(m)} \otimes _{\B _Y ^{(m)}( H _Y)} \FF ,}
  }
  \xymatrix  @R=0,3cm {
  {\widetilde{\D} ^{(m)} _{X\overset{{}_{f'}}{\rightarrow}Y}
  \otimes ^{\L} _{f ^{-1} \widetilde{\D} ^{(m)} _Y} f ^{-1} \FF [d_{X/Y}]}
  \ar[r]_(0.7){\sim}
  \ar[d] ^-{\tau _{\tilde{f},\tilde{f}'} \otimes id} _{\sim}
  &
  {\tilde{f} ^{'!} (\FF)}
  \ar[d] ^-{\tau _{\tilde{f},\tilde{f}'}} _{\sim}
  \\
  {\widetilde{\D} ^{(m)} _{X\overset{{}_f}{\rightarrow}Y}
  \otimes ^{\L} _{f ^{-1} \widetilde{\D} ^{(m)} _Y }f ^{-1} \FF[d_{X/Y}]}
  \ar[r] _(0.7){\sim}
  &
  {\tilde{f} ^{!} (\FF),}
}
\end{equation}
où $\tilde{\epsilon} _n$ désignent les $m$-PD-stratifications relatives à $\B _Y ^{(m)}( H _Y)$
(voir \cite[1.1]{caro_comparaison}).
Via un calcul utilisant la formule \cite[1.1.16.1]{caro_comparaison}, on obtient la commutativité
du diagramme de gauche. Il en résulte celle de droite.
\end{rema}

\begin{lemm}\label{lemm-comp+}
  Avec les notations de \ref{relev-comp-adj-immf0}, les isomorphismes
  \begin{gather}\notag
    \smash{\widetilde{\D}} ^{(m)} _{X \rightarrow Z} \riso
    \smash{\widetilde{\D}} ^{(m)} _{X \rightarrow Y} \otimes _{f ^{-1}\smash{\widetilde{\D}} ^{(m)} _{Y}}
    f ^{-1} \smash{\widetilde{\D}} ^{(m)} _{Y \rightarrow Z}
    \riso
    \smash{\widetilde{\D}} ^{(m)} _{X \rightarrow Y} \otimes ^\L _{f ^{-1}\smash{\widetilde{\D}} ^{(m)} _{Y}}
    f ^{-1} \smash{\widetilde{\D}} ^{(m)} _{Y \rightarrow Z}\\
    \notag
    \smash{\widetilde{\D}} ^{(m)} _{Z \leftarrow X} \riso
   f ^{-1} \smash{\widetilde{\D}} ^{(m)} _{Z \leftarrow Y}
    \otimes _{f ^{-1}\smash{\widetilde{\D}} ^{(m)} _{Y}}
    \smash{\widetilde{\D}} ^{(m)} _{Y \leftarrow X}
    \riso
       f ^{-1} \smash{\widetilde{\D}} ^{(m)} _{Z \leftarrow Y}
    \otimes ^\L _{f ^{-1}\smash{\widetilde{\D}} ^{(m)} _{Y}}
    \smash{\widetilde{\D}} ^{(m)} _{Y \leftarrow X}
  \end{gather}
  sont compatibles aux isomorphismes de la forme $\tau$.
\end{lemm}
\begin{proof}
  La construction des isomorphismes est analogue à \cite[VI.6.3]{borel}.
Pour la première ligne, la compatibilité aux isomorphismes $\tau$ est une conséquence de la remarque
\ref{rematautau} et des
formules
$ \tau _{\tilde{g} \circ \tilde{f} , \tilde{g} \circ \tilde{f}'} =
\tau _{\tilde{f},\tilde{f}'} \circ \tilde{g} ^!$
et
$\tilde{f} ^! \circ \tau _{\tilde{g} ,\tilde{g}'} =
\tau _{\tilde{g} \circ \tilde{f}, \tilde{g}'\circ \tilde{f}}$
prouvées via les diagrammes \ref{transit-crist1}. Ceux de la deuxième ligne s'en déduisent.
\end{proof}
\begin{prop}\label{prop-comp+}
  Avec les notations de \ref{relev-comp-adj-immf0}, on suppose $S=T=U$, $S$ noethérien, de
  dimension de Krull finie, $f$, $f'$, $g$ et $g'$ quasi-compacts et quasi-séparés.
  Les isomorphismes de la forme
  $\sigma $ sont alors compatibles à la composition des images directes, i.e.,
  $g _{H _Z,\,H _Y +} \circ \sigma _{f,f'} = \sigma _{g\circ f,g\circ f'}$ et
  $\sigma _{g,g'} \circ f _{H _Y,\,H _X +} =\sigma _{g\circ f,g'\circ f} $.
\end{prop}
\begin{proof}
  La commutation à la composition des foncteurs images directes se prouve de manière analogue à \cite[VI.6.4]{borel}.
  Par construction de celle-ci, cela découle du lemme \ref{lemm-comp+}.
\end{proof}
  \begin{vide}\label{relev-comp-adj-immf}
Maintenant, soient $f $ : $\X \rightarrow \Y$ un morphisme de $\V$-schémas formels lisses,
$T _Y$ un diviseur de $Y$, $T _X \supset f _0 ^{-1} (T _Y)$ un diviseur de $X$.
On se donne un deuxième morphisme $f'$ : $\X \rightarrow \Y$
tel que $f ' _0 = f _0$.
Donnons-nous une suite croissante d'entiers $(n _m) _{m \in \N}$ telle que $n _m \geq m$
et posons $\smash{\widehat{\D}} _{\X} ^{(m)}(T _X) := \B _{\X} ^{(n_m)}( T _X)
\smash{\widehat{\otimes}} _{\O _{\X} } \smash{\D} ^{(m)} _{\X}$. De même, en remplaçant $\X$ par $\Y$.
Comme pour \cite[3.2.1]{Beintro2}, on définit la catégorie des complexes de
$\smash{\widehat{\D}} _{\X} ^{(m)}(T _X)$-modules à gauche à cohomologie bornée et quasi-cohérente et qui
sera noté
$D ^\mathrm{b} _\mathrm{qc} (\smash{\widehat{\D}} _{\X} ^{(m)}(T _X) )$.
On dispose
des foncteurs images inverses extraordinaires de niveau $m$
(définis comme dans \cite[3.4.2.1]{Beintro2}) :
$f ^{!(m)} _{T_X, T _Y}$ et $f ^{'!(m)} _{T_X, T _Y}$ :
$D ^\mathrm{b} _\mathrm{qc} (\smash{\widehat{\D}} _{\Y} ^{(m)}(T _Y) )
\rightarrow
D ^\mathrm{b} _\mathrm{qc} (\smash{\widehat{\D}} _{\X} ^{(m)}(T _X) )$.
Les isomorphismes de la forme $\tau$ (voir \ref{relev-comp-adj-immf0}) commutant aux changements de base,
on obtient, pour tout
$\FF \in D ^\mathrm{b} _\mathrm{qc} (\smash{\widehat{\D}} _{\Y} ^{(m)}(T _Y) )$,
des isomorphismes
$\tau ^{ \FF (m)} _{f  ,f',T_X, T _Y}$ :
$f ^{'!(m)} _{T_X, T _Y} (\FF)\riso f ^{!(m)} _{T_X, T _Y} (\FF) $.
En procédant à deux localisations
(correspondant à la tensorisation par $\Q$ et au passage à la limite sur le niveau),
on obtient de manière analogue à \cite[4.2.1]{Beintro2} les catégories de la forme
$\smash{\underset{^{\longrightarrow}}{LD}} ^{\mathrm{b}} _{\Q ,\mathrm{qc}}
( \smash{\widehat{\D}} _{\Y} ^{(\bullet)}(T _Y))$ (ces dernières ne dépendent pas du choix de
la suite d'entiers $(n _m) _{m\in \N}$).
Il découle des isomorphismes $\tau ^{ \FF (m)} _{f  ,f',T_X, T _Y}$ avec $m$ variable,
pour tout objet $\FF$ de
$\smash{\underset{^{\longrightarrow}}{LD}} ^{\mathrm{b}} _{\Q ,\mathrm{qc}}
( \smash{\widehat{\D}} _{\Y} ^{(\bullet)}(T _Y))$,
le suivant
$\tau _{f  ,f',T_X, T _Y} ^{\FF} \ :\ f ^{'!} _{T_X, T _Y}(\FF) \riso f ^! _{T_X, T _Y}(\FF),$
où $f ^! _{T_X, T _Y}$ et $f ^{'!} _{T_X, T _Y}$ sont
les foncteurs images inverses extraordinaires
définis comme pour \cite[4.3.2.1]{Beintro2}.
On pourra le noter plus simplement $\tau _{f  ,f'}$ ou $\tau ^{\FF}$ voire $\tau$.
De la même façon, pour tout
$\E \in \smash{\underset{^{\longrightarrow}}{LD}} ^{\mathrm{b}} _{\Q ,\mathrm{qc}}
( \smash{\widehat{\D}} _{\X} ^{(\bullet)}(T _X))$,
on construit l'isomorphisme
$f ' _{T_Y, T _X +}(\E) \riso f  _{T_Y, T _X +}(\E)$,
que l'on notera $\sigma _{f  ,f',T_X, T _Y}$
(ou $\sigma _{f,f'}$ voire $\sigma$ si aucune confusion n'est à craindre).
Les isomorphismes de la forme $\tau $ (resp. $\sigma$) vérifient aussi
la formule de transitivité et celle de leur commutation à l'image inverse extraordinaire (resp. image directe)
décrites ci-dessus dans le cas des schémas.

\end{vide}

\begin{prop}
  \label{tauhdag}
  Avec les notations de \ref{relev-comp-adj-immf}, pour tout objet $\FF$ de
$\smash{\underset{^{\longrightarrow}}{LD}} ^{\mathrm{b}} _{\Q ,\mathrm{qc}}
( \smash{\widehat{\D}} _{\Y} ^{(\bullet)})$,
  on a les diagrammes commutatifs :
  $$\xymatrix  @R=0,3cm {
  {f ^! _{T_X, T _Y}(\hdag T _Y) (\FF) }
  \ar[r] _{\sim}
  &
  {(\hdag T _X) f ^!  (\FF) }
  \\
  {f ^{\prime !} _{T_X, T _Y} (\hdag T _Y) (\FF) }
  \ar[r] _{\sim}
  \ar[u] ^-{\tau _{f,f', T_X, T _Y} ^{ (\hdag T _Y) (\FF)}} _{\sim}
  &
  {(\hdag T _X) f ^{\prime !}  (\FF), }
  \ar[u] ^-{(\hdag T _X) (\tau _{f,f'} ^{\FF})} _{\sim}
  }
  \
  \
  \
  \xymatrix  @R=0,3cm {
  {f ^! (\hdag T _Y) (\FF) }
  \ar[r] _{\sim}
  &
  {(\hdag T _X) f ^!  (\FF) }
  \\
  {f ^{\prime !}  (\hdag T _Y) (\FF) }
  \ar[r] _{\sim}
  \ar[u] ^-{\tau _{f,f' } ^{(\hdag T _Y) (\FF)}} _{\sim}
  &
  {(\hdag T _X) f ^{\prime !}  (\FF). }
  \ar[u] ^-{(\hdag T _X) (\tau _{f,f'} ^{\FF})} _{\sim}
  }$$
\end{prop}
\begin{proof}
Notons $f _i$ et $f ' _i $ : $X _i \rightarrow Y_i$ les réductions modulo $\mathfrak{m} ^{i+1}$ de $f$ et $f'$, et
$\smash{\tilde{f}} _i$ et $\smash{\tilde{f}} '_i$ les morphismes d'espaces annelés canoniques
$(X _i , \B _{X _i}^{(n_m)}( T _X))
\rightarrow (Y _i,\B _{Y _i}^{(n_m)}( T _Y))$ induits par $f$ et $ f'$.
Pour tout $\FF _i \in D ^- (\B _{Y _i}^{(n_m)}( T _Y) \otimes _{\O _{Y_i}} \smash{\D} ^{(m)} _{Y _i})$, le diagramme
$$\xymatrix  @R=0,3cm {
{\smash{\tilde{f}} _i ^! (\B _{Y _i}^{(n_m)}( T _Y) \otimes ^\L _{\O _{Y_i}} \FF _i)}
\ar[r]
&
{\B _{X _i}^{(n_m)}( T _X) \otimes ^\L _{\O _{X_i}} f _i ^! (\FF _i)}
\\
{\smash{\tilde{f}} _i ^{\prime !} (\B _{Y _i}^{(n_m)}( T _Y) \otimes ^\L _{\O _{Y_i}} \FF _i)}
\ar[r]
\ar[u] ^-{\tau _{\smash{\tilde{f}} _i ,\smash{\tilde{f}}'_i}}
&
{\B _{X _i}^{(n_m)}( T _X) \otimes ^\L _{\O _{X_i}} f _i ^{\prime !} (\FF _i)}
\ar[u] ^-{id \otimes \tau _{f _i,f' _i}}
}$$
est commutatif. Il en résulte celle du diagramme de gauche de \ref{tauhdag}.

De plus, comme chacune des flèches du morphisme composé
$f ^! (\hdag T _Y) (\FF) \riso (\hdag T _X) f ^! (\hdag T _Y) (\FF) \tilde{\leftarrow} (\hdag T _X) f ^!  (\FF) $,
qui correspond à la flèche du haut du diagramme de droite de \ref{tauhdag},
commutent par fonctorialité  aux isomorphismes $\tau$, on obtient la commutativité de celui de droite.
\end{proof}

\begin{rema}
  Avec les notations de \ref{tauhdag}, supposons $f ^{-1} (T _Y) = T _X$.
  L'isomorphisme $f ^! _{T_Y} \riso f ^!$ (on a omis d'indiquer
  les foncteurs oublis) commute alors aux isomorphismes de recollements $\tau$. En effet,
  celui-ci est construit de la façon suivante :

  $f ^! _{T _Y}  ((\hdag T _Y) \FF)  \riso
  f ^! _{T _Y}(\hdag T _Y)  ((\hdag T _Y) \FF) \riso (\hdag T _X) f ^! ((\hdag T _Y) \FF)
  \riso  f ^! (\hdag T _Y) ((\hdag T _Y) \FF) \tilde{\leftarrow} f ^! ((\hdag T _Y) \FF)$.
\end{rema}

\subsection{Morphisme d'adjonction associé à un carré}
\begin{vide}\label{not-comp-comp-adj-immf}
  Soit le diagramme de $\V$-schémas formels lisses :
\begin{equation}\label{deuxcarresadj}
  \xymatrix  @R=0,3cm {
  {\PP ''} \ar[r] ^g
  &
  {\PP '} \ar[r] ^f
  &
  {\PP}
  \\
  {\X ''} \ar[u]  ^{u''}\ar[r]^b
  &
  {\X '} \ar[u] ^-{u'} \ar[r]^a
  &
  {\X ,} \ar[u] ^u
  }
\end{equation}
où $f$, $g$, $a$ et $b$ sont lisses et où $u$, $u'$ et $u''$ sont des immersions fermées.
On suppose que le diagramme \ref{deuxcarresadj} est commutatif au niveau des fibres spéciales
(mais non nécessairement commutatif au niveau des schémas formels).
De plus, on se donne $T _P$ un diviseur de $P$ tel que $T _{P'} := f ^{-1} (T _{P})$
(resp. $T _{P''} := g ^{-1} (T _{P'})$, $T _{X} := u ^{-1} (T _{P})$, $T _{X'} := u ^{\prime -1} (T _{P'})$
et $T _{X''} := u ^{\prime \prime -1} (T _{P''})$) soit un diviseur de
$P'$ (resp. $P''$, $X$, $X'$ et $X''$).

\end{vide}

\begin{prop}\label{comp-comp-adj-immf}
  Avec les notations de \ref{not-comp-comp-adj-immf},
  on dispose d'un morphisme, dit d'{\rm adjonction}, fonctoriel en
  $\E \in  D ^{\mathrm{b}} _{\mathrm{coh}} (\smash{\D} ^\dag _{\X } (\hdag T _{X}) _\Q)$,
$\phi  (\E)$ : $u ^\prime _+ \circ a ^! (\E) \rightarrow f ^{!}\circ u _+ (\E).$
Si la réduction au niveau des fibres spéciales
du carré de droite de \ref{deuxcarresadj} est cartésien, alors
$\phi  (\E)$ est un isomorphisme.
Le morphisme d'adjonction
entre foncteurs, $u ^\prime _+ \circ a ^! \rightarrow f ^{!}\circ u _+$,
sera noté $\phi$.

(i) En notant $\phi '\ : \ u ^{''} _+ \circ b ^! \rightarrow g ^! \circ u ^\prime _+$
(resp. $\phi ''\ : \ u ^{''} _+ \circ (a \circ b ) ^! \rightarrow (f \circ g) ^!\circ u  _+$)
le morphisme d'adjonction
du carré de gauche de \ref{deuxcarresadj} (resp. du grand rectangle de \ref{deuxcarresadj}),
le diagramme
$$\xymatrix  @R=0,3cm {
{ u ''  _+ \circ (a\circ b) ^!}
\ar[r] _\sim
\ar[d] ^-{\phi ''}
&
{ u ''  _+ \circ b ^!\circ  a ^! }
\ar[d]^{( g ^! \circ \phi)\circ (\phi ' \circ a ^!)}
\\
{(f\circ g) ^! \circ  u _+  }
\ar[r] _{\sim}
&
{ g ^!  \circ f ^!\circ  u _+,}
}
$$
est alors commutatif.
Ainsi, avec des abus de notations,
on dispose de la formule de transitivité
$\phi ''=( g ^! \circ \phi)\circ (\phi ' \circ a ^!)  $ (pour le composé de deux morphismes lisses)
de l'isomorphisme de changement de base par un morphisme lisse de l'image directe par une immersion fermée.

(ii) Soit $a '$ : $ \X ' \rightarrow \X$ (resp. $f'$ : $\PP ' \rightarrow \PP$)
un morphisme dont la réduction $X ' \rightarrow \X$ (resp. $P ' \rightarrow \PP$)
coïncide avec celle de $a$ (resp. $f$). Le diagramme ci-dessous
$$\xymatrix  @R=0,3cm {
{ u ^\prime_+ a ^{ !}}
\ar[r] ^-{\phi}
&
{f ^{!}\circ u _+}
\\
{ u' _+ a ^{\prime !} }
\ar[r]^{\psi}
\ar[u] ^-{ u ^\prime_+ (\tau _{a,a'})} _{\sim}
&
{f ^{\prime !}\circ u  _+,}
\ar[u] ^-{\tau _{f,f'} u _+} _{\sim}
}$$
où $\psi$ désigne le morphisme d'adjonction du diagramme de droite de \ref{deuxcarresadj} lorsque
$a$ et $f$ ont été remplacés respectivement par $a'$ et $f'$,
est commutatif.
\end{prop}
\begin{proof}
Construisons d'abord le morphisme $\phi  (\E)$. Comme
$\E \in  D ^{\mathrm{b}} _{\mathrm{coh}} (\smash{\D} ^\dag _{\X } (\hdag T _{X}) _\Q)$,
en appliquant le foncteur $u ^\prime _+  a ^!$ au morphisme d'adjonction de $u$ en $\E$ (\cite[1.2.11]{caro_courbe-nouveau}),
on obtient :
$u ^\prime _+  a ^! (\E)   \rightarrow u ^\prime _+  a ^! u ^{ !} u  _+(\E)$.
Or, comme $(f \circ u' ) ^! \riso u ^{\prime !} f ^!$ et
$(u\circ a ) ^! \riso  a ^! u ^! $,
on a
l'isomorphisme
$u ^\prime _+  \tau _{f \circ u', u\circ a} u  _+(\E)$ : $u ^\prime _+  a ^! u ^! u  _+(\E) \riso
u ^\prime _+   u ^{\prime !} f ^! u  _+(\E)$ (notations de \ref{relev-comp-adj-immf}).
Ensuite, puisque
$f ^!u  _+(\E)\in D ^{\mathrm{b}} _{\mathrm{coh}} (\smash{\D} ^\dag _{\PP ' } (\hdag T _{P '}) _\Q)$
et
$u ^{\prime !} f ^!u  _+(\E)\in D ^{\mathrm{b}} _{\mathrm{coh}} (\smash{\D} ^\dag _{\X ' } (\hdag T _{X '}) _\Q)$,
il résulte de
\cite[1.2.12]{caro_courbe-nouveau} que l'on dispose du morphisme d'adjonction de $u'$ en
$f ^!u  _+(\E)$ :
$u ^\prime _+   u ^{\prime !} f ^!u  _+(\E)
\rightarrow f ^!u  _+(\E)$.
En composant ces trois morphismes, il vient :
$\phi  (\E):=u ^\prime _+  a ^! (\E) \rightarrow f ^{!} u _+ (\E)$.

À présent, établissons que le morphisme $\phi  (\E)$ est un isomorphisme
lorsque le diagramme de droite de \ref{deuxcarresadj} est cartésien.
D'abord, comme $u$ est une immersion fermée, la première flèche dans la construction de $\phi (\E)$
est un isomorphisme.
En outre, l'hypothèse cartésienne implique
que le faisceau
$f ^!u  _+(\E)$
est à support dans $X'$ et donc que le morphisme d'adjonction de $u'$ en
$f ^!u  _+(\E)$ est un isomorphisme (voir \cite[1.2.12]{caro_courbe-nouveau}).

Prouvons à présent la formule de transitivité de (i). Pour cela, considérons le diagramme suivant :
\begin{equation}\label{diag1-comp-comp-adj-immf}
  \xymatrix  @R=0,3cm {
 { u ''  _+ (a\circ b) ^!} \ar[r]\ar[d]^{\mathrm{adj} _u}
 &
 { u ''  _+ b ^! a ^!}
 \ar[r] ^-{\mathrm{adj} _{u '}}
 \ar[d] ^-{\mathrm{adj} _u}
 &
  { u ''  _+ b^! u ^{\prime !} u' _+ a ^!} \ar[r] ^\tau \ar[d] ^-{\mathrm{adj} _u}
  &
  { u ''  _+ u ^{\prime \prime !} g ^!  u' _+ a ^!}
  \ar[r] ^-{\mathrm{adj} _{u ''}} \ar[d] ^-{\mathrm{adj} _u}
  &
  { g ^!  u' _+ a ^!}\ar[d] ^-{\mathrm{adj} _u}\\
{ u ''  _+ (a\circ b) ^! u ^! u _+}
\ar[r]\ar[ddd]^\tau
 &
 { u ''  _+ b ^! a ^! u ^! u _+}
 \ar[r] ^-{\mathrm{adj} _{u '}}
 \ar[d] ^\tau
 &
  { u ''  _+ b^! u ^{\prime !} u' _+ a ^! u ^! u _+}
  \ar[r]^\tau
  \ar[d] ^\tau
  &
  { u ''  _+ u ^{\prime \prime !} g ^!  u' _+ a ^! u ^! u _+}
  \ar[r] ^-{\mathrm{adj} _{u ''}}
  \ar[d] ^\tau
  &
  { g ^!  u' _+ a ^! u ^! u _+}\ar[d] ^\tau\\
    { }
    &
    { u ''  _+ b ^!  u ^{\prime !} f ^! u _+}
    \ar[r] ^-{\mathrm{adj} _{u '}}
    \ar@{=}[d]
 &
  { u ''  _+ b^! u ^{\prime !} u' _+  u ^{\prime !} f ^! u _+}
  \ar[r] ^\tau  \ar[d] ^-{\mathrm{adj} _{u '}}
  &
  { u ''  _+ u ^{\prime \prime !} g ^!  u' _+  u ^{\prime !} f ^! u _+}
  \ar[r] ^-{\mathrm{adj} _{u ''}} \ar[d]^{\mathrm{adj} _{u '}}
  &
  { g ^!  u' _+  u ^{\prime !} f ^! u _+}\ar[d]^{\mathrm{adj} _{u '}}\\
{ }
    &
    { u ''  _+ b ^!  u ^{\prime !} f ^! u _+} \ar@{=}[r]\ar[d]^\tau
 &
  { u ''  _+ b^!  u ^{\prime !} f ^! u _+} \ar[r]^\tau \ar[d]^\tau
  &
  { u ''  _+ u ^{\prime \prime !} g ^!   f ^! u _+}
  \ar[r] ^-{\mathrm{adj} _{u ''}} \ar@{=}[d]
  &
  { g ^!   f ^! u _+}\ar@{=}[d]\\
    { u ''  _+ u ^{\prime \prime !} (f\circ g) ^!  u _+ }\ar[d]^{\mathrm{adj} _{u ''}}\ar[r]
    &
    { u ''  _+  u ^{\prime \prime !} g ^! f ^! u _+} \ar@{=}[r]\ar[d]^{\mathrm{adj} _{u ''}}
 &
  { u ''  _+  u ^{\prime \prime !} g ^! f ^! u _+} \ar@{=}[r]\ar[d]^{\mathrm{adj} _{u ''}}
  &
  { u ''  _+  u ^{\prime \prime !} g ^! f ^! u _+}
  \ar[r]^{\mathrm{adj} _{u ''}} \ar[d]^{\mathrm{adj} _{u ''}}
  &
  { g ^!  f ^! u _+}\ar@{=}[d]\\
    { (f\circ g) ^!  u _+ }\ar[r]
    &
    {  g ^! f ^! u _+} \ar@{=}[r]
 &
  { g ^! f ^! u _+} \ar@{=}[r]
  &
  { g ^! f ^! u _+} \ar@{=}[r]
  &
  { g ^!  f ^! u _+.}
  }
\end{equation}
La commutativité du rectangle (le {\og seul\fg} rectangle : à gauche et au milieu)
de \ref{diag1-comp-comp-adj-immf} se déduit
de toutes les propriétés, données dans \ref{relev-comp-adj-immf},
des isomorphismes de la forme $\tau$, ainsi que
du diagramme suivant
\begin{equation}
  \notag
  \xymatrix  @R=0,3cm {
   {\X''}
  \ar[r] ^-{b}
  \ar@{~>}[d] ^\tau
  &
  {\X '}
  \ar[r] ^-{a}
  \ar@{~>}[d] ^\tau
  &
  {\X}
  \ar[r] ^u
  &
  {\PP}
  \ar@{=}[d]
  \\
  {\X''}
  \ar[r] ^-{b}
  \ar@{~>}[d] ^\tau
  &
  {\X '}
  \ar[r] ^-{u'}
  &
  {\PP '}
  \ar[r] ^f
  \ar@{=}[d]
  &
  {\PP}
  \ar@{=}[d]
  \\
  {\X''}
  \ar[r] ^-{u''}
  &
  {\PP ''}
  \ar[r] ^g
  &
  {\PP '}
  \ar[r] ^f
  &
  {\PP.}
}
\end{equation}
De plus, on remarque que le morphisme composé
$u ^{\prime !}\overset{\mathrm{adj _{u'}}}{\longrightarrow} u ^{\prime !} u' _+  u ^{\prime !}
\overset{\mathrm{adj _{u'}}}{\longrightarrow} u ^{\prime !}$ est l'identité.
En effet,
cela résulte de l'isomorphisme d'adjonction de bifoncteurs $\theta$ :
$\mathrm{Hom} _{\smash{\D} ^\dag _{\PP '} (\hdag T _{P'})} ( u '_+ (-),\,-)
\riso
\mathrm{Hom} _{\smash{\D} ^\dag _{\X'} (\hdag T _{X'})} (-,\,  u ^{\prime !} (-))$
(on n'utilisera que la fonctorialité à droite),
et du fait que le morphisme identité
$u'_+ u ^{\prime !} \overset{Id}{\longrightarrow} u'_+ u ^{\prime !}$ s'envoie via $\theta$
sur $u ^{\prime !} \overset{\mathrm{adj _{u'}}}{\longrightarrow} u ^{\prime !}u'_+ u ^{\prime !}$
tandis que
$u'_+ u ^{\prime !} \overset{\mathrm{adj _{u'}}}{\longrightarrow} Id$ s'envoie sur
$u ^{\prime !} \overset{Id}{\longrightarrow}u ^{\prime !}$.
Il en résulte la commutativité du carré de gauche de la troisième ligne de \ref{diag1-comp-comp-adj-immf}.

On vérifie ensuite, par définition ou par fonctorialité,
la commutativité des autres carrés de \ref{diag1-comp-comp-adj-immf}.
Ce diagramme est donc commutatif.

Or, on constate que le morphisme composé de gauche de \ref{diag1-comp-comp-adj-immf},
$u ^{''} _+ \circ (a \circ b ) ^! \riso (f \circ g) ^!\circ u  _+$,
n'est autre que $\phi ''$, tandis que celui construit
en prenant le chemin qui passe par le haut puis par la droite du contour de \ref{diag1-comp-comp-adj-immf},
$ u ''  _+ \circ b ^!\circ  a ^! \riso
 g ^!  \circ f ^!\circ  u _+$,
 correspond à $( g ^! \circ \phi)\circ (\phi ' \circ a ^!)$. D'où (i).

Démontrons maintenant (ii).
Grâce à \ref{relev-comp-adj-immf}, le diagramme ci-après
$$\xymatrix  @R=0,3cm@C=2,5cm {
{ u ^\prime_+ a ^{ !} (\E)}
\ar[r]^{\mathrm{adj} _u}
&
{u ^\prime _+  a ^! u ^{ !} u  _+(\E)}
\ar[r]^{u ^\prime _+ \tau _{f \circ u', u\circ a} u _+}
&
{u ^\prime _+   u ^{\prime !} f ^!u  _+(\E)}
\ar[r]^{\mathrm{adj} _{u'}}
&
{f ^!u  _+(\E)}
\\
{ u' _+ a ^{\prime ! (\E)}}
\ar[r]_{\mathrm{adj} _u}
\ar[u] ^-{ u ^\prime_+ (\tau _{a,a'})}
&
{u ^\prime _+  a ^{\prime !} u ^{ !} u  _+(\E)}
\ar[r]_{u ^\prime _+ \tau _{f '\circ u', u\circ a'} u _+}
\ar[u]^{u ^\prime _+ \tau _{u\circ a,u\circ a'} u _+}
&
{u ^\prime _+   u ^{\prime !} f ^{\prime !}u  _+(\E)}
\ar[r]_{\mathrm{adj} _{u'}}
\ar[u] _{u ^\prime _+ \tau _{f\circ u ',f'\circ u'} u _+}
&
{f ^{\prime !}u  _+(\E).}
\ar[u] ^-{\tau _{f,f'} u _+}
}$$
est commutatif. On conclut en remarquant que son contour correspond au diagramme de (ii).
\end{proof}

\subsection{Isomorphismes de recollement : cas rigide}
\begin{vide}
\label{notation-rig-form}

Considérons les diagrammes commutatifs
\begin{equation}
  \label{notation-rig-form-diag0}
  \xymatrix  @R=0,3cm {
{Y ' }
\ar[r] ^-{j '}
\ar[d]^b
&
{X '}
\ar[r] ^-{i '}
\ar[d]^a
&
{\PP '}
\ar[d]^u
\\
{Y }
\ar[r]^{j }
&
{X }
\ar[r]^{i }
&
{\PP ,}
}
\hfill
\xymatrix  @R=0,3cm {
{Y '' }
\ar[r] ^-{j ''}
\ar[d]^{b'}
&
{X ''}
\ar[r] ^-{i ''}
\ar[d]^{a'}
&
{\PP ''}
\ar[d]^{u'}
\\
{Y '}
\ar[r]^{j '}
&
{X '}
\ar[r]^{i '}
&
{\PP ',}
}
\end{equation}
où $u$ (resp. $u'$) est un morphisme de $\V$-schémas formels lisses, $j$, $j'$ et $j''$ sont des immersions ouvertes
de $k$-schémas, $i $, $i'$ et $i''$ sont des immersions fermées.
En notant $\phi =(u,\,a,\,b)$, on définit le foncteur {\it image inverse extraordinaire par
$\phi _K$} d'un $j ^\dag  \smash{\D} _{]X[ _{\PP}}$-module $E $ en posant
$\phi _K ^! (E ) :=
j ^{\prime \dag} \smash{\D} _{]X'[ _{\PP'}} \otimes ^\L _{ u ^{-1} _K j ^\dag  \smash{\D} _{]X[ _{\PP}}} u ^{-1} _K E  [d _{P'/P}]$,
où $d _{P'/P}$ est la dimension relative de $P '$ sur $P $.
Si aucune confusion n'est à craindre, on écrira $u _K ^!$ au lieu de $\phi _K ^!$.
En désignant par $u ^* _K$, l'image inverse
usuelle (notations de \cite[2.2.16]{Berig}), on obtient la relation :
$u _K ^! (E ) \riso u_K ^* (E ) [d _{P'/P}]$.

D'après \cite[2.2.17.(i)]{Berig}, si $v$ : $ \PP ' \rightarrow \PP $ est un morphisme
dont la restriction à $X '$ se factorise par $a$, alors on dispose
d'un isomorphisme $\epsilon _{u,\,v}$ : $v ^{*} _K \riso u _K ^*$ de foncteurs de la catégorie
des $ j ^\dag \O _{]X[ _{\PP}}$-modules à connexion intégrable et surconvergente
dans celle des $ j ^{\prime \dag} \O _{]X'[ _{\PP'}}$-modules à connexion intégrable et surconvergente, tel que
$\epsilon _{u,\,u}= Id$, et que,
si $w$ : $ \PP ' \rightarrow \PP $ est un troisième morphisme coïncidant avec $u$ et $v$ sur $X'$,
on ait la condition de transitivité
$\epsilon_{u,\, w} = \epsilon_{u,\,v}\circ \epsilon_{v,\,w}$.
Si $v'$ est un morphisme
dont la restriction à $X '$ se factorise par $a'$,
avec un raisonnement analogue à \ref{tauhdag}, on prouve
les formules $ \epsilon _{u \circ u', v \circ u'} = u ^{\prime *} _K \circ \epsilon _{u,v}$
et
$\epsilon _{u' ,v'} \circ u ^* _K =  \epsilon _{u \circ u', u \circ v'}$.

Dans le diagramme (de gauche par défaut) \ref{notation-rig-form-diag0}, supposons
que l'on ne dispose que d'un morphisme $u _0$ : $P '\rightarrow P$ rendant commutatif
le diagramme \ref{notation-rig-form-diag0} où $u$ a été remplacé par $u _0$
(on suppose toujours
que $P$ et $P'$ se relèvent en des $\V$-schémas formels lisses $\PP $ et $\PP'$).
De manière analogue à \cite[2.1.6]{Be2}, on construit alors, par recollement,
le foncteur $u _{0K} ^*$, de la catégorie des
$ j ^\dag \O _{]X[ _{\PP} }$-modules à connexion intégrable et surconvergente
dans celle des
$ j ^{\prime \dag} \O _{]X'[ _{\PP'}}$-modules à connexion intégrable et surconvergente.
Précisons sa construction.
Choisissons une application surjective $\rho$ : $\Lambda ' \rightarrow \Lambda$,
deux recouvrements ouverts affines $(\PP _\alpha) _{\alpha \in \Lambda}$
de $\PP$ et $(\PP '_{\alpha'}) _{\alpha '\in \Lambda '}$ de $\PP'$
  tels que $u _0$ se factorise par $P '_{\alpha '} \rightarrow P _{\rho(\alpha ')}$.
  On note alors $X _\alpha := X \cap P _{\alpha}$,
  $X _{\alpha \beta} := X \cap P _{\alpha}\cap P _{\beta}$ et de même en rajoutant des primes.
De plus, choisissons des relèvements $ u _{\alpha '}$ :
$\PP '_{\alpha '} \rightarrow \PP _{\rho(\alpha ')}$
des factorisations induites par $u _0$.
Si $E$ est un $ j ^\dag \O _{]X[ _{\PP} }$-module à connexion intégrable et surconvergente, on lui associe
l'objet
$(( u ^* _{\alpha ' K} E |_{]X _{\rho (\alpha ')}[ _{\PP _{\rho (\alpha ')} }}) _{\alpha ' \in \Lambda '} ,
(\eta '_{\alpha ' \beta '}) _{\alpha ', \beta ' \in \Lambda '})$, où
$\eta '_{\alpha ' \beta '}$ :
$( u ^* _{\beta ' K} E |_{]X _{\rho (\beta ')}[ _{\PP _{\rho (\beta ')} }})
|_{]X' _{\alpha ' \beta '}[ _{\PP' _{\alpha '} \cap \PP' _{\beta '}}}
\riso
( u ^* _{\alpha ' K} E |_{]X _{\rho (\alpha ')}[ _{\PP _{\rho (\alpha ')} }})
|_{]X' _{\alpha ' \beta '}[ _{\PP' _{\alpha '} \cap \PP' _{\beta '}}}$
sont les isomorphismes
de la forme $\epsilon$. Ils résultent de la condition de transitivité des isomorphismes de la forme $\epsilon$
ainsi que de leur comptabilité à l'image inverse que
la famille
$( u ^* _{\alpha ' K} E |_{]X _{\rho (\alpha ')}[ _{\PP _{\rho (\alpha ')} }}) _{\alpha ' \in \Lambda '} ,
(\eta '_{\alpha ' \beta '}) _{\alpha ', \beta ' \in \Lambda '})$
se recolle en un
$ j ^{\prime \dag} \O _{]X'[ _{\PP'}}$-module à connexion intégrable et surconvergente,
celui-ci étant par définition $u _{0K} ^* (E)$. On vérifie ensuite que ce foncteur est
indépendant des choix effectués.

Lorsque $u _0$ se relève en un morphisme
$u$ : $\PP ' \rightarrow \PP$, $u _{0K} ^*$ est canoniquement isomorphe à $u ^*_K$.
Enfin, si $ u' _0$ : $ P'' \rightarrow P'$ est un morphisme rendant commutatif le diagramme de droite
de \ref{notation-rig-form-diag0} lorsque $u'$ a été remplacé par $u' _0$,
on dispose d'un isomorphisme canonique
$(u _{0} \circ u ' _0) _K ^* \riso u _{0K} ^{\prime *} \circ u _{0K} ^*$, celui-ci
étant transitif (pour la composition de diagrammes de la forme \ref{notation-rig-form-diag0}).

Soit $Y _1$ (resp. $Y ' _1 $) un ouvert de $Y $ (resp. $Y '$) tel que
$b (Y ' _1) \subset Y _1$. En notant
$j  _1$ : $ Y  _1 \rightarrow X $ et $j  '_1$ : $ Y ' _1 \rightarrow X '$
les immersions ouvertes déduites,
pour tout $ j ^\dag \O _{]X [ _{\PP} }$-module $E $
à connexion intégrable et surconvergente,
on a le diagramme commutatif de gauche :
\begin{equation}
  \label{epsilonjdagdual}
  \xymatrix  @R=0,3cm {
  {u ^* _K j _1^{ \dag} E }
  \ar[r]
  &
  {j _1 ^{\prime \dag} u ^* _K  E }
  \\
  {v ^{*} _ K j _1^{ \dag} E _1}
  \ar[r]
  \ar[u] ^-{\epsilon _{u , v}}
  &
  {j _1 ^{\prime \dag} v ^{*} _K  E ,}
  \ar[u] ^-{j _1^{\prime \dag} \epsilon _{u , v}}
  }
  \ \ \ \ \ \
  \xymatrix  @R=0,3cm {
  {u ^* _K ( E  ^\vee)}
  \ar[r]
  &
  {( u ^* _K (E )) ^\vee}
  \\
  {v ^{*} _ K ( E  ^\vee )}
  \ar[r]
  \ar[u] ^-{\epsilon _{u , v}}
  &
  {(v ^{ *} _K ( E ))^\vee .}
  \ar[u] ^-{\epsilon ^{\vee -1} _{u , v}}
  }
\end{equation}
De plus, on vérifie par construction celui de droite.
En effet, l'isomorphisme $ p _2 ^* ( E ^\vee _1) \riso p _1 ^* ( E ^\vee _1)$, modulo
les isomorphismes (transitifs par rapport à la composée de deux morphismes) de commutation de l'image inverse au dual,
est par définition l'inverse du dual de $ p _2 ^* ( E  _1) \riso p _1 ^* ( E  _1)$.

\end{vide}

\subsection{Comparaison des isomorphismes formels et rigides de recollements}

Les deux propositions ci-après seront ensuite étendues, en \ref{spjdag} et \ref{sp+f*}, au cas non relevable.

La première partie de la proposition qui suit est due à Noot-Huyghe (\cite[1.5.3]{thesehuyghe}).
Cependant, nous donnons ici une démonstration
plus formelle qui fait intervenir la quasi-cohérence sur les schémas formels et nous
nous assurons de sa transitivité.
\begin{prop}\label{spcommup*}
  Soient $f $ : $\X '\rightarrow \X$ un morphisme de $\V$-schémas formels lisses,
  $H$ un diviseur de $X$, $H' \supset  f _0 ^{-1} (H)$ un diviseur de
  $X'$, $j $ : $X \setminus H \hookrightarrow X$ (resp. $j '$ : $X '\setminus H' \hookrightarrow X'$)
  l'immersion ouverte, $E$ un $j^\dag \O _{\X _K}$-module
  cohérent muni d'une connexion surconvergente $\nabla$ et $\E:= \sp _* (E)$,
  le $\smash{\D} _{\X,\,\Q} ^\dag (\hdag H)$-module cohérent associé.
  On désignera par $f ^* _K$ le foncteur image inverse de la catégorie des
$j ^\dag \O _{\X _K}$-modules cohérents à connexion surconvergente dans celle
des $j ^{\prime \dag} \O _{\X '_K}$-modules cohérents à connexion surconvergente.

  Il existe un isomorphisme $j ^{\prime \dag} \smash{\D} _{\X' _K}$-linéaire,
  $f _K ^* \sp ^* (\E) \riso \sp ^* f ^! _{H',H} (\E) [-d _{X'/X}]$, et un second
  $\smash{\D} ^\dag _{\X'} (\hdag H ') _{\Q}$-linéaire,
  $\sp _* f _K ^* (E)  \riso f ^! _{H',H} \sp _* (E)[-d _{X'/X}]$.

  En outre, ces derniers sont transitifs, i.e., si
  $g $ : $ \X '' \rightarrow \X '$ est un second morphisme de $\V$-schémas formels lisses,
   si $H '' \supset g _0 ^{-1} (H')$ est un diviseur de $X''$, alors le diagramme
   \begin{equation}\label{diag0spcommup*}
     \xymatrix  @R=0,3cm {
     { \sp _* g _K^* f _K^* (E) }
     \ar[r] _(0.35){\sim}
     \ar[d] _{\sim}
     &
     {g ^! _{H'',H'} \sp _* f _K^* (E) [- d_{X ''/X'}]}
     \ar[r] _(0.47){\sim}
     &
     { g ^! _{H'',H'} f ^! _{H',H}\sp _* (E)[- d_{X ''/X}]}
     \ar[d] _{\sim}
     \\
     { \sp _*  (f \circ g) _K^* (E ) }
     \ar[rr] _(0.4){\sim}
     &
     &
     { (f \circ g) ^! _{H'',H}  \sp _*   (E)[- d_{X ''/X}] }
     }
   \end{equation}
   est commutatif.
   De même, on obtient un second diagramme commutatif en remplaçant dans \ref{diag0spcommup*}
   $\sp _*$ par $\sp ^*$, $E$ par $\E$, $-$ par $+$, en inversant le sens des flèches
   et en intervertissant les symboles en haut $*$ et $!$.
\end{prop}
\begin{proof}
Nous noterons abusivement $f ^!$, $g ^!$ et $(f \circ g ) ^!$ à la place
de $f ^! _{H',H}$, $g ^! _{H'',H'}$ et $(f \circ g) ^! _{H'',H}$.
On pose
$f ^* (\E) := \O _{\X'} (\hdag H') _\Q
  \otimes _{f ^{-1} \O _{\X} (\hdag H) _\Q}   f ^{-1} \E$.
  Le foncteur
  $E \mapsto f _K ^* (E) \riso   j^{'\dag } \O _{\X '_K} \otimes _{f _K ^{-1} j^{\dag } \O _{\X_K}} f _K ^{-1} (E)$
  est une image inverse de sites annelés, de même pour $f ^*$ et $\sp ^*$.
  Grâce à la transitivité de l'isomorphisme de commutation à la composition
  des images inverses de sites annelés, on obtient
  l'isomorphisme $f _K ^* \sp ^* (\E) \riso \sp ^* f ^* (\E)$
   commutant aux stratifications respectives (induites par celle de $\E$ qui provient de sa structure de
    $\O _{\X} (\hdag H ) _\Q \otimes _{\O _{\X,\,\Q}} \smash{\D} _{\X,\,\Q} $-module).

Soit $\E ^{(\bullet)} \in \smash{\underset{^{\longrightarrow}}{LD}} ^{\mathrm{b}} _{\Q
,\mathrm{coh}} ( \smash{\widehat{\D}} _{\PP} ^{(\bullet)}(T))$
tel que $\underset{\longrightarrow}{\lim} (\E ^{(\bullet)}) \riso \E$.
On dispose, de manière analogue à \cite[4.3.2]{Beintro2},
de l'isomorphisme canonique :
$$f ^! (\E) \riso
\underset{\longrightarrow}{\lim} (
(\widehat{\B} ^{(m)} _{\X'} ( H')  \smash{\widehat{\otimes}} ^\L
_{f ^{-1}\widehat{\B} ^{(m)}  _{\X} ( H) } f ^{-1}\E ^{(m)}) _{m\in \N} )[d_{X '/X}].$$
  De plus, puisque que $\E$ est un
  $\O _{\X} (\hdag H) _\Q$-module cohérent, grâce au lemme \cite[I.7.1]{HaRD} sur les foncteurs
  way-out,
  le morphisme canonique $\O _{\X'} (\hdag H ') _\Q \otimes _{\O _{\X',\,\Q}} \smash{\D} _{\X ',\,\Q} $-linéaire
  $$\O _{\X'} (\hdag H') _\Q
  \otimes _{f ^{-1} \O _{\X} (\hdag H) _\Q}
  f ^{-1} \E
  \rightarrow
  \underset{\longrightarrow}{\lim} (
(\widehat{\B} ^{(m)} _{\X'} ( H')  \smash{\widehat{\otimes}} ^\L
_{f ^{-1}\widehat{\B} ^{(m)}  _{\X} ( H) } f ^{-1}\E ^{(m)}) _{m\in \N} )$$
  est un isomorphisme.
  Il en résulte un isomorphisme $f _K ^* \sp ^* (\E) \riso \sp ^* f ^! (\E) [-d_{X '/X}]$.
  En outre, la transitivité de cet isomorphisme est immédiate.
  En effet, il suffit d'invoquer
  la transitivité de l'isomorphisme de commutation à la composition des images inverses de sites annelés
 (et d'ajouter les isomorphismes de la forme
$ f ^* (\E) \riso f ^!(\E) [-d _{X'/X}]$).

Prouvons à présent le deuxième isomorphisme.
  D'après \cite[4.4.2]{Be1}, les foncteurs $\sp _*$ et $\sp ^*$ sont des équivalences quasi-inverses de la catégorie
  des $j^{\dag } \O _{\X_K}$-modules cohérents munis d'une connexion intégrable
  dans celle des $\O _{\X} (\hdag H) _\Q$-modules cohérents munis d'une connexion intégrable
  (de même en rajoutant des primes). Il existe donc un unique isomorphisme
  $\O _{\X'} (\hdag H ') _\Q \otimes _{\O _{\X',\,\Q}} \smash{\D} _{\X',\,\Q} $-linéaire,
  $\sp _* f _K ^* (E) \riso  f ^* \sp _* (E)$, s'inscrivant dans le diagramme canonique
  \begin{equation}\label{diag1spcommup*}
  \xymatrix  @R=0,3cm {
    {\sp ^* f ^* \sp _* (E)}
    &
    {f ^* _K \sp ^*  \sp _* (E) }
    \ar[r] _(0.6){\sim}
    \ar[l] ^-{\sim}
    &
    {f ^* _K (E)}
    \\
    {\sp ^* \sp _* f _K ^* (E)}
    \ar@{.>}[u]
    \ar[rr] _{\sim}
    &
    &
    {f ^* _K (E).}
    \ar@{=}[u]
    }
  \end{equation}
Via l'isomorphisme $\O _{\X'} (\hdag H ') _\Q \otimes _{\O _{\X',\,\Q}} \smash{\D} _{\X',\,\Q} $-linéaire,
$f ^* (\sp _* (E)) \riso  f ^!(\sp _* (E))[-d _{X'/X}]$,
on construit le suivant
 $\sp _* f _K ^* (E)\riso  f ^! (\sp _* (E))[-d _{X'/X}]$.
 Ce dernier, étant un morphisme $\O _{\X'} (\hdag H ') _\Q \otimes _{\O _{\X',\,\Q}} \smash{\D} _{\X',\,\Q} $-linéaire
 entre deux isocristaux surconvergents (ou plutôt $\smash{\D} ^\dag _{\X'} (\hdag H ') _{\Q}$-modules cohérents
 dans l'image essentielle du foncteur $\sp _*$ de \cite[4.4.5]{Be1}), est $\smash{\D} ^\dag _{\X'} (\hdag H ') _{\Q}$-linéaire.

Vérifions maintenant sa transitivité.
Le diagramme \ref{diag0spcommup*} est commutatif, si et seulement si le suivant (qui correspond globalement
à l'image de \ref{diag0spcommup*} par le foncteur $\sp ^*$) l'est.
   \begin{equation}\label{diag2spcommup*}
     \xymatrix  @R=0,3cm {
     &
     {\sp ^*g ^! \sp _* f _K^* (E) [-d _{X''/X'}]}
     \ar[r] _{\sim}
     &
     {\sp ^* g ^!f ^!\sp _* (E)[-d _{X''/X}]}
     \ar@{=}[rd]
     \\
     {\sp ^* \sp _* g _K^* f _K^* (E) }
     \ar[r] _{\sim}
     \ar[d] _{\sim}
     &
     {\sp ^*g ^* \sp _* f _K^* (E)}
     \ar[r] _{\sim}
     \ar[u] _{\sim}
     &
     {\sp ^* g ^* f ^* \sp _* (E)}
     \ar[u] _{\sim}
     \ar[r] _(0.4){\sim}
     \ar[d] _{\sim}
     &
     {\sp ^* g ^! f ^!\sp _* (E)[-d _{X''/X}]}
     \ar[d] _{\sim}
     \\
     {\sp ^* \sp _*  (f \circ g) _K ^* (E) }
     \ar[rr] _{\sim}
     &
     &
     { \sp ^*(f \circ g)^* \sp _*   (E) }
     \ar[r] _(0.4){\sim}
     &
     { \sp ^*(f \circ g)^!\sp _*   (E) [-d _{X''/X}].}
     }
   \end{equation}
   La commutativité des deux carrés et du triangle est aisée.
Il reste à démontrer celle du rectangle en bas à gauche de \ref{diag2spcommup*}.
Or, celui-ci s'inscrit (à gauche) dans le suivant :
\begin{equation}\label{diag3spcommup*}
     \xymatrix  @R=0,3cm {
     {\sp ^* \sp _* (f \circ g ) _K^* (E)}
     \ar[dd] _{\sim}
     &
     {\sp ^* \sp _*  g _K ^* f _K^* (E)}
     \ar[drr] _{\sim}
     \ar[l] ^-{\sim}
     \ar[d] _{\sim}
     \\
     &
     {\sp ^*  g ^* \sp _*  f _K^* (E)}
     \ar[d] _{\sim}
     &
     {g _K ^* \sp ^*  \sp _*  f _K^* (E)}
     \ar[r] _{\sim}
     \ar[l] ^-{\sim}
     \ar[d] _{\sim}
     &
     {g ^* _Kf ^* _K(E)}
     \ar[r] _{\sim}
     &
     {(f \circ g ) _K^* (E)}
     \\
     {\sp ^* (f \circ g )^* \sp _*  (E)}
     &
     {\sp ^*  g ^* f ^* \sp _*   (E)}
     \ar[l] ^-{\sim}
     &
     {g ^*_K \sp ^*   f ^* \sp _*   (E)}
     \ar[l] ^-{\sim}
     &
     {g _K^* f _K^* \sp ^*    \sp _*   (E)}
     \ar[r] _{\sim}
     \ar[u] _{\sim}
     \ar[l] ^-{\sim}
     &
     {(f \circ g) _K^* \sp ^*    \sp _*   (E).}
     \ar[u] _{\sim}
     }
   \end{equation}
Le deuxième carré de droite et le triangle de \ref{diag3spcommup*}
sont commutatifs par construction (voir \ref{diag1spcommup*}),
tandis que les deux autres carrés sont fonctoriels.
Or, la flèche composée du bas de \ref{diag3spcommup*}
est l'isomorphisme canonique
$\sp ^* (f \circ g) ^* \sp _* (E)  \riso (f \circ g) _K^* \sp ^*  \sp _*  (E) $.
De plus, par fonctorialité de $\sp ^* \sp _* \rightarrow Id$, le morphisme composé
du haut de \ref{diag3spcommup*} correspond au morphisme canonique
$\sp ^* \sp _*  (f \circ g) _K ^*(E)  \riso \sp ^* (f \circ g) ^* \sp _* (E)  $.
Il en découle que le contour de
\ref{diag3spcommup*} est le diagramme canonique ci-dessous
$$\xymatrix  @R=0,3cm {
{\sp ^* \sp _*  (f \circ g) _K ^*(E) }
\ar[r] _{\sim}
\ar[d] _{\sim}
&
{(f \circ g) _K^*(E) }
\\
{\sp ^* (f \circ g) ^* \sp _* (E)  }
&
{(f \circ g) _K^* \sp ^*  \sp _*  (E) ,}
\ar[u] _{\sim}
\ar[l] ^-{\sim}
}
$$
qui est commutatif d'après \ref{diag1spcommup*}.
Le rectangle de gauche de \ref{diag3spcommup*} est donc aussi commutatif.
\end{proof}

\begin{rema}\label{spcommup*2}
  Avec les notations de \ref{spcommup*},
  le diagramme canonique suivant
  $$\xymatrix  @R=0,3cm {
  {\sp _* \sp ^* f ^! _{H',H}\E}
  &
  {\sp _* f _K ^* \sp ^* (\E) [ d_{X '/X}]}
  \ar[l] ^-{\sim}
  \ar[r] _{\sim}
  &
  {f ^! _{H',H}\sp _*  \sp ^* (\E)}
  \\
  { f ^! _{H',H}(\E)}
  \ar@{=}[rr]
  \ar[u] _{\sim}
  &
  &
  {f ^! _{H',H}(\E)}
  \ar[u] _{\sim}
  }
  $$
  est commutatif.
\end{rema}
\begin{proof}
Considérons le diagramme suivant
\begin{equation}
  \label{diag1-spcommup*2}
  \xymatrix  @R=0,3cm {
  {\sp ^* f ^! _{H',H}(\E)}
  \ar@{=}[d]
  \ar[r] _(0.45){\sim}
  ^(0.45){\mathrm{adj}}
  &
  {\sp ^* \sp _* \sp ^* f ^!_{H',H} (\E)}
  &
  {\sp ^* \sp _* f _K ^* \sp ^* (\E)[ d_{X '/X}]}
  \ar[l] ^-{\sim}
  \ar[r] _{\sim}
  ^{\mathrm{adj}}
  \ar[d] _{\sim}
  &
  {f _K ^* \sp ^* (\E)[ d_{X '/X}]}
  \\
  {\sp ^* f ^!_{H',H} (\E)}
  \ar[rr] _{\sim} ^{\mathrm{adj}}
  &
  &
  {\sp ^* f ^!_{H',H}  \sp _* \sp ^* (\E)}
  &
  {f _K ^* \sp ^* \sp _* \sp ^* (\E)[ d_{X '/X}],}
  \ar[l] ^-{\sim}
  \ar[u] _{\sim} ^{\mathrm{adj}}
  }
\end{equation}
où le rectangle à gauche est l'image par $\sp ^*$ de celui de \ref{spcommup*2}.
Grâce à  que \ref{diag1spcommup*},
le carré est commutatif.
De plus, on dispose du diagramme commutatif :
\begin{equation}
  \notag
  \xymatrix  @R=0,3cm {
  {\sp ^* f ^!_{H',H} (\E)}
  \ar@{=}[d]
  \ar[r]^(0.38){\mathrm{adj}} _(0.38){\sim}
  &
  {\sp ^* \sp _* \sp ^* f ^!_{H',H} (\E)}
  \ar[d]^{\mathrm{adj}} _{\sim}
  &
  {\sp ^* \sp _* f _K ^* \sp ^* (\E)[ d_{X '/X}]}
  \ar[l] ^-{\sim}
  \ar[d]^{\mathrm{adj}} _{\sim}
  \\
  {\sp ^* f ^!_{H',H} (\E)}
  \ar@{=}[r]
  &
  {\sp ^* f ^!_{H',H} (\E)}
  &
  {f _K ^* \sp ^* (\E)[ d_{X '/X}].}
  \ar[l] ^-{\sim}
  }
\end{equation}
Il en résulte que l'isomorphisme composé du haut de
\ref{diag1-spcommup*2},
$\sp ^* f ^! _{H',H}(\E) \tilde{\leftarrow} f _K ^* \sp ^* (\E)[ d_{X '/X}]$, est le morphisme canonique.
De la même manière, en utilisant le diagramme commutatif suivant
$$\xymatrix  @R=0,3cm {
{f _K ^* \sp ^* (\E) }
\ar@{=}[d]
&
{f _K ^* \sp ^* \sp _* \sp ^* (\E) }
\ar[r] _-{\sim}
\ar[l] ^-{\sim}
_(0.55){\mathrm{adj}}
&
{\sp ^* f ^!_{H',H} \sp _* \sp ^* (\E) [-d_{X'/X}]}
\\
{f _K ^* \sp ^* (\E) }
\ar@{=}[r]
&
{f _K ^* \sp ^* (\E) }
\ar[r] _(0.4){\sim}
\ar[u] _{\sim}
^{\mathrm{adj}}
&
{\sp ^* f ^!_{H',H} (\E) [-d_{X'/X}],}
\ar[u] _{\sim}
^{\mathrm{adj}}
}
$$
on établit que le morphisme $f _K ^* \sp ^* (\E)[ d_{X '/X}]
\riso \sp ^* f ^!_{H',H} (\E)$
de \ref{diag1-spcommup*2} passant par la droite puis en bas est le morphisme canonique de commutation (construit
en \ref{spcommup*}).
Le diagramme \ref{diag1-spcommup*2}
est donc commutatif.
\end{proof}

\begin{prop}
  \label{spcommjdag}
Soient $f$: $\X' \rightarrow \X$ un morphisme de $\V$-schémas formels lisses, $H _1\subset H_2$ deux diviseurs de $X$,
$H' _1 \subset H' _2$ deux diviseurs de $X'$ tels que, pour $i =1,\,2$, $H ' _i \supset f _0 ^{-1}(H _i)$.
On note, pour $i =1,\,2$, $j _i$ : $ X \setminus H _i \hookrightarrow X$ et
$j _i'$ : $ X '\setminus H '_i\hookrightarrow X'$ les immersions ouvertes.
Pour tout isocristal $E _1$ sur $X \setminus H _1$ surconvergent le long de $H _1$, en notant
et $\E _1 = \sp _* (E _1)$,
on dispose des isomorphismes
$j _2^{ \dag } \sp ^* (\E _1) \riso \sp ^*  (\hdag H _2) (\E _1)$ et
$\sp _* j _2^{ \dag }(E _1)  \riso (\hdag H '_2) \circ \sp _* (E _1)$.
Ceux-ci commutent aux images inverses extraordinaires, i.e.,
le diagramme canonique
$$\xymatrix  @R=0,3cm {
{(\hdag H' _2) \sp _* f ^* _K (E _1)[d _{X'/X}]}
\ar[r] _{\sim}
&
{(\hdag H' _2) f ^! _{H' _1,H_1} \sp _*  (E _1)}
\ar[r] _{\sim}
&
{ f ^! _{H' _2,H_2} (\hdag H _2) \sp _*  (E _1)}
\\
{\sp _* j ^{\prime \dag} _2 f ^* _K (E _1)[d _{X'/X}]}
\ar[r] _{\sim}
\ar[u] _{\sim}
&
{\sp _*  f ^* _K j ^{\dag} _2 (E _1)[d _{X'/X}]}
\ar[r] _{\sim}
&
{f ^! _{H' _2,H_2} \sp _* j ^{\dag} _2 (E _1),}
\ar[u] _{\sim}
}$$
est commutatif.
\end{prop}
\begin{proof}
Cela découle de la transitivité de \ref{spcommup*}. En effet,
les foncteurs de la forme $j ^\dag $ (resp. $(\hdag T)$) sont aussi des images inverses
(resp. des images inverses extraordinaires).
\end{proof}

Afin de prouver que les isomorphismes de la forme
$\epsilon$ et $\tau$ se correspondent via les foncteurs quasi-inverses $\sp _*$ et $\sp ^*$ (voir \ref{sp-eps-tau}),
nous aurons besoin de quelques compléments à propos de la notion {\og topologiquement nilpotent\fg}.
\begin{lemm}\label{lemmquasi-nilp}
  Soient $f$ : $ X \rightarrow S$ un morphisme lisse, $\B$ une $\O _X$-algèbre
  munie d'une structure compatible de $\smash{\D} _{X/S} ^{(m)}$-module et
  $\E$ un $\B \otimes _{\O _X} \smash{\D} ^{(m)} _{X/S}$-module. On suppose que $p$ est localement nilpotent sur
  $S$, et qu'il existe sur $X$ un système de coordonnées locales $t _1,\dots,t_d$.
  Les conditions suivantes sont équivalentes :

  (i) Pour toute section $e$ de $\E$, il existe localement un entier $N$ tel que,
  pour tout $\underline{k}\in \N ^d$ tel que $|\underline{k}| \geq N$,
  $\underline{\partial} ^{<\underline{k}>}e =0$.

  (ii) La condition (i) est satisfaite pour tout système de coordonnées locales sur un ouvert de $X$.

  (iii) Il existe un isomorphisme de $\B \otimes _{\O _X} \PP _{X/S,(m)} $-algèbres
  $$\epsilon  \ : \ (\B \otimes _{\O _X} \PP _{X/S,(m)} )\otimes _{\B} \E
  \riso
   \E \otimes _{\B}(\B \otimes _{\O _X} \PP _{X/S,(m)})$$
   vérifiant la condition de cocycle, et induisant par réduction sur les
   $\B \otimes _{\O _X} \PP _{X/S,(m)}^n$ la $m$-PD-stratification de $\E$ relative à $\B$
   (voir \cite[1.1.15]{caro_comparaison}).
   De plus, cet isomorphisme est alors unique et déterminé,
   pour toute section locale $e$ de $\E$, par la relation
   $\epsilon ( (1\otimes 1) \otimes e ) =
   \sum _{\underline{k}\geq 0}
   \underline{\smash{\widetilde{\partial}}} ^{<\underline{k}>} e \otimes
   \underline{\smash{\widetilde{\tau}}} ^{\{ \underline{k}\}}$,
   où $\underline{\smash{\widetilde{\tau}}} ^{\{ \underline{k}\}} =1 \otimes \underline{\tau} ^{\{ \underline{k}\}}$
   et
   $\underline{\smash{\widetilde{\partial}}} ^{<\underline{k}>} := 1\otimes \underline{\partial} ^{<\underline{k}>}$.
\end{lemm}
\begin{proof}
  On calque \cite[2.3.7]{Be1}.
\end{proof}

\begin{vide}\label{espilonhat}
  On étend, de manière analogue à la suite de \cite[2.3.7]{Be1}, les notions de
{\it quasi-nilpotence}, {\it nilpotence} et {\it topologiquement nilpotence}.
Soient $\X$ un $\V$-schéma formel, $H$ un diviseur de $X$, $n_m \geq m$ deux entiers et $\E$ un
$\widehat{\B} ^{(n_m)} _{\X} (H) \widehat{\otimes} _{\O _{\X}} \widehat{\D} ^{(m)} _{\X}$-module
qui soit $\widehat{\B} ^{(n_m)} _{\X} (H)$-cohérent. Notons $X _i$ la réduction modulo $\pi ^{i+1}$ de $\X$
et $\E _i :=\O _{X _i} \otimes _{\O _{\X}} \E$.
On déduit par complétion de \ref{lemmquasi-nilp} que les affirmations suivantes sont équivalentes :

(i) Le faisceau $\E$ est topologiquement nilpotent ;

(ii) Il existe un isomorphisme de $\widehat{\B} ^{(n_m)} _{\X} (H) \widehat{\otimes} _{\O _{\X}} \PP _{X,(m)}$-algèbres
$$\widehat{\epsilon}  \ : \ (\widehat{\B} ^{(n_m)} _{\X} (H)
\widehat{\otimes} _{\O _{\X}} \PP _{\X,(m)} )\otimes _{\widehat{\B} ^{(n_m)} _{\X} (H)} \E
  \riso
   \E \otimes _{\widehat{\B} ^{(n_m)} _{\X} (H)}
   (\widehat{\B} ^{(n_m)} _{\X} (H) \widehat{\otimes} _{\O _{\X}} \PP _{\X,(m)})$$
   vérifiant la condition de cocycle, et induisant par réduction sur les
   $\B ^{(n_m)} _{X _i} (H) \otimes _{\O _{X _i} }\PP _{X _i,(m)}^n$ (pour tous entiers $i $ et $n$ )
   la $m$-PD-stratification de $\E _i $ relative à $\B ^{(n_m)} _{X _i} (H)$.
   De plus, cet isomorphisme est alors unique et satisfait, au dessus d'un ouvert possédant des coordonnées locales
   $t _1,\dots ,t_d$
   et pour toute section locale $e$ de $\E$,
   la relation
   $\widehat{\epsilon} ( (1\otimes 1) \otimes e ) =
   \sum _{\underline{k}\geq 0} \underline{\smash{\widetilde{\partial}}} ^{<\underline{k}>} e \otimes
   \underline{\smash{\widetilde{\tau}}} ^{\{ \underline{k}\}}$. Cette relation caractérise
   $\widehat{\epsilon}$.
\end{vide}

\begin{theo}\label{sp-eps-tau}
Avec les notations de \ref{spcommup*},
  pour tout morphisme $f'$ : $\X ' \rightarrow \X$ tel que $f' _0 =f _0$,
les diagrammes suivants
$$\xymatrix  @R=0,3cm@C=2cm {
{ \sp _* f ^{\prime *} _K (E )[d _{X'/X}]}
\ar[r] ^-{\sp _*(\epsilon _{f,\,f'})} _{\sim}
\ar[d] _{\sim}
&
{ \sp _* f ^{*} _K (E )[d _{X'/X}]}
\ar[d] _{\sim}
\\
{ f ^{\prime !} _{H',H} \sp _*  (E )}
\ar[r] ^-{\tau _{f ,f',H',H  }} _{\sim}
&
{ f ^{!} _{H',H}  \sp _*  (E ),}
}
\xymatrix  @R=0,3cm@C=2cm {
{ f ^{\prime *} _K \sp ^* (\E )[d _{X'/X}]}
\ar[r] ^-{\epsilon _{f,\,f'}} _{\sim}
\ar[d] _{\sim}
&
{ f ^{*} _K   \sp ^*(\E )[d _{X'/X}]}
\ar[d] _{\sim}
\\
{\sp ^*  f ^{\prime !}   _{H',H} (\E )}
\ar[r] ^-{\sp ^* (\tau _{f ,f',H',H  })} _{\sim}
&
{ f ^{!} _{H',H} \sp _*  (\E ),}
}
$$
où les isomorphismes verticaux résultent de \ref{spcommup*} et
ceux horizontaux sont construits dans les sections
\ref{relev-comp-adj-immf} et de \ref{notation-rig-form},
sont commutatifs.
\end{theo}
\begin{proof}
L'assertion étant locale en $\X$ et $\X'$, on se ramène au cas où
$\X =\Spf A$ et $\X '=\Spf A'$ sont affines et possèdent des coordonnées locales
et où $H$ (resp. $H'$) est l'ensemble des zéros modulo $\mathfrak{m}$ d'un élément $g$ de $A$ (resp. $g '$ de $ A'$).
  Nous prenons les notations de \cite[1.1.8 et 1.2]{Berig}. Pour tout $1 >\lambda \geq |\pi |$, on a ainsi
  $U _{\lambda} := ] X  [ _{\X } - ] H [ _{\X \, \lambda} = \{ x \in \X _K ; | g (x)| \geq \lambda \} $,
  $U '_{\lambda} := ] X  '[ _{\X' } - ] H '[ _{\X '\, \lambda} = \{ x '\in \X '_K ; | g ' (x')| \geq \lambda \} $
  (l'hypothèse $\lambda \geq |\pi |$ implique que ceux-ci ne dépendent pas du choix des sections $g$ et
  $g'$ définissant $H$ et $H'$).
  Pour tout $m \in \N$, on écrira $\lambda _m$ ou $\mu _m$ pour $ p ^{-1 /p ^{m+1}}$.
  On remarque que
  $\Gamma (\X,\,\widehat{\B} ^{(m)} _{\X} (H)) _\Q
  =\Gamma ( U _{\lambda _m},\, \O _{\X_K})$.

  Soit $m _0$ suffisamment grand tel
  qu'il existe un $\O _{\lambda _{m _0}}$-module cohérent $E _0$, muni d'une connexion intégrable
  $\nabla _0$, et un isomorphisme $(E,\,\nabla) \riso j ^\dag (E _0,\, \nabla _0)$.
  Pour tout $m \geq m _0$, on note $(E _m,\nabla _m)$ la restriction de $(E _0,\nabla _0)$
  à $U _{\lambda _m}$, $\E _0 = \sp _* (E _0)$
  et $\E ^{(m)} = \sp _* (E _m)$, où $\sp $ est le composé
  $U _{\lambda _m} \subset \X _K \overset{\sp}{\longrightarrow} \X$ induit par le morphisme de spécialisation.

Notons respectivement $p_1$ et $p _2$ les projections $] X[_{\X ^2 } \rightarrow \X _{K}$ à gauche et à droite,
$\I$ l'idéal définissant l'immersion diagonale
$\X _K \hookrightarrow \X _K ^2$,
$\PP ^n = \O _{\X ^2 _K} / \I ^{n+1}$
et
$\epsilon _{0 n}$ :
$\PP ^n |_{U _{\lambda _{m_0}}} \otimes _{\O _{U _{\lambda _{m_0}}}} E _{0} \riso
E _{0}
\otimes _{\O _{U _{\lambda _{m_0}}}} \PP ^n |_{U _{\lambda _{m_0}}}$,
les isomorphismes de la stratification associée à $\nabla _0$.
Comme la connexion $\nabla$ est surconvergente, d'après \cite[2.2.6]{Berig}, il existe un voisinage strict $V'$
de $]X \setminus H[ _{\X ^2 }$ dans $]X[_{\X^2}$, contenu dans
$p _1 ^{-1} (U _{\lambda _{m_0}} )\cap p _2 ^{-1} (U _{\lambda _{m_0}} )$,
tel que $\epsilon $ soit de la forme
$\epsilon = j ^{\dag} (\epsilon  _0)$,
où $\epsilon  _0$ : $p _2 ^* (E _{m_0}) |_{V'} \riso p _1^* (E _{m_0}) |_{V'}$ est un isomorphisme
de $\O _{V'}$-modules, induisant pour tout $n$, par réduction modulo $\I ^{n+1}$, la restriction à $V' \cap \X _K$ des
isomorphismes $\epsilon _{0 ,n}$.

Comme $f _0 ^{-1} (H) \subset H'$, il existe des sections $a' $ et $b'$ de $A'$ telles que
  $g' = a' f ^* (g) + \pi b'$. Si $\lambda > |\pi |$,
  il en résulte, pour tout $x ' \in U ' _{\lambda}$,
  que $\lambda \leq |f _K ^* (g) (x')|=|g (f  _K (x'))|$. Ainsi, lorsque $\lambda > |\pi |$,
  $f _K$ induit le morphisme
  $U '_{\lambda} \rightarrow U _{\lambda}$. De même en remplaçant $f$ par $f'$.

D'après \cite[4.4]{Be1}, il existe une suite d'entiers $( n _m) _{m\in \N}$
  telle que $n _m >\max \{m,\, m _0 \}$ et telle que $\E ^{(n _m)}$ soit muni d'une structure de
  $\widehat{\B} ^{(n_m)} _{\X} (H) \widehat{\otimes} \widehat{\D} ^{(m)} _{\X,\,\Q}$-module cohérent
  topologiquement nilpotent. Via \cite[4.4.7]{Be1}, il existe un
  $\widehat{\B} ^{(n_m)} _{\X} (H) \widehat{\otimes} \widehat{\D} ^{(m)} _{\X}$-module cohérent
  $\smash{\overset{_{\circ}}{\E}} ^{(n_m)}$, cohérent sur $\widehat{\B} ^{(n_m)} _{\X} (H)$
  et sans $p$-torsion, et un isomorphisme
  $\widehat{\B} ^{(n_m)} _{\X} (H) \widehat{\otimes} \widehat{\D} ^{(m)} _{\X,\,\Q}$-linéaire
  $\smash{\overset{_{\circ}}{\E}} ^{(n_m)} _\Q \riso \E ^{(n_m)}$.

  Comme $1 >\lambda _{n_m}> \eta _m > |\pi |$,
  $[X ] _{\X ^2 \,\eta _{m}} \cap p _1 ^{-1} (U _{\lambda _{n _m}} )  =
  [X ] _{\X ^2 \,\eta _m} \cap p _2 ^{-1} (U _{\lambda _{n _m}}) =:V _{\eta _m \,\lambda _{n _m}}$.
Notons $h = (f,\,f')$ : $ \X ' \rightarrow \X \times _\S \X$.
Puisque $f$ et $f '$ coïncident sur $X'$, $h ^* ( \tau _i) \in \pi A'$. On en déduit que pour
tout $x' \in \X '_K$, $h _K (x') \in [X ] _{\X ^2 \,|\pi |} \subset [X ] _{\X ^2 \,\eta _m}$.
Il en découle que le morphisme
$h _K$ induit la flèche $ U ' _{\lambda _{n _m}} \rightarrow V _{\eta _m \,\lambda _{n _m}}$,
notée abusivement $h _k$,
qui s'insère dans le diagramme
commutatif
$$\xymatrix  @R=0,3cm@C=2cm {
{U '_{\lambda _{n _m}}}
\ar@<1ex>[d] ^-{f  _K}
\ar@<-1ex>[d] _{f ' _K}
\ar[rd] ^-{h _K}
\\
{U _{\lambda _{n _m}}}
&
{V _{\eta _m \,\lambda _{n _m}}.}
\ar@<1ex>[l] ^-{p _1}
\ar@<-1ex>[l] _{p _2}
}
$$
Via \cite[1.2.2]{Berig}, quitte à augmenter $n _m$, on peut supposer $V _{\eta _m \,\lambda _{n _m}}\subset V '$.
On note alors $\epsilon _m$ :
$p _2 ^* (E _{n_m} )\riso p ^* _1  (E _{n_m})$ l'isomorphisme $\epsilon _0 |_{V _{\eta _m \,\lambda _{n _m}}}$.
Pour tout $n$, la restriction à $U _{\lambda _{n _m}}$ des réduction modulo $\I ^{n+1}$ des
isomorphismes $\epsilon _{m}$ est
$\epsilon _{0 n} |_{U _{\lambda _{n_m}}}$ et correspond donc
aux isomorphismes de la stratification de $E _{n_m}$ induite par $\nabla _{n_m}$.

En considérant $\X ^2$ comme $\X$-schéma via la projection à gauche,
notons $\phi$ : $\X ^2 \rightarrow \X \times \widehat{\A} ^d _{\V}$ le morphisme de $\X$-schémas formels
défini par $\tau _1 ,\dots , \tau _d$.
Puisque $\X$ possède des coordonnées locales, il résulte du théorème de fibration fort \cite[1.3.7]{Berig}
qu'il existe des voisinages stricts $V'$ de $]X \setminus H[ _{\X^2}$ dans $]X[ _{\X^2}$
et $V''$ de $]X \setminus H[ _{\X \times \widehat{\A} ^d _{\V}}$ dans $]X[ _{\X \times \widehat{\A} ^d _{\V}}$ tels que
$\phi _K$ induise l'isomorphisme $V ' \riso V''$.
Par construction de $\phi$,
$\phi _K ^{-1} (U _{\lambda _{n_m} } \times D ( 0,\eta _m {} ^+ )  )
= V _{\eta _m \,\lambda _{n _m}}$.
De plus, quitte à accroître $n_m$ et grâce à \cite[1.2.2]{Berig}
(et à \cite[1.2.3.(iii)]{Berig} pour la deuxième inclusion), on peut supposer
$U _{\lambda _{n_m} } \times D ( 0,\eta _m {} ^+ ) \subset V''$
et
$V _{\eta _m \,\lambda _{n _m}} \subset V'$.
L'isomorphisme $\phi _K$ : $ V' \riso V''$ induit alors le suivant
$$V _{\eta _m \,\lambda _{n _m}} \riso U _{\lambda _{n _m}} \times D (0,\eta _m {} ^+).$$
Grâce à \cite[2.2.12]{Berig}, il en dérive que l'isomorphisme $\epsilon _m$
est déterminé, pour toute section $e$ de $E _{n_m}$, par la formule
$\epsilon _m (1 \otimes e ) = \sum _{\underline{k}\geq 0} \underline{\partial} ^{[\underline{k}]} e \otimes
\underline{\tau} ^{ \underline{k}}$.
En passant à la limite sur $m$ l'isomorphisme
$h _K ^* (\epsilon _m)$ : $ f _K ^{\prime *} (E _{n_m}) \riso f _K ^{ *} (E _{n_m}) $
(ou en lui appliquant le foncteur $j ^\dag$),
on obtient $\epsilon _{f,\,f'}$.

Notons
$\smash{\widehat{\PP}} _{\X,(m)} = \underset{\longleftarrow}{\lim} _i \, \PP _{X _i\,(m)}$,
$\smash{\widehat{P}} _{\X,(m)} = \Gamma (\X,\, \smash{\widehat{\PP}} _{\X,(m)})$,
$\smash{\widehat{\PP}} '_{\X,(m)} =
\underset{\longleftarrow}{\lim} _i \,\widehat{\B} ^{(n_m)} _{X _i} (H) \otimes _{\O _{X _i}} \PP _{X _i\,(m)}$
et
$\smash{\widehat{P}} '_{\X,(m)}=\Gamma (\X, \, \smash{\widehat{\PP}} '_{\X,(m)})$.
Désignons par $\tilde{f} _i$ et $\tilde{f} ' _i$ les morphismes d'espaces annelés
$(X' _i , \B ^{(n_m)} _{X ' _i} (H'))\rightarrow (X _i , \B ^{(n_m)} _{X  _i} (H))$
induits par $f _i$ et $f'_i$, et par
$\delta _{\tilde{f} _i, \tilde{f} ' _i} $ :
$(X' _i , \B ^{(n_m)} _{X ' _i} (H'))\rightarrow \widetilde{\Delta}  _{X _i(m)}$
le morphisme déduit (notations de \ref{tauquasinilp}).
En passant à la limite inductive, on obtient le morphisme d'espaces annelés
$\delta _{\tilde{f} , \tilde{f} ' }  $ : $(\X ', \, \widehat{\B} ^{(n_m)} _{\X '} (H'))\rightarrow
\widetilde{\Delta}  _{\X(m)}$ (où $\widetilde{\Delta}  _{\X(m)}$ désigne la limite inductive des
$\widetilde{\Delta}  _{X _i(m)}$, i.e., l'espace topologique est celui des $\widetilde{\Delta}  _{X _i(m)}$
tandis que le faisceau d'anneaux correspond à la limite projective des faisceaux structuraux des
$\widetilde{\Delta}  _{X _i(m)}$)
 s'inscrivant dans le diagramme commutatif :
$$\xymatrix  @R=0,3cm {
{(\X ', \, \widehat{\B} ^{(n_m)} _{\X '} (H')) }
\ar@<1ex>[d] ^-{\tilde{f}  }
\ar@<-1ex>[d] _{\tilde{f}' }
\ar[rd] ^-{\delta _{\tilde{f} , \tilde{f} ' }}
\\
{(\X , \, \widehat{\B} ^{(n_m)} _{\X } (H))}
&
{\widetilde{\Delta}  _{\X(m)}.}
\ar@<1ex>[l] ^-{\tilde{p} _1}
\ar@<-1ex>[l] _-{\tilde{p} _2}
}$$
Il dérive de \ref{espilonhat} que
la structure de $\widehat{\B} ^{(n_m)} _{\X} (H) \widehat{\otimes} \widehat{\D} ^{(m)} _{\X}$-module
topologiquement nilpotent
de $\smash{\overset{_{\circ}}{\E}} ^{(n_m)}$ induit un isomorphisme
$\widehat{\epsilon} _m$ : $\tilde{p} _2 ^* (\smash{\overset{_{\circ}}{\E}} ^{(n_m)}) \riso
\tilde{p} _1 ^* (\smash{\overset{_{\circ}}{\E}} ^{(n_m)})$.
En lui appliquant $\delta _{\tilde{f} , \tilde{f} ' } ^{*}$,
il en découle le suivant :
$\smash{\tilde{f}}  ^{\prime *} (\smash{\overset{_{\circ}}{\E}} ^{(n_m)}) \riso
\smash{\tilde{f}} ^* (\smash{\overset{_{\circ}}{\E}} ^{(n_m)})$.
En tensorisant par $\Q$, puis
en passant à la limite inductive sur $m$, on obtient (modulo le décalage $[d _{X'/X}]$)
l'isomorphisme $ \tau _{\smash{\tilde{f}},\,\smash{\tilde{f}}'}$.

Or, on dispose d'un homomorphisme
canonique d'anneaux $\Gamma ([X ] _{\X ^2 \,\eta _m},\, \O _{[X ] _{\X ^2 \,\eta _{m}}}) \rightarrow P _{\X,(m)\Q}$
(voir la deuxième ligne de la preuve de \cite[3.1.2]{Be0}).
Ce dernier est $A\otimes _{\V} K $-linéaire pour les structures droite et gauche et envoie $\tau _i$ sur $\tau _i$.
Il en résulte un morphisme
$\Gamma (V _{\eta _{m} \,\lambda_{n_m}},\, \O _{V _{\eta _{m} \,\lambda_{n_m}}})\rightarrow P _{\X,(m)\Q}$.
Via cette extension,
les $\epsilon _m$ :
$p _2 ^* (E _{n_m} )\riso p ^* _1  (E _{n_m})$
induisent alors (modulo les foncteurs quasi-inverses $\sp _*$ et $\sp ^*$)
les isomorphismes
$\eta _m$ : $p _2 ^* (\smash{\overset{_{\circ}}{\E}} ^{(n_m)}) _\Q \riso
p _1 ^* (\smash{\overset{_{\circ}}{\E}} ^{(n_m)}) _\Q$.
Puisque $\epsilon _m (1 \otimes e ) =
\sum _{\underline{k}\geq 0} \underline{\partial} ^{[\underline{k}]} e \otimes
\underline{\tau} ^{\underline{k}}$ et que
$\widehat{\epsilon} _{m \Q}$ est caractérisé par la relation
$\widehat{\epsilon} _{m \Q} ( (1\otimes 1) \otimes e ) =
   \sum _{\underline{k}\geq 0} \underline{\smash{\widetilde{\partial}}} ^{<\underline{k}>} e \otimes
   \underline{\smash{\widetilde{\tau}}} ^{\{ \underline{k}\}}$
   (cela découle de \ref{espilonhat}),
$\eta _m =\widehat{\epsilon} _{m \Q}$.
Il en dérive que les deux isomorphismes
$ f _K ^{\prime *} (E _{n_m}) \riso f _K ^{ *} (E _{n_m}) $
et
$\smash{\tilde{f}}  ^{\prime *} (\smash{\overset{_{\circ}}{\E}} ^{(n_m)})_\Q \riso
\smash{\tilde{f}} ^* (\smash{\overset{_{\circ}}{\E}} ^{(n_m)})_\Q$
se correspondent (modulo les foncteurs quasi-inverses $\sp _*$ et $\sp ^*$).
On conclut en passant à  la limite sur $m$.
\end{proof}

\subsection{Construction de $\sp _+$}
D'après une remarque de Berthelot, la procédure de recollement de
la section \cite[2]{caro_unite} est incorrecte.
Avec les notations de \cite[2]{caro_unite},
cela vient du fait que les schémas formels de la forme
$\X _\alpha \times _{\PP} \X _\beta$ ne sont pas plats en général.
Nous donnerons ici une procédure améliorée et nous vérifierons par la suite avec celle-ci
(voir \ref{commsp+}) les résultats de la section \cite[2]{caro_unite}.

\begin{vide}
  \label{notat-construc}
 Nous garderons, sauf
 mention explicite du contraire, les suivantes :
  on se donne $\PP$ un $\V$-schéma formel séparé et lisse, $X$ un sous-schéma
fermé $k$-lisse de $P$ et $T$ un diviseur de $P$ tel que $T _X :=T
\cap X$ soit un diviseur de $X$. On note $\U := \PP \setminus T$,
$Y:= X\setminus T _X$ et $j$ : $Y \hookrightarrow X$ l'immersion ouverte.
On fixe de plus $(\PP _{\alpha}) _{\alpha \in \Lambda}$ un recouvrement d'ouverts de $\PP$.
On note $\PP _{\alpha \beta}:= \PP _\alpha \cap \PP _\beta$,
$\PP _{\alpha \beta \gamma}:= \PP _\alpha \cap \PP _\beta \cap \PP _\gamma$,
$X _\alpha := X \cap P _\alpha$,
$X_{\alpha \beta } := X _\alpha \cap X _\beta$ et
$X_{\alpha \beta \gamma } := X _\alpha \cap X _\beta \cap X _\gamma $.
De plus, on notera $Y _\alpha $ l'ouvert de $X
_\alpha$ complémentaire de $T$,
$Y _{\alpha \beta} := Y _\alpha \cap Y _\beta$,
$Y _{\alpha \beta \gamma} := Y _\alpha \cap Y _\beta \cap Y _\gamma $,
$j _\alpha$ : $ Y _\alpha
\hookrightarrow X _\alpha$,
$j _{\alpha \beta} $ :
$Y _{\alpha \beta}  \hookrightarrow  X _{\alpha \beta}$
et
$j _{\alpha \beta \gamma} $ :
$Y _{\alpha \beta \gamma }  \hookrightarrow  X _{\alpha \beta \gamma} $
les immersions ouvertes canoniques.
On suppose de plus que pour tout $\alpha\in \Lambda$, $X _\alpha$ est affine
(par exemple lorsque le recouvrement $(\PP _{\alpha}) _{\alpha \in \Lambda}$ est affine).
Comme $P$ est  séparé, pour tous $\alpha,\beta ,\gamma \in \Lambda$,
$X_{\alpha \beta }$ et $X_{\alpha \beta \gamma }$ sont donc affines.

Pour tout triplet $(\alpha, \, \beta,\, \gamma)\in \Lambda ^3$, choisissons
$\X _\alpha$ (resp. $\X _{\alpha \beta}$, $\X _{\alpha \beta \gamma}$)
des $\V$-schémas formels lisses relevant $X _\alpha$
(resp. $X _{\alpha \beta}$, $X _{\alpha \beta \gamma}$),
$p _1 ^{\alpha \beta}$ :
$\X  _{\alpha \beta} \rightarrow \X _{\alpha}$
(resp. $p _2 ^{\alpha \beta}$ :
$\X  _{\alpha \beta} \rightarrow \X _{\beta}$)
des relèvements de
$X  _{\alpha \beta} \rightarrow X _{\alpha}$
(resp. $X  _{\alpha \beta} \rightarrow X _{\beta}$).
Rappelons que grâce à Elkik (\cite{elkik} de tels relèvements existent bien.

De même, pour tout triplet $(\alpha,\,\beta,\,\gamma )\in \Lambda ^3$, on choisit des relèvements
$p _{12} ^{\alpha \beta \gamma}$ : $\X  _{\alpha \beta \gamma} \rightarrow \X  _{\alpha \beta} $,
$p _{23} ^{\alpha \beta \gamma}$ : $\X  _{\alpha \beta \gamma} \rightarrow \X  _{\beta \gamma} $,
$p _{13} ^{\alpha \beta \gamma}$ : $\X  _{\alpha \beta \gamma} \rightarrow \X  _{\alpha \gamma} $,
$p _1 ^{\alpha \beta \gamma}$ : $\X  _{\alpha \beta \gamma} \rightarrow \X  _{\alpha} $,
$p _2 ^{\alpha \beta \gamma}$ : $\X  _{\alpha \beta \gamma} \rightarrow \X  _{\beta} $,
$p _3 ^{\alpha \beta \gamma}$ : $\X  _{\alpha \beta \gamma} \rightarrow \X  _{\gamma} $,
$u _{\alpha}$ : $\X _{\alpha } \hookrightarrow \PP _{\alpha }$,
$u _{\alpha \beta}$ : $\X _{\alpha \beta} \hookrightarrow \PP _{\alpha \beta}$
et
$u _{\alpha \beta \gamma}$ : $\X _{\alpha \beta \gamma } \hookrightarrow \PP _{\alpha \beta \gamma}$
induisant les morphismes canoniques au niveau des fibres spéciales.
Sauf mention du contraire, tous ces relèvements seront supposés fixés par la suite.

Via les isomorphismes de la forme $\tau $ (\ref{relev-comp-adj-immf}), on remarque que l'on dispose
des diagrammes commutatifs de foncteurs suivants
\begin{equation}
  \label{diag-transp12p1}
  \xymatrix  @R=0,3cm@C=1,5cm {
  {p _1 ^{\alpha \beta \gamma!}}
  \ar[r] ^-{\tau } _-{\sim}
  \ar@{=}[d]
  &
  {p _{12} ^{\alpha \beta \gamma!} \circ p _1 ^{\alpha \beta !}}
  \ar[d]^\tau _{\sim}
  \\
  {p _1 ^{\alpha \beta \gamma!}}
  \ar[r] ^-{\tau } _-{\sim}
  &
  {p _{13} ^{\alpha \beta \gamma!} \circ p _1 ^{\alpha \gamma !},}
  }
\xymatrix  @R=0,3cm {
  {p _2 ^{\alpha \beta \gamma!}}
  \ar[r] ^-{\tau } _-{\sim}
  \ar@{=}[d]
  &
  {p _{12} ^{\alpha \beta \gamma!} \circ p _2 ^{\alpha \beta !}}
  \ar[d] ^-{\tau } _{\sim}
  \\
  {p _2 ^{\alpha \beta \gamma!}}
  \ar[r] ^-{\tau } _-{\sim}
  &
  {p _{23} ^{\alpha \beta \gamma!} \circ p _1 ^{\beta \gamma !},}
  }
  \xymatrix  @R=0,3cm {
  {p _3 ^{\alpha \beta \gamma!}}
  \ar[r] ^-{\tau } _-{\sim}
  \ar@{=}[d]
  &
  {p _{13} ^{\alpha \beta \gamma!} \circ p _2 ^{\alpha \gamma !}}
  \ar[d] ^-{\tau } _{\sim}
  \\
  {p _3 ^{\alpha \beta \gamma!}}
  \ar[r] ^-{\tau } _-{\sim}
  &
  {p _{23} ^{\alpha \beta \gamma!} \circ p _2 ^{\beta \gamma !}.}
  }
\end{equation}
\end{vide}

\begin{defi}\label{defindonnederecol}
Pour tout $\alpha \in \Lambda$, donnons-nous $\E _\alpha$,
un $\smash{\D} ^{\dag} _{\X _{\alpha} } (\hdag T  \cap X _{\alpha}) _{\Q}$-module cohérent.
On appelle \textit{donnée de recollement} sur $(\E _{\alpha})_{\alpha \in \Lambda}$,
la donnée pour tous $\alpha,\,\beta \in \Lambda$, d'un isomorphisme
$\smash{\D} ^{\dag} _{\X _{\alpha \beta} }(\hdag T  \cap X _{\alpha \beta}) _{ \Q}$-linéaire
$ \theta _{  \alpha \beta} \ : \  p _2  ^{\alpha \beta !} (\E _{\beta}) \riso p  _1 ^{\alpha \beta !} (\E _{\alpha}),$
ceux-ci vérifiant la condition de cocycle :
$\theta _{13} ^{\alpha \beta \gamma }=
\theta _{12} ^{\alpha \beta \gamma }
\circ
\theta _{23} ^{\alpha \beta \gamma }$,
où $\theta _{12} ^{\alpha \beta \gamma }$, $\theta _{23} ^{\alpha \beta \gamma }$
et $\theta _{13} ^{\alpha \beta \gamma }$ sont définis par les diagrammes commutatifs
\begin{equation}
  \label{diag1-defindonnederecol}
\xymatrix  @R=0,3cm {
{  p _{12} ^{\alpha \beta \gamma !} p  _2 ^{\alpha \beta !}  (\E _\beta )}
\ar[r] ^-{\tau} _{\sim}
\ar[d] ^-{p _{12} ^{\alpha \beta \gamma !} (\theta _{\alpha \beta})} _{\sim}
&
{p _2 ^{\alpha \beta \gamma!}  (\E _\beta )}
\ar@{.>}[d] ^-{\theta _{12} ^{\alpha \beta \gamma }}
\\
{ p _{12} ^{\alpha \beta \gamma !}  p  _1 ^{\alpha \beta !}  (\E _\alpha)}
\ar[r]^{\tau} _{\sim}
&
{p _1 ^{\alpha \beta \gamma!}(\E _\alpha),}
}
%
%
\xymatrix  @R=0,3cm {
{  p _{23} ^{\alpha \beta \gamma !} p  _2 ^{\beta \gamma!}  (\E _\gamma )}
\ar[r] ^-{\tau} _{\sim}
\ar[d] ^-{p _{23} ^{\alpha \beta \gamma !} (\theta _{ \beta \gamma})} _{\sim}
&
{p _3 ^{\alpha \beta \gamma!}  (\E _\gamma )}
\ar@{.>}[d] ^-{\theta _{23} ^{\alpha \beta \gamma }}
\\
{ p _{23} ^{\alpha \beta \gamma !}  p  _1 ^{ \beta \gamma !}  (\E _\beta)}
\ar[r]^{\tau} _{\sim}
&
{p _2 ^{\alpha \beta \gamma!}(\E _\beta),}
}
%
%
\xymatrix  @R=0,3cm {
{  p _{13} ^{\alpha \beta \gamma !} p  _2 ^{\alpha \gamma !}  (\E _\gamma )}
\ar[r] ^-{\tau} _{\sim}
\ar[d] ^-{p _{13} ^{\alpha \beta \gamma !} (\theta _{\alpha \gamma})} _{\sim}
&
{p _3 ^{\alpha \beta \gamma!}  (\E _\gamma )}
\ar@{.>}[d]^{\theta _{13} ^{\alpha \beta \gamma }}
\\
{ p _{13} ^{\alpha \beta \gamma !}  p  _1 ^{\alpha \gamma !}  (\E _\alpha)}
\ar[r]^{\tau} _{\sim}
&
{p _1 ^{\alpha \beta \gamma!}(\E _\alpha).}
}
\end{equation}

On construit la catégorie $\mathrm{Coh} (X,\, (\PP _\alpha) _{\alpha \in \Lambda},\, T)$ de la manière
suivante :

- un objet est une famille $(\E _\alpha) _{\alpha \in \Lambda}$
de $\smash{\D} ^{\dag} _{\X _{\alpha} } (\hdag T  \cap X _{\alpha}) _{\Q}$-modules cohérents, $\E _\alpha$,
munie d'une donnée de recollement $ (\theta _{\alpha\beta}) _{\alpha ,\beta \in \Lambda}$,

- un morphisme
$((\E _{\alpha})_{\alpha \in \Lambda},\, (\theta _{\alpha\beta}) _{\alpha ,\beta \in \Lambda})
\rightarrow
((\E ' _{\alpha})_{\alpha \in \Lambda},\, (\theta '_{\alpha\beta}) _{\alpha ,\beta \in \Lambda})$
est une famille de morphismes $f _\alpha$ : $\E _\alpha \rightarrow \E '_\alpha$
commutant aux données de recollement, i.e., telle que le diagramme suivant soit commutatif :
\begin{equation}
  \label{diag2-defindonnederecol}
\xymatrix  @R=0,3cm {
{ p _2  ^{\alpha \beta !} (\E _{\beta}) }
\ar[d] ^-{p _2  ^{\alpha \beta !} (f _{\beta}) }
\ar[r] ^-{\theta _{\alpha\beta}} _{\sim}
&
{  p  _1 ^{\alpha \beta !} (\E _{\alpha}) }
\ar[d] ^-{p  _1 ^{\alpha \beta !} (f _{\alpha})}
\\
{p _2  ^{\alpha \beta !} (\E '_{\beta})  }
\ar[r]^{\theta '_{\alpha\beta}} _{\sim}
&
{ p  _1 ^{\alpha \beta !} (\E '_{\alpha})  .}
}
\end{equation}
Lorsque $T$ est vide, on omettra comme d'habitude de l'indiquer.
\end{defi}
\begin{rema}\label{rem-defindonnederecol}
  Pour tous $\alpha ,\beta \in \Lambda$, soient
$f _\alpha$ : $\E _\alpha \rightarrow \E '_\alpha$ un morphisme
de $\smash{\D} ^{\dag} _{\X _{\alpha} } (\hdag T  \cap X _{\alpha}) _{\Q}$-modules cohérents,
$ \theta _{  \alpha \beta} \ : \  p _2  ^{\alpha \beta !} (\E _{\beta}) \riso p  _1 ^{\alpha \beta !} (\E _{\alpha})$
et
$ \theta '_{  \alpha \beta} \ :
\  p _2  ^{\alpha \beta !} (\E '_{\beta}) \riso p  _1 ^{\alpha \beta !} (\E' _{\alpha})$
des isomorphismes
$\smash{\D} ^{\dag} _{\X _{\alpha \beta} }(\hdag T  \cap X _{\alpha \beta}) _{ \Q}$-linéaires.
On suppose en outre que
les morphismes $f _\alpha$ et les isomorphismes $\theta _{\alpha \beta}$
et $\theta ' _{\alpha \beta}$ induisent le diagramme commutatif \ref{diag2-defindonnederecol}.

Alors, les isomorphismes $\theta _{\alpha\beta}$ satisfont à la condition de cocycle
si et seulement s'il en est de même des isomorphismes $\theta '_{\alpha\beta}$. En effet,
en transformant, via \ref{diag2-defindonnederecol},
les carrés \ref{diag1-defindonnederecol} en trois cubes commutatifs,
on obtient les trois carrés commutatifs suivants :
\begin{equation}
  \notag
  \xymatrix  @R=0,3cm@C=2cm {
{p _2 ^{\alpha \beta \gamma!}  (\E _\beta )}
\ar[d] ^-{\theta _{12} ^{\alpha \beta \gamma }} _{\sim}
\ar[r] ^-{p _2 ^{\alpha \beta \gamma!}  (f _\beta )}
&
{p _2 ^{\alpha \beta \gamma!}  (\E' _\beta )}
\ar[d] ^-{\theta _{12} ^{\prime \alpha \beta \gamma }} _{\sim}
\\
{p _1 ^{\alpha \beta \gamma!}(\E _\alpha)}
\ar[r] ^-{p _1 ^{\alpha \beta \gamma!}(f _\alpha)}
&
{p _1 ^{\alpha \beta \gamma!}(\E '_\alpha),}
}
  \xymatrix  @R=0,3cm@C=2cm {
{p _3 ^{\alpha \beta \gamma!}  (\E _\gamma )}
\ar[d] ^-{\theta _{23} ^{\alpha \beta \gamma }} _{\sim}
\ar[r] ^-{p _3 ^{\alpha \beta \gamma!}  (f _\gamma )}
&
{p _3 ^{\alpha \beta \gamma!}  (\E' _\gamma )}
\ar[d] ^-{\theta _{23} ^{\prime \alpha \beta \gamma }} _{\sim}
\\
{p _2 ^{\alpha \beta \gamma!}(\E _\beta)}
\ar[r] ^-{p _2 ^{\alpha \beta \gamma!}(f _\beta)}
&
{p _2 ^{\alpha \beta \gamma!}(\E '_\beta),}
}
  \xymatrix  @R=0,3cm@C=2cm {
{p _3 ^{\alpha \beta \gamma!}  (\E _\gamma )}
\ar[d] ^-{\theta _{13} ^{\alpha \beta \gamma }} _{\sim}
\ar[r] ^-{p _3 ^{\alpha \beta \gamma!}  (f _\gamma )}
&
{p _3 ^{\alpha \beta \gamma!}  (\E' _\gamma )}
\ar[d] ^-{\theta _{13} ^{\prime \alpha \beta \gamma }} _{\sim}
\\
{p _1 ^{\alpha \beta \gamma!}(\E _\alpha)}
\ar[r] ^-{p _1 ^{\alpha \beta \gamma!}(f _\alpha)}
&
{p _1 ^{\alpha \beta \gamma!}(\E '_\alpha).}
}
\end{equation}
\end{rema}

\begin{prop}\label{prop1}
En notant $\mathrm{Coh} (X,\, \PP,\, T)$,
la catégorie des $\smash{\D} ^{\dag} _{\PP }(\hdag T ) _{\Q}$-modules
cohérents à support dans $X$,
on dispose d'une équivalence de catégories entre
$\mathrm{Coh} (X,\, \PP,\, T)$
et $\mathrm{Coh} (X,\, (\PP _\alpha) _{\alpha \in \Lambda},\, T)$.
\end{prop}
\begin{proof}
Construisons d'abord le foncteur canonique
$\mathcal{L}oc$ : $\mathrm{Coh} (X,\, \PP ,\, T)
\rightarrow
\mathrm{Coh} (X,\, (\PP _\alpha) _{\alpha \in \Lambda},\, T)$.
Pour tout $\smash{\D} ^{\dag} _{\PP , \Q}(\hdag T )$-module cohérent $\E$ à support dans $X$,
on définit l'isomorphisme
$\theta _{\alpha \beta}$ : $ p  _2 ^{\alpha \beta !}  u _{\beta } ^! (\E |_{\PP _\beta})
\riso
 p  _1 ^{\alpha \beta !}  u _{\alpha } ^! (\E |_{\PP _\alpha}) $,
 comme étant l'unique flèche rendant commutatif le diagramme suivant
\begin{equation}
  \label{coh->cohloc}
  \xymatrix  @R=0,3cm {
{ p  _2 ^{\alpha \beta !}  u _{\beta } ^! (\E |_{\PP _\beta}) }
\ar[r]^{\tau} _{\sim}
\ar@{.>}[d] ^-{\theta _{\alpha \beta}}
&
{u ^! _{\alpha \beta} ( (\E |_{\PP _\beta}) |_{\PP _{\alpha \beta}}) }\ar@{=}[d] \\
{ p  _1 ^{\alpha \beta !}  u _{\alpha } ^! (\E |_{\PP _\alpha}) }
\ar[r]^{\tau} _{\sim}
&
{u ^! _{\alpha \beta} ( (\E |_{\PP _\alpha}) |_{\PP _{\alpha \beta}}) .}
}
\end{equation}
Via l'isomorphisme $\tau $ :
$ p _{12} ^{\alpha \beta \gamma !} u ^! _{\alpha \beta} ( (\E |_{\PP _\alpha}) |_{\PP _{\alpha \beta}})
\riso
 u ^! _{\alpha \beta \gamma} ((\E |_{\PP _\beta}) |_{\PP _{\alpha \beta \gamma }})$,
 en appliquant le foncteur $p _{12} ^{\alpha \beta \gamma !} $ au carré \ref{coh->cohloc}
 et avec \ref{diag1-defindonnederecol},
 on obtient le diagramme commutatif :
$$\xymatrix  @R=0,3cm {
{p _2 ^{\alpha \beta \gamma!}  (u _{\beta } ^! (\E |_{\PP _\beta}) )}
\ar[d] ^-{\theta _{12} ^{\alpha \beta \gamma }} _{\sim}
\ar[r] ^\tau _{\sim}
&
{u ^! _{\alpha \beta \gamma} ((\E |_{\PP _\beta}) |_{\PP _{\alpha \beta \gamma }})}
\ar@{=}[d]
\\
{p _1 ^{\alpha \beta \gamma!}( u _{\alpha } ^! (\E |_{\PP _\alpha}))}
\ar[r] ^\tau _{\sim}
&
{u ^! _{\alpha \beta \gamma} ((\E |_{\PP _\alpha}) |_{\PP _{\alpha \beta \gamma }}),}
}
$$
où les isomorphismes horizontaux sont de la forme $\tau$, grâce à la formule de transitivité et
à la commutation aux images inverses extraordinaires des isomorphismes de la forme $\tau$ (\ref{relev-comp-adj-immf}).
De même, on construit les deux autres diagrammes analogues.
Avec ces trois diagrammes, on vérifie que
$\mathcal{L}oc (\E) :=
((u _{\alpha } ^! (\E |_{\PP _\alpha})) _{\alpha \in \Lambda},\, (\theta _{\alpha \beta})_{\alpha , \beta \in \Lambda})$
satisfait à la condition de cocycle et est ainsi un objet de
$\mathrm{Coh} (X,\, (\PP _\alpha) _{\alpha \in \Lambda},\, T)$.

En outre, si $f$ : $\E \rightarrow \E'$ est un morphisme
de $\smash{\D} ^{\dag} _{\PP }(\hdag T ) _{\Q}$-modules
cohérents à support dans $X$, alors, par fonctorialité en $\E$
de \ref{coh->cohloc} (on transforme le carré \ref{coh->cohloc} en cube),
la famille $(u _\alpha ^! (f|_{\PP _\alpha})) _{\alpha \in \Lambda}$ commute aux données de recollement.

\medskip

Construisons à présent un foncteur quasi-inverse canonique $\mathcal{R}ecol$ :
$\mathrm{Coh} (X,\, (\PP _\alpha) _{\alpha \in \Lambda},\, T)
\rightarrow
\mathrm{Coh} (X,\, \PP ,\, T)$.

Soit $(\E _{\alpha}) _{\alpha \in \Lambda}$
une famille de $\smash{\D} ^{\dag} _{\X _{\alpha} , \Q}(\hdag T  \cap X _{\alpha})$-modules cohérents
munie d'une donnée de recollement
$(\theta _{\alpha \beta})_{\alpha , \beta \in \Lambda}$.
Prouvons que la famille $( u _{\alpha +} ( \E _{\alpha} ))_{\alpha \in \Lambda}$
se recolle
en un $\smash{\D} ^{\dag} _{\PP , \Q}(\hdag T )$-module cohérent
à support dans $X$.

Pour cela, notons $\phi _1 ^{\alpha \beta}$ (resp. $\phi _2 ^{\alpha \beta}$) le morphisme d'adjonction
(voir \ref{comp-comp-adj-immf})
du carré de gauche (resp. de droite)
\begin{equation}\label{carre0eqcatrecolD-mod}
\xymatrix  @R=0,3cm {
{\PP  _{\alpha \beta}} \ar[rr]
&
&{ \PP  _{\alpha}}
\\
{\X  _{\alpha \beta}} \ar[rr] ^-{p _1 ^{\alpha \beta}} \ar[u] ^-{u _{\alpha \beta}}
&
&
{ \X _{\alpha},} \ar[u] ^-{u _{\alpha}}
}
\hspace{3cm}
\xymatrix  @R=0,3cm{
{\PP  _{\alpha \beta}} \ar[rr] &  &{ \PP  _{\beta}} \\
{\X  _{\alpha \beta}} \ar[rr] ^-{p _2 ^{\alpha \beta}} \ar[u] ^-{u _{\alpha \beta}}
&  &{ \X _{\beta}.} \ar[u] ^-{u _{\beta}}
}
\end{equation}
On construit pour tous $\alpha ,\beta \in \Lambda$ un isomorphisme
$\tau _{\alpha \beta}\ : \  {(u _{\beta +} ( \E _{\beta} ))} |_{\PP _{\alpha \beta}}
\riso {(u _{\alpha +} ( \E _{\alpha} ))} |_{\PP _{\alpha \beta}}$
comme
étant l'unique morphisme rendant commutatif le diagramme suivant
\begin{equation}\label{carre1eqcatrecolD-mod}
\xymatrix  @R=0,3cm {
{u _{\alpha \beta +} \circ p  _1 ^{\alpha \beta !} (\E _{\alpha})}
\ar[rr] ^-{\phi _1 ^{\alpha \beta}(\E _{\alpha})} _{\sim}
& &
 {{(u _{\alpha +} ( \E _{\alpha} ))} |_{ \PP _{\alpha \beta}}   } \\
{ u _{\alpha \beta +} \circ p  _2 ^{\alpha \beta !} (\E _{\beta}) }
 \ar[rr] ^-{\phi _2 ^{\alpha \beta}(\E _{\beta})} _{\sim}
 \ar[u] ^-{ u _{\alpha \beta +}(\theta _{  \alpha \beta}) } _{\sim}
 &  &
{ {(u _{\beta +} ( \E _{\beta} )) }|_{\PP _{\alpha \beta}}. }
\ar@{.>}[u] _{\tau _{\alpha \beta}}
}
\end{equation}
Il reste maintenant à établir que les isomorphismes $\tau _{\alpha \beta}$ vérifient la condition de
recollement.
À cette fin,
notons $\phi _{12} ^{\alpha \beta \gamma}$ (resp. $\phi _{23} ^{\alpha \beta \gamma}$
et $\phi _{13} ^{\alpha \beta \gamma}$)
le morphisme d'adjonction (toujours \ref{comp-comp-adj-immf})
du carré de gauche (resp. du centre et de droite) suivant
\begin{equation}\label{carre2eqcatrecolD-mod}
\xymatrix  @R=0,3cm {
{\PP _{\alpha \beta \gamma } }         \ar[rr]
&
&
{\PP _{\alpha \beta } }
\\
{\X _{\alpha \beta \gamma } }
    \ar[rr] ^-{p _{12} ^{\alpha \beta \gamma}} \ar[u] ^-{u _{\alpha \beta \gamma}}
&
&
{\X _{\alpha \beta },} \ar[u] ^-{u _{\alpha \beta }}
}
\xymatrix  @R=0,3cm {
{\PP _{\alpha \beta \gamma } }         \ar[rr]
&
&
{\PP _{\beta \gamma} }
\\
{\X _{\alpha \beta \gamma } }
    \ar[rr] ^-{p _{23} ^{\alpha \beta \gamma}} \ar[u] ^-{u _{\alpha \beta \gamma}}
&
&
{\X _{ \beta \gamma},} \ar[u] ^-{u _{\beta \gamma}}
}
\xymatrix  @R=0,3cm {
{\PP _{\alpha \beta \gamma } }         \ar[rr]
&
&
{\PP _{\alpha \gamma} }
\\
{\X _{\alpha \beta \gamma } }
    \ar[rr] ^-{p _{13} ^{\alpha \beta \gamma}} \ar[u] ^-{u _{\alpha \beta \gamma}}
&
&
{\X _{ \alpha \gamma},} \ar[u] ^-{u _{\alpha \gamma}}
}
\end{equation}
et $\phi _{1} ^{\alpha \beta \gamma}$
(resp. $\phi _{2} ^{\alpha \beta \gamma}$ et $\phi _{3} ^{\alpha \beta \gamma}$)
celui du diagramme
de gauche (resp. du centre et de droite) :
\begin{equation}\label{carre4eqcatrecolD-mod}
\xymatrix  @R=0,3cm {
{\PP _{\alpha \beta \gamma } }         \ar[rrr]
&
&
&
{\PP _{\alpha } }
\\
{\X _{\alpha \beta \gamma } }
    \ar[rrr]^{p  _{1} ^{\alpha \beta \gamma}}  \ar[u] ^-{u _{\alpha \beta \gamma}}
&
&
&
{\X _{\alpha }, } \ar[u] ^-{u _{\alpha  }}
}
\xymatrix  @R=0,3cm {
{\PP _{\alpha \beta \gamma } }         \ar[rrr]
&
&
&
{\PP _{\beta } }
\\
{\X _{\alpha \beta \gamma } }
    \ar[rrr]^{p  _{2} ^{\alpha \beta \gamma}}  \ar[u] ^-{u _{\alpha \beta \gamma}}
&
&
&
{\X _{\beta }, } \ar[u] ^-{u _{\beta  }}
}
\xymatrix  @R=0,3cm {
{\PP _{\alpha \beta \gamma } }         \ar[rrr]
&
&
&
{\PP _{\gamma } }
\\
{\X _{\alpha \beta \gamma } }
    \ar[rrr]^{p  _{3} ^{\alpha \beta \gamma}}  \ar[u] ^-{u _{\alpha \beta \gamma}}
&
&
&
{\X _{\gamma }. } \ar[u] ^-{u _{\gamma  }}
}
\end{equation}

Considérons le diagramme commutatif suivant
\begin{equation}\label{carre3eqcatrecolD-mod}
\xymatrix  @R=0,3cm@C=2cm {
{u _{\alpha \beta \gamma +}p _1 ^{\alpha \beta \gamma!}(\E _\alpha)}
\ar[r] ^-{u _{\alpha \beta \gamma +} (\tau)} _{\sim}
&
{u _{\alpha \beta \gamma +}\circ p _{12} ^{\alpha \beta \gamma !} ( p  _1 ^{\alpha \beta !} (\E _{\alpha}))}
\ar[r] ^-{\phi _{12} ^{\alpha \beta \gamma}( p  _1 ^{\alpha \beta !} (\E _{\alpha})  )}
_{\sim}
&
{u _{\alpha \beta +} ( p  _1 ^{\alpha \beta !} (\E _{\alpha}) )|_{ \PP _{\alpha \beta \gamma} }}
\\
{u _{\alpha \beta \gamma +}p _2 ^{\alpha \beta \gamma!}  (\E _\beta )}
\ar[u]  ^{u _{\alpha \beta \gamma +} (\theta _{12} ^{\alpha \beta \gamma })}
_{\sim}
\ar[r]^{u _{\alpha \beta \gamma +} (\tau)} _{\sim}
&
{u _{\alpha \beta \gamma +}\circ p _{12} ^{\alpha \beta \gamma !} ( p  _2 ^{\alpha \beta !} (\E _{\beta})) }
\ar[r]  ^{\phi _{12} ^{\alpha \beta \gamma}( p  _2 ^{\alpha \beta !} (\E _{\beta}) )}
_{\sim}
\ar[u] ^-{u _{\alpha \beta \gamma +}\circ p _{12} ^{\alpha \beta \gamma !} ( \theta _{\alpha \beta } ) }
_{\sim}
&
{ u _{\alpha \beta +} ( p  _2 ^{\alpha \beta !} (\E _{\beta}))|_{ \PP _{\alpha \beta \gamma} }.}
\ar[u] _{ u _{\alpha \beta +}(\theta _{  \alpha \beta}) |_{ \PP _{\alpha \beta \gamma} }}
^{\sim}
}
\end{equation}
Grâce à \ref{comp-comp-adj-immf}.i) et \ref{comp-comp-adj-immf}.ii), on obtient alors
\begin{gather}
(\phi _1 ^{\alpha \beta}(\E _{\alpha}) |_{ \PP _{\alpha \beta \gamma} } )\circ
\phi _{12} ^{\alpha \beta \gamma}( p  _1 ^{\alpha \beta !} (\E _{\alpha})  )
\circ u _{\alpha \beta \gamma +} (\tau)=
\phi _{1} ^{\alpha \beta \gamma} (\E _{\alpha}), \label{gather1eqcatrecolD-mod}\\
(\phi _2 ^{\alpha \beta}(\E _{\beta}) |_{ \PP _{\alpha \beta \gamma} } )\circ
\phi _{12} ^{\alpha \beta \gamma}( p  _2 ^{\alpha \beta !} (\E _{\beta})  )
\circ u _{\alpha \beta \gamma +} (\tau)=
\phi _{2} ^{\alpha \beta \gamma} (\E _{\beta}), \label{gather2eqcatrecolD-mod}\\
(\phi _1 ^{\beta \gamma}(\E _{\beta}) |_{ \PP _{\alpha \beta \gamma} } )\circ
\phi _{23} ^{\alpha \beta \gamma}( p  _1 ^{\beta \gamma!} (\E _{\beta})  )
\circ u _{\alpha \beta \gamma +} (\tau)=
\phi _{2} ^{\alpha \beta \gamma} (\E _{\beta}), \label{gather3eqcatrecolD-mod}\\
(\phi _2 ^{ \beta \gamma}(\E _{\gamma}) |_{ \PP _{\alpha \beta \gamma} } )\circ
\phi _{23} ^{\alpha \beta \gamma}( p  _2 ^{ \beta \gamma !} (\E _{\gamma})  )
\circ u _{\alpha \beta \gamma +} (\tau)=
\phi _{3} ^{\alpha \beta \gamma} (\E _{\gamma}), \label{gather4eqcatrecolD-mod}\\
(\phi _1 ^{\alpha \gamma}(\E _{\alpha}) |_{ \PP _{\alpha \beta \gamma} } )\circ
\phi _{13} ^{\alpha \beta \gamma}( p  _1 ^{\alpha \gamma !} (\E _{\alpha})  )
\circ u _{\alpha \beta \gamma +} (\tau)=
\phi _{1} ^{\alpha \beta \gamma} (\E _{\alpha}), \label{gather5eqcatrecolD-mod}\\
(\phi _2 ^{ \alpha \gamma}(\E _{\gamma}) |_{ \PP _{\alpha \beta \gamma} } )\circ
\phi _{13} ^{\alpha \beta \gamma}( p  _2 ^{ \alpha \gamma!} (\E _{\gamma})  )
\circ u _{\alpha \beta \gamma +} (\tau)=
\phi _{3} ^{\alpha \beta \gamma} (\E _{\gamma}).\label{gather6eqcatrecolD-mod}
\end{gather}

En composant \ref{carre1eqcatrecolD-mod} restreint à $\PP _{\alpha \beta \gamma}$
et
\ref{carre3eqcatrecolD-mod}, via les égalités
\ref{gather1eqcatrecolD-mod} et \ref{gather2eqcatrecolD-mod},
on obtient le carré commutatif :
\begin{equation}\label{carre5eqcatrecolD-mod}
\xymatrix  @R=0,3cm@C=2cm {
{u _{\alpha \beta \gamma +}p _1 ^{\alpha \beta \gamma!}(\E _\alpha)}
\ar[r] ^-{\phi _{1} ^{\alpha \beta \gamma}( \E _{\alpha})  }
_{\sim}
&
{(u _{\alpha +} ( \E _{\alpha} )) |_{ \PP _{\alpha \beta \gamma} }   }
\\
{u _{\alpha \beta \gamma +}p _2 ^{\alpha \beta \gamma!}  (\E _\beta )}
\ar[r]  ^{\phi _{2} ^{\alpha \beta \gamma}(\E _{\beta} )}
_{\sim}
\ar[u]^{u _{\alpha \beta \gamma +} (\theta _{12} ^{\alpha \beta \gamma })}
_{\sim}
&
{ (u _{\beta +} ( \E _{\beta} )) |_{ \PP _{\alpha \beta \gamma} }. }
\ar[u] _{ {\tau _{\alpha \beta }} |_{ \PP _{\alpha \beta \gamma} }}
^{\sim}
}
\end{equation}
De façon analogue, en utilisant \ref{gather3eqcatrecolD-mod} et
\ref{gather4eqcatrecolD-mod} (resp. \ref{gather5eqcatrecolD-mod} et
\ref{gather6eqcatrecolD-mod}) on obtient les diagrammes commutatifs suivants :
\begin{equation}\label{carre6eqcatrecolD-mod}
\xymatrix  @R=0,3cm@C=1,5cm {
{u _{\alpha \beta \gamma +}\circ p _{2} ^{\alpha \beta \gamma !} (\E _{\beta})}
\ar[r] ^-{\phi _{2} ^{\alpha \beta \gamma}( \E _{\beta})  }
_{\sim}
&
{(u _{\beta +} ( \E _{\beta} )) |_{ \PP _{\alpha \beta \gamma} }   }
\\
{u _{\alpha \beta \gamma +}\circ p _{3} ^{\alpha \beta \gamma !} (\E _{\gamma}) }
\ar[r]  ^{\phi _{3} ^{\alpha \beta \gamma}(\E _{\gamma} )}
_{\sim}
\ar[u] ^-{u _{\alpha \beta \gamma +} (\theta _{23} ^{\alpha \beta \gamma })}
_{\sim}
&
{ (u _{\gamma +} ( \E _{\gamma} )) |_{ \PP _{\alpha \beta \gamma} },}
\ar[u] _{ {\tau _{ \beta \gamma }} |_{ \PP _{\alpha \beta \gamma} }}
^{\sim}
}
\xymatrix  @R=0,3cm@C=1,5cm {
{u _{\alpha \beta \gamma +}\circ p _{1} ^{\alpha \beta \gamma !} (\E _{\alpha})}
\ar[r] ^-{\phi _{1} ^{\alpha \beta \gamma}( \E _{\alpha})  }
_{\sim}
&
{(u _{\alpha +} ( \E _{\alpha} )) |_{ \PP _{\alpha \beta \gamma} }   }
\\
{u _{\alpha \beta \gamma +}\circ p _{3} ^{\alpha \beta \gamma !} (\E _{\gamma}) }
\ar[r]  ^{\phi _{3} ^{\alpha \beta \gamma}(\E _{\gamma} )}
_{\sim}
\ar[u] ^-{u _{\alpha \beta \gamma +}(\theta _{13} ^{\alpha \beta \gamma } )}
_{\sim}
&
{ (u _{\gamma +} ( \E _{\gamma} )) |_{ \PP _{\alpha \beta \gamma} }. }
\ar[u] _{{ \tau _{\alpha \gamma }} |_{ \PP _{\alpha \beta \gamma} }}
^{\sim}
}
\end{equation}
De ces trois derniers diagrammes, comme le foncteur $u _{\alpha \beta \gamma +}$
est (pleinement) fidèle (pour les modules cohérents), il en dérive que
les isomorphismes $\theta _{\alpha \beta}$ vérifient la condition de
cocycle si et seulement si les isomorphismes $\tau _{\alpha
\beta}$ se recollent.
\medskip

Soit
$f = (f _\alpha)_{\alpha \in \Lambda}$ :
$((\E _{\alpha}) _{\alpha \in \Lambda}, (\theta _{\alpha \beta }) _{\alpha, \beta \in \Lambda})
\rightarrow
((\E ' _{\alpha}) _{\alpha \in \Lambda}, (\theta '_{\alpha \beta }) _{\alpha, \beta \in \Lambda})$,
un morphisme de
$\mathrm{Coh} (X,\, (\PP _\alpha) _{\alpha \in \Lambda},\, T)$.
On lui associe la famille
$\mathcal{R}ecol (f)$ :
$(u _{\alpha +} (f _\alpha)) _{\alpha \in \Lambda}$.
En notant $\tau _{\alpha \beta}$ (resp. $\tau '_{\alpha \beta}$) l'isomorphisme
rendant commutatif \ref{carre1eqcatrecolD-mod} pour $\theta _{\alpha \beta}$
(resp. $\theta '_{\alpha \beta}$),
on obtient le cube
\begin{equation}\label{carre7eqcatrecolD-mod}
  \xymatrix  @R=0,3cm@C=2cm  {
   &
   {u _{\alpha \beta +} \circ p  _1 ^{\alpha \beta !} (\E '_{\alpha})}
        \ar[rr] ^-{\phi _1 ^{\alpha \beta}(\E '_{\alpha})}
   & &
   {{(u _{\alpha +} ( \E ' _{\alpha} ))} |_{ \PP _{\alpha \beta}}   }
   \\
   {u _{\alpha \beta +} \circ p  _1 ^{\alpha \beta !} (\E _{\alpha})}
       \ar[rr] ^(0.6){\phi _1 ^{\alpha \beta}(\E _{\alpha})}
       \ar[ur] ^-{{u _{\alpha \beta +} \circ p  _1 ^{\alpha \beta !} (f _{\alpha})}}
       & &
   {{(u _{\alpha +} ( \E _{\alpha} ))} |_{ \PP _{\alpha \beta}}   }
   \ar[ur] _{{{(u _{\alpha +} ( f _{\alpha} ))} |_{ \PP _{\alpha \beta}}   }}
   \\
   & { u _{\alpha \beta +} \circ p  _2 ^{\alpha \beta !} (\E '_{\beta}) }
   \ar'[u]^(0.7){ u _{\alpha \beta +}(\theta '_{  \alpha \beta}) }[uu]
   \ar'[r]_(0.8){\phi _2 ^{\alpha \beta}(\E '_{\beta})}[rr]
   & &
   { {(u _{\beta +} ( \E '_{\beta} )) }|_{\PP _{\alpha \beta}} }
    \ar[uu] _{\tau '_{\alpha \beta}}
   \\
    { u _{\alpha \beta +} \circ p  _2 ^{\alpha \beta !} (\E _{\beta}) }
       \ar[rr] _{\phi _2 ^{\alpha \beta}(\E _{\beta})}
       \ar[uu] ^-{ u _{\alpha \beta +}(\theta _{  \alpha \beta}) }
       \ar[ur] ^-{ u _{\alpha \beta +} \circ p  _2 ^{\alpha \beta !} (f _{\beta}) }
       & &
    { {(u _{\beta +} ( \E _{\beta} )) }|_{\PP _{\alpha \beta}} ,}
    \ar[uu] _(0.4){\tau _{\alpha \beta}}
    \ar[ur] _{ {(u _{\beta +} ( f _{\beta} )) }|_{\PP _{\alpha \beta}} }
}
\end{equation}
dont les carrés de devant, de derrière, du bas et du haut sont commutatifs par fonctorialité ou grâce à
\ref{carre1eqcatrecolD-mod}. Comme celui de gauche l'est (via \ref{diag2-defindonnederecol}),
il en résulte qu'il en est de même du carré de droite. Les morphismes $u _{\alpha +} (f _\alpha)$ se recollent donc.

On a donc construit le foncteur $\mathcal{R}ecol$.
Prouvons maintenant que celui-ci est quasi-inverse de $\mathcal{L}oc$.
Soit $((\E _{\alpha}) _{\alpha \in \Lambda}, (\theta _{\alpha \beta }) _{\alpha, \beta \in \Lambda})$
un objet de $\mathrm{Coh} (X,\, (\PP _\alpha) _{\alpha \in \Lambda},\, T)$.
Dans un premier temps, il s'agit d'établir un isomorphisme fonctoriel
$\mathcal{L}oc \circ \mathcal{R}ecol
((\E _{\alpha}) _{\alpha \in \Lambda}, (\theta _{\alpha \beta }) _{\alpha, \beta \in \Lambda})
\riso
((\E _{\alpha}) _{\alpha \in \Lambda}, (\theta _{\alpha \beta }) _{\alpha, \beta \in \Lambda})$.

Notons
$\E : =
\mathcal{R} ecol ((\E _{\alpha}) _{\alpha \in \Lambda}, (\theta _{\alpha \beta }) _{\alpha, \beta \in \Lambda})$,
$\mathcal{L}oc (\E) = ((u _{\alpha} ^! ( \E |_{\PP _{\alpha}}) )
_{\alpha \in \Lambda}, (\theta '' _{\alpha \beta }) _{\alpha, \beta \in \Lambda})$ et
$ \tau _\alpha$ : $u _{\alpha +} (\E _\alpha) \riso \E |_{\PP _\alpha}$ les
isomorphismes canoniques vérifiant
$\tau _{\alpha \beta } = \tau _\alpha ^{-1}|_{\PP _{\alpha \beta}} \circ \tau _\beta |_{\PP _{\alpha \beta}}$,
où $\tau _{\alpha \beta}$
a été défini via \ref{carre1eqcatrecolD-mod}.
Considérons le diagramme suivant
\begin{equation}\label{carre8eqcatrecolD-mod}
  \xymatrix  @R=0,3cm@C=1,5cm  {
   &
   { p  _1 ^{\alpha \beta !}\circ u^! _\alpha ( \E |_{\PP _\alpha}) }
        \ar[rr] ^-{\tau}
   & &
   { u ^! _{\alpha \beta } ((\E |_{\PP _\alpha}) |_{ \PP _{\alpha \beta} })}
   \\
   { p  _1 ^{\alpha \beta !}\circ u^! _\alpha ( u _{\alpha +}  (\E _{\alpha})) }
       \ar[rr] ^(0.45){\tau}
       \ar[ur] ^-{ p  _1 ^{\alpha \beta !}\circ u^! _\alpha ( \tau _{\alpha}) }
       & &
   { u ^! _{\alpha \beta } (u _{\alpha +}  (\E _{\alpha}) |_{ \PP _{\alpha \beta} })}
   \ar[ur] _{ u ^! _{\alpha \beta } ({\tau _\alpha} |_{ \PP _{\alpha \beta} })}
   \\
   &
   { p  _2 ^{\alpha \beta !}\circ u^! _\beta ( \E |_{\PP _\beta}) }
        \ar'[r]^{\tau}[rr]
        \ar'[u] ^-{\theta '' _{\alpha \beta}}[uu]
   & &
   { u ^! _{\alpha \beta } ((\E |_{\PP _\beta}) |_{ \PP _{\alpha \beta} })}
   \ar@{=}[uu]
   \\
   { p  _2 ^{\alpha \beta !}\circ u^! _\beta ( u _{\beta +}  (\E _{\beta})) }
       \ar[rr] ^-{\tau}
       \ar[ur] ^-{ p  _2 ^{\alpha \beta !}\circ u^! _\beta ( \tau _{\beta}) }
       \ar@{.>}[uu] ^-{\theta ' _{\alpha \beta}}
       & &
   { u ^! _{\alpha \beta } (u _{\beta +}  (\E _{\beta}) |_{ \PP _{\alpha \beta} }),}
   \ar[ur] _{ u ^! _{\alpha \beta } ({\tau _\beta} |_{ \PP _{\alpha \beta} })}
   \ar[uu] ^(0.4){ u ^! _{\alpha \beta } (\tau _{\alpha \beta})}
}
\end{equation}
où la flèche $\theta '_{\alpha \beta}$ est par définition celle rendant commutatif le carré
de devant. Les carrés du fond et de droite sont commutatifs par définition et ceux du haut et du bas
le sont par fonctorialité.
Grâce la remarque \ref{rem-defindonnederecol}, on en déduit que
les isomorphismes $\theta '_{\alpha \beta}$ vérifie la condition de cocycle et
on se ramène à prouver que l'isomorphisme d'adjonction
$\E _\alpha \riso  u^! _\alpha \circ u _{\alpha +} (\E _\alpha)$
est compatible aux données de recollement respectives, i.e., que le carré de gauche suivant
\begin{equation}
  \label{locrecol<->id}
  \xymatrix  @R=0,3cm@C=2cm {
  { p  _1 ^{\alpha \beta !} (\E _{\alpha}) }
  \ar[r] ^(0.4){\mathrm{adj} _{u _\alpha}} _(0.4){\sim}
  &
  { p  _1 ^{\alpha \beta !}\circ u^! _\alpha \circ u _{\alpha +}  (\E _{\alpha}) }
  \ar[r] ^-{\tau} _{\sim}
  &
  { u ^! _{\alpha \beta } (u _{\alpha +}  (\E _{\alpha}) |_{ \PP _{\alpha \beta} })}
  \\
  {p  _2 ^{\alpha \beta !} (\E _{\beta})   }
  \ar[r]  ^(0.4){\mathrm{adj} _{u _\beta}} _(0.4){\sim}
  \ar[u] ^-{\theta _{\alpha \beta}} _{\sim}
  &
  { p  _2 ^{\alpha \beta !} \circ u^! _\beta \circ u _{\beta +} (\E _{\beta})  }
  \ar[r]  ^{\tau} _{\sim}
  \ar[u] ^-{\theta '_{\alpha \beta}} _{\sim}
  &
  {u ^! _{\alpha \beta } (u _{\beta +}  (\E _{\beta}) |_{ \PP _{\alpha \beta} })  }
  \ar[u] ^-{u ^! _{\alpha \beta }(\tau _{\alpha \beta} )  } _{\sim}
  }
\end{equation}
est commutatif. Or, en appliquant le foncteur
$u_{\alpha \beta +}$ au diagramme \ref{locrecol<->id} et en composant avec le diagramme commutatif
$$\xymatrix  @R=0,3cm@C=2cm {
 {u_{\alpha \beta +} u ^! _{\alpha \beta } (u _{\alpha +}  (\E _{\alpha}) |_{ \PP _{\alpha \beta} })}
 \ar[r]  ^(0.6){\mathrm{adj} _{u _{\alpha \beta}}} _(0.6){\sim}
 &
 {u _{\alpha +}  (\E _{\alpha}) |_{ \PP _{\alpha \beta} }}
 \\
 {u_{\alpha \beta +} u ^! _{\alpha \beta } (u _{\beta +}  (\E _{\beta}) |_{ \PP _{\alpha \beta} })}
 \ar[u]^{u_{\alpha \beta +} u ^! _{\alpha \beta } (\tau _{\alpha \beta})} _{\sim}
 \ar[r]   ^(0.6){\mathrm{adj} _{u _{\alpha \beta}}} _(0.6){\sim}
 &
 {u _{\beta +}  (\E _{\beta}) |_{ \PP _{\alpha \beta} },}
 \ar[u] ^-{\tau _{\alpha \beta}} _{\sim}
 }
 $$
on obtient \ref{carre1eqcatrecolD-mod} (par construction de la flèche d'adjonction de la proposition
\ref{comp-comp-adj-immf}), qui est commutatif. Le foncteur $u _{\alpha \beta +}$ étant fidèle,
on démontre ainsi la commutativité du contour de \ref{locrecol<->id}.
Puisque le carré de droite de \ref{locrecol<->id} est commutatif (il correspond au carré de devant
de \ref{carre8eqcatrecolD-mod}), il en découle celle du carré gauche.

Réciproquement, si $\E$ est un
$\smash{\D} ^{\dag} _{\PP }(\hdag T ) _{\Q}$-module cohérent à support dans $X$,
vérifions que l'on dispose d'un isomorphisme
$\mathcal{R}ecol \circ \mathcal{L}oc (\E) \riso \E$ fonctoriel en $\E$.
Notons $( \theta _{\alpha \beta} )_{\alpha \beta \in \Lambda}$, la donnée de recollement
de $(u ^! _\alpha (\E |_{\PP _\alpha} ))_{\alpha\in \Lambda}$ définie dans \ref{coh->cohloc}
et $( \tau _{\alpha \beta} )_{\alpha \beta \in \Lambda}$, la donnée de recollement
de $(u _{\alpha +} u ^! _\alpha (\E |_{\PP _\alpha} ))_{\alpha\in \Lambda}$
déduite de celle de $(u ^! _\alpha (\E |_{\PP _\alpha} ))_{\alpha\in \Lambda}$ via
\ref{carre1eqcatrecolD-mod}.
Prouvons maintenant que les isomorphismes d'adjonction
$u _{\alpha +} u ^! _\alpha (\E |_{\PP _\alpha}) \rightarrow
\E |_{\PP _\alpha}$ sont compatibles aux données de recollement respectives, i.e.,
que le carré de droite suivant (la commutativité des deux autres est tautologique)
\begin{equation}
  \label{recolloc<->id1}
\xymatrix  @R=0,3cm {
{u _{\alpha \beta +} \circ u _{\alpha \beta } ^! (\E |_{ \PP _{\alpha \beta }})}
\ar[r] ^-{u _{\alpha \beta +} (\tau)} _{\sim}
&
{u _{\alpha \beta +} \circ p  _1 ^{\alpha \beta !}( u^! _\alpha (\E |_{ \PP _\alpha}))}
\ar[rr] ^-{\phi _1 ^{\alpha \beta} (u^! _\alpha (\E |_{ \PP _\alpha}) )}
_{\sim}
&
&
{(u _{\alpha +}  u^! _\alpha (\E |_{ \PP _{\alpha }}))   |_{ \PP _{\alpha \beta }}}
\ar[r] ^(0.7){\mathrm{adj _{u _\alpha}} |_{ \PP _{\alpha \beta }}} _(0.7){\sim}
&
{\E  |_{ \PP _{\alpha \beta }}}
\\
{u _{\alpha \beta +} \circ u _{\alpha \beta } ^! (\E |_{ \PP _{\alpha \beta }})}
\ar[r] ^-{u _{\alpha \beta +} (\tau)} _{\sim}
\ar@{=}[u]
&
{u _{\alpha \beta +} \circ p  _2 ^{\alpha \beta !}(u^! _\beta (\E |_{ \PP _\beta}))}
\ar[rr]^{\phi _2 ^{\alpha \beta} (u^! _\beta (\E |_{ \PP _\beta}) )}
_{\sim}
\ar[u] ^-{u _{\alpha \beta +} \theta _{\alpha \beta}}
_{\sim}
&
&
{(u _{\beta +} u^! _\beta (\E |_{ \PP _{\beta }}))   |_{ \PP _{\alpha \beta }}}
\ar[r] ^(0.7){\mathrm{adj _{u _\beta}} |_{ \PP _{\alpha \beta }}} _(0.7){\sim}
\ar[u] ^-{\tau _{\alpha \beta} } _{\sim}
&
{\E  |_{ \PP _{\alpha \beta }}}
\ar@{=}[u]
}
\end{equation}
est commutatif. Or, on dispose du diagramme commutatif suivant :
\begin{equation}
  \label{recolloc<->id2}
  \xymatrix  @R=0,3cm {
  {u _{\alpha \beta +}  p  _1 ^{\alpha \beta !} u^! _\alpha (\E |_{ \PP _\alpha})}
  \ar[r]^(.4){\mathrm{adj _{u _\alpha}}}_(0.4){\sim}
  \ar@{=}[d]
  &
  {u _{\alpha \beta +}  p  _1 ^{\alpha \beta !} u^! _\alpha  u _{\alpha +} u^! _\alpha (\E |_{ \PP _\alpha})}
  \ar[r]^{u _{\alpha \beta +} \tau} _{\sim}
  \ar[d]^{\mathrm{adj _{u _\alpha}}} _{\sim}
  &
  {u _{\alpha \beta +}  u _{\alpha \beta } ^! (u _{\alpha +} u^! _\alpha (\E |_{ \PP _{\alpha }}))
  |_{ \PP _{\alpha \beta }}}
  \ar[d]^{\mathrm{adj _{u _\alpha}}} _{\sim}
  \ar[r] ^(0.6){\mathrm{adj _{u _{\alpha \beta }}}} _(0.6){\sim}
  &
  {u _{\alpha +} u^! _\alpha (\E |_{ \PP _{\alpha }})  |_{ \PP _{\alpha \beta }}}
  \ar[d] ^-{\mathrm{adj} _{u_\alpha}} _{\sim}
  \\
  {u _{\alpha \beta +}  p  _1 ^{\alpha \beta !} u^! _\alpha (\E |_{ \PP _\alpha})}
  \ar@{=}[r]
  &
  {u _{\alpha \beta +}  p  _1 ^{\alpha \beta !} u^! _\alpha (\E |_{ \PP _\alpha})}
  \ar[r]^{u _{\alpha \beta +} \tau} _{\sim}
  &
  {u _{\alpha \beta +}  u _{\alpha \beta } ^! (\E |_{ \PP _{\alpha \beta }})}
    \ar[r] ^(0.6){\mathrm{adj _{u _{\alpha \beta }}}} _(0.6){\sim}
  &
  {\E |_{ \PP _{\alpha \beta }}  .}
  }
\end{equation}
En effet, les deux carrés de droite de \ref{recolloc<->id2} sont commutatifs par fonctorialité,
tandis que celui de gauche l'est pour les mêmes raisons que
celle du deuxième carré de gauche de la troisième ligne de \ref{diag1-comp-comp-adj-immf}.
Or, le morphisme
$\phi _1 ^{\alpha \beta} (u^! _\alpha (\E |_{ \PP _\alpha}) )$
est, par construction (voir la preuve de \ref{comp-comp-adj-immf}),
égal au morphisme composé horizontal du haut de
\ref{recolloc<->id2}. Via la commutativité du diagramme \ref{recolloc<->id2},
il en découle que le morphisme composé horizontal
du haut de \ref{recolloc<->id1} est le morphisme d'adjonction
$u _{\alpha \beta +}  u _{\alpha \beta } ^! (\E |_{ \PP _{\alpha \beta }})
\rightarrow
\E |_{ \PP _{\alpha \beta }}$.
De même, on vérifie que morphisme composé horizontal du bas de \ref{recolloc<->id1}
est égal au morphisme d'adjonction par $u _{\alpha \beta }$.
\end{proof}

\begin{rema}
  Nous avons en fait prouvé mieux que \ref{prop1} :
  nous avons construit un foncteur {\it canonique} $\mathcal{R}ecol$
  quasi-inverse du foncteur canonique
  $\mathcal{L}oc$ : $\mathrm{Coh} (X,\, \PP ,\, T)
\rightarrow
\mathrm{Coh} (X,\, (\PP _\alpha) _{\alpha \in \Lambda},\, T)$.
Cela correspond à une extension de l'analogue $p$-adique de Berthelot du théorème de Kashiwara
(voir \cite[5.3.3]{Beintro2})
qui correspond au cas où $X \hookrightarrow P$
se relève en un morphisme $u$ : $\X \hookrightarrow \PP$
de $\V$-schémas formels lisses. Le foncteur $\mathcal{L}oc$ (resp. $\mathcal{R}ecol$) étend d'une certaine manière
$u ^!$ (resp. $u _+$).
\end{rema}
\begin{defi}
  On définit la catégorie
  $\mathrm{Isoc} ^\dag (Y,\,X,\, (\PP _\alpha) _{\alpha \in \Lambda}/K)$ de la manière suivante :
  les objets sont les familles $(E _\alpha) _{\alpha \in \Lambda}$,
  de $j _\alpha ^\dag \O _{\X _\alpha K}$-modules cohérents $E _\alpha$
  possédant une connexion intégrable surconvergente,
  ces familles étant munies d'une {\it donnée de recollement}, i.e., d'isomorphismes,
$ \eta _{  \alpha \beta} \ : \
p _{2 K}  ^{\alpha \beta !} (E _{\beta})
\riso
p  _{1 K} ^{\alpha \beta !} (E _{\alpha})$,
$j _{\alpha \beta} ^{\dag} \smash{\D} _{\X _{\alpha \beta K} }$-linéaires
et vérifiant la condition de cocycle :
$\eta _{13} ^{\alpha \beta \gamma }=
\eta _{12} ^{\alpha \beta \gamma }
\circ
\eta _{23} ^{\alpha \beta \gamma }$,
où $\eta _{12} ^{\alpha \beta \gamma }$, $\eta _{23} ^{\alpha \beta \gamma }$
et $\eta _{13} ^{\alpha \beta \gamma }$ sont définis par les diagrammes commutatifs
\begin{equation}
  \label{diag1-defindonnederecolK}
\xymatrix  @R=0,3cm {
{  p _{12K} ^{\alpha \beta \gamma !} p  _{2K} ^{\alpha \beta !}  (E _\beta )}
\ar[r] ^-{\epsilon} _{\sim}
\ar[d] ^-{p _{12K} ^{\alpha \beta \gamma !} (\eta _{\alpha \beta})} _{\sim}
&
{p _{2K} ^{\alpha \beta \gamma!}  (E _\beta )}
\ar@{.>}[d] ^-{\eta _{12} ^{\alpha \beta \gamma }}
\\
{ p _{12K} ^{\alpha \beta \gamma !}  p  _{1K} ^{\alpha \beta !}  (E _\alpha)}
\ar[r]^{\epsilon} _{\sim}
&
{p _{1K} ^{\alpha \beta \gamma!}(E _\alpha),}
}
%
%
\xymatrix  @R=0,3cm {
{  p _{23K} ^{\alpha \beta \gamma !} p  _{2K} ^{\beta \gamma!}  (E _\gamma )}
\ar[r] ^-{\epsilon} _{\sim}
\ar[d] ^-{p _{23K} ^{\alpha \beta \gamma !} (\eta _{ \beta \gamma})} _{\sim}
&
{p _{3K} ^{\alpha \beta \gamma!}  (E _\gamma )}
\ar@{.>}[d] ^-{\eta _{23} ^{\alpha \beta \gamma }}
\\
{ p _{23K} ^{\alpha \beta \gamma !}  p  _{1K} ^{ \beta \gamma !}  (E _\beta)}
\ar[r]^{\epsilon} _{\sim}
&
{p _{2K} ^{\alpha \beta \gamma!}(E _\beta),}
}
%
%
\xymatrix  @R=0,3cm {
{  p _{13K} ^{\alpha \beta \gamma !} p  _{2K} ^{\alpha \gamma !}  (E _\gamma )}
\ar[r] ^-{\epsilon} _{\sim}
\ar[d] ^-{p _{13K} ^{\alpha \beta \gamma !} (\eta _{\alpha \gamma})} _{\sim}
&
{p _{3K} ^{\alpha \beta \gamma!}  (E _\gamma )}
\ar@{.>}[d]^{\eta _{13} ^{\alpha \beta \gamma }}
\\
{ p _{13K} ^{\alpha \beta \gamma !}  p  _{1K} ^{\alpha \gamma !}  (E _\alpha)}
\ar[r]^{\epsilon} _{\sim}
&
{p _{1K} ^{\alpha \beta \gamma!}(E _\alpha).}
}
\end{equation}
Les morphismes
$f= (f _\alpha)_{\alpha \in \Lambda}$ :
$((E _{\alpha})_{\alpha \in \Lambda},\, (\eta _{\alpha\beta}) _{\alpha ,\beta \in \Lambda})
\rightarrow
((E '_{\alpha})_{\alpha \in \Lambda},\, (\eta '_{\alpha\beta}) _{\alpha ,\beta \in \Lambda})$
sont les familles
de morphismes $f _\alpha$ : $ E _\alpha \rightarrow E ' _\alpha$
commutant aux données de recollement.

On remarque enfin que, comme $p _{1K}$ est une immersion ouverte,
les foncteurs
$p  _{1 K} ^{\alpha \beta !} $ et $p  _{1 K} ^{\alpha \beta *} $ sont égaux.
De même, pour $p  _{2 K} ^{\alpha \beta } $,
$p _{12 K}^{\alpha \beta \gamma }$,
$p _{23 K}^{\alpha \beta \gamma }$
et $p _{13 K}^{\alpha \beta \gamma }$.
De plus, lorsque $\Lambda $ est réduit à un élément,
$\mathrm{Isoc} ^\dag (Y,\,X,\, \PP /K)$ est égale à
la catégorie des isocristaux sur $Y$ surconvergent le long de $T \cap X$, habituellement noté
$\mathrm{Isoc} ^\dag (Y,\,X/K)$.
\end{defi}

\begin{prop}\label{eqcat-iso-reco}
  Il existe un foncteur canonique
  $\mathcal{L}oc$ :
  $\mathrm{Isoc} ^\dag (Y,\,X/K)
  \rightarrow
  \mathrm{Isoc} ^\dag (Y,\,X,\, (\PP _\alpha) _{\alpha \in \Lambda}/K)$
  induisant une équivalence de catégorie.
\end{prop}
\begin{proof}
 On associe à chaque $j  ^\dag \O _{]X[ _{\PP}}$-module cohérent, $E $,
  muni d'une connexion intégrable surconvergente, la famille
  $ (u ^* _{\alpha K} ( E |_{]X _\alpha[ _{\PP _\alpha}})) _{\alpha \in \Lambda}$,
 où $u ^* _{\alpha K}$ est l'image inverse correspondant au diagramme commutatif suivant
 (voir les notations de \ref{notation-rig-form}) :
 $$\xymatrix  @R=0,3cm {
 { Y _\alpha}
 \ar[r]
 \ar@{=}[d]
 &
 {X _\alpha}
 \ar[r]
 \ar@{=}[d]
 &
 {\PP _\alpha}
\\
 { Y _\alpha}
 \ar[r]
 &
 {X _\alpha}
 \ar[r]
 &
 {\X _\alpha .}
 \ar[u] ^-{u _\alpha}
 }
 $$
De manière analogue au début de la preuve de \ref{prop1} (on remplace $\tau $ par $\epsilon$),
cette famille est munie d'une structure canonique de donnée de recollement.
  De même, si $f$ : $ E \rightarrow E'$ est un morphisme de $\mathrm{Isoc} ^\dag (Y,\,X/K)$,
la famille de morphismes $(u ^* _{\alpha K}( f |_{]X _\alpha[ _{\PP _\alpha}})) _{\alpha \in \Lambda}$
 commute aux données de recollement.
  On a ainsi construit un foncteur $\mathcal{L}oc$ :
    $\mathrm{Isoc} ^\dag (Y,\,X/K)
  \rightarrow
  \mathrm{Isoc} ^\dag (Y,\,X,\, (\PP _\alpha) _{\alpha \in \Lambda}/K)$.

  Grâce à \cite[2.3.1]{Berig} (qui implique que les foncteurs $u _{\alpha K} ^*$ et  $u _{\alpha \beta K} ^*$ sont
  pleinement fidèles pour les isocristaux surconvergents) on vérifie que ce foncteur est pleinement fidèle.

  Prouvons à présent son essentielle surjectivité.
  Soit $((E _{\alpha})_{\alpha \in \Lambda},\, (\eta _{\alpha\beta}) _{\alpha ,\beta \in \Lambda})
  \in \mathrm{Isoc} ^\dag (Y,\,X,\, (\PP _\alpha) _{\alpha \in \Lambda}/K)$.
  D'après \cite[2.3.1]{Berig}, il existe un
  $j _\alpha ^\dag \O _{]X _\alpha[ _{\PP _\alpha}}$-module cohérent $E '_\alpha$
  et un isomorphisme $\iota _\alpha$ : $ u ^* _{\alpha K} (E '_\alpha) \riso E _\alpha$.
  Il existe alors un unique isomorphisme $\tau ' _{\alpha \beta}$ :
  ${E '_\beta} |_{]X _{\alpha \beta}[ _{\PP _{\alpha \beta}}}
  \riso {E '_\alpha} |_{]X _{\alpha \beta}[ _{\PP _{\alpha \beta}}}$ s'inscrivant dans le diagramme
  commutatif :
  \begin{equation}
    \label{diag2-eqcat-iso-reco}
    \xymatrix  @R=0,3cm {
  {u ^* _{\alpha \beta K} ({E '_\beta} |_{]X _{\alpha \beta}[ _{\PP _{\alpha \beta}}})}
  \ar[r] ^-{\epsilon} _{\sim}
  \ar@{.>}[d] ^-{u ^* _{\alpha \beta K}(\tau ' _{\alpha \beta} )}
  &
  { p _{2 K}  ^{\alpha \beta *} u ^* _{\beta K} (E '_{\beta}) }
  \ar[rr] ^-{ p _{2 K}  ^{\alpha \beta *}(\iota _\beta)} _{\sim}
  &
  &
  { p _{2 K}  ^{\alpha \beta *} (E _{\beta}) }
    \ar[d] ^-{\eta  _{\alpha \beta} } _{\sim}
  \\
   {u ^* _{\alpha \beta K}({E '_\alpha} |_{]X _{\alpha \beta}[ _{\PP _{\alpha \beta}}})}
   \ar[r]^{\epsilon} _{\sim}
   &
  {p  _{1 K} ^{\alpha \beta *}u ^* _{\alpha K}  (E ' _{\alpha})}
  \ar[rr] ^-{p  _{1 K} ^{\alpha \beta *}(\iota _\alpha)} _{\sim}
  &
  &
  {p  _{1 K} ^{\alpha \beta *} (E _{\alpha}).}
    }
  \end{equation}
    Vérifions à présent que les isomorphismes $\tau ' _{\alpha \beta}$ se recollent.
    Considérons le diagramme commutatif suivant
    \begin{equation}
      \label{diag3-eqcat-iso-reco}\xymatrix  @R=0,3cm {
    {u ^* _{\alpha \beta \gamma K} ({E '_\beta} |_{]X _{\alpha \beta\gamma}[ _{\PP _{\alpha \beta\gamma}}})}
  \ar[r] ^-{\epsilon} _{\sim}
  \ar[d]
  ^{u ^* _{\alpha \beta \gamma K}({\tau ' _{\alpha \beta}} |_{]X _{\alpha \beta \gamma}[ _{\PP _{\alpha \beta\gamma}}})}
  _{\sim}
  &
  {p _{12K} ^{\alpha \beta \gamma *}u ^* _{\alpha \beta K} ({E '_\beta} |_{]X _{\alpha \beta}[ _{\PP _{\alpha \beta}}})}
  \ar[r] ^-{\epsilon} _{\sim}
  \ar[d] ^-{p _{12K} ^{\alpha \beta \gamma *} u ^* _{\alpha \beta K}(\tau ' _{\alpha \beta} )}
  _{\sim}
  &
  {p _{12K} ^{\alpha \beta \gamma *} p _{2 K}  ^{\alpha \beta *} u ^* _{\beta K} (E '_{\beta}) }
  \ar[r] ^-{\iota _\beta} _{\sim}
  &
  {p _{12K} ^{\alpha \beta \gamma *} p _{2 K}  ^{\alpha \beta *} (E _{\beta}) }
  \ar[d] ^-{p _{12K} ^{\alpha \beta \gamma *} (\eta _{\alpha \beta})} _{\sim}
  \ar[r] ^-{\epsilon} _{\sim}
&
{p _{2K} ^{\alpha \beta \gamma*}  (E _\beta ),}
\ar[d] ^-{\eta _{12} ^{\alpha \beta \gamma }} _{\sim}
\\
  {u ^* _{\alpha \beta \gamma K} ({E '_\alpha} |_{]X _{\alpha \beta \gamma}[ _{\PP _{\alpha \beta\gamma}}})}
  \ar[r] ^-{\epsilon} _{\sim}
  &
   {p _{12K} ^{\alpha \beta \gamma *} u ^* _{\alpha \beta K}({E '_\alpha} |_{]X _{\alpha \beta}[ _{\PP _{\alpha \beta}}})}
   \ar[r]^{\epsilon} _{\sim}
   &
  {p _{12K} ^{\alpha \beta \gamma *} p  _{1 K} ^{\alpha \beta *}u ^* _{\alpha K}  (E ' _{\alpha})}
  \ar[r] ^-{\iota _\alpha} _{\sim}
  &
  {p _{12K} ^{\alpha \beta \gamma *} p  _{1 K} ^{\alpha \beta *} (E _{\alpha})}
  \ar[r]^{\epsilon} _{\sim}
&
{p _{1K} ^{\alpha \beta \gamma*}(E _\alpha),}
   }
    \end{equation}
    où le rectangle du milieu se déduit de \ref{diag2-eqcat-iso-reco}
    par application du foncteur $ p _{12K} ^{\alpha \beta \gamma *}$.
   On remarque que la flèche composée du bas de \ref{diag3-eqcat-iso-reco}
   est indépendant du chemin suivi, i.e.,
le diagramme suivant est commutatif :
$$\xymatrix  @R=0,3cm@C=0,5cm {
  {u ^* _{\alpha \beta \gamma K} ({E '_\alpha} |_{]X _{\alpha \beta \gamma}[ _{\PP _{\alpha \beta\gamma}}})}
  \ar[r] ^(0.45){\epsilon} _(0.45){\sim}
  &
   {p _{12K} ^{\alpha \beta \gamma *} u ^* _{\alpha \beta K}({E '_\alpha} |_{]X _{\alpha \beta}[ _{\PP _{\alpha \beta}}})}
   \ar[r]^{\epsilon} _{\sim}
   &
  {p _{12K} ^{\alpha \beta \gamma *} p  _{1 K} ^{\alpha \beta *}u ^* _{\alpha K}  (E ' _{\alpha})}
  \ar[r] ^-{\iota _\alpha} _{\sim}
  &
  {p _{12K} ^{\alpha \beta \gamma *} p  _{1 K} ^{\alpha \beta *} (E _{\alpha})}
  \ar[r]^{\epsilon} _{\sim}
&
{p _{1K} ^{\alpha \beta \gamma*}(E _\alpha)}
\\
  {u ^* _{\alpha \beta \gamma K} ({E '_\alpha} |_{]X _{\alpha \beta \gamma}[ _{\PP _{\alpha \beta \gamma}}})}
  \ar[r] ^(0.45){\epsilon} _(0.45){\sim}
  \ar@{=}[u]
  &
   {p _{13K} ^{\alpha \beta \gamma *} u ^* _{\alpha \gamma K}({E '_\alpha} |_{]X _{\alpha \gamma}[ _{\PP _{\alpha \gamma}}})}
   \ar[r]^{\epsilon} _{\sim}
   \ar[u] ^\epsilon _{\sim}
   &
  {p _{13K} ^{\alpha \beta \gamma *} p  _{1 K} ^{\alpha \gamma *}u ^* _{\alpha K}  (E ' _{\alpha})}
  \ar[r] ^-{\iota _\alpha} _{\sim}
  \ar[u] ^\epsilon _{\sim}
  &
  {p _{13K} ^{\alpha \beta \gamma *} p  _{1 K} ^{\alpha \gamma *} (E _{\alpha})}
  \ar[r]^{\epsilon} _{\sim}
  \ar[u] ^\epsilon _{\sim}
&
{p _{1K} ^{\alpha \beta \gamma*}(E _\alpha).}
\ar@{=}[u]
}$$
De même pour la flèche du haut de \ref{diag3-eqcat-iso-reco}.
En écrivant les deux autres diagrammes analogues à \ref{diag3-eqcat-iso-reco}, on établit
    que la famille $(E ' _{\alpha}) _{\alpha \in \Lambda}$ se recolle en un
    $j  ^\dag \O _{]X _{\alpha }[ _{\PP _{\alpha}}}$-module cohérent $E '$,
  muni d'une connexion intégrable surconvergente. En outre, les isomorphismes $\iota _\alpha$ induisent
  $\mathcal{L}oc (E ') \riso ((E _{\alpha})_{\alpha \in \Lambda},\, (\eta _{\alpha\beta}) _{\alpha ,\beta \in \Lambda})$.
\end{proof}
\begin{defi}
  Avec les notations de \ref{eqcat-iso-reco}, soient
  $((E _{\alpha})_{\alpha \in \Lambda},\, (\eta _{\alpha\beta}) _{\alpha ,\beta \in \Lambda})
  \in \mathrm{Isoc} ^\dag (Y,\,X,\, (\PP _\alpha) _{\alpha \in \Lambda}/K)$ et
  $E \in \mathrm{Isoc} ^\dag (Y,\,X/K)$.
  On dira que $((E _{\alpha})_{\alpha \in \Lambda},\, (\eta _{\alpha\beta}) _{\alpha ,\beta \in \Lambda})$
  {\it se recolle} en $E$,
  s'il existe un isomorphisme
  $\mathcal{L}oc (E) \riso ((E _{\alpha})_{\alpha \in \Lambda},\, (\eta _{\alpha\beta}) _{\alpha ,\beta \in \Lambda})$.
\end{defi}

\begin{prop}\label{pre-sp+plfid}
  Les foncteurs $\sp _*$ et $\sp ^*$ induisent des équivalences quasi-inverses
  entre la catégorie
  $\mathrm{Isoc} ^\dag (Y,\,X,\, (\PP _\alpha) _{\alpha \in \Lambda}/K)$ et la sous-catégorie
  pleine de $\mathrm{Coh} (X,\, (\PP _\alpha) _{\alpha \in \Lambda},\, T)$ des objets
  $((\E _{\alpha})_{\alpha \in \Lambda},\, (\theta _{\alpha\beta}) _{\alpha ,\beta \in \Lambda})$
  tels que, pour tout $\alpha \in \Lambda$,
  $\E _{\alpha}$ soit $\O _{\X _\alpha,\,\Q} (\hdag T \cap X _\alpha)$-cohérent.
\end{prop}
\begin{proof}
i) Construisons d'abord le foncteur
$\sp _*$ : $\mathrm{Isoc} ^\dag (Y,\,X,\, (\PP _\alpha) _{\alpha \in \Lambda}/K)
  \rightarrow \mathrm{Coh} (X,\, (\PP _\alpha) _{\alpha \in \Lambda},\, T)$.

  Si $((E _{\alpha})_{\alpha \in \Lambda},\, (\eta _{\alpha\beta}) _{\alpha ,\beta \in \Lambda})$
est un objet de
  $\mathrm{Isoc} ^\dag (Y,\,X,\, (\PP _\alpha) _{\alpha \in \Lambda}/K)$,
  démontrons alors que
  $( (\sp _* E _{\alpha})_{\alpha \in \Lambda},\, (\theta _{\alpha\beta}) _{\alpha ,\beta \in \Lambda})$,
    où $\theta _{\alpha \beta}$ est l'unique morphisme rendant commutatif le diagramme suivant
  \begin{equation}
    \label{defdonneesp*}
    \xymatrix  @R=0,3cm
{
    {\sp _*  p _{1 K}  ^{\alpha \beta !} (E _{\alpha})}
    \ar[r] _{\sim}
    &
    {p _{1 }  ^{\alpha \beta !} \sp _*    (E _{\alpha})}
    \\
    {\sp _*  p _{2 K}  ^{\alpha \beta !} (E _{\beta})}
    \ar[r] _{\sim}
    \ar[u] ^-{\sp _* \eta _{\alpha \beta}} _{\sim}
    &
    {p _{2 }  ^{\alpha \beta !} \sp _*    (E _{\beta})}
    \ar@{.>}[u]^{\theta _{\alpha \beta}}
}
  \end{equation}
dont les isomorphismes horizontaux découlent de \ref{spcommup*},
est un objet de $\mathrm{Coh} (X,\, (\PP _\alpha) _{\alpha \in \Lambda},\, T)$, i.e.,
que les isomorphismes $\theta _{\alpha \beta}$
vérifient la condition de cocycle. Considérons le diagramme commutatif suivant :
\begin{equation}
  \label{diag1-pre-sp+plfid}
  \xymatrix  @R=0,3cm {
  {\sp _*  p _{1K} ^{\alpha \beta \gamma!}(E _\alpha)}
  \ar[r] ^(0.45){\sp _* (\epsilon)} _(0.45){\sim}
  &
  {\sp _* p _{12 K} ^{\alpha \beta \gamma *}  p _{1 K}  ^{\alpha \beta *} (E _{\alpha})}
  \ar[r] _{\sim}
  &
  {p _{12 } ^{\alpha \beta \gamma !}  \sp _* p _{1 K}  ^{\alpha \beta *} (E _{\alpha})}
  \ar[r] _{\sim}
  &
  {p _{12 } ^{\alpha \beta \gamma !}  p _{1 }  ^{\alpha \beta !}  \sp _* (E _{\alpha})}
  \ar[r] ^-{\tau} _{\sim}
  &
  {p _1 ^{\alpha \beta \gamma!}(\sp _* (E _{\alpha}))}
  \\
  {\sp _* p _{2K} ^{\alpha \beta \gamma!}  (E _\beta )}
  \ar[r] ^(0.45){\sp _* (\epsilon)} _(0.45){\sim}
  \ar[u] ^-{\sp _* (\eta _{12} ^{\alpha \beta \gamma })} _{\sim}
  &
  {\sp _* p _{12 K} ^{\alpha \beta \gamma *}  p _{2 K}  ^{\alpha \beta *} (E _{\beta})}
  \ar[r] _{\sim}
  \ar[u] ^-{\sp _* p _{12 K} ^{\alpha \beta \gamma *} (\eta _{\alpha \beta})} _{\sim}
  &
  {p _{12 } ^{\alpha \beta \gamma !}  \sp _* p _{2 K}  ^{\alpha \beta *} (E _{\beta})}
  \ar[r] _{\sim}
  \ar[u]^{p _{12 } ^{\alpha \beta \gamma !} \sp _* (\eta _{\alpha \beta})} _{\sim}
  &
  {p _{12 } ^{\alpha \beta \gamma !}  p _{2 }  ^{\alpha \beta !}  \sp _* (E _{\beta})}
  \ar[u] ^-{p _{12 } ^{\alpha \beta \gamma !} (\theta _{\alpha \beta}) } _{\sim}
  \ar[r] ^-{\tau} _{\sim}
  &
  {p _2 ^{\alpha \beta \gamma!}  ( \sp _* (E _{\beta})).}
  \ar[u] ^-{\theta _{12} ^{\alpha \beta \gamma }} _{\sim}
  }
\end{equation}
Or, il découle de \ref{spcommup*} et de \ref{sp-eps-tau} que
 le morphisme composé du haut (resp. du bas) de \ref{diag1-pre-sp+plfid} est
 l'isomorphisme canonique
 $\sp _*  (p _{1 K}  ^{\alpha \beta \gamma } ) ^* (E _{\alpha})
 \riso
  p _{1 }  ^{\alpha \beta \gamma!} \sp _* (E _{\alpha})$
  (resp.  $\sp _*  (p _{2 K}  ^{\alpha \beta \gamma } ) ^* (E _{\beta})
 \riso
  p _{2 }  ^{\alpha \beta \gamma!}  \sp _* (E _{\beta})$).
  En écrivant les deux autres diagrammes analogues, on conclut que la condition
  de cocycle est validée.

En outre, si $f = (f _\alpha)_{\alpha \in \Lambda}$ :
$((E _{\alpha})_{\alpha \in \Lambda},\, (\eta _{\alpha\beta}) _{\alpha ,\beta \in \Lambda})
\rightarrow ((E '_{\alpha})_{\alpha \in \Lambda},\, (\eta '_{\alpha\beta}) _{\alpha ,\beta \in \Lambda})$
est un morphisme de la catégorie $\mathrm{Isoc} ^\dag (Y,\,X,\, (\PP _\alpha) _{\alpha \in \Lambda}/K)$,
on vérifie, par fonctorialité de \ref{defdonneesp*}, que la famille
$\sp _* (f): = (\sp _* (f _\alpha))_{\alpha \in \Lambda}$ commute aux données de recollement.

ii) Réciproquement, on construit un foncteur, noté $\sp ^*$, de la manière suivante :
soit $((\E _{\alpha})_{\alpha \in \Lambda},\, (\theta _{\alpha\beta}) _{\alpha ,\beta \in \Lambda})$
est un objet de
$\mathrm{Coh} (X,\, (\PP _\alpha) _{\alpha \in \Lambda},\, T)$
tel que
$\E _{\alpha}$ soit $\O _{\X _\alpha,\,\Q} (\hdag T \cap X _\alpha)$-cohérent.
En notant $\eta _{\alpha \beta}$ l'unique morphisme rendant commutatif le diagramme suivant
  \begin{equation}
    \label{defdonneesp*2}
    \xymatrix  @R=0,3cm
{
    {\sp ^*  p _{1 }  ^{\alpha \beta !} (\E _{\alpha})}
    &
    {p _{1 K}  ^{\alpha \beta *} \sp ^*    (\E _{\alpha})}
    \ar[l] ^-{\sim}
    \\
    {\sp ^*  p _{2 }  ^{\alpha \beta !} (\E _{\beta})}
    \ar[u] ^-{\sp ^* \theta _{\alpha \beta}} _{\sim}
    &
    {p _{2 K}  ^{\alpha \beta *} \sp ^*    (\E _{\beta}),}
    \ar@{.>}[u]^{\eta _{\alpha \beta}}
    \ar[l] ^-{\sim}
}
  \end{equation}
  on pose
$\sp ^* ((\E _{\alpha})_{\alpha \in \Lambda},\, (\theta _{\alpha\beta}) _{\alpha ,\beta \in \Lambda})
:=
((\sp ^* \E _{\alpha})_{\alpha \in \Lambda},\, (\eta _{\alpha\beta}) _{\alpha ,\beta \in \Lambda})$.
De manière analogue à l'étape i) (en utilisant \ref{spcommup*} et \ref{sp-eps-tau}),
on établit que
$((\sp ^* \E _{\alpha})_{\alpha \in \Lambda},\, (\eta _{\alpha\beta}) _{\alpha ,\beta \in \Lambda}) \in
\mathrm{Isoc} ^\dag (Y,\,X,\, (\PP _\alpha) _{\alpha \in \Lambda}/K)$.
  De plus, si $f = (f _\alpha)_{\alpha \in \Lambda}$ :
$((\E _{\alpha})_{\alpha \in \Lambda},\, (\theta _{\alpha\beta}) _{\alpha ,\beta \in \Lambda})
\rightarrow ((\E '_{\alpha})_{\alpha \in \Lambda},\, (\theta '_{\alpha\beta}) _{\alpha ,\beta \in \Lambda})$
est un morphisme, par fonctorialité de \ref{defdonneesp*2}, la famille
$\sp ^* (f) := (\sp ^* (f _\alpha))_{\alpha \in \Lambda}$ commute aux données de recollement.

iii) Soit $((E _{\alpha})_{\alpha \in \Lambda},\, (\eta _{\alpha\beta}) _{\alpha ,\beta \in \Lambda})$
un objet de
  $\mathrm{Isoc} ^\dag (Y,\,X,\, (\PP _\alpha) _{\alpha \in \Lambda}/K)$.
  Notons $\theta _{\alpha \beta}$ (resp. $\eta ' _{\alpha \beta}$) les isomorphismes
  de recollement canoniques induits sur $\sp _* (E _{\alpha})$ (resp. $\sp ^* \sp _* (E _{\alpha})$).
  Démontrons maintenant que les isomorphismes d'adjonctions
$\sp ^* \sp _* (E _\alpha) \riso E _\alpha$ commutent aux données de recollement, i.e., que
le carré de droite du diagramme suivant
\begin{equation}\label{diag2-pre-sp+plfid}
  \xymatrix  @R=0,3cm {
    {p _{1 K}  ^{\alpha \beta !} (E _{\alpha})}
    &
    {\sp ^* \sp _*  p _{1 K}  ^{\alpha \beta !} (E _{\alpha})}
    \ar[r] _{\sim}
    \ar[l] ^-{\sim}
    &
    {\sp ^* p _{1 }  ^{\alpha \beta !} \sp _*    (E _{\alpha})}
    &
    {p _{1 K}  ^{\alpha \beta *} \sp ^*  \sp _*   (E _{\alpha})}
    \ar[r] _{\sim}
    \ar[l] ^-{\sim}
    &
    {p _{1 K}  ^{\alpha \beta *}  (E _{\alpha})}
    \\
    { p _{2 K}  ^{\alpha \beta !} (E _{\beta})}
    \ar[u] ^-{\eta _{\alpha \beta}} _{\sim}
    &
    {\sp ^* \sp _*  p _{2 K}  ^{\alpha \beta !} (E _{\beta})}
    \ar[r] _{\sim}
    \ar[l] ^-{\sim}
    \ar[u] ^-{\sp ^* \sp _* \eta _{\alpha \beta}} _{\sim}
    &
    {\sp ^*  p _{2 }  ^{\alpha \beta !} \sp _*    (E _{\beta})}
    \ar[u]^{\sp ^* (\theta _{\alpha \beta})} _{\sim}
    &
    {p _{2 K}  ^{\alpha \beta *} \sp ^*  \sp _*   (E _{\beta})}
    \ar[u]^{\eta ' _{\alpha \beta}} _{\sim}
    \ar[r] _{\sim}
    \ar[l] ^-{\sim}
    &
    {p _{2 K}  ^{\alpha \beta *} (E _{\beta}),}
    \ar[u] ^-{\eta _{\alpha \beta}} _{\sim}
    }
\end{equation}
est commutatif.
La commutativité des deux carrés du milieu découlent respectivement de \ref{defdonneesp*} et \ref{defdonneesp*2},
tandis que celle du carré de gauche se vérifie par fonctorialité.
Or, il dérive de \ref{diag1spcommup*} que les morphismes composés du haut et du bas
de \ref{diag2-pre-sp+plfid} sont les morphismes identités.
Comme toutes les flèches sont des isomorphismes, le diagramme
\ref{diag2-pre-sp+plfid} est commutatif.

iv) Soit $((\E _{\alpha})_{\alpha \in \Lambda},\, (\theta _{\alpha\beta}) _{\alpha ,\beta \in \Lambda})$
un objet de
$\mathrm{Coh} (X,\, (\PP _\alpha) _{\alpha \in \Lambda},\, T)$
tel que
$\E _{\alpha}$ soit $\O _{\X _\alpha,\,\Q} (\hdag T \cap X _\alpha)$-cohérent.
  On vérifie de façon similaire à l'étape iii) (on utilise \ref{spcommup*2}
  au lieu de \ref{diag1spcommup*}), que les isomorphismes d'adjonctions
$ \E _\alpha \riso \sp _*  \sp ^*(\E _\alpha)$ commutent aux données de recollement.
\end{proof}

\begin{theo}\label{sp+plfid}
  On dispose d'un foncteur canonique pleinement fidèle $\sp _{X \hookrightarrow \PP,\,T \,+}$ :
  $\mathrm{Isoc} ^\dag (Y,\,X/K) \rightarrow \mathrm{Coh} (X,\, \PP ,\, T)$.
  Son image essentielle est constituée par les
  $\smash{\D} ^{\dag}  _{\PP}(\hdag T) _{\Q}$-modules cohérents $\E$ à support dans $X$ satisfaisant l'une des conditions
  équivalentes suivantes :

  (*) pour tout ouvert $\PP '$ de $\PP$ tel que l'immersion fermée $X \cap P ' \hookrightarrow P'$
  se relève en un morphisme $v$ : $ \X' \hookrightarrow \PP '$ de $\V$-schémas formels lisses,
  le faisceau $v ^! (\E |_{\PP'}) $ est $\O _{\X',\,\Q} (\hdag T \cap X')$-cohérent,

  (**) pour tout ouvert $\PP '$ affine de $\PP$, pour tout relèvement
  $v$ : $ \X' \hookrightarrow \PP '$ de $X \cap P ' \hookrightarrow P'$,
  le faisceau $v ^! (\E |_{\PP'}) $ est $\O _{\X',\,\Q} (\hdag T \cap X')$-cohérent.

Si $T$ est vide, on le notera $\sp _{X \hookrightarrow \PP\,+}$.
  Si aucune confusion n'est à craindre, on écrira simplement $\sp _+$.
\end{theo}
\begin{proof}
  Il découle de \ref{eqcat-iso-reco}, \ref{prop1} et \ref{pre-sp+plfid} et avec leurs notations,
  la construction d'un foncteur canonique pleinement fidèle $\sp _{X \hookrightarrow \PP,\,T \,+}$ :
  $\mathrm{Isoc} ^\dag (Y,\,X/K) \rightarrow \mathrm{Coh} (X,\, \PP ,\, T)$ défini en posant
  $\sp _{X \hookrightarrow \PP,\,T \,+} := \mathcal{R}ecol \circ \sp _* \circ \mathcal{L} oc$.
  À présent, caractérisons son image essentielle. La condition (*) implique (**). De plus, soit
  $\E$ un $\smash{\D} ^{\dag}  _{\PP}(\hdag T) _{\Q}$-module cohérent à support dans $X$
  satisfaisant la condition (**).
  Avec les notations de la preuve de \ref{prop1},
  soit $\mathcal{L}oc (\E) :=
((u _{\alpha } ^! (\E |_{\PP _\alpha})) _{\alpha \in \Lambda},\, (\theta _{\alpha \beta})_{\alpha , \beta \in \Lambda})$.
Via \ref{pre-sp+plfid}, comme $u _{\alpha } ^! (\E |_{\PP _\alpha})$ est
$\O _{\X _\alpha,\,\Q} (\hdag T \cap X _\alpha)$-cohérent, il en résulte que $\E$ est dans l'image
essentielle de $\sp _{X \hookrightarrow \PP,\,T \,+}$. Réciproquement,
donnons-nous un objet $((\E _{\alpha})_{\alpha \in \Lambda},\, (\theta _{\alpha\beta}) _{\alpha ,\beta \in \Lambda})$
de
$\mathrm{Coh} (X,\, (\PP _\alpha) _{\alpha \in \Lambda},\, T)$
  tel que, pour tout $\alpha \in \Lambda$,
  $\E _{\alpha}$ soit $\O _{\X _\alpha,\,\Q} (\hdag T \cap X _\alpha)$-cohérent.
  Posons $\E := \mathcal{R}ecol
  ((\E _{\alpha})_{\alpha \in \Lambda},\, (\theta _{\alpha\beta}) _{\alpha ,\beta \in \Lambda})$
  et prouvons que $\E$ vérifie la condition (*).
  Soit $\PP '$ un ouvert de $\PP$ tel que l'immersion fermée $X \cap P ' \hookrightarrow P'$
  se relève en un morphisme $v$ : $ \X' \hookrightarrow \PP '$ de $\V$-schémas formels lisses.
  Il s'agit de prouver que $v ^! (\E |_{\PP '}) $ est $\O _{\X',\,\Q} (\hdag T \cap X')$-cohérent,
  ce qui est local en $\X'$ et $\PP'$.
  On peut donc supposer que $\X'$ est affine et qu'il existe $\alpha _0 \in \Lambda$ tel que
  $\PP '  \subset \PP _{\alpha _0}$.
  Choisissons $\rho _\alpha$ : $\X '  \hookrightarrow \X _{\alpha _0}$
 un relèvement de $X '  \subset X _{\alpha _0}$.
  On conclut via l'isomorphisme $  \xymatrix  @R=0,3cm  {
  {v ^!  (\E |_{\PP '})}
  \ar[r] ^\tau _{\sim}
  &
  {\rho _\alpha ^! u _{\alpha } ^! (\E |_{\PP _\alpha})}
  }
  $.
\end{proof}
\begin{rema}\label{rema-sp+plfid}
  Avec les notations de \ref{sp+plfid}, grâce à la description de Berthelot
des isocristaux surconvergents énoncée dans \cite[2.2.12]{caro_courbe-nouveau},
  $v ^! (\E |_{\PP'}) $ est $\O _{\X',\,\Q} (\hdag T \cap X')$-cohérent
  si et seulement si, en notant $\Y '$ l'ouvert de $\X '$ complémentaire
  de $T \cap X'$,  $v ^! (\E |_{\PP'}) | _{\Y '}$ est $\O _{\Y',\,\Q}$-cohérent.
\end{rema}

\begin{prop}\label{sp+indPalpha}
  Le foncteur $\sp _{X \hookrightarrow \PP,\,T \,+}$
  ne dépend que de $X \hookrightarrow \PP$ et de $T$.
\end{prop}
\begin{proof}
  Il s'agit de prouver que le foncteur $\sp _{X \hookrightarrow \PP,\,T \,+}$
  est indépendant du choix des recouvrements ouverts de $P$, ni de celui des relèvements qu'il induit.
  Donnons-nous un deuxième recouvrement d'ouverts affines $\PP = \cup _{\alpha '\in \Lambda'} \PP _{\alpha'}$.
  En remplaçant respectivement $\alpha, \,\beta,\,\gamma \in \Lambda$ par
  $\alpha', \,\beta',\,\gamma'\in \Lambda '$,
  on prend pour ce recouvrement les notations analogues à \ref{notat-construc}.

  Dans un premier temps, supposons qu'il existe une application surjective $\rho$ : $\Lambda ' \rightarrow \Lambda$
  telle que, $\PP _{\alpha '} \subset \PP _{\rho(\alpha ')}$ et
  $\PP _\alpha = \cup _{\alpha ' \in \rho ^{-1} (\alpha)}\PP _{\alpha '}$.
  Pour tout $(\alpha ',\,\beta ',\,\gamma ') \in \Lambda ^{'3}$, on choisit alors des relèvements
  $\epsilon _{\alpha '}$ : $\X _{\alpha '} \rightarrow \X _{\rho(\alpha')}$,
  $\epsilon _{\alpha '\beta '}$ : $\X _{\alpha ' \beta '} \rightarrow \X _{\rho(\alpha')\rho(\beta ')}$
  et
  $\epsilon _{\alpha '\beta '\gamma '}$ : $\X _{\alpha ' \beta '\gamma '}
  \rightarrow \X _{\rho(\alpha')\rho(\beta ')\rho(\gamma ')}$.
  On construit alors un foncteur
  $$\mathcal{L}oc _{(\PP _{\alpha'}) _{\alpha '\in \Lambda'},\,(\PP _\alpha) _{\alpha \in \Lambda}}\ : \
  \mathrm{Coh} (X,\, (\PP _\alpha) _{\alpha \in \Lambda},\, T)
  \rightarrow \mathrm{Coh} (X,\, (\PP _{\alpha'}) _{\alpha '\in \Lambda'},\, T)$$
  de la façon suivante :
pour tout objet $((\E _{\alpha}) _{\alpha \in \Lambda}, (\theta _{\alpha \beta }) _{\alpha, \beta \in \Lambda})$
de la catégorie $\mathrm{Coh} (X,\, (\PP _\alpha) _{\alpha \in \Lambda},\, T)$, on définit
$\mathcal{L}oc _{(\PP _{\alpha'}) _{\alpha '\in \Lambda'},\,(\PP _\alpha) _{\alpha \in \Lambda}}
((\E _{\alpha}) _{\alpha \in \Lambda}, (\theta _{\alpha \beta }) _{\alpha, \beta \in \Lambda})
=((\E _{\alpha'}) _{\alpha '\in \Lambda'}, (\theta _{\alpha '\beta '}) _{\alpha', \beta' \in \Lambda'})$,
où $\E _{\alpha'} :=\epsilon _{\alpha'} ^! (\E _{\rho (\alpha')})$ et
$\theta _{\alpha '\beta '}$ est l'unique morphisme induisant le diagramme commutatif :
\begin{equation}\label{sp+indPalpha-diag1}
  \xymatrix  @R=0,3cm {
{p _1 ^{\alpha '\beta '!} \epsilon _{\alpha'} ^! (\E _{\rho (\alpha')})}
\ar[r]^(0.4){\tau} _(0.4){\sim}
&
{ \epsilon _{\alpha'\beta '} ^! p _1 ^{\rho(\alpha ')\rho(\beta ')!}(\E _{\rho (\alpha')})}
\\
{p _2 ^{\alpha '\beta '!} \epsilon _{\beta'} ^! (\E _{\rho (\beta')})}
\ar[r]^(0.4){\tau} _(0.4){\sim}
\ar@{.>}[u] ^-{\theta _{\alpha '\beta '}}
&
{ \epsilon _{\alpha'\beta '} ^! p _2 ^{\rho(\alpha ')\rho(\beta ')!}(\E _{\rho (\beta')}).}
\ar[u] ^-{ \epsilon _{\alpha'\beta '} ^! (\theta _{\rho(\alpha ') \rho (\beta ')}) }
_{\sim}
}
\end{equation}
On vérifie que les isomorphismes $\theta _{\alpha '\beta '}$ satisfont à la condition de
cocycle. Enfin, en transformant le carré \ref{sp+indPalpha-diag1} en cube,
on obtient la fonctorialité de
$\mathcal{L}oc _{(\PP _{\alpha'}) _{\alpha '\in \Lambda'},\,(\PP _\alpha) _{\alpha \in \Lambda}}$.

De manière analogue (en remplaçant $\E$ par $E$, $\theta $ par $\eta$ et
$\tau $ par $\epsilon$), on construit un foncteur
  $$\mathcal{L}oc _{(\PP _{\alpha'}) _{\alpha '\in \Lambda'},\,(\PP _\alpha) _{\alpha \in \Lambda}}
  \ : \
  \mathrm{Isoc} ^\dag (Y,\,X,\, (\PP _\alpha) _{\alpha \in \Lambda}/K)
\rightarrow
\mathrm{Isoc} ^\dag (Y,\,X,\, (\PP _{\alpha'}) _{\alpha '\in \Lambda'}/K).$$
Pour ces deux foncteurs, on dispose d'un isomorphisme canonique
$\mathcal{L}oc _{(\PP _{\alpha'}) _{\alpha '\in \Lambda'},\,(\PP _\alpha) _{\alpha \in \Lambda}}
\circ \mathcal{L}oc \riso \mathcal{L}oc$.
Il reste alors à prouver la commutation suivante :
$\sp ^* \circ \mathcal{L}oc _{(\PP _{\alpha'}) _{\alpha '\in \Lambda'},\,(\PP _\alpha) _{\alpha \in \Lambda}}
\riso \mathcal{L}oc _{(\PP _{\alpha'}) _{\alpha '\in \Lambda'},\,(\PP _\alpha) _{\alpha \in \Lambda}}
\circ \sp ^*$,
où $\sp ^*$
désigne le foncteur construit dans \ref{pre-sp+plfid}.
Soit
$((\E _{\alpha})_{\alpha \in \Lambda},\, (\theta _{\alpha\beta}) _{\alpha ,\beta \in \Lambda})$
un objet de $\mathrm{Coh} (X,\, (\PP _\alpha) _{\alpha \in \Lambda},\, T)$ tel que,
pour tout $\alpha \in \Lambda$,
  $\E _{\alpha}$ soit $\O _{\X _\alpha,\,\Q} (\hdag T \cap X _\alpha)$-cohérent.
  Avec les notations analogues à \ref{defdonneesp*2},
  considérons le diagramme
\begin{equation}
  \notag
  \xymatrix  @R=0,3cm @C=-1,5cm @L=0,1cm  {
 &
 {p _{1K} ^{\alpha '\beta '*} \sp ^* \epsilon _{\alpha'} ^! (\E _{\rho (\alpha')})}
 \ar[rr]
 &&
 {\sp ^* p _1 ^{\alpha '\beta '!} \epsilon _{\alpha'} ^! (\E _{\rho (\alpha')})}
 \ar[rr]^{\tau}
 &&
 {\sp ^* \epsilon _{\alpha'\beta '} ^! p _1 ^{\rho(\alpha ')\rho(\beta ')!}(\E _{\rho (\alpha')})}
 \\
 {p _{2K} ^{\alpha '\beta '*} \sp ^* \epsilon _{\beta'} ^! (\E _{\rho (\beta')})}
 \ar[rr]
 \ar[ur] ^(0.4){\eta _{\alpha ' \beta '}}
 &&
 {\sp ^* p _2 ^{\alpha '\beta '!} \epsilon _{\beta'} ^! (\E _{\rho (\beta')})}
 \ar[rr]^{\tau}
 \ar[ur]  ^(0.4){\sp ^*(\theta _{\alpha '\beta '})}
 &&
 {\sp ^* \epsilon _{\alpha'\beta '} ^! p _2 ^{\rho(\alpha ')\rho(\beta ')!}(\E _{\rho (\beta')})}
 \ar[ur] ^(0.4){\sp ^* \epsilon _{\alpha'\beta '} ^! (\theta _{\rho(\alpha ') \rho (\beta ')}) }
 \\
   &
 {p _{1K} ^{\alpha '\beta '*}  \epsilon _{\alpha'K} ^* \sp ^*(\E _{\rho (\alpha')})}
 \ar[rr] ^-{\epsilon}
 \ar[uu]
 &&
 { \epsilon _{\alpha'\beta 'K} ^*  p _{1K} ^{\rho(\alpha ')\rho(\beta ')*} \sp ^*(\E _{\rho (\alpha')})}
 \ar[rr]
 &&
 {\epsilon _{\alpha'\beta 'K} ^* \sp ^* p _1 ^{\rho(\alpha ')\rho(\beta ')*}(\E _{\rho (\alpha')})}
 \ar[uu]
 \\
 {p _{2K} ^{\alpha '\beta '*}  \epsilon _{\beta'K} ^* \sp ^* (\E _{\rho (\beta')})}
 \ar[rr]^{\epsilon}
 \ar[uu]
 \ar@{.>}[ur] _(0.6){\eta' _{\alpha ' \beta '}}
 &&
 { \epsilon _{\alpha'\beta 'K} ^* p _{2K} ^{\rho(\alpha ')\rho(\beta ')*} \sp ^* (\E _{\rho (\beta')})}
 \ar[rr]
 \ar[ur]  _(0.6){ \epsilon _{\alpha'\beta 'K} ^* (\eta _{\rho(\alpha ')\rho (\beta ')})}
 &&
 { \epsilon _{\alpha'\beta 'K} ^* \sp ^* p _2 ^{\rho(\alpha ')\rho(\beta ')*}(\E _{\rho (\beta')}),}
 \ar[ur] _(0.6){ \epsilon _{\alpha'\beta 'K} ^* \sp ^* (\theta _{\rho(\alpha ') \rho (\beta ')}) }
 \ar[uu]
}
\end{equation}
où les carrés horizontaux sont commutatifs par définition. Il découle
de \ref{spcommup*} et de \ref{sp-eps-tau} que les rectangles de devant et de
derrière le sont aussi. Comme il en est de même par fonctorialité du carré de droite,
il en dérive que le carré de gauche est commutatif, i.e.,
que l'isomorphisme canonique
$ \sp ^* \epsilon _{\alpha'} ^! (\E _{\rho (\alpha')})
\riso \epsilon _{\alpha'K} ^* \sp ^*(\E _{\rho (\alpha')})$
commute aux données de recollement respectives.

Passons à présent au cas général.
En utilisant le troisième recouvrement
$\PP = \cup _{(\alpha,\, \alpha')\in \Lambda \times \Lambda '} \PP _\alpha \cap \PP _{\alpha '}$,
on se ramène au premier cas.
\end{proof}

\section{Comparaison des foncteurs duaux des isocristaux surconvergents
: compatibilité aux images inverses (extraordinaires) par une immersion ouverte}

On désignera par
$f$ : $\X' \rightarrow \X$ une immersion ouverte de $\V$-schémas formels lisses,
$T$ un diviseur de $X$ et $T ' := f ^{-1} (T)$ le diviseur de $X'$ correspondant.
On note $\Y$ l'ouvert de $\X$ complémentaire de $T$ et on désigne par
$E$ un isocristal sur $Y$ surconvergent le long de $T$.
Le dual de $E$ sera noté $E ^\vee$.
De même pour tout
$\O _{\X } (\hdag T) _\Q$-module $\E$, $\E ^\vee$ sera son dual $\O _{\X } (\hdag T) _\Q$-linéaire.

Nous avons construit dans \cite[2.2.12]{caro_comparaison}, l'isomorphisme canonique de commutation
aux foncteurs duaux suivant :
$\sp _* ( E ^\vee) \riso \DD ^\dag _T \sp _* (E)$.
Nous vérifions dans cette section que celui-ci commute aux images inverses (extraordinaires)
respectives par $f$, i.e., que le diagramme \ref{spEdualf!diagnew} est commutatif.

Comme le cas de la restriction à un ouvert est immédiat,
sauf pour \ref{corospEdualf!},
on supposera que $f$ est un isomorphisme.

\subsection{Images directes et images inverses extraordinaires dans le cas d'un isomorphisme}

\begin{vide}\label{defrhoisoann}
    On note $\smash{\D} _{\X } (\hdag T ) _\Q:= \O  _{\X } (\hdag T ) _\Q \otimes _{\O  _{\X ,\Q}}\smash{\D} _{\X ,\Q} $
    (de même avec des primes).
    Le morphisme $\smash{\D} _{\X '} (\hdag T ') _\Q$-linéaire à gauche canonique
  $\smash{\D} _{\X '} (\hdag T ') _\Q  \rightarrow f ^* \smash{\D} _{\X } (\hdag T) _\Q$ (voir \cite[2.1.3.1]{Be2})
  est un isomorphisme. On note $\iota$, l'isomorphisme composé :
$f ^{-1} \smash{\D} _{\X } (\hdag T) _\Q \riso f ^* \smash{\D} _{\X } (\hdag T) _\Q \liso \smash{\D} _{\X '} (\hdag T ') _\Q $.
\end{vide}

\begin{prop}\label{rhoisoann}
  L'isomorphisme
  $\iota$ : $f ^{-1} \smash{\D} _{\X } (\hdag T) _\Q \riso \smash{\D} _{\X '} (\hdag T ') _\Q $
  est un isomorphisme d'anneaux.
\end{prop}
\begin{proof}
Pour tout entier $n$, on écrit $\PP _{\X,n}:= \O _{\X \times \X} / \I ^{n+1}$, où
$\I$ désigne l'idéal de l'immersion diagonale de $\X$. On posera
$\widetilde{\PP} _{\X,n}:= \O _{\X} (\hdag T) _{\Q}  \otimes _{\O _{\X}} \PP _{\X,n}$
et $\widetilde{\D} _{\X,n}$ son dual $\O _{\X} (\hdag T) _{\Q}$-linéaire pour la structure gauche
(de même en rajoutant des primes).

On remarque d'abord que l'isomorphisme $\iota$
envoie $f ^{-1}\widetilde{\D} _{\X,n}$ sur $\widetilde{\D} _{\X',n}$.
En effet, cela dérive du diagramme commutatif :
\begin{equation}
  \label{preP1PP'}
  \xymatrix  @R=0,3cm  {
{f ^{-1} \smash{\D} _{\X } (\hdag T) _\Q }
\ar[r] _-\sim
&
{f ^{*} \smash{\D} _{\X } (\hdag T) _\Q }
&
{\smash{\D} _{\X '} (\hdag T ') _\Q }
\ar[l] ^-\sim
\\
{f ^{-1} \widetilde{\D} _{\X,n}}
\ar[r] _-\sim
\ar@{^{(}->}[u]
&
{f ^{*} \widetilde{\D} _{\X,n}}
\ar@{^{(}->}[u]
&
{\widetilde{\D} _{\X',n},}
\ar[l]  _-\sim      \ar@{^{(}->}[u]
}
\end{equation}
dont la flèche en bas à droite découle par dualité du morphisme canonique
$f^* \widetilde{\PP} _{\X,n} \rightarrow \widetilde{\PP} _{\X',n}$.

Soient $n$, $n'$ deux entiers,
$P \in f ^{-1} \widetilde{\D} _{\X ,n} $, $Q \in f ^{-1} \widetilde{\D} _{\X ,n'} $.
En notant $1\otimes P$ : $f ^* \widetilde{\PP} _{\X,n} \rightarrow \O _{\X'} (\hdag T') _{\Q}$
l'image de $P$ par
$f ^{-1}\widetilde{\D} _{\X,n} \riso f ^{*}\widetilde{\D} _{\X,n} \riso
\mathcal{H} om _{\O _{\X'} (\hdag T') _{\Q}} ( f ^* \widetilde{\PP} _{\X,n},\O _{\X'} (\hdag T') _{\Q})$,
via \ref{preP1PP'},
on obtient tautologiquement le diagramme commutatif :
\begin{equation}\label{P1PP'}
\xymatrix  @R=0,3cm   {
  {\widetilde{\PP} _{\X',n}} \ar[rd] _{\iota (P)}
  &
  { f ^* \widetilde{\PP} _{\X,n}} \ar[d] ^-{1 \otimes P} \ar[l] ^-{\sim}
  &
  { f ^{-1} \widetilde{\PP} _{\X,n}} \ar[d] ^-{P} \ar[l] ^-{\sim}
  \\
  &
  {\O _{\X'} (\hdag T') _{\Q}}
  &
  {f ^{-1} \O _{\X} (\hdag T) _{\Q}.} \ar[l] ^\sim
  }
\end{equation}
Considérons maintenant le diagramme suivant :
\begin{equation}
  \label{fisoDD'ann}
  \xymatrix  @R=0,3cm   {
{f ^{-1} \widetilde{\PP} _{\X ,n+n'}}
\ar[d] ^\wr \ar[r] ^-{\smash{\widetilde{\delta} } ^{n,n'}}
&
{f ^{-1} \widetilde{\PP} _{\X ,n} \otimes _{f ^{-1}\O _{\X} (\hdag T) _{\Q}}f ^{-1} \widetilde{\PP} _{\X ,n'}}
\ar[r] ^-{Q} \ar[d] ^\wr
&
{f ^{-1} \widetilde{\PP} _{\X ,n} }
\ar[r] ^-{P} \ar[d] ^\wr
&
{f ^{-1} \O _{\X} (\hdag T) _{\Q}}
\ar[d] ^\wr
\\
{\widetilde{\PP} _{\X ' ,n+n'}}
\ar[r] ^-{\smash{\widetilde{\delta} } ^{n,n'}}
&
{\widetilde{\PP} _{\X' ,n} \otimes _{\O _{\X'} (\hdag T') _{\Q}}\widetilde{\PP} _{\X ',n'}}
\ar[r] ^-{\iota (Q)}
&
{ \widetilde{\PP} _{\X ',n} }
\ar[r] ^-{\iota (P)}
&
{\O _{\X'} (\hdag T') _{\Q},}
}
\end{equation}
dont les isomorphismes $\smash{\widetilde{\delta} } ^{n,n'}$
sont ceux analogues à \cite[1.1.10]{caro_comparaison}.
Par définition, le composé du haut (resp. du bas) donne $P\cdot Q$ (resp. $\iota (P) \cdot \iota (Q)$).
Ainsi, par \ref{P1PP'},
il s'agit d'établir la commutativité de \ref{fisoDD'ann}.
On déduit de \ref{P1PP'} la commutativité des carrés du milieu et de droite.
Quant à celle du carré de gauche, cela dérive d'un calcul immédiat
(on utilise la caractérisation de $\widetilde{\delta} ^{n,n'}$ :
$\widetilde{\delta} ^{n,n'} (a \otimes b)= (a \otimes 1) \otimes (1 \otimes b)$).
D'où le résultat.
\end{proof}

\begin{vide}\label{f*=f!}
On remarque que le diagramme
$$\xymatrix  @R=0,3cm  {
{f ^{-1} \O _{\X } (\hdag T) _\Q}
\ar[r] _-\sim \ar[d]
&
{\O _{\X '} (\hdag T') _\Q}
\ar[d]
\\
{f ^{-1}\smash{\D}  _{\X } (\hdag T) _\Q}
\ar[r] _-\sim ^{\iota}
&
{\smash{\D} _{\X '} (\hdag T') _\Q}
}$$
est commutatif.
On dispose ainsi,
pour tout $\E \in D _\mathrm{coh} ^\mathrm{b} (\overset{^\mathrm{g}}{} \smash{\D}  _{\X } ( \hdag T ) _{\Q})$,
de l'isomorphisme canonique
$f ^* (\E) =
\O _{\X '} (\hdag T') _\Q \otimes _{f ^{-1} \O _{\X } (\hdag T) _\Q} \E
\riso
\smash{\D} _{\X '} (\hdag T') _\Q \otimes _{f ^{-1} \smash{\D} _{\X } (\hdag T) _\Q} \E$.
De même en remplaçant $\smash{\D}$ par $\smash{\D} ^\dag$.

En voyant $\smash{\D} _{\X '} (\hdag T ') _\Q $ comme un
$(\smash{\D} _{\X '} (\hdag T') _\Q , f ^{-1} \smash{\D} _{\X } (\hdag T) _\Q )$-bimodule, on obtient
l'isomorphisme canonique
$\smash{\D} _{\X '} (\hdag T') _\Q \riso \smash{\D} _{\X '\rightarrow \X} (\hdag T) _\Q$ de
$(\smash{\D} _{\X '} (\hdag T') _\Q , f ^{-1} \smash{\D} _{\X } (\hdag T) _\Q )$-bimodules.
Dans la suite, on pourra identifier canoniquement
$\smash{\D} _{\X '\rightarrow \X} (\hdag T) _\Q $
et $\smash{\D} _{\X '} (\hdag T') _\Q $.
\end{vide}

\begin{prop}\label{f-1DDomega}
Les isomorphismes canoniques
\begin{gather}\label{f-1DDomegaiso1}
  f ^{-1} \omega _{\X} \otimes _{f ^{-1} \O _{\X}} \smash{\D} _{\X '} (\hdag T') _\Q
\riso
\omega _{\X'} \otimes _{\O _{\X'}} \smash{\D} _{\X '} (\hdag T') _\Q ,\\
\label{f-1DDomegaiso2}
\smash{\D} _{\X '} (\hdag T') _\Q\otimes _{f ^{-1} \O _{\X}} f ^{-1} \omega _{\X} ^{-1}
\riso
\smash{\D} _{\X '} (\hdag T') _\Q\otimes _{\O _{\X'}} \omega _{\X'}  ^{-1} ,\\
\label{f-1DDomegaiso3}
\omega _{\X '} \otimes _{\O _{\X '}} \smash{\D} _{\X '} (\hdag T') _\Q\otimes _{\O _{\X'}} \omega _{\X'}  ^{-1}
\riso \smash{\D} _{\X \leftarrow \X'} (\hdag T) _\Q,
\end{gather}
sont
$f ^{-1} \smash{\D} _{\X } (\hdag T) _\Q$-bilinéaires pour les deux premiers et
$(f ^{-1} \smash{\D} _{\X } (\hdag T) _\Q , \smash{\D} _{\X '} (\hdag T') _\Q )$-bilinéaire
pour le dernier.
\end{prop}
\begin{proof}
La proposition a bien un sens. En effet,
en considérant le faisceau $\smash{\D} _{\X '} (\hdag T') _\Q$ comme un $(f ^{-1} \smash{\D} _{\X } (\hdag T) _\Q , \smash{\D} _{\X '} (\hdag T') _\Q )$-bimodule,
$f ^{-1} \omega _{\X} \otimes _{f ^{-1} \O _{\X}} \smash{\D} _{\X '} (\hdag T') _\Q$ devient
un $(f ^{-1} \smash{\D} _{\X } (\hdag T) _\Q , \smash{\D} _{\X '} (\hdag T') _\Q )$-bimodule à droite et donc,
grâce à $\iota$,
un $\smash{\D} _{\X '} (\hdag T') _\Q $-bimodule à droite. De même, pour les deux autres isomorphismes
\ref{f-1DDomegaiso2} et \ref{f-1DDomegaiso3}.

La proposition est locale.
Supposons donc $\X$ affine muni de coordonnées locales $t_1,\dots, t_d$.
Notons $t' _1,\dots ,t' _d$ les coordonnées locales
sur $\X'$ correspondantes, $\partial _1,\dots, \partial_d$ et $\partial _1',\dots, \partial_d '$
les dérivations associées.
On peut alors identifer $\omega _{\X}$ à $\O _{\X}$ grâce à la base $d t _1 \wedge \dots \wedge d t _d$,
et
$\omega _{\X} ^{-1}$ à $\O _{\X}$ grâce à la base duale.
De même en remplaçant $\X$ par $\X '$ et $t$ par $t'$. Ces identifications sont compatibles
car l'image de $d t _1 \wedge \dots \wedge d t _d$
par l'isomorphisme canonique
$f ^{-1} \omega _{\X} \otimes _{f ^{-1} \O _{\X}} \O _{\X'} \riso \omega _{\X'}$
est $d t' _1 \wedge \dots \wedge d t' _d$ (de même en passant au dual).

Avec ces identifications,
l'action par la structure tordue
se calcule avec l'{\it adjoint} (voir \cite[1.2.2.1]{Be2}).
En effet, cela découle de la formule essentielle \cite[1.2.3]{Be2}.

La flèche $\iota$ envoie $\underline{\partial} ^{[\underline{k}]}:= \partial _1 ^{k_1} \cdots \partial _d ^{k_d}$
sur $\underline{\partial} ^{\prime [\underline{k}]}:= \partial  _1 ^{\prime k_1} \cdots \partial _d ^{\prime k_d}$
et $\sum _{\underline{k}} a _{\underline{k}} \underline{\partial} ^{[\underline{k}]}$
sur $\sum _{\underline{k}} (1\otimes a _{\underline{k}} ) \underline{\partial} ^{'[\underline{k}]}$,
où $1\otimes a _{\underline{k}}$ est l'image canonique de $a _{\underline{k}} \in f ^{-1} \O _{\X } (\hdag T ) _\Q$ dans
$\O _{\X '} (\hdag T ') _\Q$. Puisque $\iota $ est un homomorphisme d'anneaux, on calcule que pour tout
$P \in f ^{-1} \smash{\D} _{\X } (\hdag T) _\Q $, $\overset{\mathrm{t}}{} \iota (P) = \iota (\overset{\mathrm{t}}{} P)$,
où {\og $\mathrm{t}$ \fg} désigne l'adjoint.

\end{proof}

\begin{vide}
  De manière analogue à \ref{rhoisoann}, on prouve que l'isomorphisme composé
  $f ^{-1} \smash{\D} ^{(m)} _{X _i} (T) \riso f ^* \smash{\D} ^{(m)} _{X _i} (T) \liso \smash{\D} ^{(m)} _{X '_i} (T')$
  est un isomorphisme d'anneaux. Il en résulte, avec les notations
  de \cite[1.1.2]{caro_courbe-nouveau}, par complétion et passage à la limite sur le niveau,
  que
  $f ^{-1} \smash{\D} ^\dag _{\X } (\hdag T) _\Q \riso \smash{\D} ^\dag _{\X '\rightarrow \X} (\hdag T) _\Q
  \liso \smash{\D} ^\dag _{\X '} (\hdag T ') _\Q $
  est un isomorphisme d'anneaux.
  Les isomorphismes de \ref{f-1DDomegaiso1} sont encore valables en remplaçant
  $\smash{\D}$ par $\smash{\D} ^\dag$.
  Ainsi, l'isomorphisme canonique
$\omega _{\X '} \otimes _{\O _{\X '}} \smash{\D} ^\dag _{\X '} (\hdag T') _\Q\otimes _{\O _{\X'}} \omega _{\X'}  ^{-1}
\riso \smash{\D} ^\dag _{\X \leftarrow \X'} (\hdag T) _\Q$ est
$(f ^{-1} \smash{\D} ^\dag _{\X } (\hdag T) _\Q,\smash{\D} ^\dag _{\X '} (\hdag T') _\Q)$-bilinéaire.
  Par la suite, on identifiera
  $\smash{\D} ^\dag _{\X '\rightarrow \X } (\hdag T) _\Q$ à $\smash{\D} ^\dag _{\X '} (\hdag T ') _\Q $
  et
  $\smash{\D} ^\dag _{\X \leftarrow \X'} (\hdag T) _\Q$
à
$\omega _{\X '} \otimes _{\O _{\X '}} \smash{\D} ^\dag _{\X '} (\hdag T') _\Q\otimes _{\O _{\X'}} \omega _{\X'}  ^{-1}$.
\end{vide}

\begin{nota}\label{identification}
Pour alléger les notations, on pose
$\widetilde{\O}  _{\X}:= \O  _{\X} (\hdag T ) _{\Q}$
$\widetilde{\omega} _{\X}  := \omega _{\X} \otimes _{\O _{\X}} \widetilde{\O}  _{\X}$,
$\smash{\widetilde{\D}}  _{\X} : = \smash{\D} ^\dag _{\X} (\hdag T) _{\Q}$
(par défaut) ou $\smash{\D}  _{\X} (\hdag T) _{\Q}$,
$\widetilde{\D}  _{\X' \rightarrow \X }:=
\underset{\longrightarrow}{\lim} _m \widehat{\B} ^{(m)}  _{\X'} (T')
\widehat{\otimes} _{\O _{\X'}}\widehat{\D} ^{(m)} _{\X' \rightarrow \X,\Q} $
(par défaut) ou $\O   _{\X'} (\hdag T') _{\Q}\otimes _{\O _{\X',\Q}} \smash{\D}_{\X' \rightarrow \X,\Q}  $.
De même en rajoutant des primes.

On notera
  $(\widetilde{\D} _{\X} \otimes _{\widetilde{\O} _{\X}} \widetilde{\omega} _{\X } ^{-1}) _\mathrm{t}$
  le bimodule déduit de $\widetilde{\D} _{\X} \otimes _{\widetilde{\O} _{\X}} \widetilde{\omega} _{\X } ^{-1}$
  via l'isomorphisme de transposition $\beta$ de \cite[1.3.4.3]{Be2} (en d'autres termes, la structure gauche est
  la structure droite et vis versa).
De même,
$(\widetilde{\omega} _{\X } \otimes _{\widetilde{\O} _{\X}} \widetilde{\D} _{\X} ) _\mathrm{t}$
désignera le bimodule déduit de
$\widetilde{\omega} _{\X } \otimes _{\widetilde{\O} _{\X}} \widetilde{\D} _{\X} $
  via l'isomorphisme de transposition $\delta$ de \cite[1.3.4.1]{Be2}.
De cette façon, les isomorphismes
$    \beta \  : \
 \widetilde{\D} _{\X} \otimes _{\widetilde{\O} _{\X}} \widetilde{\omega} _{\X } ^{-1}
  \riso
  (\widetilde{\D} _{\X} \otimes _{\widetilde{\O} _{\X}} \widetilde{\omega} _{\X } ^{-1}) _\mathrm{t}$
  et
$\delta \ :\
\widetilde{\omega} _{\X } \otimes _{\widetilde{\O} _{\X}} \widetilde{\D} _{\X}
\riso
(\widetilde{\omega} _{\X } \otimes _{\widetilde{\O} _{\X}} \widetilde{\D} _{\X} ) _\mathrm{t}$
sont bilinéaires.
\end{nota}

\begin{vide}
Soient $\M \in D _\mathrm{coh} ^\mathrm{b} ( \smash{\D} ^{\dag} _{\X } ( \hdag T ) _{\Q} \overset{^{\mathrm{d}}}{})$,
$\M '    \in D _\mathrm{coh} ^\mathrm{b} ( \smash{\D} ^{\dag} _{\X '} ( \hdag T ') _{\Q}\overset{^{\mathrm{d}}}{})$,
$\E \in D _\mathrm{coh} ^\mathrm{b} (\overset{^\mathrm{g}}{} \smash{\D} ^{\dag} _{\X } ( \hdag T ) _{\Q})$,
$\E '    \in D _\mathrm{coh} ^\mathrm{b} (\overset{^\mathrm{g}}{} \smash{\D} ^{\dag} _{\X '} ( \hdag T ') _{\Q} )$.
Via les deux identifications de \ref{identification}, on obtient :
\begin{align}
  \notag
f ^! _T (\E) = &\smash{\D} ^\dag _{\X '} (\hdag T') _\Q \otimes _{f ^{-1} \smash{\D} ^\dag _{\X } (\hdag T) _\Q} f ^{-1} \E,
\\ \notag
f ^! _T (\M) = &f ^{-1} \M \otimes _{f ^{-1} \smash{\D} ^\dag _{\X } (\hdag T) _\Q}
(\omega _{\X '} \otimes _{\O _{\X '}} \smash{\D} ^\dag _{\X '} (\hdag T') _\Q\otimes _{\O _{\X'}} \omega _{\X'}  ^{-1}),
\\ \notag
f _{T,+} (\E') = & f _* (
(\omega _{\X '} \otimes _{\O _{\X '}} \smash{\D} ^\dag _{\X '} (\hdag T') _\Q\otimes _{\O _{\X'}} \omega _{\X'}  ^{-1})
\otimes _{\smash{\D} ^\dag _{\X '} (\hdag T') _\Q}\E'),
\\ \notag
f _{T,+} (\M') = &f _* (\M').
\end{align}

On dispose des isomorphismes :
\begin{gather}\label{f!dg}
  \omega _{\X '} \otimes _{\O _{\X '}} f ^! _T (\E) \riso
f ^! _T (\omega _{\X } \otimes _{\O _{\X }}\E), \
 f ^! _T (\M) \otimes _{\O _{\X '}} \omega _{\X '} ^{-1}
\riso
f ^! _T ( \M   \otimes _{\O _{\X }} \omega _{\X } ^{-1}   ),
\\
\label{f+dg}
\omega _{\X } \otimes _{\O _{\X }}f _{T,+} (\E')
\riso
f _{T,+} (\omega _{\X '} \otimes _{\O _{\X '}} \E'),\
f _{T,+} (\M ') \otimes _{\O _{\X }} \omega _{\X } ^{-1}
\riso
f _{T,+} ( \M ' \otimes _{\O _{\X '}} \omega _{\X '} ^{-1} )
\end{gather}
(pour le passage de droite à gauche, on pourra consulter \cite[I.2]{virrion}).
Décrivons localement ces isomorphismes.
Supposons donc $\X$ affine muni de coordonnées locales $t_1,\dots, t_d$.
Notons $t' _1,\dots ,t' _d$ les coordonnées locales
sur $\X'$ correspondantes et identifions
$\omega _{\X}$ à $\O _{\X}$ grâce à la base $d t _1 \wedge \dots \wedge d t _d$,
et
$\omega _{\X} ^{-1}$ à $\O _{\X}$ grâce à la base duale et de même avec des primes.

Avec ces identifications,
on calcule d'abord que l'isomorphisme canonique (voir sa construction dans \cite[I.2.1.2]{virrion})
$(\widetilde{\D} _{\X} \otimes _{\widetilde{\O} _{\X}} \widetilde{\omega} _{\X } ^{-1} ) _\mathrm{t}
\otimes _{\widetilde{\D} _{\X}}
(\widetilde{\omega} _{\X } \otimes _{\widetilde{\O} _{\X}} \widetilde{\D} _{\X} ) _\mathrm{t}
\riso \widetilde{\D} _{\X}$
envoie
$P \otimes Q$ sur $Q  P$,
avec $P$ et $Q$ des sections locales de $\widetilde{\D} _{\X}$.
En particulier $1 \otimes 1$ s'envoie sur $1$.
Il en découle
que l'isomorphisme canonique de \cite[I.2.2.(i)]{virrion}
$\M \otimes _{\widetilde{\D} _{\X}} \E \riso
(\widetilde{\omega} _{\X } \otimes _{\widetilde{\O} _{\X}} \E)
\otimes _{\widetilde{\D} _{\X}}
(\M \otimes _{\widetilde{\O} _{\X}} \widetilde{\omega} _{\X } ^{-1} )$
expédie $x \otimes e$ sur $e \otimes x$, où $x$ et $e$ sont des sections locales
respectives de $\M$ et $\E$.
Or, l'isomorphisme
$\omega _{\X '} \otimes _{\O _{\X '}} f ^! _T (\E) \riso
f ^! _T (\omega _{\X } \otimes _{\O _{\X }}\E)$ est le composé suivant :
\begin{gather}
  \notag
  \omega _{\X '} \otimes _{\O _{\X '}}
(\widetilde{\D} _{\X'} \otimes _{f ^{-1} \widetilde{\D} _{\X}}  f ^{-1} \E)
\riso
\omega _{\X '} \otimes _{\O _{\X '}}
(f ^{-1} (\widetilde{\omega} _{\X } \otimes _{\widetilde{\O} _{\X}} \E)
\otimes ^{\mathrm{d} } _{f ^{-1} \widetilde{\D} _{\X}}
(\widetilde{\D} _{\X'} \otimes _{\widetilde{\O} _{\X'}} \widetilde{\omega} _{\X '} ^{-1} ))
\riso
\\ \label{isoomegaf!dh}
\riso
f ^{-1} (\widetilde{\omega} _{\X } \otimes _{\widetilde{\O} _{\X}} \E)
\otimes  _{f ^{-1} \widetilde{\D} _{\X}}
(\omega _{\X '} \otimes _{\O _{\X '}}
\widetilde{\D} _{\X'} \otimes _{\widetilde{\O} _{\X'}} \widetilde{\omega} _{\X '} ^{-1} )
,
\end{gather}
où le symbole {\og $\mathrm{d}$ \fg} signifie que l'on a pris,
pour calculer le produit tensoriel, la structure droite (i.e. la structure tordue)
de $\widetilde{\D} _{\X'} \otimes _{\widetilde{\O} _{\X'}} \widetilde{\omega} _{\X '} ^{-1} $,
et où le premier isomorphisme découle de
\cite[I.2.2.(i)]{virrion}.
Le composé
\ref{isoomegaf!dh} envoie
$P' \otimes x$ sur $x \otimes P'$, où $P'$ et $x$ sont des sections locales respectives de
$\widetilde{\D} _{\X'} $ et $\E$. On bénéficie de descriptions analogues pour
les trois autres isomorphismes de \ref{f!dg} et \ref{f+dg}.

De plus, on définit les opérations analogues en remplaçant {\og $\smash{\D} ^\dag$\fg} par {\og $\smash{\D} $\fg}.
Si un risque de confusion est à craindre,
on notera $f ^{! ^\dag} _T$ (resp. $f _{T,+} ^\dag$) pour $f ^! _T $ (resp. $f _{T,+}$).
Pour tout $\E \in D _\mathrm{coh} ^\mathrm{b} (\overset{^\mathrm{g}}{} \smash{\D}  _{\X } ( \hdag T ) _{\Q})$,
on dispose d'un isomorphisme canonique
\begin{equation}
  \label{extf!}
  f ^{! ^\dag} _T ( \smash{\D} ^{\dag} _{\X } ( \hdag T ) _{\Q}
\otimes _{\smash{\D}  _{\X } ( \hdag T ) _{\Q}} \E ) \riso
\smash{\D} ^{\dag} _{\X' } ( \hdag T' ) _{\Q}
\otimes _{\smash{\D}  _{\X '} ( \hdag T' ) _{\Q}} f ^! _T (\E).
\end{equation}

\begin{vide}
  \label{fdomegabeta}
  On dispose des isomorphismes canoniques
  $f ^! _{T,\mathrm{g}} (\widetilde{\D} _{\X} \otimes _{\widetilde{\O} _{\X}} \widetilde{\omega} _{\X } ^{-1})
  \riso
  \widetilde{\D} _{\X'} \otimes _{\widetilde{\O} _{\X '}} \widetilde{\omega} _{\X '} ^{-1}$
  et
  $f ^! _{T,\mathrm{d}} ((\widetilde{\D} _{\X} \otimes _{\widetilde{\O} _{\X}} \widetilde{\omega} _{\X } ^{-1}) _{\mathrm{t}})
  \riso
  (\widetilde{\D} _{\X'} \otimes _{\widetilde{\O} _{\X '}} \widetilde{\omega} _{\X '} ^{-1} ) _{\mathrm{t}} $.
 De plus, on bénéficie du composé suivant :
  \begin{gather}
\notag
f ^! _{T,\mathrm{d}} (\widetilde{\D} _{\X} \otimes _{\widetilde{\O} _{\X}} \widetilde{\omega} _{\X } ^{-1})
=
 \widetilde{\D} _{\X'}
 \otimes ^{\mathrm{d}\mathrm{d}} _{ f^{-1} \widetilde{\D} _{\X} }
f^{-1} (\widetilde{\D} _{\X} \otimes _{\widetilde{\O} _{\X}} \widetilde{\omega} _{\X } ^{-1})
\riso
 \\
 \label{fdomegabetaiso1}
 \riso
f^{-1} \widetilde{\D} _{\X}
 \otimes ^{\mathrm{d} \mathrm{d}} _{ f^{-1} \widetilde{\D} _{\X} }
 ( \widetilde{\D} _{\X'} \otimes _{\widetilde{\O} _{\X'}} \widetilde{\omega} _{\X '} ^{-1})
 \riso
( \widetilde{\D} _{\X'} \otimes _{\widetilde{\O} _{\X'}} \widetilde{\omega} _{\X '} ^{-1} ) _{\mathrm{t}},
  \end{gather}
où les symboles {\og $\mathrm{d}$ \fg} signifient que l'on choisit dans le calcul des produits tensoriels
les structures droites des bimodules respectifs, et où le premier isomorphisme découle de
\cite[I.2.2.(i)]{virrion}.
De manière analogue à \ref{fdomegabetaiso1}, on construit
$f ^! _{T,\mathrm{g}} (  (  \widetilde{\D} _{\X} \otimes _{\widetilde{\O} _{\X}} \widetilde{\omega} _{\X } ^{-1})_{\mathrm{t}})
\riso
\widetilde{\D} _{\X'} \otimes _{\widetilde{\O} _{\X'}} \widetilde{\omega} _{\X '} ^{-1} $.

Par un calcul local, on vérifie que les diagrammes
\begin{equation}
    \label{fdomegabetadiag1}
    \xymatrix  @R=0,3cm  {
{f ^! _{T,\mathrm{g}} (   \widetilde{\D} _{\X} \otimes _{\widetilde{\O} _{\X}} \widetilde{\omega} _{\X } ^{-1} )}
    \ar[r] _-\sim ^{\beta} \ar[d] _-\sim
    &
    {f ^! _{T,\mathrm{g}} ((\widetilde{\D} _{\X} \otimes _{\widetilde{\O} _{\X}} \widetilde{\omega} _{\X } ^{-1})_{\mathrm{t}})}
    \ar[d] _-\sim
    \\
    {\widetilde{\D} _{\X'} \otimes _{\widetilde{\O} _{\X'}} \widetilde{\omega} _{\X '} ^{-1} }
\ar@{=}[r]
&
    { \widetilde{\D} _{\X'} \otimes _{\widetilde{\O} _{\X'}} \widetilde{\omega} _{\X '} ^{-1} ,}
}
\
    \xymatrix  @R=0,3cm  {
{f ^! _{T,\mathrm{d}} (  (  \widetilde{\D} _{\X} \otimes _{\widetilde{\O} _{\X}} \widetilde{\omega} _{\X } ^{-1})_{\mathrm{t}})}
    \ar[r] _-\sim ^{\beta} \ar[d] _-\sim
    &
    {f ^! _{T,\mathrm{d}} (\widetilde{\D} _{\X} \otimes _{\widetilde{\O} _{\X}} \widetilde{\omega} _{\X } ^{-1})}
    \ar[d] _-\sim
    \\
    {( \widetilde{\D} _{\X'} \otimes _{\widetilde{\O} _{\X'}} \widetilde{\omega} _{\X '} ^{-1} ) _{\mathrm{t}}}
\ar@{=}[r]
&
    {( \widetilde{\D} _{\X'} \otimes _{\widetilde{\O} _{\X'}} \widetilde{\omega} _{\X '} ^{-1} ) _{\mathrm{t}},}
}
\end{equation}
sont commutatifs.
Lorsque $f $ est l'identité, les flèches de gauche des deux carrés sont égaux à l'identité et celles de droite
s'identifient à $\beta$.

\end{vide}

\end{vide}

\subsection{Sur la commutation des foncteurs duaux aux images inverses par un isomorphisme}

\begin{nota}\label{nota-morptraceOdagTpre}
Commençons par quelques notations et rappels sur les foncteurs duaux.

  Pour tout $\E \in D _{\mathrm{parf}} ( \overset{^\mathrm{g}}{} \smash{\D}  _{\X} ( \hdag T ) _{\Q} )$,
le complexe dual $\smash{\D}  _{\X} ( \hdag T ) _{\Q}$-linéaire de $\E$
au sens de Virrion (voir \cite[I.3.2]{virrion}) est défini en posant :
  \begin{gather}\label{DDpasdag}
 \DD  _{\X , T } (\E):= \R \mathcal{H} om _{\smash{\D}  _{\X} ( \hdag T ) _{\Q}}
  (\E , \smash{\D}  _{\X} ( \hdag T ) _{\Q} \otimes _{\O _{\X}} \omega _{\X } ^{-1}  ) [d _{X}].
  \end{gather}

De même, pour tous $\E \in D _{\mathrm{parf}} ( \overset{^\mathrm{g}}{} \smash{\D} ^{\dag} _{\X  ,\Q} ( \hdag T ) )$,
  $\M \in D _{\mathrm{parf}} ( \smash{\D} ^{\dag} _{\X  ,\Q} ( \hdag T ) \overset{^{\mathrm{d}}}{})$,
  leurs complexes duaux $\smash{\D} ^{\dag} _{\X  ,\Q} ( \hdag T )$-linéaires sont définis comme suit :
  \begin{gather}\notag
 \DD ^\dag _{\X , T } (\E):= \R \mathcal{H} om _{\widetilde{\D} _{\X}}
  (\E , \widetilde{\D} _{\X} \otimes _{\widetilde{\O} _{\X}} \widetilde{\omega} _{\X } ^{-1}  ) [d _{X}],
  \\
  \notag
  \DD ^{\prime \dag} _{\X , T } (\E):=
  \R \mathcal{H} om _{\widetilde{\D} _{\X}}
  (\E , (\widetilde{\D} _{\X} \otimes _{\widetilde{\O} _{\X}} \widetilde{\omega} _{\X } ^{-1} ) _{\mathrm{t}}  ) [d _{X}],
  \\
  \notag
  \DD ^\dag _{\X , T } (\M):=
  \widetilde{\omega} _{\X } \otimes _{\widetilde{\O} _{\X}}  \R \mathcal{H} om _{\widetilde{\D} _{\X}}
  (\M , \widetilde{\D} _{\X} ) [d _{X}]
  \riso
\R \mathcal{H} om _{\widetilde{\D} _{\X}}
  (\M , (\widetilde{\omega} _{\X } \otimes _{\widetilde{\O} _{\X}}   \widetilde{\D} _{\X} ) _{\mathrm{t}}) [d _{X}],
  \\
  \notag
  \DD ^{\prime \dag}  _{\X , T } (\M):=
\R \mathcal{H} om _{\widetilde{\D} _{\X}}
  (\M , \widetilde{\omega} _{\X } \otimes _{\widetilde{\O} _{\X}}   \widetilde{\D} _{\X} ) [d _{X}].
  \end{gather}
  Le foncteur dual commutant à l'extension des scalaires (\cite{virrion}),
  pour tout $\E \in D _{\mathrm{parf}} ( \overset{^\mathrm{g}}{} \smash{\D}  _{\X} ( \hdag T ) _{\Q} )$,
\begin{equation}
  \label{extdual}
  \DD ^\dag _{\X,T} ( \smash{\D}  ^\dag _{\X} ( \hdag T ) _{\Q} \otimes _{\smash{\D}  _{\X} ( \hdag T ) _{\Q}} \E)
  \riso \smash{\D}  ^\dag _{\X} ( \hdag T ) _{\Q} \otimes _{\smash{\D}  _{\X} ( \hdag T ) _{\Q}} \DD  _{\X,T} (\E).
\end{equation}
  Si aucune confusion avec \ref{DDpasdag} n'est à craindre,
  on omettra d'indiquer le symbole {\og $\dag$\fg}.

  On dispose des isomorphismes canoniques :
  \begin{gather}\label{DDdg}
    \DD ' _{\X , T } (\M) \otimes _{\widetilde{\O} _{\X}} \widetilde{\omega} _{\X } ^{-1}
  \riso
  \R \mathcal{H} om _{\overset{^\mathrm{g}}{} \widetilde{\D} _{\X}}
  (\M \otimes _{\widetilde{\O} _{\X}} \widetilde{\omega} _{\X } ^{-1}, \widetilde{\D} _{\X} )
  \otimes _{\widetilde{\O} _{\X}} \widetilde{\omega} _{\X } ^{-1} [d _{X}]
  \riso
  \DD _{\X , T } (\M \otimes _{\widetilde{\O} _{\X}} \widetilde{\omega} _{\X } ^{-1}),
  \\
  \label{DDgd}
  \widetilde{\omega} _{\X }  \otimes _{\widetilde{\O} _{\X}} \DD ' _{\X , T } (\E)
  \riso
  \widetilde{\omega} _{\X }  \otimes _{\widetilde{\O} _{\X}}
  \R \mathcal{H} om _{\widetilde{\D} _{\X} \overset{^\mathrm{d}}{}}
  (\widetilde{\omega} _{\X }  \otimes _{\widetilde{\O} _{\X}} \E , \widetilde{\D} _{\X} ) [d _{X}]
  \riso
  \DD _{\X , T } (\widetilde{\omega} _{\X }  \otimes _{\widetilde{\O} _{\X}} \E ).
  \end{gather}
  Les isomorphismes \ref{DDdg} et \ref{DDgd} n'utilisent pas celui de transposition,
  ce qui facilite les calculs locaux.
  En effet, supposons $\X$ affine et muni de coordonnées locales $t_1,\dots, t_d$.
Notons $t' _1,\dots ,t' _d$ les coordonnées locales
sur $\X'$ correspondantes et identifions
$\omega _{\X}$ à $\O _{\X}$ grâce à la base $d t _1 \wedge \dots \wedge d t _d$,
et
$\omega _{\X} ^{-1}$ à $\O _{\X}$ grâce à la base duale et de même avec des primes.
Les isomorphismes \ref{DDdg} et \ref{DDgd} s'identifient alors à l'identité.

  Enfin, on remarque que dans le lemme \ref{lemmisoduareliso}
  l'isomorphisme de transposition n'apparaît pas non plus.

\end{nota}
\begin{nota}\label{nota-morptraceOdagTpre2}
$\bullet$ Pour tout $\M '\in D _{\mathrm{parf}} ( \smash{\D} ^{\dag} _{\X'  ,\Q} ( \hdag T ') \overset{^{\mathrm{d}}}{})$,
l'isomorphisme de dualité relative (\cite[1.2.7]{caro_courbe-nouveau})
 $f _{T , +} \circ \DD _{\X ' , T '} (\M ')\riso \DD _{\X , T}\circ f_{T , +} (\M ' )$
 sera noté $\chi$. En fait, via les isomorphismes de transposition $\delta$ de \cite[1.3.4.1]{Be2},
 nous avions d'abord construit l'isomorphisme
 $f _{T , +} \circ \DD '_{\X ' , T '} (\M ')\riso \DD '_{\X , T}\circ f_{T , +} (\M ' )$,
 que l'on désignera encore par $\chi$.

$\bullet$
Pour tout $\E \in D _{\mathrm{parf}} (\overset{^\mathrm{g}}{} \smash{\D} ^{\dag} _{\X  ,\Q} ( \hdag T ) )$,
on notera :
\begin{align}
  \notag
  \mathrm{proj} '\ :\
  &
  f _+ ( \widetilde{\omega} _{\X'} ) \otimes  _{\widetilde{\O} _{\X} } \E
=  f _* ( \widetilde{\omega} _{\X'} ) \otimes  _{\widetilde{\O} _{\X} } \E
  \riso
  f _* ( \widetilde{\omega} _{\X'} ) \otimes  _{\widetilde{\O} _{\X} }
  f _* ( \widetilde{\D}  _{\X' }\otimes _{f ^{-1} \widetilde{\D}  _{\X}} f ^{-1} \E)
  \\
  &
  \riso
  f _* ( \widetilde{\omega} _{\X'}  \otimes  _{\widetilde{\O} _{\X'} }
   (\widetilde{\D}  _{\X' }\otimes _{f ^{-1} \widetilde{\D}  _{\X}} f ^{-1} \E))
   =
   f _+ ( \widetilde{\omega} _{\X'} \otimes  _{\widetilde{\O} _{\X'} } f ^! _T (\E)).
   \label{projiso}
\end{align}
On obtient en particulier
$\mathrm{proj} '$ :
$ f _* ( \widetilde{\omega} _{\X'} ) \otimes _{\widetilde{\O} _{\X} }\widetilde{\D}  _{\X}
  \riso
   f _* ( \widetilde{\omega} _{\X'} \otimes _{\widetilde{\O} _{\X'} } \widetilde{\D}  _{\X'})$,
 ce dernier étant par fonctorialité $\widetilde{\D}  _{\X}$-bilinéaire.
On vérifie par un calcul que cet isomorphisme $\mathrm{proj} '$ correspond bien
à l'isomorphisme canonique lorsque $f$ est un morphisme quasi-compact et quasi-séparé quelconque
(pour la construction dans le cas général, voir \cite[1.2.27]{caro_surcoherent} ou \cite[1.4.1]{caro-frobdualrel}).

$\bullet$
Le morphisme trace : $ f _* ( \widetilde{\omega} _{\X'} ) =
f _+ ( \widetilde{\omega} _{\X'} ) \rightarrow \widetilde{\omega} _{\X}$
sera noté $\mathrm{Tr}$.

 \end{nota}

\begin{lemm}
\label{lemmisoduareliso}
Pour tout $\M '\in D _{\mathrm{parf}} ( \smash{\D} ^{\dag} _{\X'  ,\Q} ( \hdag T ') \overset{^{\mathrm{d}}}{})$,
avec les notations \ref{nota-morptraceOdagTpre} et \ref{nota-morptraceOdagTpre2},
le diagramme ci-après
\begin{equation}
\label{isoduareliso}
\xymatrix  @R=0,3cm    {
{f _* \R \mathcal{H} om _{\widetilde{\D} _{\X'}}
  (\M ', \widetilde{\omega} _{\X '} \otimes _{\widetilde{\O} _{\X'}} \widetilde{\D} _{\X'} ) [d _{X'}]}
  \ar[d] _-\sim
  \ar@{=}[r]
  &
  { f _{T,+}\circ   \DD ' _{\X ' , T '} (\M')}
  \ar[dddd] ^-\chi _-\sim
  \\
  {f _*  \R \mathcal{H} om _{f ^{-1}\widetilde{\D} _{\X}}
  (\M ', \widetilde{\omega} _{\X '} \otimes _{\widetilde{\O} _{\X'}} \widetilde{\D} _{\X'} ) [d _{X'}] }
  \ar[d] _-\sim ^-{f _*}
  \\
  {\R \mathcal{H} om _{\widetilde{\D} _{\X}}
  (f _* \M ', f _* (\widetilde{\omega} _{\X '} \otimes _{\widetilde{\O} _{\X'}} \widetilde{\D} _{\X'}) ) [d _{X'}] }
  \\
 {\R \mathcal{H} om _{\widetilde{\D} _{\X}}
  (f _* \M ', f _* (\widetilde{\omega} _{\X '} )\otimes _{\widetilde{\O} _{\X}} \widetilde{\D} _{\X}) ) [d _{X'}] }
  \ar[u] _-\sim ^{\mathrm{proj} '}
  \ar[d] _-\sim ^{\mathrm{Tr} }
  \\
    {\R \mathcal{H} om _{\widetilde{\D} _{\X}}
  ( f _* \M ' ,
  \widetilde{\omega}  _{\X} \otimes _{\widetilde{\O}  _{\X}} \widetilde{\D}  _{\X})[d _{X}] }
  \ar@{=}[r]
  &
  { \DD ' _{\X , T }\circ  f _{T,+}(\M') ,}
  }
\end{equation}
où l'isomorphisme à droite en haut dérive de
l'identification entre $\widetilde{\D} _{\X'} $ et $\widetilde{\D} _{\X'\rightarrow \X}$,
est commutatif.
\end{lemm}

\begin{proof}
Modulo l'identification entre $\widetilde{\D} _{\X'} $ et $\widetilde{\D} _{\X'\rightarrow \X}$,
le morphisme \cite[1.2.7.8]{caro_courbe-nouveau}
est égal au morphisme $\mathrm{Tr}   \circ \mathrm{proj} ^{\prime -1} $ :
$f _* ( \widetilde{\omega} _{\X'} \otimes _{\widetilde{\O} _{\X'} } \widetilde{\D} _{\X'})
\riso
\widetilde{\omega} _{\X} \otimes _{\widetilde{\O} _{\X} } \widetilde{\D} _{\X}$.
De plus, le composé de \cite[1.2.7.2]{caro_courbe-nouveau}
correspond à
$f _* \R \mathcal{H} om _{\widetilde{\D} _{\X'}}
(\M ', \widetilde{\omega} _{\X '} \otimes _{\widetilde{\O} _{\X'}} \widetilde{\D} _{\X'} ) [d _{X'}]
\riso
\R \mathcal{H} om _{\widetilde{\D} _{\X}}
 (f _* \M ', f _* (\widetilde{\omega} _{\X '} \otimes _{\widetilde{\O} _{\X'}} \widetilde{\D} _{\X'}) ) [d _{X'}]$
 de \ref{isoduareliso}. D'où le résultat par construction de l'isomorphisme $\chi $ de
\cite[1.2.7]{caro_courbe-nouveau}.

\end{proof}

\begin{vide}\label{sssection125}
Dans ce paragraphe \ref{sssection125}, supposons $f $ seulement propre.
Pour la commodité du lecteur,
rappelons que l'on dispose,
pour tous $\M '\in D _{\mathrm{parf}} ( \smash{\D} ^{\dag} _{\X ' ,\Q} ( \hdag T ') \overset{^{\mathrm{d}}}{})$,
$\E \in D _{\mathrm{parf}} (\overset{^\mathrm{g}}{} \smash{\D} ^{\dag} _{\X  ,\Q} ( \hdag T ) )$,
  de l'isomorphisme composé :
  \begin{gather}\notag
    \R f _*  \R \mathcal{H} om _{\widetilde{\D} _{\X'}}
  (\M ', f ^! _T (\widetilde{\omega} _{\X } \otimes _{\widetilde{\O} _{\X}} \E ))
  \riso
  \R f _*  \R \mathcal{H} om _{\widetilde{\D} _{\X'}}
  (\M ', \widetilde{\omega} _{\X '} \otimes _{\widetilde{\O} _{\X'}} f ^! _T ( \E ))
  \\ \notag
  =
  \R f _*  \R \mathcal{H} om _{\widetilde{\D} _{\X'}}
  (\M ', \widetilde{\omega} _{\X '} \otimes _{\widetilde{\O} _{\X'}}
  \widetilde{\D} _{\X'\rightarrow \X} \otimes ^\L _{f ^{-1}\widetilde{\D} _{\X}} f ^{-1} \E ) [d _{X'/X}]
  \\ \notag
  \underset{\mathrm{proj}_f}{\liso}
    (\R f _*  \R \mathcal{H} om _{\widetilde{\D} _{\X'}}
  (\M ', \widetilde{\omega} _{\X '} \otimes _{\widetilde{\O} _{\X'}}
  \widetilde{\D} _{\X'\rightarrow \X} )) \otimes ^\L _{\widetilde{\D} _{\X}}  \E [d _{X'/X}]
  \\ \notag
  \underset{\otimes}{\liso}
  (\R f _* ( \R \mathcal{H} om _{\widetilde{\D} _{\X'}}
  (\M ', \widetilde{\omega} _{\X '} \otimes _{\widetilde{\O} _{\X'}} \widetilde{\D} _{\X'})
  \otimes _{\widetilde{\D} _{\X'}} ^\L \widetilde{\D} _{\X'\rightarrow \X} ))
\otimes ^\L _{\widetilde{\D} _{\X}}  \E [d _{X'/X}]
\\ \notag
\underset{\chi}{\riso}
  \R \mathcal{H} om _{\widetilde{\D} _{\X}}
  (f _{T,+}  (\M ' ),
  \widetilde{\omega} _{\X } \otimes _{\widetilde{\O} _{\X}} \widetilde{\D} _{\X'})
\otimes ^\L _{\widetilde{\D} _{\X}}  \E
\\\label{adjdagt}
\underset{\otimes}{\riso}
\R \mathcal{H} om _{\widetilde{\D} _{\X}}
  (f _{T,+}  (\M ' ),   \widetilde{\omega} _{\X } \otimes _{\widetilde{\O} _{\X}}   \E).
\end{gather}
En lui appliquant $H ^0 \circ \Gamma (\X,-)$,
on obtient l'isomorphisme d'adjonction suivant :
\begin{equation}
  \label{adjdagt2}
\mathrm{adj}\ :\  \mathrm{Hom} _{\widetilde{\D} _{\X'}}
  (\M ', f ^! _T (\widetilde{\omega} _{\X } \otimes _{\widetilde{\O} _{\X}} \E ))
  \riso
 \mathrm{Hom} _{\widetilde{\D} _{\X}}
  (f _{T,+}  (\M ' ),   \widetilde{\omega} _{\X } \otimes _{\widetilde{\O} _{\X}}   \E).
\end{equation}

Lorsque $f$ est un isomorphisme,
on notera
$\mathrm{adj} _{\widetilde{\omega} _{\X} \otimes _{\widetilde{\O} _{\X} } \E}$ :
$f _{T+}  f _T ^!  ( \widetilde{\omega} _{\X} \otimes _{\widetilde{\O} _{\X} } \E )
\rightarrow \widetilde{\omega} _{\X} \otimes _{\widetilde{\O} _{\X} } \E$,
l'image de l'identité de
$f _T ^!  ( \widetilde{\omega} _{\X} \otimes _{\widetilde{\O} _{\X} } \E )$
par \ref{adjdagt2}.
Comme $f _{T+}  f _T ^! ( \widetilde{\omega} _{\X} \otimes _{\widetilde{\O} _{\X} } \E )
\riso \widetilde{\omega} _{\X} \otimes _{\widetilde{\O} _{\X} } f _{T+}  f _T ^! (  \E )$
(voir \ref{f!dg} et \ref{f+dg}),
on obtient le morphisme d'adjonction
$\mathrm{adj} _{\E}$ : $f _{T+}  f _T ^! ( \E ) \rightarrow \E$.
\end{vide}

\begin{lemm}  \label{def-alpha}
Pour tout $\E \in D _{\mathrm{parf}} (\overset{^\mathrm{g}}{} \smash{\D} ^{\dag} _{\X  ,\Q} ( \hdag T ) )$,
le morphisme composé
\begin{equation}\label{def-alphadiag}
f _{T+}  f _T ^!  ( \widetilde{\omega} _{\X} \otimes _{\widetilde{\O} _{\X} } \E ) \riso
 f _* ( \widetilde{\omega} _{\X'} \otimes _{\widetilde{\O} _{\X'} } f ^! _T (\E)) \underset{\mathrm{proj} '}{\liso}
  f _* ( \widetilde{\omega} _{\X'} ) \otimes _{\widetilde{\O} _{\X} } \E
  \underset{\mathrm{Tr}}{\riso}
\widetilde{\omega} _{\X} \otimes _{\widetilde{\O} _{\X} } \E,
\end{equation}
où le premier isomorphisme est \ref{f!dg},
est égal à $\mathrm{adj} _{\widetilde{\omega} _{\X} \otimes _{\widetilde{\O} _{\X} } \E}$.
\end{lemm}
\begin{proof}
Pour tous $\M '\in D _{\mathrm{parf}} ( \smash{\D} ^{\dag} _{\X ' ,\Q} ( \hdag T ') \overset{^{\mathrm{d}}}{})$,
$\E \in D _{\mathrm{parf}} (\overset{^\mathrm{g}}{} \smash{\D} ^{\dag} _{\X  ,\Q} ( \hdag T ) )$,
le diagramme ci-après
\begin{equation}
  \label{alphs=tr-diag}
  \xymatrix  @R=0,3cm     {
  {f _*  \R \mathcal{H} om _{\widetilde{\D} _{\X'}}
  (\M ',  f ^! _T  ( \widetilde{\omega} _{\X} \otimes _{\widetilde{\O} _{\X} } \E ))}
  \ar[r] _-\sim ^{f _*} \ar[d] _-\sim
  &
  { \R \mathcal{H} om _{\widetilde{\D} _{\X}}
  (f _*  \M ',f _*  f ^! _T  ( \widetilde{\omega} _{\X} \otimes _{\widetilde{\O} _{\X} } \E ))}
  \ar[d] _-\sim
  \\
  {f _*  \R \mathcal{H} om _{\widetilde{\D} _{\X'}}
  (\M ',  \widetilde{\omega} _{\X' } \otimes _{\widetilde{\O} _{\X'}} f ^! _T (\E ))}
  \ar[r]_-\sim ^{f _*}
  &
  { \R \mathcal{H} om _{\widetilde{\D} _{\X}}
  (f _*  \M ',f _*  (\widetilde{\omega} _{\X' } \otimes _{\widetilde{\O} _{\X'}} f ^! _T (\E )))}
  \\
  {(f _*  \R \mathcal{H} om _{\widetilde{\D} _{\X'}}
  (\M ', \widetilde{\omega} _{\X '} \otimes _{\widetilde{\O} _{\X'}}
  \widetilde{\D} _{\X'} )) \otimes ^\L _{\widetilde{\D} _{\X}}  \E}
  \ar[u]_-\sim ^-{\mathrm{proj} _f}
  \ar[d]_-\sim
  &
  { \R \mathcal{H} om _{\widetilde{\D} _{\X}}
  (f _*  \M ',f _* ( \widetilde{\omega} _{\X'} ) \otimes _{\widetilde{\O} _{\X} } \E)}
  \ar[u]_-\sim ^-{\mathrm{proj} '}
  \ar[d]_-\sim ^-{\mathrm{Tr}}
  \\
  { \R \mathcal{H} om _{\widetilde{\D} _{\X}}
  (f _*  \M ', f _*  (\widetilde{\omega} _{\X '} \otimes _{\widetilde{\O} _{\X'}}  \widetilde{\D} _{\X'}) )
  \otimes ^\L _{\widetilde{\D} _{\X}}  \E}
  &
  { \R \mathcal{H} om _{\widetilde{\D} _{\X}}
  (f _*  \M ',\widetilde{\omega} _{\X} \otimes _{\widetilde{\O} _{\X} } \E)}
  \ar@{=}[d]
  \\
{ \R \mathcal{H} om _{\widetilde{\D} _{\X}}
  (f _*  \M ', f _*  (\widetilde{\omega} _{\X '} )\otimes _{\widetilde{\O} _{\X}}  \widetilde{\D} _{\X})
  \otimes ^\L _{\widetilde{\D} _{\X}}  \E}
  \ar[u]_-\sim ^-{\mathrm{proj} '}
  \ar[r]_-\sim ^-{\otimes \circ \mathrm{Tr}}
  &
  { \R \mathcal{H} om _{\widetilde{\D} _{\X}}
  (f _*  \M ',\widetilde{\omega} _{\X} \otimes _{\widetilde{\O} _{\X} } \E),}
}
\end{equation}
où le morphisme de projection $\mathrm{proj} _f$ en haut à droite est l'isomorphisme de la troisième
ligne de \ref{adjdagt}
et dont les deux isomorphismes $\mathrm{proj}'$ sont induits par \ref{projiso},
est commutatif. En effet, cela est local en $\X$.
On peut donc supposer que $\M'$ possède une résolution bornée
par des $\widetilde{\D} _{\X'}$-modules localement projectifs et de type fini. En résolvant platement
$\E$, la commutativité du rectangle du bas de \ref{alphs=tr-diag}
résulte d'un calcul. Enfin, celle du carré se vérifie par fonctorialité.

Or, grâce à \ref{isoduareliso} (et via l'identification $f _* =f _{T+}$), la flèche
$$(f _*  \R \mathcal{H} om _{\widetilde{\D} _{\X'}}
  (\M ', \widetilde{\omega} _{\X '} \otimes _{\widetilde{\O} _{\X'}}
  \widetilde{\D} _{\X'} )) \otimes ^\L _{\widetilde{\D} _{\X}}  \E
  \riso
  \R \mathcal{H} om _{\widetilde{\D} _{\X}}
  (f _*  \M ',\widetilde{\omega} _{\X} \otimes _{\widetilde{\O} _{\X} } \E)$$
  de \ref{alphs=tr-diag}
  (passant par la gauche puis par le bas)
  est $\otimes \circ \chi$ de \ref{adjdagt} (qui s'identifie d'ailleurs à
$\otimes \circ \chi \circ \otimes $).
Il en résulte que l'isomorphisme
$$f _*  \R \mathcal{H} om _{\widetilde{\D} _{\X'}}
  (\M ',  f ^! _T  ( \widetilde{\omega} _{\X} \otimes _{\widetilde{\O} _{\X} } \E ))
  \riso
  \R \mathcal{H} om _{\widetilde{\D} _{\X}}
  (f _*  \M ',\widetilde{\omega} _{\X} \otimes _{\widetilde{\O} _{\X} } \E)$$
  passant par la gauche puis par le bas de \ref{alphs=tr-diag} est le composé de \ref{adjdagt}.
Pour $\M ' = f ^! _T (\widetilde{\omega} _{\X } \otimes _{\widetilde{\O} _{\X}}\E )$,
en appliquant $H ^0 \circ \R \Gamma (\X,-)$ à \ref{alphs=tr-diag}, il en dérive
le diagramme commutatif
\begin{equation}
  \label{alphs=tr-diagfin}
  \xymatrix  @R=0,3cm    {
{\mathrm{Hom} _{\widetilde{\D} _{\X'}}
  (f ^! _T (\widetilde{\omega} _{\X } \otimes _{\widetilde{\O} _{\X}}\E ),
  f ^! _T (\widetilde{\omega} _{\X } \otimes _{\widetilde{\O} _{\X}}\E ))}
  \ar[r] _-\sim\ar[rd] _-{\mathrm{adj}} ^-\sim
  &
{\mathrm{Hom} _{\widetilde{\D} _{\X}}
  (f _{T+} f ^! _T (\widetilde{\omega} _{\X } \otimes _{\widetilde{\O} _{\X}}\E ),
  f _{T+} f ^! _T (\widetilde{\omega} _{\X } \otimes _{\widetilde{\O} _{\X}}\E ))}
  \ar[d]  _-\sim
  \\
  {}
  &
  {\mathrm{Hom} _{\widetilde{\D} _{\X}}
  (f _{T+} f ^! _T (\widetilde{\omega} _{\X } \otimes _{\widetilde{\O} _{\X}}\E ),
  \widetilde{\omega} _{\X } \otimes _{\widetilde{\O} _{\X}}\E )},
  }
\end{equation}
  dont le morphisme de droite se déduit du composé de droite de \ref{alphs=tr-diag}
  (qui découle fonctoriellement de \ref{def-alphadiag}).
On conclut en calculant l'image de l'identité de
$f ^! _T (\widetilde{\omega} _{\X } \otimes _{\widetilde{\O} _{\X}}\E )$
via les deux chemins de \ref{alphs=tr-diagfin}.
\end{proof}

\begin{vide}
  \label{predefDf!=f!D}
  Soit $\FF$ un $\widetilde{\D} _{\X}$-bimodule à gauche.
On obtient le $(f ^{-1} \widetilde{\D} _{\X}, \widetilde{\D} _{\X'})$-bimodule
$f ^! _{T,\mathrm{d}} (\FF)$. Via l'isomorphisme canonique
$\iota$ : $f ^{-1} \widetilde{\D} _{\X}
\riso
\widetilde{\D} _{\X'}$,
$f ^! _{T,\mathrm{d}} (\FF)$ est ainsi muni d'une structure canonique
de $\widetilde{\D} _{\X'}$-bimodule à gauche.
De même, en remplaçant {\og à gauche \fg} par {\og à droite \fg}
ou le symbole {\og $d$ \fg} par {\og $g$ \fg} ou
en passant aux complexes.

Pour tous $\E \in D  (\overset{^\mathrm{g}}{} \smash{\D} ^{\dag} _{\X  ,\Q} ( \hdag T ) )$
et
$\FF \in D   ^+  (\overset{^\mathrm{g}}{} \smash{\D} ^{\dag} _{\X  ,\Q} ( \hdag T ),
\overset{^\mathrm{g}}{} \smash{\D} ^{\dag} _{\X  ,\Q} ( \hdag T )  )$,
on obtient ainsi les deuxièmes isomorphisme des composés :
\begin{gather}\label{f!hom=homf!f!}
  f ^! _T \R \mathcal{H} om _{\widetilde{\D} _{\X}}
  ( \E ,\FF) \riso
\R \mathcal{H} om _{f ^{-1}\widetilde{\D} _{\X}}
  (f ^{-1} \E ,f ^! _{T,\mathrm{d}} \FF) \riso
  \R \mathcal{H} om _{\widetilde{\D} _{\X'}} (f ^! _T (\E), f ^! _{T,\mathrm{d}} (\FF)),
  \\
  \label{f!hom=homf!f!2}
f ^! _{T,\mathrm{d}} \R \mathcal{H} om _{\widetilde{\O} _{\X}} ( \E ,\FF)
  \riso
\R \mathcal{H} om _{f ^{-1}\widetilde{\O} _{\X}}
  (f ^{-1} \E ,f ^! _{T,\mathrm{d}} \FF) \riso
  \R \mathcal{H} om _{\widetilde{\O} _{\X'}} (f ^! _T (\E), f ^! _{T,\mathrm{d}} (\FF)).
\end{gather}
De même, pour tous $\E \in D ^- (\overset{^\mathrm{g}}{} \smash{\D} ^{\dag} _{\X  ,\Q} ( \hdag T ) )$,
$\FF \in D   ^-  (\overset{^\mathrm{g}}{} \smash{\D} ^{\dag} _{\X  ,\Q} ( \hdag T ),
\overset{^\mathrm{g}}{} \smash{\D} ^{\dag} _{\X  ,\Q} ( \hdag T )  )$,
\begin{equation}
  \label{predefDf!=f!Dgath2}
f ^! _{T,\mathrm{d}}  (\E \otimes _{\widetilde{\O} _{\X}} ^\L \FF)
\riso
  f ^{-1}  (\E) \otimes _{f ^{-1} \widetilde{\O} _{\X}} ^\L   f ^! _{T,\mathrm{d}} (\FF)
 \riso
f ^! _{T}  (\E) \otimes _{\widetilde{\O} _{\X'}} ^\L  f ^! _{T,\mathrm{d}} (\FF)
\end{equation}
On remarque que la structure gauche de $\E \otimes _{\widetilde{\O} _{\X}} ^\L \FF$
étant celle du produit tensoriel,
$f ^! _{T,\mathrm{g}}  (\E \otimes _{\widetilde{\O} _{\X}} ^\L \FF)
\riso
f ^! _{T}  ( \E )\otimes _{\widetilde{\O} _{\X'}} ^\L  f ^! _{T,\mathrm{g}} (\FF)$
(pour ce dernier isomorphisme, l'hypothèse que $f$ est un isomorphisme est superflu, mais il faut
ajouté le décalage $[d _{X'/X}]$).
\end{vide}

\begin{vide}\label{defDf!=f!D}
Soit $\E \in D _{\mathrm{parf}} (\overset{^\mathrm{g}}{} \smash{\D} ^{\dag} _{\X  ,\Q} ( \hdag T ) )$.
On note $\xi$ : $f ^! _T \DD '_{T} (\E) \riso \DD '_{T'}  f ^! _T (\E)$ l'isomorphisme défini comme suit
\begin{gather}\notag
\xi \ :\   f ^! _T \DD ' _{T} (\E)
\riso
\R \mathcal{H} om _{\widetilde{\D} _{\X'}}
(f ^! _T (\E) ,
f ^! _{T,\mathrm{d}}((\widetilde{\D} _{\X}\otimes _{\widetilde{\O} _{\X}} \widetilde{\omega} _{\X } ^{-1}) _{\mathrm{t}} ))[d _X]
\\ \label{defDf!=f!D1}
\riso
\R \mathcal{H} om _{\widetilde{\D} _{\X'}}
(f ^! _T \E,
(\widetilde{\D} _{\X'} \otimes  _{\widetilde{\O} _{\X'}} \widetilde{\omega} _{\X '} ^{-1}) _{\mathrm{t}})[d _{X'}]
=  \DD '_{T'}  f ^! _T (\E),
\end{gather}
où l'isomorphisme de la première ligne (resp. deuxième ligne) est \ref{f!hom=homf!f!} (resp. \ref{fdomegabeta}).
Il en découle l'isomorphisme :
$ f ^! _T \DD _{T} (\E) \underset{\beta}{\riso} f ^! _T \DD '_{T} (\E)
\underset{\xi}{\riso}
\DD '_{T'}  f ^! _T (\E) \underset{\beta}{\riso} \DD _{T'}  f ^! _T (\E) $,
que l'on notera à nouveau $\xi$.

Grâce à \ref{f!dg} et \ref{DDgd}, on obtient le composé :
\begin{gather}\notag
  f ^! _T \circ \DD _T (\widetilde{\omega} _{\X } \otimes _{\widetilde{\O} _{\X}} \E)
  \riso
  f ^! _T (\widetilde{\omega} _{\X } \otimes _{\widetilde{\O} _{\X}}  \DD ' _T (\E))
  \riso
  \widetilde{\omega} _{\X '} \otimes _{\widetilde{\O} _{\X'}}   f ^! _T (\DD' _T (\E))
  \overset{\xi}{\underset{\sim}{\longrightarrow}}
  \widetilde{\omega} _{\X '} \otimes _{\widetilde{\O} _{\X'}} \DD' _{T'} (f ^! _T(\E))
  \\
  \riso
  \DD _{T'} (\widetilde{\omega} _{\X '} \otimes _{\widetilde{\O} _{\X'}} f ^! _T(\E))
  \riso
  \DD _{T'}  \circ f ^! _T  (\widetilde{\omega} _{\X } \otimes _{\widetilde{\O} _{\X}} \E)
  \label{xid}
\end{gather}
encore noté $\xi$. Pour tout
$\M \in D _{\mathrm{parf}} ( \smash{\D} ^{\dag} _{\X  ,\Q} ( \hdag T ) \overset{^{\mathrm{d}}}{})$,
on a ainsi l'isomorphisme $\xi $ : $ f ^! _T \DD _T (\M) \riso \DD _{T'} f ^! _T  (\M)$.
\end{vide}

\begin{rema}
Grâce au diagramme de droite de \ref{fdomegabetadiag1}, avec les notations \ref{defDf!=f!D}, on aurait pu définir
$\xi $ comme suit :
\begin{gather}\notag
\xi \ :\   f ^! _T \DD  _{T} (\E)
\riso
\R \mathcal{H} om _{\widetilde{\D} _{\X'}}
(f ^{!} _T ( \E),
f ^! _{T,\mathrm{d}}   (\widetilde{\D} _{\X}  \otimes _{\widetilde{\O} _{\X}} \widetilde{\omega} _{\X } ^{-1}))[d _X]
\\
\riso
\R \mathcal{H} om _{\widetilde{\D} _{\X'}}
(f ^{!} _T ( \E),
\widetilde{\D} _{\X'} \otimes  _{\widetilde{\O} _{\X'}} \widetilde{\omega} _{\X '} ^{-1} ) _\mathrm{t})[d _X]
\underset{\beta}{\riso}
\DD _{T'}  f ^! _T (\E).
\label{defDf!=f!D1bis}
\end{gather}
Avec la remarque de \ref{fdomegabeta},
on vérifie que lorsque $f$ est l'identité $\xi$ est l'identité.
\end{rema}

\begin{prop}\label{xichitrcomp}
  Soit $\M \in D _{\mathrm{parf}} (\overset{^*}{} \smash{\D} ^{\dag} _{\X  ,\Q} ( \hdag T ) )$. Le diagramme
\begin{equation}
\label{xichitrcompdiag}
\xymatrix  @R=0,3cm   {
{ f _{T,+} f ^! _T \DD _T (\M) }
\ar[r] _-\sim ^-{\xi} \ar[d]  _-{\mathrm{adj} _{\DD _T (\M) }}    ^-\sim
&
{ f _{T,+} \DD _{T'} f ^! _T  (\M) }
\ar[d] ^\chi _-\sim
\\
{\DD _T (\M) }
\ar[r] _-{\DD _T  \mathrm{adj} _{\M}} ^-\sim
&
{\DD _T f _{T,+} f ^! _T  (\M) }
}
\end{equation}
  est commutatif. De même, en remplaçant $\DD$ par $\DD '$.
\end{prop}
\begin{proof}
Comme les carrés de droite, de gauche, du haut et du bas du diagramme
\begin{equation}
  \xymatrix  @R=0,3cm   {
&{ f _{T,+} f ^! _T \DD '_T (\M) }
\ar[rr]  ^-{\xi} \ar'[d][dd]  _-(0.75){\mathrm{adj} _{\DD '_T (\M) }}
&&
{ f _{T,+} \DD '_{T'} f ^! _T  (\M) }
\ar[dd] ^\chi
\\
{ f _{T,+} f ^! _T \DD _T (\M) }
\ar[rr]  ^-(0.65){\xi} \ar[dd]  _-{\mathrm{adj} _{\DD _T (\M) }}
\ar[ur] ^-{\beta}
&&
{ f _{T,+} \DD _{T'} f ^! _T  (\M) }
\ar[dd] ^(0.4){\chi}
\ar[ur] ^-{\beta}
\\
&
{\DD '_T (\M) }
\ar'[r]_-{\DD ' _T \mathrm{adj} _{\M}}[rr]
&&
{\DD '_T f _{T,+} f ^! _T  (\M) }
\\
{\DD _T (\M) }
\ar[rr] _-{\DD _T \mathrm{adj} _{\M}}
\ar[ur] ^-{\beta}
&&
{\DD _T f _{T,+} f ^! _T  (\M) }
\ar[ur] ^-{\beta}
}
\end{equation}
sont commutatif, il suffit de prouver la proposition pour $\DD '$.
De même, grâce à \ref{f!dg}, \ref{f+dg}, \ref{DDdg} et \ref{DDgd},
on vérifie que le cas où $* =g$ se déduit du cas $*=d$. Traitons donc ce dernier.
Pour tout $\FF \in D _{\mathrm{parf}} (\overset{^\mathrm{g}}{} \smash{\D} ^{\dag} _{\X  ,\Q} ( \hdag T ) )$,
notons $\alpha$ le morphisme composé
$f _* ( \widetilde{\omega} _{\X'} \otimes _{\widetilde{\O} _{\X'} } f ^! _T (\FF)) \underset{\mathrm{proj}'}{\liso}
  f _* ( \widetilde{\omega} _{\X'} ) \otimes _{\widetilde{\O} _{\X} } \FF
  \underset{\mathrm{Tr}}{\riso}
\widetilde{\omega} _{\X} \otimes _{\widetilde{\O} _{\X} } \FF$, où
$\mathrm{proj}'$ a été défini en \ref{projiso}.

Soit $\E \in D _{\mathrm{parf}} (\overset{^\mathrm{g}}{} \smash{\D} ^{\dag} _{\X  ,\Q} ( \hdag T ) )$
tel que $\M = \widetilde{\omega} _{\X} \otimes _{\widetilde{\O} _{\X} }  \E$.
D'après \ref{def-alpha}, par fonctorialité en l'isomorphisme
$\DD ' _T(\widetilde{\omega} _{\X} \otimes _{\widetilde{\O} _{\X} } \E)
\riso
\widetilde{\omega} _{\X} \otimes _{\widetilde{\O} _{\X} } \DD  _T ( \E)$,
le composé de gauche du diagramme ci-dessous
\begin{equation}
  \label{xichitrcompdiag1}
\xymatrix  @R=0,3cm   {
{f _* f ^! _T \DD' _T(\widetilde{\omega} _{\X} \otimes _{\widetilde{\O} _{\X} } \E)}
\ar[rrr] _-\sim ^\xi
\ar[d] _-\sim
&&&
{f _*  \DD' _T  f ^! _T (\widetilde{\omega} _{\X} \otimes _{\widetilde{\O} _{\X} }  \E)}
\ar@{=}[dd]
\\
{f _* f ^! _T (\widetilde{\omega} _{\X} \otimes _{\widetilde{\O} _{\X} } \DD _T ( \E))}
\ar[d] _-\sim
&
\\
{f _* ( \widetilde{\omega} _{\X'} \otimes _{\widetilde{\O} _{\X'} } f ^! _T \DD _T(\E))}
\ar[d] ^-{\alpha} _-\sim
\ar[r] ^-{\xi} _-\sim
&
{f _* ( \widetilde{\omega} _{\X'} \otimes _{\widetilde{\O} _{\X'} } \DD _T f ^! _T (\E))}
\ar[r] _-\sim
&
{f _*  \DD' _T (\widetilde{\omega} _{\X'} \otimes _{\widetilde{\O} _{\X'} }  f ^! _T (\E))}
\ar[dd] _-\sim ^\chi
\ar[r] _-\sim
&
{f _*  \DD' _T  f ^! _T (\widetilde{\omega} _{\X} \otimes _{\widetilde{\O} _{\X} }  \E)}
\ar[dd] _-\sim ^\chi
\\
{\widetilde{\omega} _{\X} \otimes _{\widetilde{\O} _{\X} } \DD _T ( \E)}
\ar[d] _-\sim
\\
{\DD' _T (\widetilde{\omega} _{\X} \otimes _{\widetilde{\O} _{\X} }   \E)}
\ar[rr] _-{\DD ' _T \alpha} ^-\sim
&
&
{\DD' _T f _*  (\widetilde{\omega} _{\X'} \otimes _{\widetilde{\O} _{\X'} }  f ^! _T (\E))}
\ar[r] _-\sim
&
{\DD' _T f _* f ^! _T (\widetilde{\omega} _{\X} \otimes _{\widetilde{\O} _{\X} }   \E)}
}
\end{equation}
est égal à $\mathrm{adj} _{\DD' _T (\widetilde{\omega} _{\X} \otimes _{\widetilde{\O} _{\X} }   \E)}$.
Avec \ref{def-alpha}, le composé du bas est
$\DD' _T  \mathrm{adj} _{\widetilde{\omega} _{\X} \otimes _{\widetilde{\O} _{\X} }   \E}$.
Ainsi, le contour de \ref{xichitrcompdiag1} correspond à celui de \ref{xichitrcompdiag}.
Le rectangle du haut est commutatif par définition (\ref{xid}),
celui en bas à droite l'est par fonctorialité.
Enfin, via \ref{isoduareliso} et \ref{defDf!=f!D1bis},
on établit par un calcul local (où $f _*$ est le foncteur oubli et
où on identifie $\widetilde{\omega} _{\X}$ à $\widetilde{\O} _{\X}$
grâce à la base $d t _1 \wedge \cdots \wedge d t _d$ etc.)
celle du rectangle en bas à sgauche.

\end{proof}

\begin{rema}
  Avec les notations de \ref{xichitrcomp} et dans le cas {\it à gauche},
  il résulte de \ref{xichitrcomp} et de \ref{def-dg!-g!d2}
  que l'isomorphisme $\xi$ est égal à l'isomorphisme $\theta _{f,T,\M}$
  défini en \ref{videdef-dg!-g!d}.
\end{rema}

\subsection{Preuve de la compatibilité aux images inverses par une immersion ouverte}
\begin{vide}[Convention]
\label{convention}
On pose ici $\widetilde{\O}  _{\X}:= \O  _{\X} (\hdag T ) _{\Q}$
$\widetilde{\omega} _{\X}  := \omega _{\X} \otimes _{\O _{\X}} \widetilde{\O}  _{\X}$,
$\smash{\widetilde{\D}}  _{\X} : =\smash{\D}  _{\X} (\hdag T) _{\Q}$.
De même en rajoutant des primes.

  Soient $\E$ un $\smash{\widetilde{\D}}  _{\X}$-module à gauche
  et $\FF$ un $\smash{\widetilde{\D}}  _{\X}$-bimodule à gauche.

Pour calculer le produit tensoriel $\G _1 :=\FF \otimes _{\widetilde{\O} _{\X}} \E$,
on prendra la structure droite de $\FF$.
Par fonctorialité, $\G _1$ est muni d'une structure de
$\smash{\widetilde{\D}}  _{\X}$-bimodule à gauche.
Par convention,
la structure gauche de $\G _1$ sera celle induite fonctoriellement par
la structure gauche de $\FF$ et la structure droite
sera celle induite par le produit tensoriel.

Lorsque l'on choisira d'inverser ces deux structures, on notera
$\G _2:= \FF \otimes ^{\mathrm{t}} _{\widetilde{\O} _{\X}} \E$.
Ainsi, $\G _2 =\G _1$ sauf que la structure droite de $\G _1$ est la structure gauche de $\G _2$
et vice versa.

Pour calculer $\G _3 :=\E \otimes _{\widetilde{\O} _{\X}} \FF$, on prendra la structure
gauche de $\FF$. Par convention,
la structure gauche de $\G _3$ sera la structure induite par le produit tensoriel et
la structure droite sera celle déduite de la structure droite de $\FF$.

On définit de même $\E \otimes ^{\mathrm{t}} _{\widetilde{\O} _{\X}} \FF$ comme étant le
$\smash{\widetilde{\D}}  _{\X}$-bimodule à gauche égale à $\G _3$
modulo un inversement de ses structures gauche et droite.

\end{vide}
Afin de prouver la proposition \ref{spEdualf!}, établissons d'abord les lemmes qui suivent.

\begin{lemm}
\label{lemmhomDB}
Soient
$\E  \in D ^- ( \overset{^\mathrm{g}}{} \smash{\widetilde{\D}} _{\X }  ) $,
$\FF \in D^{+} (\overset{^\mathrm{g}}{} \smash{\widetilde{\D}} _{\X }  ,
\overset{^\mathrm{g}}{}\smash{\widetilde{\D}} _{\X }  )$
et
$\G \in D  _{\mathrm{tdf}} (\overset{^\mathrm{g}}{} \smash{\widetilde{\D}} _{\X }  ) $.
Le diagramme
\begin{equation}
  \label{lemmhomDBdiag}
  \xymatrix  @R=0,3cm   {
  {   f ^! _T ( \R \mathcal{H}om _{ \smash{\widetilde{\D}} _{\X } }
( \E   , \FF   ) \otimes ^\L  _{\widetilde{\O} _{\X }}    \G           )}
\ar[r]  _-\sim
\ar[d] _-\sim
&
{   f ^! _T ( \R \mathcal{H}om _{ \smash{\widetilde{\D}} _{\X } }
( \E   , \FF   \otimes ^\L  _{\widetilde{\O} _{\X }}    \G           ))  }
\ar[d] _-\sim
\\
  {   f ^! _T ( \R \mathcal{H}om _{ \smash{\widetilde{\D}} _{\X } } ( \E   , \FF   ) )
\otimes ^\L  _{\widetilde{\O} _{\X '}}   f ^! _T ( \G      )     }
\ar[d] _-\sim
&
{   \R \mathcal{H}om _{ \smash{\widetilde{\D}} _{\X '} }
(  f ^! _T  \E   ,   f ^! _{T,\mathrm{d}} (\FF   \otimes ^\L  _{\widetilde{\O} _{\X }}    \G )  )  }
\ar[d] _-\sim
\\
 { \R \mathcal{H}om _{ \smash{\widetilde{\D}} _{\X '} } ( f ^! _T \E   , f ^! _{T,\mathrm{d}} \FF   )
\otimes ^\L  _{\widetilde{\O} _{\X '}}   f ^! _T ( \G      )     }
\ar[r]  _-\sim
&
{   \R \mathcal{H}om _{ \smash{\widetilde{\D}} _{\X '} }
(  f ^! _T  \E   ,   f ^! _{T,\mathrm{d}} (\FF  ) \otimes ^\L  _{\widetilde{\O} _{\X '}}   f ^! _T ( \G )  )  ,}
}
\end{equation}
dont les isomorphismes $\smash{\widetilde{\D}} _{\X }  $-linéaires horizontaux
sont les isomorphismes (\cite[2.1.26]{caro_comparaison}) et
dont deux des flèches verticales sont \ref{f!hom=homf!f!},
est commutatif.
\end{lemm}
\begin{proof}
On résout $\FF$ injectivement et $\G$ platement. La commutativité de \ref{lemmhomDBdiag}
découle alors d'un calcul local.
\end{proof}

\begin{lemm}
\label{lemmhomDcartanB}
Soient
$\E  \in D ^- ( \overset{^\mathrm{g}}{} \smash{\widetilde{\D}} _{\X }  ) $,
$\FF  \in D^{-} (\overset{^\mathrm{g}}{} \smash{\widetilde{\D}} _{\X }  )$
et
$\boldsymbol{\mathcal{G}} \in D^{+} (\overset{^\mathrm{g}}{}
\smash{\widetilde{\D}} _{\X }  ,\overset{^*}{}\smash{\widetilde{\D}} _{\X }  )$.
Le diagramme
\begin{equation}
  \label{lemmhomDcartanBdiag}
  \xymatrix  @R=0,3cm   {
  {f ^! _T \R \mathcal{H}om _{ \smash{\widetilde{\D}} _{\X } }
( \E  \otimes ^\L  _{\widetilde{\O} _{\X }} \FF  , \G   )}
\ar[r] ^-{\mathrm{Ca}} _-\sim
\ar[d] _-\sim
&
{f ^! _T \R \mathcal{H}om _{ \smash{\widetilde{\D}} _{\X } }
( \E  , \R \mathcal{H}om _{\widetilde{\O} _{\X }} ( \FF  , \G )) }
\ar[d] _-\sim
\\
{\R \mathcal{H}om _{ \smash{\widetilde{\D}} _{\X '} }
( f ^! _T ( \E  \otimes ^\L  _{\widetilde{\O} _{\X }}  \FF ) , f ^! _{T,\mathrm{d}}   \G   )}
\ar[d] _-\sim
&
{\R \mathcal{H}om _{ \smash{\widetilde{\D}} _{\X' } }
( f ^! _T  \E  , f ^! _{T,\mathrm{d}}  \R \mathcal{H}om _{\widetilde{\O} _{\X }} ( \FF  , \G )) }
\ar[d] _-\sim
\\
{\R \mathcal{H}om _{ \smash{\widetilde{\D}} _{\X '} }
( f ^! _T  \E  \otimes ^\L  _{\widetilde{\O} _{\X '}} f ^! _T   \FF  , f ^! _{T,\mathrm{d}}   \G   )}
\ar[r] ^-{\mathrm{Ca}} _-\sim
&
{\R \mathcal{H}om _{ \smash{\widetilde{\D}} _{\X' } }
( f ^! _T  \E  , \R \mathcal{H}om _{\widetilde{\O} _{\X '}} (f ^! _T  \FF  , f ^! _{T,\mathrm{d}} \G )) ,}
}
\end{equation}
où les isomorphismes $\smash{\widetilde{\D}} _{\X }  $-linéaires horizontaux
sont les isomorphismes de Cartan (\cite[2.1.34]{caro_comparaison}),
où les isomorphismes verticaux du haut (resp. en bas à droite) découlent de \ref{f!hom=homf!f!}
(resp. \ref{f!hom=homf!f!2}),
est commutatif.
\end{lemm}
\begin{proof}
On résout $\FF$ par $\smash{\widetilde{\D}} _{\X }$-modules à gauche plat
et $\G$ par des $\smash{\widetilde{\D}} _{\X }$-bimodules à gauche injectifs.
\end{proof}

\begin{lemm}\label{lemmgammabeta}
Soit $\E $ un $\smash{\widetilde{\D}}  _{\X}$-module à gauche.
Avec les notations de \ref{convention}, on dispose du diagramme commutatif :
  \begin{equation}\label{lemmgammabetadiag}
  \xymatrix  @R=0,3cm   {
{(\smash{\widetilde{\D}} _{\X } \otimes _{\widetilde{\O} _{\X}} \widetilde{\omega} _{\X } ^{-1})
\otimes _{\widetilde{\O} _{\X}} \E}
\ar[r] _-\sim ^-{\gamma _{\E} \otimes \omega ^{-1}}
\ar[d] _-\sim ^-{\beta}
&
{\E \otimes _{\widetilde{\O} _{\X}}
(  \smash{\widetilde{\D}} _{\X } \otimes _{\widetilde{\O} _{\X}} \widetilde{\omega} _{\X } ^{-1} )   }
\ar[d] _-\sim ^-{\beta}
\\
{(\smash{\widetilde{\D}} _{\X } \otimes _{\widetilde{\O} _{\X}} \widetilde{\omega} _{\X } ^{-1} )_\mathrm{t}
\otimes _{\widetilde{\O} _{\X}} \E }
\ar@{=}[d]
&
{\E \otimes _{\widetilde{\O} _{\X}}
(\smash{\widetilde{\D}} _{\X } \otimes _{\widetilde{\O} _{\X}} \widetilde{\omega} _{\X } ^{-1} )_\mathrm{t}}
\ar@{=}[d]
\\
{\E \otimes _{\widetilde{\O} _{\X}} ^{\mathrm{t}}
(\smash{\widetilde{\D}} _{\X } \otimes _{\widetilde{\O} _{\X}}\widetilde{\omega} _{\X } ^{-1}   )  }
&
{(\smash{\widetilde{\D}} _{\X } \otimes _{\widetilde{\O} _{\X}} \widetilde{\omega} _{\X } ^{-1}  )
\otimes ^{\mathrm{t}} _{\widetilde{\O} _{\X}} \E.}
\ar[l] ^-\sim _-{\gamma _{\E} \otimes \omega ^{-1}}
}
\end{equation}

\end{lemm}
\begin{proof}
  En utilisant les caractérisations de \cite[1.3.1 et 1.3.3]{Be1} des isomorphismes de transposition,
  on calcule que l'image de $(1 \otimes \omega) \otimes e$, avec $e$ (resp. $\omega$) une section locale de $\E$
  (resp. $\widetilde{\omega} _{\X } ^{-1}$), ne dépend pas du chemin suivi.
  On conclut par $\smash{\widetilde{\D}} _{\X }$-linéarité.
\end{proof}

Afin d'établir le diagramme commutatif \ref{coroeveef=homefdiag}, nous aurons besoin des
deux lemmes ci-dessous.
\begin{lemm}
\label{gammaomegcompiso}
  Soit $\G \in D ( \overset{^\mathrm{g}}{} \smash{\widetilde{\D}} _{\X }  ) $.
  Le diagramme ci-après
\begin{equation}
  \label{gammaomegcompisodiag}
  \xymatrix  @R=0,3cm  {
{ f ^! _{T, \mathrm{d}} (
  (\smash{\widetilde{\D}} _{\X } \otimes \widetilde{\omega} _{\X } ^{-1} )_\mathrm{t}
  \otimes _{\widetilde{\O} _{\X}} \G )}
  \ar[d] _-\sim
  &
  { f ^! _{T, \mathrm{d}} (
  \G   \otimes _{\widetilde{\O} _{\X}}
  (\smash{\widetilde{\D}} _{\X } \otimes \widetilde{\omega} _{\X } ^{-1} )_\mathrm{t} )   }
\ar[l] ^-{\gamma _{\G} \otimes \omega ^{-1}} _-\sim
\ar[d] _-\sim
\\
{ f ^! _{T, \mathrm{d}} (   (\smash{\widetilde{\D}} _{\X } \otimes \widetilde{\omega} _{\X } ^{-1} )_\mathrm{t})
  \otimes _{\widetilde{\O} _{\X'}} f ^! _{T} (\G )}
    \ar[d] _-\sim
&
{ f ^! _{T} (   \G  ) \otimes _{\widetilde{\O} _{\X'}}
 f ^! _{T, \mathrm{d}}( (\smash{\widetilde{\D}} _{\X } \otimes \widetilde{\omega} _{\X } ^{-1} )_\mathrm{t} )   }
 \ar[d] _-\sim
\\
  {(\smash{\widetilde{\D}} _{\X '} \otimes \widetilde{\omega} _{\X '} ^{-1} )_\mathrm{t}
  \otimes _{\widetilde{\O} _{\X'}} f ^! _{T} (\G )  }
&
{ f ^! _{T} (   \G  ) \otimes _{\widetilde{\O} _{\X'}}
 (\smash{\widetilde{\D}} _{\X '} \otimes \widetilde{\omega} _{\X '} ^{-1} )_\mathrm{t}    }
  \ar[l] ^-{\gamma _{f ^! _{T} (   \G  ) } \otimes \omega ^{-1}} _-\sim
}
\end{equation}
dont la flèche en haut à droite est \ref{predefDf!=f!Dgath2}
et dont les flèches horizontales dérivent fonctoriellement de celle du bas de \ref{lemmgammabetadiag}, est commutatif.
\end{lemm}
\begin{proof}
Cela se vérifie par un calcul local.
\end{proof}

\begin{lemm}
  \label{eveef=homef}
Soient
$\E  \in D ( \overset{^\mathrm{g}}{} \smash{\widetilde{\D}} _{\X }  ) $,
$\FF \in D^{+} _{\mathrm{tdf},.}
(\overset{^\mathrm{g}}{} \smash{\widetilde{\D}} _{\X }, \overset{^\mathrm{g}}{}\smash{\widetilde{\D}} _{\X }  )$.
On bénéficie du diagramme commutatif ci-dessous :
\begin{equation}
    \label{eveef=homefdiag}
    \xymatrix  @R=0,3cm  {
      {   f ^! _{T,\mathrm{d}} (\E ^\vee \otimes ^\L  _{\widetilde{\O} _{\X }} \FF   )  }
      \ar[r] _-\sim
      \ar[d] _-\sim
      &
      {   f ^! _{T,\mathrm{d}} ( \R \mathcal{H}om _{ \smash{\widetilde{\O}} _{\X } } ( \E   , \FF    )) }
      \ar[dd] _-\sim
      \\
      { f ^! _{T} (\E     ^\vee ) \otimes ^\L  _{\widetilde{\O} _{\X' }} f ^! _{T,\mathrm{d}} (\FF )    }
      \ar[d] _-\sim
      \\
      {   f ^! _{T} (\E)  ^\vee \otimes ^\L  _{\widetilde{\O} _{\X' }} f ^! _{T,\mathrm{d}} (\FF )    }
      \ar[r] _-\sim
      &
      {  \R \mathcal{H}om _{ \smash{\widetilde{\O}} _{\X } } ( f ^!  _T (\E )  , f ^! _{T,\mathrm{d}} (\FF)    ),}
      }
\end{equation}
dont la flèche en haut à gauche est \ref{predefDf!=f!Dgath2} et celle de droite est \ref{f!hom=homf!f!2}.
\end{lemm}
\begin{proof}
  On résout $\FF$ platement et $\widetilde{\O} _{\X}$ injectivement puis on conclut via
  un calcul local immédiat.
\end{proof}

\begin{vide}
Soit $\E  \in D ( \overset{^\mathrm{g}}{} \smash{\widetilde{\D}} _{\X }  ) $.
  En appliquant \ref{gammaomegcompisodiag} à $\G = \E ^\vee$,
  on obtient le carré supérieur du diagramme :
\begin{equation}
    \label{coroeveef=homefdiagpre}
    \xymatrix  @R=0,3cm  {
      { f ^! _{T, \mathrm{d}} (
  (\smash{\widetilde{\D}} _{\X } \otimes \widetilde{\omega} _{\X } ^{-1} )_\mathrm{t}
  \otimes _{\widetilde{\O} _{\X}} \E ^\vee )}
      \ar[d] _-\sim
      &
  { f ^! _{T, \mathrm{d}} (
  \E ^\vee   \otimes _{\widetilde{\O} _{\X}}
  (\smash{\widetilde{\D}} _{\X } \otimes \widetilde{\omega} _{\X } ^{-1} )_\mathrm{t} )   }
    \ar[d] _-\sim
    \ar[l] _-\sim ^{\gamma \otimes \omega ^{-1}}
        \\
    {(\smash{\widetilde{\D}} _{\X '} \otimes \widetilde{\omega} _{\X '} ^{-1} )_\mathrm{t}
  \otimes _{\widetilde{\O} _{\X'}} f ^! _{T} (\E   ^\vee ) }
  \ar[d] _-\sim
      &
{ f ^! _{T}   (   \E ^\vee  ) \otimes _{\widetilde{\O} _{\X'}}
(\smash{\widetilde{\D}} _{\X' } \otimes \widetilde{\omega} _{\X' } ^{-1} )_\mathrm{t}   }
    \ar[d] _-\sim
    \ar[l] _-\sim ^{\gamma \otimes \omega ^{-1}}
    \\
    {(\smash{\widetilde{\D}} _{\X '} \otimes \widetilde{\omega} _{\X '} ^{-1} )_\mathrm{t}
  \otimes _{\widetilde{\O} _{\X'}} f ^! _{T} (\E  ) ^\vee  }
      &
{ f ^! _{T}   (   \E ) ^\vee   \otimes _{\widetilde{\O} _{\X'}}
(\smash{\widetilde{\D}} _{\X '} \otimes \widetilde{\omega} _{\X '} ^{-1} )_\mathrm{t}    .}
    \ar[l] _-\sim ^{\gamma \otimes \omega ^{-1}}
    }
\end{equation}
Le carré du bas est commutatif par fonctorialité.
Le diagramme \ref{eveef=homefdiag} utilisé pour
$\FF =(\smash{\widetilde{\D}} _{\X } \otimes \widetilde{\omega} _{\X } ^{-1} )_\mathrm{t}$ donne le carré
du haut de :
\begin{equation}
  \label{coroeveef=homefdiagpre2}
\xymatrix  @R=0,3cm   {
{   f ^! _{T,\mathrm{d}} (\E ^\vee \otimes ^\L  _{\widetilde{\O} _{\X }}
(\smash{\widetilde{\D}} _{\X } \otimes \widetilde{\omega} _{\X } ^{-1} )_\mathrm{t}   )  }
      \ar[r] _-\sim
      \ar[d] _-\sim
  &
   {   f ^! _{T,\mathrm{d}} ( \R \mathcal{H}om _{ \smash{\widetilde{\O}} _{\X } } ( \E   ,
   (\smash{\widetilde{\D}} _{\X } \otimes \widetilde{\omega} _{\X } ^{-1} )_\mathrm{t}   )     ) }
   \ar[d] _-\sim
   \\
   {   f ^! _{T} (\E)  ^\vee \otimes ^\L  _{\widetilde{\O} _{\X' }} f ^! _{T,\mathrm{d}}
   ((\smash{\widetilde{\D}} _{\X } \otimes \widetilde{\omega} _{\X } ^{-1} )_\mathrm{t} )    }
      \ar[r] _-\sim \ar[d] _-\sim
   &
   {  \R \mathcal{H}om _{ \smash{\widetilde{\O}} _{\X } } ( f ^!  _T (\E )  ,
   f ^! _{T,\mathrm{d}} ((\smash{\widetilde{\D}} _{\X } \otimes \widetilde{\omega} _{\X } ^{-1} )_\mathrm{t})    )}
   \ar[d] _-\sim
   \\
   {   f ^! _{T} (\E)  ^\vee \otimes ^\L  _{\widetilde{\O} _{\X' }}
   (\smash{\widetilde{\D}} _{\X '} \otimes \widetilde{\omega} _{\X '} ^{-1} )_\mathrm{t}     }
      \ar[r] _-\sim
   &
   {  \R \mathcal{H}om _{ \smash{\widetilde{\O}} _{\X } } ( f ^!  _T (\E )  ,
   (\smash{\widetilde{\D}} _{\X '} \otimes \widetilde{\omega} _{\X' } ^{-1} )_\mathrm{t}    ).}
   }
\end{equation}
On bénéficie par fonctorialité de la commutativité du carré du bas et donc
celle \ref{coroeveef=homefdiagpre2}.
En composant \ref{coroeveef=homefdiagpre} et \ref{coroeveef=homefdiagpre2},
on obtient ainsi le diagramme commutatif :
\begin{equation}
    \label{coroeveef=homefdiag}
    \xymatrix  @R=0,3cm  {
      { f ^! _{T, \mathrm{d}} (
  (\smash{\widetilde{\D}} _{\X } \otimes \widetilde{\omega} _{\X } ^{-1} )_\mathrm{t}
  \otimes _{\widetilde{\O} _{\X}} \E ^\vee )}
 \ar[r] _-\sim
      \ar[d] _-\sim
      &
{   f ^! _{T,\mathrm{d}} ( \R \mathcal{H}om _{ \smash{\widetilde{\O}} _{\X } } ( \E    ,
(\smash{\widetilde{\D}} _{\X } \otimes \widetilde{\omega} _{\X } ^{-1} )_\mathrm{t}   )) }
    \ar[d] _-\sim
        \\
    {(\smash{\widetilde{\D}} _{\X '} \otimes \widetilde{\omega} _{\X '} ^{-1} )_\mathrm{t}
  \otimes _{\widetilde{\O} _{\X'}} f ^! _{T} (\E  ) ^\vee  }
      \ar[r] _-\sim
      &
    {  \R \mathcal{H}om _{ \smash{\widetilde{\O}} _{\X } } ( f ^!  _T (\E )  ,
    (\smash{\widetilde{\D}} _{\X '} \otimes \widetilde{\omega} _{\X '} ^{-1} )_\mathrm{t}  ).}
    }
\end{equation}

\end{vide}

\begin{vide}\label{notaisoc}
On rappelle que $E$ est un isocristal sur $Y$ surconvergent le long de $T$.
On notera par la suite $\E := \sp _* (E)$, $f ^* (E)$ l'image inverse de $E$ par $f$
(voir \cite[2.3.2.2]{Berig})
et $f ^* ( \E) :=\O _{\X '} (\hdag T') _\Q \otimes _{f ^{-1} \O _{\X } (\hdag T) _\Q} f ^{-1}\E$.
On vérifie à la main la commutativité du diagramme ci-dessous :
\begin{equation}
  \label{mumu}
  \xymatrix  @R=0,3cm    {
  {\sp _* f ^* (E ^\vee)}
  \ar[r] _-\sim        \ar[d]  _-\sim
  &
  {f ^*  \sp _* (E ^\vee)}
  \ar[r]  _-\sim
  &
  {f ^*  (\sp _* E )^\vee}
  \ar@{=}[r]
  &
  {f ^*  (\E ^\vee )}
  \ar[d]  _-\sim
  \\
  {\sp _* (f ^* E )^\vee}
  \ar[r]  _-\sim
  &
  {(\sp _* f ^* E )^\vee}
  \ar[r] _-\sim
  &
  {(f ^* \sp _*  E )^\vee}
  \ar@{=}[r]
  &
  {(f ^* \E )^\vee.}
}
\end{equation}

Avec \ref{f*=f!}, on a les isomorphismes canoniques
$f ^* ( \E)
\riso
f ^{! } _T (\E )
\riso
f ^{! \dag } _T (\E )$.
Il résulte alors de l'isomorphisme canonique $f ^* (\E ^\vee) \riso f ^* (\E) ^\vee$
les suivants
$f ^{! } _T (\E ^\vee) \riso  f ^{! } _T (\E )^\vee$
et
$\mu$ : $f ^{! ^\dag } _T (\E ^\vee) \riso  f ^{! ^\dag } _T (\E )^\vee$.
\end{vide}

\begin{prop}\label{spEdualf!}
Avec les notations \ref{notaisoc}, le diagramme
\begin{equation}
  \label{spEdualf!diag}
\xymatrix  @R=0,3cm  {
{f ^{! ^\dag} _T (\DD ^\dag  _T ( \O _{\X } (\hdag T) _\Q) \otimes _{\O _{\X } (\hdag T) _\Q} \E ^\vee)}
\ar[r] _-\sim \ar[d] _-\sim
&
{f ^{! ^\dag} _T \DD ^\dag _T ( \E )   }
\ar[dd] _-\sim ^-\xi
\\
{f ^{! ^\dag} _T  \DD ^\dag  _T ( \O _{\X } (\hdag T) _\Q)
\otimes _{\O _{\X '} (\hdag T ') _\Q}
f ^{! ^\dag} _T( \E ^\vee)}
\ar[d] _-\sim ^-{\xi \otimes \mu}
\\
{\DD ^\dag _{T'} ( \O _{\X '} (\hdag T') _\Q) \otimes _{\O _{\X '} (\hdag T ') _\Q} f ^{! ^\dag} _T (\E) ^\vee}
\ar[r] _-\sim
&
{ \DD ^\dag _{T'}  f ^{! ^\dag} _T (\E)     ,}
}
\end{equation}
où les isomorphismes horizontaux dérivent de \cite[2.2.1]{caro_comparaison},
est commutatif.
\end{prop}

\begin{proof}
Le carré
$$\xymatrix  @R=0,3cm   {
{f ^{! ^\dag} _T \DD ^\dag _T ( \smash{\D}  ^\dag _{\X} ( \hdag T ) _{\Q} \otimes _{\smash{\D}  _{\X} ( \hdag T ) _{\Q}} \E )}
\ar[d] _-\sim ^\xi \ar[r] _-\sim
&
{ \smash{\D}  ^\dag _{\X'} ( \hdag T' ) _{\Q} \otimes _{\smash{\D}  _{\X'} ( \hdag T ') _{\Q}} f ^{!} _T \DD  _T ( \E )}
\ar[d] _-\sim ^\xi
\\
{ \DD ^\dag _{T'}  f ^{! ^\dag} _T (\smash{\D}  ^\dag _{\X} ( \hdag T ) _{\Q} \otimes _{\smash{\D}  _{\X} ( \hdag T ) _{\Q}}\E)}
\ar[r] _-\sim
&
{\smash{\D}  ^\dag _{\X'} ( \hdag T' ) _{\Q} \otimes _{\smash{\D}  _{\X'} ( \hdag T ') _{\Q}}  \DD _{T'}  f ^{!} _T (\E),}
}$$
où l'isomorphisme de droite se construit de manière analogue à \ref{defDf!=f!D}, est commutatif.
On se ramène ainsi à prouver la commutativité de \ref{spEdualf!diag} sans $\dag$.
Rappelons que les isomorphismes horizontaux de \ref{spEdualf!diag} sans $\dag$
sont construits dans \cite[2.2.1]{caro_comparaison}.

D'après \ref{lemmhomDB},
on bénéficie de la commutativité du deuxième carré du haut du diagramme :
\begin{equation}
  \label{spEdualf!diag1}
\xymatrix  @R=0,3cm   {
{f ^! _T (\R \mathcal{H} om _{\smash{\widetilde{\D}} _{\X }}
( \widetilde{\O} _{\X} ,
\smash{\widetilde{\D}} _{\X } \otimes \widetilde{\omega} _{\X } ^{-1} )\otimes _{\widetilde{\O} _{\X}} \E ^\vee)}
\ar[r] _-\sim
\ar[d] _-\sim ^\beta
&
{f ^! _T \R \mathcal{H} om _{\smash{\widetilde{\D}} _{\X }}
( \widetilde{\O} _{\X} ,
\smash{\widetilde{\D}} _{\X } \otimes \widetilde{\omega} _{\X } ^{-1} \otimes _{\widetilde{\O} _{\X}} \E ^\vee)}
\ar[d] _-\sim^\beta
\\
{f ^! _T (\R \mathcal{H} om _{\smash{\widetilde{\D}} _{\X }}
( \widetilde{\O} _{\X} ,
(\smash{\widetilde{\D}} _{\X } \otimes \widetilde{\omega} _{\X } ^{-1} )_\mathrm{t})
\otimes _{\widetilde{\O} _{\X}} \E ^\vee)}
\ar[r] _-\sim
\ar[d] _-\sim
&
{f ^! _T \R \mathcal{H} om _{\smash{\widetilde{\D}} _{\X }}
( \widetilde{\O} _{\X} ,
(\smash{\widetilde{\D}} _{\X } \otimes \widetilde{\omega} _{\X } ^{-1} )_\mathrm{t}
\otimes _{\widetilde{\O} _{\X}} \E ^\vee)}
\ar[d] _-\sim
\\
{\R \mathcal{H} om _{\smash{\widetilde{\D}} _{\X ' }}
( f ^! _T \widetilde{\O} _{\X } ,
f ^! _{T,\mathrm{d}} ( (\smash{\widetilde{\D}} _{\X } \otimes \widetilde{\omega} _{\X } ^{-1} )_\mathrm{t})
\otimes _{\widetilde{\O} _{\X'}} f ^! _T (\E   ^\vee)}
\ar[r] _-\sim
\ar[d] _-\sim
&
{\R \mathcal{H} om _{\smash{\widetilde{\D}} _{\X ' }}
( f ^! _T \widetilde{\O} _{\X } ,
f ^! _{T,\mathrm{d}} ( (\smash{\widetilde{\D}} _{\X } \otimes \widetilde{\omega} _{\X } ^{-1} )_\mathrm{t})
\otimes _{\widetilde{\O} _{\X'}} f ^! _T (\E ^\vee ) )}
\ar[d] _-\sim
\\
{\R \mathcal{H} om _{\smash{\widetilde{\D}} _{\X ' }}
( \widetilde{\O} _{\X '} ,
(\smash{\widetilde{\D}} _{\X '} \otimes \widetilde{\omega} _{\X '} ^{-1} )_\mathrm{t})
\otimes _{\widetilde{\O} _{\X'}} f ^! _T (\E)  ^\vee}
\ar[r] _-\sim
\ar[d] _-\sim ^\beta
&
{\R \mathcal{H} om _{\smash{\widetilde{\D}} _{\X ' }}
( \widetilde{\O} _{\X '} ,
(\smash{\widetilde{\D}} _{\X '} \otimes \widetilde{\omega} _{\X '} ^{-1} )_\mathrm{t}
\otimes _{\widetilde{\O} _{\X'}} f ^! _T (\E)  ^\vee)}
\ar[d] _-\sim ^\beta
\\
{\R \mathcal{H} om _{\smash{\widetilde{\D}} _{\X ' }}
( \widetilde{\O} _{\X '} ,
\smash{\widetilde{\D}} _{\X '} \otimes \widetilde{\omega} _{\X' } ^{-1} )
\otimes _{\widetilde{\O} _{\X'}} f ^! _T (\E )^\vee}
\ar[r] _-\sim
&
{\R \mathcal{H} om _{\smash{\widetilde{\D}} _{\X ' }}
( \widetilde{\O} _{\X '} ,
\smash{\widetilde{\D}} _{\X' } \otimes \widetilde{\omega} _{\X '} ^{-1} \otimes _{\widetilde{\O} _{\X'}} f ^! _T (\E )^\vee).}
}
\end{equation}
Celle des autres carrés se vérifient par fonctorialité
(pour le deuxième du bas, on dispose des isomorphismes :
$ f _{T,\mathrm{d}} ^! ((\smash{\widetilde{\D}} _{\X } \otimes \widetilde{\omega} _{\X } ^{-1}) _\mathrm{t})
\riso (\smash{\widetilde{\D}} _{\X '} \otimes \widetilde{\omega} _{\X '} ^{-1}) _\mathrm{t}$
et $f ^! _T  (\E ^\vee )\riso f ^! _T (\E )^\vee$).
Ainsi, \ref{spEdualf!diag1} est commutatif.

D'après \ref{lemmgammabeta}, le carré de gauche du diagramme :
\begin{equation}  \label{spEdualf!prediag2}
  \xymatrix  @R=0,3cm   {
{\smash{\widetilde{\D}} _{\X } \otimes \widetilde{\omega} _{\X } ^{-1} \otimes _{\widetilde{\O} _{\X}} \E ^\vee}
\ar[r] _-\sim ^-{\gamma \otimes \omega ^{-1}}
\ar[d] _-\sim ^-{\beta}
&
{\E ^\vee\otimes _{\widetilde{\O} _{\X}}
\smash{\widetilde{\D}} _{\X } \otimes \widetilde{\omega} _{\X } ^{-1}}
\ar[r] _-\sim ^-\sim
\ar[d] _-\sim ^-{\beta}
&
{\R \mathcal{H} om _{\widetilde{\O} _{\X}} (\E, \smash{\widetilde{\D}} _{\X } \otimes \widetilde{\omega} _{\X } ^{-1})}
\ar[d] _-\sim ^-{\beta}
\\
{(\smash{\widetilde{\D}} _{\X } \otimes \widetilde{\omega} _{\X } ^{-1} )_\mathrm{t}
\otimes _{\widetilde{\O} _{\X}} \E ^\vee}
&
{\E ^\vee\otimes _{\widetilde{\O} _{\X}}
(\smash{\widetilde{\D}} _{\X } \otimes \widetilde{\omega} _{\X } ^{-1}} )_\mathrm{t}
\ar[r] _-\sim ^-\sim
\ar[l] _-\sim ^-{\gamma \otimes \omega ^{-1}}
&
{\R \mathcal{H} om _{\widetilde{\O} _{\X}} (\E,
(\smash{\widetilde{\D}} _{\X } \otimes \widetilde{\omega} _{\X } ^{-1})_\mathrm{t})}
}
\end{equation}
est commutatif. Comme le carré de droite de \ref{spEdualf!prediag2} est commutatif par fonctorialité,
il en résulte la commutativité de \ref{spEdualf!prediag2}. On dispose du diagramme analogue à
\ref{spEdualf!prediag2}
en remplaçant $\E$ par $ f ^! _T (\E)$ et $\X$ et $\X'$.
On construit fonctoriellement à partir de \ref{spEdualf!prediag2} (resp. de son analogue), le carré du haut
(resp. du bas) du diagramme :
\begin{equation}
  \label{spEdualf!diag2}
\xymatrix  @R=0,3cm   {
{f ^! _T \R \mathcal{H} om _{\smash{\widetilde{\D}} _{\X }}
( \widetilde{\O} _{\X} ,
\smash{\widetilde{\D}} _{\X } \otimes \widetilde{\omega} _{\X } ^{-1} \otimes _{\widetilde{\O} _{\X}} \E ^\vee)}
\ar[r] _-\sim
\ar[d] _-\sim ^\beta
&
{f ^! _T \R \mathcal{H} om _{\smash{\widetilde{\D}} _{\X }}
( \widetilde{\O} _{\X} ,
\R \mathcal{H} om _{\widetilde{\O} _{\X}} (\E, \smash{\widetilde{\D}} _{\X } \otimes \widetilde{\omega} _{\X } ^{-1}))}
\ar[d] _-\sim^\beta
\\
{f ^! _T \R \mathcal{H} om _{\smash{\widetilde{\D}} _{\X }}
( \widetilde{\O} _{\X} ,
(\smash{\widetilde{\D}} _{\X } \otimes \widetilde{\omega} _{\X } ^{-1} )_\mathrm{t}
\otimes _{\widetilde{\O} _{\X}} \E ^\vee)}
\ar[r] _-\sim
\ar[d] _-\sim
&
{f ^! _T \R \mathcal{H} om _{\smash{\widetilde{\D}} _{\X }}
( \widetilde{\O} _{\X} ,
\R \mathcal{H} om _{\widetilde{\O} _{\X}}
(\E, (\smash{\widetilde{\D}} _{\X } \otimes \widetilde{\omega} _{\X } ^{-1}) _\mathrm{t}))}
\ar[d] _-\sim
\\
{\R \mathcal{H} om _{\smash{\widetilde{\D}} _{\X '}}
( \widetilde{\O} _{\X'} ,
f ^! _{T,\mathrm{d}}  ((\smash{\widetilde{\D}} _{\X } \otimes \widetilde{\omega} _{\X } ^{-1} )_\mathrm{t}
\otimes _{\widetilde{\O} _{\X}} \E ^\vee))}
\ar[r] _-\sim
\ar[d] _-\sim
&
{\R \mathcal{H} om _{\smash{\widetilde{\D}} _{\X '}}
( \widetilde{\O} _{\X'} ,
f ^! _{T,\mathrm{d}} (\R \mathcal{H} om _{\widetilde{\O} _{\X}}
(\E, (\smash{\widetilde{\D}} _{\X } \otimes \widetilde{\omega} _{\X } ^{-1}) _\mathrm{t})))}
\ar[d] _-\sim
\\
{\R \mathcal{H} om _{\smash{\widetilde{\D}} _{\X ' }}
( \widetilde{\O} _{\X '} ,
(\smash{\widetilde{\D}} _{\X '} \otimes \widetilde{\omega} _{\X '} ^{-1} )_\mathrm{t}
\otimes _{\widetilde{\O} _{\X'}} f ^! _T (\E)  ^\vee)}
\ar[r] _-\sim
\ar[d] _-\sim ^\beta
&
{\R \mathcal{H} om _{\smash{\widetilde{\D}} _{\X ' }}
( \widetilde{\O} _{\X '} ,
\R \mathcal{H} om _{\widetilde{\O} _{\X '}}
(f ^! _T \E, (\smash{\widetilde{\D}} _{\X ' } \otimes \widetilde{\omega} _{\X ' } ^{-1}) _\mathrm{t}))}
\ar[d] _-\sim ^\beta
\\
{\R \mathcal{H} om _{\smash{\widetilde{\D}} _{\X ' }}
( \widetilde{\O} _{\X '} ,
\smash{\widetilde{\D}} _{\X '}
\otimes \widetilde{\omega} _{\X '} ^{-1}\otimes _{\widetilde{\O} _{\X'}} f ^! _T (\E )^\vee)}
\ar[r] _-\sim
&
{\R \mathcal{H} om _{\smash{\widetilde{\D}} _{\X ' }}
( \widetilde{\O} _{\X '} ,
\R \mathcal{H} om _{\widetilde{\O} _{\X '}}
(f ^! _T \E, \smash{\widetilde{\D}} _{\X ' } \otimes \widetilde{\omega} _{\X ' } ^{-1}))}.
}
\end{equation}
On déduit de \ref{coroeveef=homefdiag} que le deuxième carré du bas
de \ref{spEdualf!diag2} est commutatif. Comme les deux autres carrés le sont par fonctorialité,
il s'en suit la commutativité de \ref{spEdualf!diag2}.

Comme $ f _{T,\mathrm{d}} ^! ((\smash{\widetilde{\D}} _{\X } \otimes \widetilde{\omega} _{\X } ^{-1}) _\mathrm{t})
\riso (\smash{\widetilde{\D}} _{\X '} \otimes \widetilde{\omega} _{\X '} ^{-1}) _\mathrm{t}$
et $f _T ^!   (\widetilde{\O} _{\X}) \riso \widetilde{\O} _{\X'}$,
il dérive de \ref{lemmhomDcartanB} le carré du milieu de
\begin{equation}
  \label{spEdualf!diag3}
\xymatrix  @R=0,3cm   {
{f ^! _T \R \mathcal{H} om _{\smash{\widetilde{\D}} _{\X }}
( \widetilde{\O} _{\X} ,
\R \mathcal{H} om _{\widetilde{\O} _{\X}} (\E, \smash{\widetilde{\D}} _{\X } \otimes \widetilde{\omega} _{\X } ^{-1}))}
\ar[d] _-\sim ^\beta
&
{f ^! _T \R \mathcal{H} om _{\smash{\widetilde{\D}} _{\X }}
( \E, \smash{\widetilde{\D}} _{\X } \otimes \widetilde{\omega} _{\X } ^{-1})}
\ar[d] _-\sim^\beta
\ar[l] _-\sim ^-{\mathrm{Ca}}
\\
{f ^! _T \R \mathcal{H} om _{\smash{\widetilde{\D}} _{\X }}
( \widetilde{\O} _{\X} ,
\R \mathcal{H} om _{\widetilde{\O} _{\X}}
(\E, (\smash{\widetilde{\D}} _{\X } \otimes \widetilde{\omega} _{\X } ^{-1}) _\mathrm{t}))}
\ar[d] _-\sim
&
{f ^! _T \R \mathcal{H} om _{\smash{\widetilde{\D}} _{\X }}
( \E, (\smash{\widetilde{\D}} _{\X } \otimes \widetilde{\omega} _{\X } ^{-1}) _\mathrm{t})}
\ar[d] _-\sim
\ar[l] _-\sim ^-{\mathrm{Ca}}
\\
{\R \mathcal{H} om _{\smash{\widetilde{\D}} _{\X ' }}
( \widetilde{\O} _{\X '} ,
\R \mathcal{H} om _{\widetilde{\O} _{\X '}}
(f ^! _T \E, (\smash{\widetilde{\D}} _{\X ' } \otimes \widetilde{\omega} _{\X ' } ^{-1}) _\mathrm{t}))}
\ar[d] _-\sim ^\beta
&
{\R \mathcal{H} om _{\smash{\widetilde{\D}} _{\X ' }}
(f ^! _T  \E, (\smash{\widetilde{\D}} _{\X ' } \otimes \widetilde{\omega} _{\X ' } ^{-1}) _\mathrm{t})}
\ar[d] _-\sim ^\beta
\ar[l] _-\sim ^-{\mathrm{Ca}}
\\
{\R \mathcal{H} om _{\smash{\widetilde{\D}} _{\X ' }}
( \widetilde{\O} _{\X '} ,
\R \mathcal{H} om _{\widetilde{\O} _{\X '}}
(f ^! _T \E, \smash{\widetilde{\D}} _{\X ' } \otimes \widetilde{\omega} _{\X ' } ^{-1}))}
&
{\R \mathcal{H} om _{\smash{\widetilde{\D}} _{\X ' }}
(f ^! _T  \E, \smash{\widetilde{\D}} _{\X ' } \otimes \widetilde{\omega} _{\X ' } ^{-1}),}
\ar[l] _-\sim ^-{\mathrm{Ca}}
}
\end{equation}
où les flèches horizontales sont les isomorphismes de Cartan (\cite[2.1.34]{caro_comparaison}).
Comme les carrés du haut et du bas sont commutatifs par fonctorialité,
le diagramme \ref{spEdualf!diag3} l'est donc.

On constate que les isomorphismes du haut et du bas de \ref{spEdualf!diag1}
(resp. \ref{spEdualf!diag2}, resp. \ref{spEdualf!diag3})
correspondent à \cite[2.2.1.5]{caro_comparaison}
(resp. aux compositions de \cite[2.2.1.2]{caro_comparaison}
avec \cite[2.2.1.4]{caro_comparaison}, resp. à \cite[2.2.1.3]{caro_comparaison}).
On aboutit alors,
en composant \ref{spEdualf!diag1}, \ref{spEdualf!diag2} et \ref{spEdualf!diag3}
et en ajoutant le décalage $[d _X]$,
au diagramme \ref{spEdualf!diag} sans $\dag$. D'où le résultat.
\end{proof}

\begin{vide}
Soient $\E \in D (\O _{\X,\Q})$ et
$\FF \in D   ^+  (\overset{^\mathrm{g}}{} \smash{\D} ^{\dag} _{\X  ,\Q} ,
\overset{^\mathrm{g}}{} \smash{\D} ^{\dag} _{\X  ,\Q}   )$.
De manière analogue à \ref{predefDf!=f!D}, on dispose des isomorphismes :
\begin{gather}
\notag
  f ^! _T (\R \mathcal{H} om _{\O _{\X,\Q}} (\E , \FF))
  \riso
  \R \mathcal{H} om _{f ^{-1} \O _{\X,\Q}  } (f ^{-1} \E ,f^! _{T,\mathrm{d}}( \FF))
  \riso
  \R \mathcal{H} om _{\O _{\X',\Q}} (f ^{*} \E ,f^! _{T,\mathrm{d}}( \FF))
  \\
f ^!  (\smash{\D} ^{\dag} _{\X  ,\Q} \otimes _{\O _{\X,\Q} } \E)
\riso
f ^! _{\mathrm{g}} (\smash{\D} ^{\dag} _{\X  ,\Q}) \otimes _{f ^{-1} \O _{\X,\Q}} f ^{-1} \E
\riso
\smash{\D} ^{\dag} _{\X ' ,\Q} \otimes _{\O _{\X',\Q}} f ^* (\E).
  \label{homood2}
\end{gather}
De même, pour tous $\E \in D ^- (\O _{\X,\Q})$,
$\FF \in D   ^-  (\overset{^\mathrm{g}}{} \smash{\D} ^{\dag} _{\X  ,\Q} ( \hdag T ),
\overset{^\mathrm{g}}{} \smash{\D} ^{\dag} _{\X  ,\Q} ( \hdag T )  )$,
\begin{equation}
  \label{otimesood}
  f ^! _T (\E \otimes _{\O _{\X,\Q}} ^\L \FF)
  \riso
   f ^{-1} \E \otimes _{f ^{-1} \O _{\X,\Q}} ^\L  f ^! _{T,\mathrm{d}}  (\FF)
  \riso
 f ^{*} \E \otimes _{\O _{\X',\Q}} ^\L f ^! _{T,\mathrm{d}}  (\FF).
\end{equation}
\end{vide}

On dispose d'un isomorphisme canonique
$\DD _T ( \O _{\X} (\hdag T) _\Q ) \riso  \O _{\X} (\hdag T) _\Q $.
La proposition qui suit signifie que celui-ci commute au foncteur
$f ^! _T$. Nous aurons besoin pour sa preuve du lemme ci-dessous.
\begin{lemm}
\label{DOEF}
  Soient $\E \in D (\O _{\X,\Q})$ et
$\FF \in D   ^\mathrm{b}  (\overset{^\mathrm{g}}{} \smash{\D} ^{\dag} _{\X  ,\Q},
\overset{^\mathrm{g}}{} \smash{\D} ^{\dag} _{\X  ,\Q} )$.
   Le diagramme canonique
\begin{equation}
  \label{DOEFdiag}
  \xymatrix  @R=0,3cm  {
  { f ^!  \R \mathcal{H} om _{\smash{\D} ^{\dag} _{\X  ,\Q}} (\smash{\D} ^{\dag} _{\X  ,\Q} \otimes _{\O _{\X,\Q}} \E, \FF)}
  \ar[d] _-\sim
  &
  { f ^!  (\R \mathcal{H} om _{\O _{\X  ,\Q}} (\E,\O _{\X  ,\Q}) \otimes_{\O _{\X  ,\Q}} \FF)}
  \ar[d] _-\sim
   \ar[l] _-\sim
  \\
{\R \mathcal{H} om _{\smash{\D} ^{\dag} _{\X ' ,\Q}} (f ^!  (\smash{\D} ^{\dag} _{\X  ,\Q} \otimes _{\O _{\X,\Q}} \E) ,
f ^! _{\mathrm{d}}(\FF)    ) }
 \ar[d] _-\sim
 &
{ f ^* (\R \mathcal{H} om _{\O _{\X  ,\Q}} (\E,\O _{\X  ,\Q}) )\otimes_{\O _{\X ' ,\Q}} f ^! _{\mathrm{d}} ( \FF)}
\ar[d] _-\sim
\\
{\R \mathcal{H} om _{\smash{\D} ^{\dag} _{\X ' ,\Q}} (\smash{\D} ^{\dag} _{\X  ',\Q} \otimes _{\O _{\X',\Q}} f ^* (\E) ,
f ^! _{\mathrm{d}}  (\FF)  )       }
&
{ \R \mathcal{H} om _{\O _{\X ' ,\Q}} (f ^* ( \E) ,\O _{\X  ',\Q}) \otimes_{\O _{\X ' ,\Q}} f ^! _{\mathrm{d}} ( \FF)}
 \ar[l] _-\sim
}
\end{equation}
est commutatif.

\end{lemm}

\begin{proof}
Pour construire trois des flèches verticales de \ref{DOEFdiag},
on a utilisé respectivement \ref{otimesood}, \ref{homood2} et \ref{f!hom=homf!f!}.
Le diagramme \ref{DOEFdiag} est le composé des deux diagrammes suivants :
\begin{equation}
  \label{DOEFdiag1}
  \xymatrix  @R=0,3cm  {
  { f ^!  \R \mathcal{H} om _{\smash{\D} ^{\dag} _{\X  ,\Q}} (\smash{\D} ^{\dag} _{\X  ,\Q} \otimes _{\O _{\X,\Q}} \E, \FF)}
  \ar[d] _-\sim
  &
  { f ^!  \R \mathcal{H} om _{\O _{\X  ,\Q}} (\E, \FF)}
   \ar[l] _-\sim
   \ar[dd] _-\sim
  \\
{\R \mathcal{H} om _{\smash{\D} ^{\dag} _{\X ' ,\Q}} (f ^!  (\smash{\D} ^{\dag} _{\X  ,\Q} \otimes _{\O _{\X,\Q}} \E) ,
f ^! _{\mathrm{d}}(\FF)    ) }
 \ar[d] _-\sim
\\
{\R \mathcal{H} om _{\smash{\D} ^{\dag} _{\X ' ,\Q}} (\smash{\D} ^{\dag} _{\X  ',\Q} \otimes _{\O _{\X',\Q}} f ^* (\E) ,
f ^! _{\mathrm{d}}  (\FF)  )       }
&
{\R \mathcal{H} om _{\O _{\X ' ,\Q}} (f ^*(\E),  f ^! _{\mathrm{d}}(\FF)),}
 \ar[l] _-\sim
}
\end{equation}
\begin{equation}
  \label{Homeootf}
  \xymatrix  @R=0,3cm     {
  { f ^!  \R \mathcal{H} om _{\O _{\X  ,\Q}} (\E, \FF)}
  \ar[dd] _-\sim
  &
  { f ^!  (\R \mathcal{H} om _{\O _{\X  ,\Q}} (\E,\O _{\X  ,\Q}) \otimes_{\O _{\X  ,\Q}} \FF)}
  \ar[l] _-\sim
    \ar[d] _-\sim
\\
&
{ f ^* (\R \mathcal{H} om _{\O _{\X  ,\Q}} (\E,\O _{\X  ,\Q}) )\otimes_{\O _{\X ' ,\Q}} f ^! _{\mathrm{d}} ( \FF)}
\ar[d] _-\sim
\\
{\R \mathcal{H} om _{\O _{\X ' ,\Q}} (f ^*(\E),  f ^! _{\mathrm{d}}(\FF))}
&
{ \R \mathcal{H} om _{\O _{\X ' ,\Q}} (f ^* ( \E) ,\O _{\X  ',\Q}) \otimes_{\O _{\X ' ,\Q}} f ^! _{\mathrm{d}} ( \FF),}
 \ar[l] _-\sim
}
\end{equation}
dont la commutativité se vérifie, après avoir résolu injectivement $\FF$, par un calcul.
\end{proof}

\begin{prop}\label{DO=Ociso}
  Le diagramme canonique suivant
$$\xymatrix  @R=0,3cm   {
{f _T ^! \DD _T ( \O _{\X} (\hdag T) _\Q )}
\ar[r] _-\sim ^\xi
\ar[d] _-\sim
&
{\DD _{T'}f _T ^! (\O _{\X} (\hdag T) _{\Q})}
\ar[r] _-\sim
&
{\DD _{T'}( \O _{\X'} (\hdag T') _\Q)}
\ar[d] _-\sim
\\
{f _T ^! (\O _{\X} (\hdag T) _{\Q})}
\ar[rr] _-\sim
&&
{\O _{\X '} (\hdag T') _{\Q}}
}$$
est commutatif.
\end{prop}
\begin{proof}
Par \cite[4.3.12]{Be1}, il suffit de prouver la proposition lorsque $T$ est vide.
  On notera $\smash{\widetilde{\D}} ^\dag _{\X,\Q} :=\smash{\D} ^\dag _{\X,\Q}\otimes _{\O _{\X}} \omega _{\X} ^{-1}$
  et $(\smash{\widetilde{\D}} ^\dag _{\X,\Q} )_\mathrm{t} :=
  (\smash{\D} ^\dag _{\X,\Q}\otimes _{\O _{\X}} \omega _{\X} ^{-1}) _\mathrm{t}$.
  Notons $\T _\X$ le faisceau tangent de $\X$ et
  $\smash{\D} ^\dag _{\X,\Q} \otimes _{\O _{\X}} \wedge ^\bullet \T _\X$
  le complexe de Spencer
  \begin{equation}
    \label{spencerdef}
  \smash{\D} ^\dag _{\X,\Q} \otimes _{\O _{\X}} \wedge ^{\mathrm{d}} \T _\X
  \rightarrow \cdots \rightarrow
  \smash{\D} ^\dag _{\X,\Q} \otimes _{\O _{\X}} \T _\X
  \rightarrow
  \smash{\D} ^\dag _{\X,\Q}
  \end{equation}
  (voir \cite[4.3.1]{Be2}).
L'application canonique $\smash{\D} ^\dag _{\X,\Q}\rightarrow \O _{\X,\Q}$ envoyant un opérateur $P$ sur $P \cdot 1$
induit un quasi-isomorphisme $\smash{\D} ^\dag _{\X,\Q} \otimes _{\O _{\X}} \wedge ^\bullet \T _\X \riso \O _{\X,\Q}$.
  L'isomorphisme canonique $f ^* \Omega _{\X} \riso \Omega _{\X'}$ induit par dualité
  le suivant $\T _{\X'} \riso f ^* \T _{\X}$.
 Avec \ref{homood2}, on en déduit, pour tout entier $r$, le composé :
  $f ^!  (\smash{\D} ^\dag _{\X,\Q} \otimes _{\O _{\X}} \wedge ^r \T _\X)
  \riso
  \smash{\D} ^\dag _{\X',\Q} \otimes _{\O _{\X'}} f ^* (\wedge ^r \T _{\X})
\liso
  \smash{\D} ^\dag _{\X',\Q} \otimes _{\O _{\X'}} \wedge ^r \T _{\X'}$.
  On construit ainsi le diagramme suivant :
  \begin{equation}
    \label{DO=Ociso-diag1}
    \xymatrix  @R=0,3cm   {
    {f ^! (\smash{\D} ^\dag _{\X,\Q} \otimes _{\O _{\X}} \wedge ^{\mathrm{d}} \T _\X)}
    \ar[r] \ar[d] _-\sim
    &
    {\cdots}
    \ar[r]
    &
    {f ^! (\smash{\D} ^\dag _{\X,\Q} \otimes _{\O _{\X}} \T _\X)}
    \ar[r]  \ar[d] _-\sim
    &
    {f ^! \smash{\D} ^\dag _{\X,\Q}}
    \ar[r] \ar[d] _-\sim
    &
    {f ^!  \O _{\X,\Q}}
    \ar[d] _-\sim
    \\
    {\smash{\D} ^\dag _{\X',\Q} \otimes _{\O _{\X'}} \wedge ^{\mathrm{d}} \T _{\X'}} \ar[r]
    &
    {\cdots} \ar[r]
    &
    {\smash{\D} ^\dag _{\X',\Q} \otimes _{\O _{\X'}} \T _{\X'}} \ar[r]
    &
    {\smash{\D} ^\dag _{\X',\Q}}
    \ar[r]
    &
    {\O _{\X',\Q},}
    }
  \end{equation}
dont les complexes horizontaux sont exacts ($f ^!$ est exact).
On vérifie de plus que \ref{DO=Ociso-diag1} est commutatif.
  On obtient en particulier un isomorphisme
$f ^! (\smash{\D} ^\dag _{\X,\Q} \otimes _{\O _{\X}} \wedge ^\bullet \T _\X)
\riso
\smash{\D} ^\dag _{\X',\Q} \otimes _{\O _{\X'}} \wedge ^\bullet \T _{\X'}$.
Il en dérive aussi par dualité la commutativité du carré du bas de :
  \begin{equation}
  \label{DO=Ociso-diag2}
  \xymatrix  @R=0,3cm   {
  {f  ^! \R \mathcal{H} om _{\smash{\D} ^\dag _{\X,\Q}}( \O _{\X,\Q}, \smash{\widetilde{\D}} ^\dag _{\X,\Q})}
  \ar[d] _-\sim ^\xi \ar[r] _-\sim
  &
  {f  ^! \mathcal{H} om _{\smash{\D} ^\dag _{\X,\Q}}
  ( \smash{\D} ^\dag _{\X,\Q} \otimes _{\O _{\X}} \wedge ^\bullet \T _\X ,
  \smash{\widetilde{\D}} ^\dag _{\X,\Q})}
  \ar[d] _-\sim ^\xi
  \\
  {\R \mathcal{H} om _{\smash{\D} ^\dag _{\X',\Q}}( f ^!  \O _{\X,\Q}, \smash{\widetilde{\D}} ^\dag _{\X',\Q})}
  \ar[r] _-\sim
  &
  {\mathcal{H} om _{\smash{\D} ^\dag _{\X',\Q}}
  ( f ^! (\smash{\D} ^\dag _{\X,\Q} \otimes _{\O _{\X}} \wedge ^\bullet \T _\X),
  \smash{\widetilde{\D}} ^\dag _{\X',\Q})}
  \\
  {\R \mathcal{H} om _{\smash{\D} ^\dag _{\X',\Q}}( \O _{\X',\Q}, \smash{\widetilde{\D}} ^\dag _{\X',\Q})}
  \ar[r] _-\sim \ar[u] _-\sim
  &
  {\mathcal{H} om _{\smash{\D} ^\dag _{\X',\Q}}( \smash{\D} ^\dag _{\X',\Q} \otimes _{\O _{\X'}} \wedge ^\bullet \T _{\X'},
  \smash{\widetilde{\D}} ^\dag _{\X',\Q}).}
  \ar[u] _-\sim
    }
  \end{equation}
On établit par fonctorialité la commutativité du carré du haut de \ref{DO=Ociso-diag2}.
D'où celle de \ref{DO=Ociso-diag2}.

Le diagramme
\begin{equation}
\label{DO=Ociso-diag3}
  \xymatrix  @R=0,3cm   {
{f  ^! \mathcal{H} om _{\smash{\D} ^\dag _{\X,\Q}}
  ( \smash{\D} ^\dag _{\X,\Q} \otimes _{\O _{\X,\Q}} \wedge ^\bullet \T _{\X,\Q} ,
  \smash{\widetilde{\D}} ^\dag _{\X,\Q})}
  \ar[r] _-\sim \ar[d] _-\sim ^\beta
&
{f  ^! (\mathcal{H} om _{\O _{\X,\Q}}
  ( \wedge ^\bullet \T _{\X,\Q} , \O _{\X,\Q}) \otimes _{\O _{\X,\Q}}
  \smash{\widetilde{\D}} ^\dag _{\X,\Q})}
  \ar[d] _-\sim ^\beta
\\
{f  ^! \mathcal{H} om _{\smash{\D} ^\dag _{\X,\Q}}
  ( \smash{\D} ^\dag _{\X,\Q} \otimes _{\O _{\X,\Q}} \wedge ^\bullet \T _{\X,\Q} ,
  (\smash{\widetilde{\D}} ^\dag _{\X,\Q}) _\mathrm{t})}
    \ar[r] _-\sim \ar[d] _-\sim
&
{f  ^! (\mathcal{H} om _{\O _{\X,\Q}}
  ( \wedge ^\bullet \T _{\X,\Q} , \O _{\X,\Q}) \otimes _{\O _{\X,\Q}}
  (\smash{\widetilde{\D}} ^\dag _{\X,\Q}) _\mathrm{t})}
  \ar[d] _-\sim
\\
{\mathcal{H} om _{\smash{\D} ^\dag _{\X',\Q}}
  ( \smash{\D} ^\dag _{\X',\Q} \otimes _{\O _{\X',\Q}} f ^* (\wedge ^\bullet \T _{\X,\Q}) ,
  f^! _{\mathrm{d}} ((\smash{\widetilde{\D}} ^\dag _{\X,\Q}) _\mathrm{t}))}
    \ar[r] _-\sim \ar[d] _-\sim
&
{\mathcal{H} om _{\O _{\X',\Q}}
  ( f ^*(\wedge ^\bullet \T _{\X,\Q} ), \O _{\X',\Q}) \otimes _{\O _{\X',\Q}}
f^! _{\mathrm{d}} ((\smash{\widetilde{\D}} ^\dag _{\X,\Q}) _\mathrm{t})}
  \ar[d] _-\sim
\\
{\mathcal{H} om _{\smash{\D} ^\dag _{\X',\Q}}
  ( \smash{\D} ^\dag _{\X',\Q} \otimes _{\O _{\X',\Q}} \wedge ^\bullet \T _{\X ',\Q} ,
  (\smash{\widetilde{\D}} ^\dag _{\X',\Q}) _\mathrm{t} )}
    \ar[r] _-\sim \ar[d] _-\sim ^\beta
&
{\mathcal{H} om _{\O _{\X',\Q}}
  (\wedge ^\bullet \T _{\X '}, \O _{\X',\Q}) \otimes _{\O _{\X',\Q}}
  (\smash{\widetilde{\D}} ^\dag _{\X',\Q}) _\mathrm{t} }
  \ar[d] _-\sim ^\beta
\\
{\mathcal{H} om _{\smash{\D} ^\dag _{\X',\Q}}
  ( \smash{\D} ^\dag _{\X',\Q} \otimes _{\O _{\X',\Q} } \wedge ^\bullet \T _{\X ',\Q},
  \smash{\widetilde{\D}} ^\dag _{\X',\Q})}
    \ar[r] _-\sim
&
{\mathcal{H} om _{\O _{\X',\Q}}
  (\wedge ^\bullet \T _{\X '}, \O _{\X',\Q}) \otimes _{\O _{\X',\Q}}
  \smash{\widetilde{\D}} ^\dag _{\X',\Q}}
}
\end{equation}
est commutatif.
En effet, pour le carré du haut et les deux du bas,
cela se vérifie par fonctorialité tandis que le dernier résulte de \ref{DOEFdiag}.

Les deux carrés du haut et celui du bas du diagramme
\begin{equation}
\label{DO=Ociso-diag4}
  \xymatrix  @R=0,3cm     {
{f  ^! ( \mathcal{H} om _{\O _{\X,\Q}}
(\wedge ^\bullet \T _{\X,\Q} , \O _{\X,\Q}) \otimes _{\O _{\X,\Q}} \smash{\widetilde{\D}} ^\dag _{\X,\Q} )  }
  \ar[d] _-\sim ^\beta
  \ar[r] _-\sim
  &
{f  ^! ( \Omega ^\bullet  _{\X,\Q}  \otimes _{\O _{\X,\Q}} \smash{\widetilde{\D}} ^\dag _{\X,\Q})}
 \ar[d] _-\sim ^\beta
\\
{f  ^! (\mathcal{H} om _{\O _{\X,\Q}}
  ( \wedge ^\bullet \T _{\X,\Q} ,\O _{\X,\Q}) \otimes _{\O _{\X,\Q}}
  (\smash{\widetilde{\D}} ^\dag _{\X,\Q} )_\mathrm{t})}
  \ar[d] _-\sim
    \ar[r] _-\sim
  &
{f  ^! ( \Omega ^\bullet  _{\X,\Q}  \otimes _{\O _{\X,\Q}} (\smash{\widetilde{\D}} ^\dag _{\X,\Q}) _\mathrm{t})}
  \ar[d] _-\sim
\\
{f  ^* \mathcal{H} om _{\O _{\X,\Q}}
  ( \wedge ^\bullet \T _{\X,\Q} ,\O _{\X,\Q}) \otimes _{\O _{\X ',\Q}}
  f^! _{\mathrm{d}} ((\smash{\widetilde{\D}} ^\dag _{\X,\Q} )_\mathrm{t})}
  \ar[d] _-\sim
  \ar[r] _-\sim
  &
{f  ^* ( \Omega ^\bullet  _{\X,\Q}  )\otimes _{\O _{\X',\Q}}
f^! _{\mathrm{d}} ((\smash{\widetilde{\D}} ^\dag _{\X,\Q} )_\mathrm{t})}
  \ar[d] _-\sim
\\
{\mathcal{H} om _{\O _{\X',\Q}}
  (\wedge ^\bullet \T _{\X ',\Q}, \O _{\X',\Q}) \otimes _{\O _{\X',\Q}}
f^! _{\mathrm{d}} ((\smash{\widetilde{\D}} ^\dag _{\X,\Q} )_\mathrm{t})}
  \ar[d] _-\sim
  \ar[r] _-\sim
&
{ \Omega ^\bullet  _{\X ',\Q} \otimes _{\O _{\X',\Q}}
f^! _{\mathrm{d}} ((\smash{\widetilde{\D}} ^\dag _{\X,\Q} )_\mathrm{t})}
\ar[d] _-\sim
\\
{\mathcal{H} om _{\O _{\X',\Q}}
  (\wedge ^\bullet \T _{\X ',\Q}, \O _{\X',\Q}) \otimes _{\O _{\X',\Q}}
  \smash{\widetilde{\D}} ^\dag _{\X',\Q}}
  \ar[r] _-\sim
&
{ \Omega ^\bullet  _{\X ',\Q} \otimes _{\O _{\X',\Q}} \smash{\widetilde{\D}} ^\dag _{\X',\Q},}
  }
\end{equation}
dont les flèches verticales du bas sont induits par
$f^! _{\mathrm{d}} ((\smash{\widetilde{\D}} ^\dag _{\X,\Q} )_\mathrm{t}) \riso
(\smash{\widetilde{\D}} ^\dag _{\X',\Q} )_\mathrm{t} \underset{\beta}{\riso}
f^! _{\mathrm{d}} ((\smash{\widetilde{\D}} ^\dag _{\X,\Q} )_\mathrm{t})$,
sont commutatifs par fonctorialité. Le dernier l'est par définition.

Par un calcul local, on établit la commutativité du diagramme ci-après :
\begin{equation}
  \label{DO=Ociso-diag5}
  \xymatrix  @R=0,3cm  {
  { f ^! ( \omega _{\X,\Q} \otimes _{\O _{\X,\Q}} \smash{\D} ^\dag _{\X,\Q} \otimes _{\O _{\X,\Q}} \omega ^{-1} _{\X,\Q})}
\ar[rr]
\ar[d] _-\sim ^{\beta}
&&
{ f ^! ( \omega _{\X,\Q}  \otimes _{\O _{\X,\Q}} \omega ^{-1} _{\X,\Q})}
\ar[d] _-\sim
\\
  {f ^! (\omega _{\X,\Q} \otimes _{\O _{\X,\Q}}
  (\smash{\D} ^\dag _{\X,\Q} \otimes _{\O _{\X,\Q}} \omega ^{-1} _{\X,\Q}) _\mathrm{t})     }
  \ar[r] _-\sim \ar[d] _-\sim
  &
  { f ^! (\smash{\D} ^\dag _{\X,\Q} )}
\ar[r]
\ar[dd] _-\sim
&
{ f ^! ( \O _{\X,\Q})}
\ar[dd] _-\sim
  \\
{f ^{*} \omega _{\X,\Q} \otimes _{\O _{\X',\Q}}
f ^! _{\mathrm{d}} (   (\smash{\D} ^\dag _{\X,\Q} \otimes _{\O _{\X,\Q}} \omega ^{-1} _{\X,\Q}  )  _\mathrm{t}) }
\ar[d] _-\sim
\\
{\omega _{\X',\Q} \otimes _{\O _{\X',\Q}}
(\smash{\D} ^\dag _{\X',\Q} \otimes _{\O _{\X',\Q}} \omega ^{-1} _{\X',\Q}  )  _\mathrm{t} }
\ar[r] _-\sim
\ar[d] _-\sim ^\beta
&
{\smash{\D} ^\dag _{\X',\Q} }
\ar[r]
&
{ \O _{\X',\Q}}
\ar[d] _-\sim
\\
{ \omega _{\X',\Q} \otimes _{\O _{\X',\Q}} \smash{\D} ^\dag _{\X',\Q} \otimes _{\O _{\X',\Q}} \omega ^{-1} _{\X',\Q}}
\ar[rr]
&&
{\omega _{\X',\Q}  \otimes _{\O _{\X',\Q}} \omega ^{-1} _{\X',\Q},}
}
\end{equation}
où la flèche horizontale du haut (resp. du bas) est induite par l'action à droite
de $\smash{\D} ^\dag _{\X,\Q}$ sur $\omega _{\X,\Q}$ (resp. de $\smash{\D} ^\dag _{\X',\Q}$ sur $\omega _{\X',\Q}$)
tandis que celles du carré de droite découlent
de l'action à gauche de $\smash{\D} ^\dag _{\X,\Q}$ sur $\O _{\X,\Q}$
ou de celle de $\smash{\D} ^\dag _{\X',\Q}$ sur $\O _{\X',\Q}$.
On dispose d'un morphisme canonique entre les flèches de droite de \ref{DO=Ociso-diag4}
décalées de $[d _X]$ et celle de gauche de \ref{DO=Ociso-diag5}.
Via ce morphisme,
en composant \ref{DO=Ociso-diag2}, \ref{DO=Ociso-diag3}, \ref{DO=Ociso-diag4} décalés de $[d _X]$ avec
\ref{DO=Ociso-diag5}, on obtient le diagramme commutatif de \ref{DO=Ociso}.

\end{proof}

\begin{coro}
\label{corospEdualf!}
Avec les notations \ref{notaisoc}, on suppose ici que $f$ est une immersion ouverte.
Le diagramme
  \begin{equation}
  \label{spEdualf!diagnew}
\xymatrix  @R=0,3cm  {
{\sp _* f ^* (E ^\vee)}
  \ar[r] _-\sim        \ar[d]  _-\sim
  &
{f ^{! ^\dag} _T (\E ^\vee)}
\ar[d] _-\sim ^{\mu}
&
{f ^{! ^\dag} _T (\DD ^\dag  _T ( \O _{\X } (\hdag T) _\Q) \otimes _{\O _{\X } (\hdag T) _\Q} \E ^\vee)}
\ar[r] _-\sim
\ar[l] ^-\sim
&
{f ^{! ^\dag} _T \DD ^\dag _T ( \E )}
\ar[d] _-\sim ^\xi
\\
{\sp _* (f ^* E )^\vee}
  \ar[r]  _-\sim
  &
{ f ^{! ^\dag} _T (\E) ^\vee}
&
{\DD ^\dag _{T'} ( \O _{\X '} (\hdag T') _\Q) \otimes _{\O _{\X '} (\hdag T ') _\Q} f ^{! ^\dag} _T (\E) ^\vee}
\ar[r] _-\sim
\ar[l] ^-\sim
&
{ \DD ^\dag _{T'} ( f ^{! ^\dag} _T (\E)),}
}
\end{equation}
où $\mu $ et $\xi$ ont été défini respectivement dans \ref{notaisoc} et \ref{defDf!=f!D},
est commutatif.
\end{coro}

\begin{proof}
La commutativité du carré de gauche est tautologique (voir \ref{notaisoc}).
Il résulte par fonctorialité de \ref{DO=Ociso} que le rectangle du milieu du diagramme :
\begin{equation}
\label{spEdualf!diagnew1}
\xymatrix  @R=0,3cm  {
{f ^{! ^\dag} _T ( \O _{\X } (\hdag T) _\Q  \otimes _{\O _{\X } (\hdag T) _\Q}  \E ^\vee)    }
\ar[d] _-\sim
&
{f ^{! ^\dag} _T (\DD ^\dag  _T ( \O _{\X } (\hdag T) _\Q) \otimes _{\O _{\X } (\hdag T) _\Q} \E ^\vee) }
\ar[l] ^-\sim
\ar[d] _-\sim
\\
{f ^{! ^\dag} _T ( \O _{\X } (\hdag T) _\Q )
\otimes _{\O _{\X '} (\hdag T') _\Q}
f ^{! ^\dag} _T  (\E ^\vee )      }
\ar[dd]  _-\sim
&
{f ^{! ^\dag} _T \DD ^\dag  _T ( \O _{\X } (\hdag T) _\Q)
\otimes _{\O _{\X '} (\hdag T ') _\Q} f ^{! ^\dag} _T  (\E ^\vee) }
\ar[l] ^-\sim
\ar[d] ^-{\xi \otimes \mathrm{Id}} _-\sim
\\
&
{\DD ^\dag  _T  f ^{! ^\dag} _T  ( \O _{\X } (\hdag T) _\Q)
\otimes _{\O _{\X '} (\hdag T ') _\Q} f ^{! ^\dag} _T  (\E ^\vee) }
\ar[d] _-\sim
\\
{\O _{\X '} (\hdag T') _\Q \otimes _{\O _{\X '} (\hdag T') _\Q  }  f ^{! ^\dag} _T  (\E ^\vee )  }
\ar[d] _-\sim ^-{\mu}
&
{\DD ^\dag  _T ( \O _{\X '} (\hdag T') _\Q )
\otimes _{\O _{\X '} (\hdag T ') _\Q} f ^{! ^\dag} _T  (\E ^\vee) }
\ar[l] ^-\sim
\ar[d] _-\sim ^-{\mu}
\\
{\O _{\X '} (\hdag T') _\Q \otimes _{\O _{\X '} (\hdag T') _\Q  }  f ^{! ^\dag} _T  (\E ) ^\vee  }
&
{\DD ^\dag  _T ( \O _{\X '} (\hdag T') _\Q )
\otimes _{\O _{\X '} (\hdag T ') _\Q} f ^{! ^\dag} _T  (\E )^\vee }
\ar[l] ^-\sim
}
\end{equation}
est commutatif. Comme les carrés supérieur et inférieur le sont par fonctorialité, \ref{spEdualf!diagnew1}
est commutatif. En composant
\ref{spEdualf!diagnew1} et \ref{spEdualf!diag}, on obtient \ref{spEdualf!diagnew}.
\end{proof}

\begin{rema}
  \label{rema-corospEdualf!}
Avec les notations \ref{notaisoc}, on se note
$g$ : $\X'' \hookrightarrow \X'$ une seconde immersion ouverte et
$T '' := g _0^{-1} (T ')$.
On bénéficie par fonctorialité du diagramme commutatif :
\begin{equation}
  \label{rema-corospEdualf!diag}
  \xymatrix  @R=0,3cm  {
  {\sp _* (f \circ g) ^* E ^\vee}
  \ar[r] _-\sim
  \ar[d] _-\sim
  &
  {\sp _*  g ^* f ^* E ^\vee}
  \ar[r] _-\sim
  \ar[d] _-\sim
  &
  {\sp _* g ^* ( f ^* E )^\vee}
  \ar[r] _-\sim
  \ar[d] _-\sim
  &
  {\sp _*     ( g ^*  f ^* E )^\vee}
  \ar[r] _-\sim
  \ar[d] _-\sim
  &
  {\sp _*  ( (f\circ g ) ^* E ) ^\vee}
  \ar[d] _-\sim
  \\
  {(f \circ g) ^! _{T} \DD _{T} \sp _* E}
  \ar[r] _-\sim
  &
  {g ^! _{T'} f ^! _{T} \DD _{T} \sp _* E}
  \ar[r] _-\sim
  &
  {g ^! _{T'} \DD _{T'} f ^! _{T} \sp _* E}
  \ar[r] _-\sim
  &
  {\DD _{T''} g ^! _{T'} f ^! _{T}  \sp _* E}
  \ar[r] _-\sim
  &
  {\DD _{T''} (f\circ g) ^! _T  \sp _* E.}
}
\end{equation}
De plus, par transitivité des isomorphismes de commutation des duaux aux images inverses (extraordinaires)
le composé du haut est l'isomorphisme canonique
$\sp _* (f \circ g) ^* E ^\vee
\riso
\sp _*  ( (f\circ g ) ^* E ) ^\vee$
et celui du bas est
$(f \circ g) ^! _T \DD _T \sp _* E
\riso
\DD _{T''} (f\circ g) ^! _T \sp _* E$.

\end{rema}

\section{\label{commsp+}Commutation de $\sp _+$ aux opérations cohomologiques}

\subsection{Commutation aux foncteurs restrictions et images inverses extraordinaires}
Nous reprenons les notations de \ref{notat-construc}.
La proposition qui suit est immédiate.
\begin{prop}\label{EEhatsp+}
  Soient $E$ un isocristal sur $Y$ surconvergent le long de $T _X$,
  $\PP' $ un ouvert de $\PP$, $X ' := X \cap P'$ et $T' :=T \cap P'$.
Il existe alors un isomorphisme canonique :
$\sp _{X \hookrightarrow \PP , T\,+} (E) |_{\PP'} \riso
\sp _{X '\hookrightarrow \PP ', T'\,+} (E |_{]X'[ _{\PP'} }) $ fonctoriel en $E$.

 En particulier, en notant $\widehat{E}$ l'isocristal convergent sur $Y$ associé à $E$,
  $\sp _{X \hookrightarrow \PP , T\,+} (E) |_{\U} \riso
  \sp _{Y \hookrightarrow \U \,+} (\widehat{E})$.
\end{prop}
\begin{prop}\label{spjdag}
Soient $T' \subset T$ un second diviseur de $X$
et $E$ un isocristal sur $X \setminus T'$ surconvergent le long de $T'$.
On dispose de l'isomorphisme :
$(\hdag T ) \sp _{X \hookrightarrow \PP, T '\,+} (E)\riso \sp _{X \hookrightarrow \PP,\,T \,+} (j ^{\dag} E)$
fonctoriel en $E$.
\end{prop}
La preuve étant identique dans le cas général et afin de ne pas alourdir les notations, supposons $T' = \emptyset$.
La proposition résulte aussitôt des trois lemmes qui suivent.
\begin{lemm}\label{lem1spjdag}
  Si $((E _{\alpha})_{\alpha \in \Lambda},\, (\eta _{\alpha\beta}) _{\alpha ,\beta \in \Lambda})
  \in \mathrm{Isoc} ^\dag (X,\,X,\, (\PP _\alpha) _{\alpha \in \Lambda}/K)$
  se recolle en $E\in \mathrm{Isoc} ^\dag (X,\,X/K)$, alors
$(( j ^\dag _\alpha E _{\alpha})_{\alpha \in \Lambda},\, (\eta '_{\alpha\beta}) _{\alpha,  \beta \in \Lambda})$
se recolle en $j ^\dag E$, où $\eta '_{\alpha\beta}$ est l'unique isomorphisme rendant
commutatif le diagramme :
\begin{equation}
  \notag
  \xymatrix  @R=0,3cm {
  { p  _{1 K} ^{\alpha \beta *} ( j ^\dag _\alpha E _{\alpha})}
  \ar[r] _{\sim}
  &
  { j^\dag _{\alpha \beta} p _{1 K}  ^{\alpha \beta *} (E _{\alpha})}
  \\
  { p _{2 K}  ^{\alpha \beta *} ( j^\dag _\beta E _{\beta})}
  \ar[r] _{\sim}
  \ar@{.>}[u] ^-{\eta '_{  \alpha \beta}}
  &
  { j^\dag _{\alpha \beta} p _{2 K}  ^{\alpha \beta *} (E _{\beta}),}
  \ar[u] ^-{j^\dag _{\alpha \beta} \eta _{  \alpha \beta}} _{\sim}
  }
\end{equation}
dont les isomorphismes horizontaux sont ceux de la commutation à l'image inverse
des foncteurs
de la forme $j ^\dag$.
\end{lemm}
\begin{proof}
i) Dans un premier temps, vérifions que la famille d'isomorphismes
$\eta ' _{\alpha \beta}$ définit bien une donnée de recollement.

Soit le diagramme
\begin{equation}
  \notag
\xymatrix  @R=0,3cm {
&
{p _{12 K} ^{\alpha \beta \gamma *}  p  _{1 K} ^{\alpha \beta *} j ^\dag _\alpha (E _{\alpha})}
\ar[r]
&
{p _{12 K} ^{\alpha \beta \gamma *} j^\dag _{\alpha \beta} p _{1 K}  ^{\alpha \beta *} (E _{\alpha})}
\ar[rr]
&
&
{j^\dag _{\alpha \beta \gamma } p _{12 K} ^{\alpha \beta \gamma *}  p _{1 K}  ^{\alpha \beta *} (E _{\alpha})}
\\
{p  _{1 K} ^{\alpha \beta \gamma *}   j ^\dag _\alpha (E _{\alpha})}
\ar[ur] ^-{\epsilon}
\ar[rrr]
&
&
&
{j^\dag _{\alpha \beta \gamma }  p  _{1 K} ^{\alpha \beta \gamma * }(E _{\alpha})}
\ar[ur] _{j^\dag _{\alpha \beta \gamma } (\epsilon)}
\\
&
{p _{12 K} ^{\alpha \beta \gamma *}  p  _{2 K} ^{\alpha \beta *} j ^\dag _\beta (E _{\beta})}
\ar[r]
\ar'[u]^(0.7){p _{12 K} ^{\alpha \beta \gamma *} (\eta ' _{\alpha \beta})}[uu]
&
{p _{12 K} ^{\alpha \beta \gamma *} j^\dag _{\alpha \beta} p _{2 K}  ^{\alpha \beta *} (E _{\beta})}
\ar'[r][rr]
\ar'[u]^(0.7){p _{12 K} ^{\alpha \beta \gamma *} j^\dag _{\alpha \beta} (\eta  _{\alpha \beta})}[uu]
&
&
{j^\dag _{\alpha \beta \gamma } p _{12 K} ^{\alpha \beta \gamma *}  p _{2 K}  ^{\alpha \beta *} (E _{\beta})}
\ar[uu] _{j^\dag _{\alpha \beta \gamma } p _{12 K} ^{\alpha \beta \gamma *} (\eta  _{\alpha \beta})}
\\
{ p  _{2 K} ^{\alpha \beta \gamma *} j ^\dag _\beta (E _{\beta})}
\ar[ur]^\epsilon
\ar[uu] ^-{ \eta ^{\prime \alpha \beta \gamma } _{12}}
\ar[rrr]
&
&
&
{j^\dag _{\alpha \beta \gamma }  p  _{2 K} ^{\alpha \beta \gamma *} (E _{\beta}).}
\ar[uu]^(0.35){ j^\dag _{\alpha \beta \gamma } \eta ^{ \alpha \beta \gamma } _{12}}
\ar[ur] _{j^\dag _{\alpha \beta \gamma }( \epsilon)}
}
\end{equation}
Les carrés du fond, de droite et de gauche sont commutatifs par définition ou par fonctorialité.
De plus, grâce à \ref{epsilonjdagdual}
et
à la transitivité de l'isomorphisme de commutation des foncteurs de la forme
$j ^\dag$ aux images inverses,
il en est de même des diagrammes du haut et du bas.
Toutes les flèches étant des isomorphismes, le diagramme de devant est donc commutatif.
En écrivant les deux autres diagrammes analogues à ce dernier, on vérifie la condition de cocycle.

\medskip

  ii) Par hypothèse, il existe, pour tout $\alpha \in \Lambda$, un isomorphisme
  $\iota _\alpha$ : $E _\alpha \riso u ^* _{\alpha K} ( E |_{]X _\alpha[ _{\PP _\alpha}})$, ceux-ci étant compatibles aux données de
  recollement respectives. Notons $\smash{\widetilde{\eta}} _{\alpha \beta}$,
  la donnée de recollement canonique de $(u ^* _{\alpha K} ( E |_{]X _\alpha[ _{\PP _\alpha}})) _{\alpha \in \Lambda}$, et
  $\smash{\widetilde{\eta}} '_{\alpha \beta}$ celle que l'on déduit pour la famille
  $ (j _\alpha ^\dag u ^* _{\alpha K} ( E |_{]X _\alpha[ _{\PP _\alpha}})) _{\alpha \in \Lambda}$.
Dans le diagramme suivant
  \begin{equation}
    \notag
    \xymatrix  @R=0,3cm @C=1,5cm {
    &
  { j^\dag _{\alpha \beta} p _{1 K}  ^{\alpha \beta *} (E _{\alpha})}
  \ar[rr] ^-{ j^\dag _{\alpha \beta} p _{1 K}  ^{\alpha \beta *} (\iota _\alpha)}
  &
  &
  { j^\dag _{\alpha \beta} p _{1 K}  ^{\alpha \beta *} (u ^* _{\alpha K} ( E |_{]X _\alpha[ _{\PP _\alpha}}))}
  \\
{ p  _{1 K} ^{\alpha \beta *} ( j ^\dag _\alpha E _{\alpha})}
  \ar[ur]
  \ar[rr] ^(0.6){p  _{1 K} ^{\alpha \beta *}  j ^\dag _\alpha (\iota _\alpha)}
&&
  { p  _{1 K} ^{\alpha \beta *} ( j ^\dag _\alpha u ^* _{\alpha K} ( E |_{]X _\alpha[ _{\PP _\alpha}}))}
   \ar[ur]
   \\
  &
  { j^\dag _{\alpha \beta} p _{2 K}  ^{\alpha \beta *} (E _{\beta})}
  \ar'[u]^{j^\dag _{\alpha \beta} \eta _{  \alpha \beta}}[uu]
    \ar'[r]^{j^\dag _{\alpha \beta} p _{2 K}  ^{\alpha \beta *} (\iota _\beta)}[rr]
  &
  &
  { j^\dag _{\alpha \beta} p _{2 K}  ^{\alpha \beta *} (u ^* _{\beta K} ( E |_{]X _\beta[ _{\PP _\beta}}))}
    \ar[uu] ^-{j^\dag _{\alpha \beta} \smash{\widetilde{\eta}} _{\alpha \beta}}
  \\
  { p _{2 K}  ^{\alpha \beta *} ( j^\dag _\beta E _{\beta})}
\ar[rr] ^-{p  _{2 K} ^{\alpha \beta *}  j ^\dag _\beta (\iota _\beta)}
  \ar[uu] ^-{\eta '_{  \alpha \beta}}
  \ar[ur]
&
&
  { p  _{2 K} ^{\alpha \beta *} ( j ^\dag _\beta u ^* _{\beta K} ( E |_{]X _\beta[ _{\PP _\beta}})),}
  \ar[uu] ^(0.4){\smash{\widetilde{\eta}} '_{\alpha \beta}}
  \ar[ur]
  }
  \end{equation}
les carrés du fond, de droite et de gauche sont commutatifs par hypothèse ou définition, tandis que
ceux du haut et du bas le sont par fonctorialité.
Il en résulte que celui de devant l'est aussi.
On obtient ainsi un isomorphisme
$(( j_\alpha ^\dag E _\alpha) _{\alpha \in \Lambda},\, (\eta '_{\alpha \beta}) _{\alpha, \beta \in \Lambda})
\riso
(( j_\alpha ^\dag u ^* _{\alpha K} ( E |_{]X _\alpha[ _{\PP _\alpha}})) _{\alpha \in \Lambda},
\, (\smash{\widetilde{\eta}} '_{\alpha \beta}) _{\alpha, \beta \in \Lambda})$.
\medskip

iii)
Il suffit alors de prouver que l'isomorphisme canonique
$u ^* _{\alpha K} (  j ^\dag E |_{]X _\alpha[ _{\PP _\alpha}}) \riso
j _\alpha ^\dag u ^* _{\alpha K} ( E |_{]X _\alpha[ _{\PP _\alpha}})$ est
compatible aux données de recollement respectives.
Pour cela, il s'agit d'établir la commutativité du carré de gauche
du rectangle du fond du diagramme
\begin{equation}
  \label{lem1spjdagdiag2}
  \xymatrix  @R=0,3cm@C=0,5cm {
  &
  {p _{1 K}  ^{\alpha \beta *} u ^* _{\alpha K}  j ^\dag E |_{]X _\alpha[ _{\PP _\alpha}} }
  \ar[r]
  &
  {p _{1 K}  ^{\alpha \beta *}  j_\alpha ^\dag u ^* _{\alpha K}  E |_{]X _\alpha[ _{\PP _\alpha}}}
  \ar[rr]
  &
  &
  {j _{\alpha \beta} ^\dag p _{1 K}  ^{\alpha \beta *} u ^* _{\alpha K}  E |_{]X _\alpha[ _{\PP _\alpha}}}
  \\
  { u ^* _{\alpha \beta K} j ^\dag E |_{]X _{\alpha \beta}[ _{\PP _{\alpha\beta}}}}
  \ar[ur] ^-{\epsilon}
  \ar[rrr]
  &
  &
  &
  {j _{\alpha \beta} ^\dag u ^* _{\alpha \beta K}  E |_{]X _{\alpha \beta}[ _{\PP _{\alpha\beta}}}}
  \ar[ur] ^(0.4){j _{\alpha \beta} ^\dag  (\epsilon)}
  \\
  &
  {p _{2 K}  ^{\alpha \beta *} u ^* _{\beta K}  j ^\dag E |_{]X _\beta[ _{\PP _\beta}} }
  \ar[r]
  \ar'[u]^{\epsilon}[uu]
  &
  {p _{2 K}  ^{\alpha \beta *}  j_\beta ^\dag u ^* _{\beta K}  E |_{]X _\beta[ _{\PP _\beta}}}
  \ar'[r][rr]
  \ar'[u]^{\smash{\widetilde{\eta}} '_{\alpha \beta}}[uu]
  &
  &
  {j _{\alpha \beta} ^\dag p _{2 K}  ^{\alpha \beta *} u ^* _{\beta K}  E |_{]X _\beta[ _{\PP _\beta}}}
  \ar[uu] ^-{j ^\dag _{\alpha \beta} \smash{\widetilde{\eta}} _{\alpha \beta}}
  \\
      { u ^* _{\alpha \beta K} j ^\dag E |_{]X _{\alpha \beta}[ _{\PP _{\alpha \beta}}}}
  \ar[rrr]
  \ar@{=}[uu]
  \ar[ur] ^-{\epsilon}
  &&
  &
  {j _{\alpha \beta} ^\dag u ^* _{\alpha \beta K}  E |_{]X _{\alpha \beta}[ _{\PP _{\alpha \beta}}}.}
  \ar@{=}[uu]
  \ar[ur] ^(0.4){j _{\alpha \beta} ^\dag  (\epsilon)}
  }
\end{equation}
Le carré de droite du fond, les carrés de droite, de gauche et de devant
sont commutatifs par définition.
En outre, via \ref{epsilonjdagdual} et par transitivité de l'isomorphisme de commutation des foncteurs de la forme
$j ^\dag$ aux images inverses, il en est de même des rectangles du haut et du bas.
Comme les flèches de \ref{lem1spjdagdiag2} sont des isomorphismes, il en découle
la commutativité du carré de gauche du fond.
\end{proof}

\begin{lemm}\label{lem2spjdag}
Avec les notations \ref{defindonnederecol} et de \ref{prop1},
  soient $((\E _{\alpha})_{\alpha \in \Lambda},\, (\theta _{\alpha\beta}) _{\alpha ,\beta \in \Lambda})
  \in \mathrm{Coh} (X,\, (\PP _\alpha) _{\alpha \in \Lambda})$
  et
  $\E\in \mathrm{Coh} (X,\, \PP)$.
   Avec les notations de la preuve de \ref{prop1},
   si $\mathcal{R}ecol ((\E _{\alpha})_{\alpha \in \Lambda},\, (\theta _{\alpha\beta}) _{\alpha ,\beta \in \Lambda})
   \riso \E$, alors
   $\mathcal{R}ecol
(( (\hdag T \cap X _\alpha) \E _{\alpha})_{\alpha \in \Lambda},\,
(\theta ' _{\alpha\beta}) _{\alpha ,\beta \in \Lambda}) \riso (\hdag T ) (\E)$,
où $\theta '_{\alpha\beta}$ est l'unique isomorphisme rendant
commutatif le diagramme :
\begin{equation}
  \notag
  \xymatrix  @R=0,3cm {
  { p  _{1 } ^{\alpha \beta !} ( (\hdag T \cap X _\alpha) (\E _{\alpha}))}
  \ar[r] _{\sim}
  &
  { (\hdag T \cap X _{\alpha \beta} ) \circ p _{1 }  ^{\alpha \beta !} (\E _{\alpha})}
  \\
  { p _{2 }  ^{\alpha \beta !} ((\hdag T \cap X _\beta) (\E _{\beta}))}
  \ar[r] _{\sim}
  \ar@{.>}[u] ^-{\theta '_{  \alpha \beta}}
  &
  { (\hdag T \cap X _{\alpha \beta} )\circ  p _{2 }  ^{\alpha \beta !} (\E _{\beta}),}
  \ar[u] ^-{(\hdag T \cap X _{\alpha \beta} )( \theta _{  \alpha \beta})} _{\sim}
  }
\end{equation}
où les isomorphismes horizontaux résultent de \cite[2.2.18.1]{caro_surcoherent}.
\end{lemm}
\begin{proof}
  Analogue à \ref{lem1spjdag} en remplaçant l'utilisation de la commutativité
  du carré de gauche de \ref{epsilonjdagdual} par celle de \ref{tauhdag}.
\end{proof}

\begin{lemm}\label{lem3spjdag}
  Soit $((E _{\alpha})_{\alpha \in \Lambda},\, (\eta _{\alpha\beta}) _{\alpha ,\beta \in \Lambda})
  \in \mathrm{Isoc} ^\dag (X,\,X,\, (\PP _\alpha) _{\alpha \in \Lambda}/K)$. On note,
  via \ref{lem1spjdag} puis \ref{defdonneesp*},
  $\eta ' _{\alpha \beta}$ (resp. $\theta '' _{\alpha \beta}$) les isomorphismes canoniques de recollement
  de $(j ^\dag _\alpha E _\alpha)_{\alpha \in \Lambda}$
  (resp. $(\sp _* j ^\dag _\alpha E _\alpha)_{\alpha \in \Lambda}$), tandis que ceux de $(\sp _* E _\alpha )_{\alpha \in \Lambda}$
  (resp.
  $((\hdag T \cap X _\alpha) \sp _* E _\alpha )_{\alpha \in \Lambda}$) seront notés $\theta _{\alpha \beta}$
  (resp. $\theta '_{\alpha \beta}$).

  On a l'isomorphisme
  $( (\sp _* j ^\dag _\alpha E _\alpha )_{\alpha \in \Lambda},\, (\theta '' _{\alpha\beta}) _{\alpha, \beta \in \Lambda})
  \riso
  (((\hdag T \cap X _\alpha) (\sp _* E _\alpha))_{\alpha \in \Lambda},\,
  (\theta '_{\alpha\beta}) _{\alpha ,\beta \in \Lambda})$.
\end{lemm}
\begin{proof}
Il s'agit de vérifier la commutativité du carré de droite
du diagramme suivant :
  \begin{equation}
\label{diag1lem3spjdag}
    \xymatrix  @R=0,3cm@C=-0.3cm  {
    &
    { (\hdag T \cap X _{\alpha \beta})  \sp _* p _{1 K}  ^{\alpha \beta !} E _\alpha}
    \ar[rr]
    &&
    { (\hdag T \cap X _{\alpha \beta})  p _{1 }  ^{\alpha \beta !}\sp _* E _\alpha}
    \ar[rr]
    &&
    { p _{1 }  ^{\alpha \beta !}(\hdag T \cap X _\alpha) \sp _* E _\alpha}
    \\
    {\sp _* j _{\alpha \beta} ^\dag p _{1 K}  ^{\alpha \beta !} (E _{\alpha})}
    \ar[rr]
    \ar[ur]
    &&
    {\sp _*  p _{1 K}  ^{\alpha \beta !} j _{\alpha} ^\dag (E _{\alpha})}
    \ar[rr]
    &&
    {p _{1 }  ^{\alpha \beta !} \sp _*   j _{\alpha} ^\dag (E _{\alpha})}
    \ar[ur]
    \\
    &
    { (\hdag T \cap X _{\alpha \beta}) \sp _* p _{2 K}  ^{\alpha \beta !} E _\beta}
    \ar'[u]^{(\hdag T \cap X _{\alpha \beta})(\sp _* \eta _{\alpha \beta})}[uu]
    \ar'[r][rr]
    &&
    { (\hdag T \cap X _{\alpha \beta})  p _{2 }  ^{\alpha \beta !}\sp _* E _\beta,}
    \ar'[u]^{(\hdag T \cap X _{\alpha \beta})(\theta  _{\alpha \beta})}[uu]
    \ar'[r][rr]
    &&
    { p _{2 }  ^{\alpha \beta !}(\hdag T \cap X _\beta) \sp _* E _\beta}
    \ar[uu] ^-{\theta ' _{\alpha \beta}}
    \\
    {\sp _* j _{\alpha \beta} ^\dag p _{2 K}  ^{\alpha \beta !} (E _{\beta})}
    \ar[rr]
    \ar[uu] ^-{\sp _* j _{\alpha \beta} ^\dag \eta _{\alpha \beta}}
    \ar[ur]
    &&
    {\sp _*  p _{2 K}  ^{\alpha \beta !} j _{\beta} ^\dag (E _{\beta})}
    \ar[rr]
    \ar[uu] ^(0.4){\sp _* \eta ' _{\alpha \beta} }
    &&
    {p _{2 }  ^{\alpha \beta !} \sp _*   j _{\beta} ^\dag (E _{\beta}).}
    \ar[ur]
    \ar[uu] ^(0.4){\theta '' _{\alpha \beta}}
    }
  \end{equation}
Le carré de gauche est commutatif par fonctorialité. Ceux de devant et de derrière le sont par définition.
  De plus, grâce à \ref{spcommup*} (les foncteurs de la forme $j ^\dag$ et $(\hdag T)$
  sont respectivement des cas particuliers d'images inverses et d'images inverses extraordinaires),
  les rectangles du haut et du bas sont commutatifs.
  Le carré de droite de \ref{diag1lem3spjdag} l'est donc aussi.
\end{proof}

\begin{prop}\label{f0*=f*isoc}
  Soit le diagramme commutatif
  \begin{equation}
  \label{f0*=f*isoc-diag0}
  \xymatrix  @R=0,3cm {
{Y ' }
\ar[r] ^-{j '}
\ar[d]^h
&
{X '}
\ar[r] ^-{i '}
\ar[d]^g
&
{P '}
\ar[d]^{f_0}
\ar[r]
&
{\PP '}
\\
{Y }
\ar[r]^{j }
&
{X }
\ar[r]^{i }
&
{P }
\ar[r]
&
{\PP ,}
}
\end{equation}
où $\PP$ et $\PP'$ sont des $\V$-schémas formels lisses, $j$ et $j'$ sont des immersions ouvertes
de $k$-schémas lisses, $i $ et $i'$ sont des immersions fermées.

Le foncteur image inverse $g ^*$ :
$\mathrm{Isoc} ^\dag (Y,\,X/K) \rightarrow
\mathrm{Isoc} ^\dag (Y',\,X'/K)$, défini dans \cite[2.3.2.(iv)]{Berig},
est canoniquement isomorphisme au foncteur
$f _{0K} ^*$ de \ref{notation-rig-form}.
\end{prop}
\begin{proof}
  Dans un premier temps, supposons que $X'=X$, $Y' =Y$ et que $f _0$ soit une immersion fermée.
Notons $\delta _{f _0}$ : $P' \hookrightarrow P' \times P $ le graphe de $f _0$,
$p _1$ : $\PP ' \times _\S \PP \rightarrow \PP '$ et $p _2$ : $\PP ' \times _\S \PP \rightarrow \PP $
les projections canoniques.
Soient $E$ un
$ j ^\dag \O _{]X[ _{\PP} }$-module à connexion intégrable et surconvergente
et $E '  $ le $ j ^\dag \O _{]X[ _{\PP '} }$-module à connexion intégrable et surconvergente
correspondant (voir \cite[2.3.2(i)]{Berig}).
Il s'agit de prouver l'isomorphisme : $ E'\riso  f _{0K} ^* (E) $.
Or, par construction, $p _{1K} ^* (E ') \riso p _{2K} ^* (E )$.
Il en dérive
$\delta _{f _0 ,K} ^* p _{1K} ^* (E ')  \riso \delta _{f _0,K} ^* p _{2K} ^* (E )$, ce qui
se réécrit $E '\riso f _{0K} ^* (E) $.

De plus, le cas où le morphisme $f _0$ se relève en un morphisme
$f$ : $\PP' \rightarrow \PP$ est immédiat. Celui où
$X'=X$ et $f _0$ est une immersion est donc aussi résolu.

Passons à présent au cas général. Le morphisme $f _0$ est le composé de l'immersion fermée
$\delta _{f _0}$
suivi de la projection $ P ' \times P \rightarrow P$
dont $\PP ' \times _\S \PP \rightarrow \PP$ constitue un relèvement.
Les cas déjà traités nous permettent de conclure.
\end{proof}

Avant de prouver la proposition suivante, nous aurons besoin de la remarque ci-après.
\begin{rema}\label{exteloc}
  Le foncteur $\mathcal{L}oc$ construit dans la preuve de \ref{prop1}, s'étend, avec ses notations,
  en un foncteur de la catégorie $D ^\mathrm{b} _{\mathrm{coh}} (\smash{\D} ^\dag _{\PP} (\hdag T) _{\Q})$
  dans celle des familles
  $((\E _{\alpha})_{\alpha \in \Lambda},\, (\theta _{\alpha\beta}) _{\alpha ,\beta \in \Lambda})$,
  où $\E _{\alpha} \in D ^\mathrm{b} _{\mathrm{coh}}  (\smash{\D} ^{\dag} _{\X _{\alpha} } (\hdag T  \cap X _{\alpha}) _{\Q})$
  et
$ \theta _{  \alpha \beta} \ : \  p _2  ^{\alpha \beta !} (\E _{\beta}) \riso p  _1 ^{\alpha \beta !} (\E _{\alpha})$
sont des isomorphismes $\smash{\D} ^{\dag} _{\X _{\alpha \beta} }(\hdag T  \cap X _{\alpha \beta}) _{ \Q}$-linéaires
vérifiant la condition de cocycle
$\theta _{13} ^{\alpha \beta \gamma }=
\theta _{12} ^{\alpha \beta \gamma }
\circ
\theta _{23} ^{\alpha \beta \gamma }$,
où $\theta _{12} ^{\alpha \beta \gamma }$, $\theta _{23} ^{\alpha \beta \gamma }$
et $\theta _{23} ^{\alpha \beta \gamma }$ sont définis par les diagrammes commutatifs analogues à
  \ref{diag1-defindonnederecol}.

Si $\E \in D ^\mathrm{b} _{\mathrm{coh}} (\smash{\D} ^\dag _{\PP} (\hdag T) _{\Q})$,
le morphisme canonique $\mathcal{L}oc ( \R \underline{\Gamma} ^\dag _{X } \E)
\rightarrow \mathcal{L}oc (\E)$ est un isomorphisme.

\end{rema}

\begin{prop}\label{sp+f*}
On garde les notations et hypothèses de \ref{f0*=f*isoc}
et on suppose qu'il existe un diviseur $T$ de $P$ tel que
$T ':= f _0 ^{-1} (T)$ (resp. $T \cap X$ et $T' \cap X'$) soit un diviseur de $P'$
(resp. $X$ et $X'$) et tel que $Y = X \setminus T$ et $Y' =X' \setminus T'$.

Pour tout isocristal $E$ sur $Y$ surconvergent le long de $T \cap X$,
on dispose de l'isomorphisme :
$$\sp _{X '\hookrightarrow \PP',\,T ' \,+} f _{0K}^*  (E)
\riso
\R \underline{\Gamma} ^\dag _{X '} \, f _0  ^!\, \sp _{X \hookrightarrow \PP,\,T  \,+} (E)[-d _{X'/X}] $$
fonctoriel en $E$.
\end{prop}
\begin{proof}

Traitons d'abord le cas où $\PP$ et $\PP' $ sont affines.
Comme $f _0$ se relève et, pour tout
relèvement $i '$ : $\X' \hookrightarrow \PP'$ de l'immersion fermée $X' \hookrightarrow P'$,
grâce à l'isomorphisme $\R \underline{\Gamma} ^\dag _{X '} \riso i' _+ \circ i ^{\prime !}$,
ce cas est validé.

Puisque cela est local en $\PP$ et $\PP'$, il résulte du premier cas que
$\R \underline{\Gamma} ^\dag _{X '} \, f _0  ^!\, \sp _{X \hookrightarrow \PP,\,T  \,+} (E)[-d _{X'/X}] $
est un $\smash{\D} ^\dag _{\PP' } (\hdag T') _{\Q}$-module cohérent à support dans $X'$.
Via \ref{prop1}, il suffit alors de prouver que l'on dispose d'un isomorphisme
$\mathcal{L} oc (\sp _{X '\hookrightarrow \PP',\,T ' \,+} f _{0K}^*  (E))
\riso
\mathcal{L} oc
(\R \underline{\Gamma} ^\dag _{X '} \, f _0  ^!\, \sp _{X \hookrightarrow \PP,\,T  \,+} (E)[-d _{X'/X}] )$
(notations de \ref{prop1}).

À cette fin, choisissons une application surjective $\rho$ : $\Lambda ' \rightarrow \Lambda$,
deux recouvrements ouverts affines $(\PP _\alpha) _{\alpha \in \Lambda}$
de $\PP$ et $(\PP '_{\alpha'}) _{\alpha '\in \Lambda '}$ de $\PP'$
  tels que $f _0$ se factorise par $P '_{\alpha '} \rightarrow P _{\rho(\alpha ')}$.
  On note alors $X _\alpha := X \cap P _{\alpha}$,
  $X _{\alpha \beta} := X \cap P _{\alpha}\cap P _{\beta}$, $j _\alpha$ :
  $\PP _\alpha \subset \PP$, et de même en rajoutant des primes.
De plus, choisissons des relèvements $ f _{\alpha '}$ :
$\PP '_{\alpha '} \rightarrow \PP _{\rho(\alpha ')}$,
$i _{\alpha}$ : $\X _{\alpha } \hookrightarrow \PP _{\alpha}$,
$i ' _{\alpha'}$ : $\X ' _{\alpha } \hookrightarrow \PP '_{\alpha '}$
et
$g _{\alpha '}$ : $ \X' _{\alpha '} \rightarrow \X _{\rho(\alpha')}$
des relèvements des factorisations induites par $f _0$, $i$, $i'$ et $g$.

Grâce à \ref{exteloc}, le morphisme
$\mathcal{L} oc (\R \underline{\Gamma} ^\dag _{X '} \, f _0  ^!\, \sp _{X \hookrightarrow \PP,\,T  \,+} (E) )
\rightarrow
\mathcal{L} oc (f _0  ^!\, \sp _{X \hookrightarrow \PP,\,T  \,+} (E) )$ est un isomorphisme.
On se ramène ainsi à établir
$\mathcal{L} oc (\sp _{X '\hookrightarrow \PP',\,T ' \,+} f _{0K}^*  (E))
\riso
\mathcal{L} oc (f _0  ^!\, \sp _{X \hookrightarrow \PP,\,T  \,+} (E) [-d _{X'/X}]) $.

L'objet $\mathcal{L} oc (\sp _{X '\hookrightarrow \PP',\,T ' \,+} f _{0K}^*  (E))$
correspond à la famille
$( \sp _* i ^{\prime *} _{\alpha 'K} f ^* _{\alpha ' K} ( E |_{]X _{\rho (\alpha ')}[ _{\PP _{\rho (\alpha ')} }}))
_{\alpha ' \in \Lambda '} $,
dont les isomorphismes de recollement sont induits par ceux de la forme $\epsilon$
(voir \ref{notation-rig-form-diag0}).

L'objet
$\mathcal{L} oc (f _0  ^!\, \sp _{X \hookrightarrow \PP,\,T  \,+} (E) )$ est
isomorphe à la famille
$( i ^{\prime !} _{\alpha '} j ^{\prime !} _{\alpha ' } f _0  ^!\,
\sp _{X \hookrightarrow \PP,\,T  \,+} (E) [-d _{X'/X}] )
_{\alpha ' \in \Lambda '}$
dont les isomorphismes de recollement sont induits par ceux de la forme $\tau$.

Il résulte des isomorphismes
$ i ^{\prime !} _{\alpha '} j ^{\prime !} _{\alpha ' } f _0  ^!
\riso  i ^{\prime !} _{\alpha '} f _{\alpha '}  ^! j ^! _{\rho (\alpha ') }
\riso  g _{\alpha '}  ^!   i ^{!} _{\rho (\alpha ')} j ^! _{\rho (\alpha ') }$
que
$\mathcal{L} oc (f _0  ^!\, \sp _{X \hookrightarrow \PP,\,T  \,+} (E) )$
correspond à la famille
$( g _{\alpha '}  ^!   i ^{!} _{\rho (\alpha ')} j ^! _{\rho (\alpha ') }\,
\sp _{X \hookrightarrow \PP,\,T  \,+} (E) [-d _{X'/X}] )
_{\alpha ' \in \Lambda '}$
dont les isomorphismes de recollement sont induits par ceux de la forme $\tau$.
Puisque $\mathcal{L}oc \circ \sp _{X \hookrightarrow \PP,\,T  \,+} (E)
\riso \sp _* \mathcal{L}oc  (E)$ puis grâce à \ref{sp-eps-tau}, on obtient \newline
$g _{\alpha '}  ^!   i ^{!} _{\rho (\alpha ')} j ^! _{\rho (\alpha ') }
\, \sp _{X \hookrightarrow \PP,\,T  \,+} (E) ) [-d _{X'/X}]
\riso g _{\alpha '}  ^!   \sp _* i ^{*} _{\rho (\alpha ')K}
(E |_{]X _{\rho (\alpha ')}[ _{\PP _{\rho (\alpha ')} }}) [-d _{X'/X}]
\riso \sp _* g _{\alpha 'K}  ^*    i ^{*} _{\rho (\alpha ')K}
(E |_{]X _{\rho (\alpha ')}[ _{\PP _{\rho (\alpha ')} }})$, ceux-ci étant
compatibles aux données de recollement respectives.
On conclut ensuite via l'isomorphisme canonique
$\epsilon$ : $i ^{*} _{\rho (\alpha ')K} \riso i ^{\prime *} _{\alpha 'K} f ^* _{\alpha ' K} $.

\end{proof}
\begin{coro}\label{coro-frob-sp+}
(i) Pour tout isocristal $E$ sur $Y$ surconvergent le long de $T _X$,
on dispose d'un isomorphisme
$\sp _{X \hookrightarrow \PP,\,T  \,+} ( F ^* E)
\riso F ^* \sp _{X \hookrightarrow \PP,\,T  \,+} (E)$ fonctoriel en $E$.
On notera $F$-$\mathrm{Coh} (X,\, \PP ,\, T)$, la catégorie des
$F$-$\smash{\D} ^{\dag} _{\PP }(\hdag T ) _{\Q}$-modules cohérents à support dans $X$.

  (ii) Le foncteur $ \sp _{X \hookrightarrow \PP,\,T  \,+}$ :
  $F$-$\mathrm{Isoc} ^\dag (Y,\,X/K) \rightarrow F$-$\mathrm{Coh} (X,\, \PP ,\, T)$
  est pleinement fidèle et son image essentielle est constituée par les
  $F$-$\smash{\D} ^{\dag}  _{\PP}(\hdag T) _{\Q}$-modules cohérents $\E$ à support dans $X$ satisfaisant la condition
  suivante :

  (*) pour tout ouvert $\PP '$ de $\PP$ tel que l'immersion fermée $X \cap P ' \hookrightarrow P'$
  se relève en un morphisme $v$ : $ \X' \hookrightarrow \PP '$ de $\V$-schémas formels lisses,
  le faisceau $v ^! (\E |_{\PP'}) $ est $\O _{\X',\,\Q} (\hdag T \cap X')$-cohérent.
\end{coro}
\begin{proof}
L'isomorphisme $\sp _{X \hookrightarrow \PP,\,T  \,+} ( F ^* E)
\riso F ^* \sp _{X \hookrightarrow \PP,\,T  \,+} (E)$
résulte de \ref{f0*=f*isoc} et de \ref{sp+f*}.
En effet, pour tout $\smash{\D} ^{\dag} _{\PP }(\hdag T ) _{\Q}$-module
cohérent $\E$ à support dans $X$, le faisceau $F ^* (\E) $ est aussi à support dans $X$ et donc
$\R \underline{\Gamma} ^\dag _X F ^* (\E) \riso F ^* (\E)$.
Le foncteur $ \sp _{X \hookrightarrow \PP,\,T  \,+}$ induit donc le suivant
  $F$-$\mathrm{Isoc} ^\dag (Y,\,X/K) \rightarrow F$-$\mathrm{Coh} (X,\, \PP ,\, T)$.
La fidélité de celui-ci résulte du cas sans structure de Frobenius (\ref{sp+plfid}).
Pour la pleine fidélité, soient
$E$ et $E'$ deux objets de $F$-$\mathrm{Isoc} ^\dag (Y,\,X/K) $
et $g$ un morphisme $\sp _{X \hookrightarrow \PP,\,T  \,+} (E) \rightarrow
\sp _{X \hookrightarrow \PP,\,T  \,+} (E')$ commutant à Frobenius.
D'après \ref{sp+plfid}, il existe un morphisme
$f $ : $E\rightarrow E'$ dans $\mathrm{Isoc} ^\dag (Y,\,X/K) $
tel que $\sp _{X \hookrightarrow \PP,\,T  \,+} (f) = g$.
Puisque le foncteur $\sp _{X \hookrightarrow \PP,\,T  \,+}$ (sans structure de Frobenius)
est fidèle, le morphisme $f$ commute à Frobenius.
On a donc vérifié la pleine fidélité.
Enfin, la description de l'image essentielle découle de celle de \ref{sp+plfid}
ainsi que de la pleine fidélité de $\sp _{X \hookrightarrow \PP,\,T  \,+}$ (sans structure de Frobenius).
\end{proof}

\begin{rema}
  Les isomorphismes \ref{sp+f*} et \ref{spjdag} sont compatibles à Frobenius.
\end{rema}

\subsection{Isomorphisme de commutation des foncteurs duaux aux images inverses extraordinaires
par une immersion}
Cette section est dédiée à l'établissement du résultat suivant : soient
 $f$ : $ \X ' \hookrightarrow \X$ une immersion de $\V$-schémas formels lisses,
  $T$ un diviseur de $X$ tel que $T' :=f ^{-1} (T)$ soit un diviseur de $X'$,
  $\E \in D ^{\mathrm{b}}_{\mathrm{coh}} ( \overset{^\mathrm{g}}{} \smash{\D} ^{\dag}   _{\X}(\hdag T) _\Q)$
 à support dans l'adhérence schématique $\overline{X'}$ de $X'$ dans $X$.
 On construit, via \ref{def-thetafeq1}, un isomorphisme canonique
$\theta _{f,T,\E}$ : $f _T ^! \DD _T (\E) \riso  \DD _{T'} f ^! _T  (\E)$.
Celui-ci est en outre transitif pour
la composition d'immersions (voir \ref{trans-theta-f2}).
\begin{vide}
    \label{videdef-dg!-g!d}
Soient $g$ : $\X' \hookrightarrow \X$ une immersion fermée de $\V$-schémas formels lisses,
$T$ un diviseur de $X$ tel que $T' := T \cap X'$ soit un diviseur de $X'$.
  Pour tout
  $\E \in D ^{\mathrm{b}}_{\mathrm{coh}} ( \overset{^\mathrm{g}}{} \smash{\D} ^\dag _{\X,\Q}(\hdag T))$ à support dans $X'$,
  on définit l'isomorphisme $\theta _{g,T,\E}$ : $ g ^! _T \circ \DD _T (\E)\riso \DD _{T'} \circ  g ^! _T (\E) $
  via le diagramme commutatif suivant :
  \begin{equation}
    \label{def-dg!-g!d}
    \xymatrix  @R=0,3cm    @C=2cm {
    {g^! _T \DD _{T} g_{T,+} g^! _T (\E)}
    &
    {g ^! _T g _{T,+} \DD _{T'} g ^! _T (\E)}
    \ar[l] _-{g ^! _T \chi _g g ^! _T } ^-\sim
    \\
    {g ^! _T \DD _T (\E)}
    \ar[u] ^-{g^! _T \DD _{T} \mathrm{adj} _{g,T,\E}} _-\sim
    \ar@{.>}[r] _-{\theta _{g,T,\E}} ^-\sim
    &
    {\DD _{T'} g^! _T (\E).}
    \ar[u] _-{\mathrm{adj} _{g,T,\DD _{T'} g^! _T (\E) }} ^-\sim
    }
  \end{equation}
On pourra le noter plus simplement $\theta _{g,T}$ ou $\theta _{g,\E}$ voire $\theta _g$.
\end{vide}
\begin{vide}
\label{remdef-dg!-g!d}
Avec les notations de \ref{videdef-dg!-g!d}, on aurait pu définir d'une seconde façon,
via le diagramme commutatif \ref{def-dg!-g!d2}, l'isomorphisme $\theta _{g,T,\E}$ :
pour vérifier la commutativité de \ref{def-dg!-g!d2},
on remarque d'abord que les contours des deux diagrammes ci-dessous :
\begin{gather}
\label{def-dg!-g!d1}
    \xymatrix  @R=0,3cm     {
{g_{T,+} g_{T} ^! \DD _{T} (\E)}
\ar[r] ^-{\mathrm{adj} _g} _-\sim
\ar[d] ^-{\theta _g} _-\sim
&
{g_{T,+} g_{T} ^! \DD _{T} g_{T,+} g_{T} ^! (\E)}
\ar[r] ^-{\mathrm{adj} _g} _-\sim
\ar[d] ^-{\chi_g} _-\sim
&
{\DD _{T} g_{T,+} g_{T} ^! (\E)}
\ar[d] ^-{\chi_g} _-\sim
\\
{g_{T,+} \DD _{T'} g_{T} ^!  (\E)}
\ar[r] ^-{\mathrm{adj} _g} _-\sim
&
{g_{T,+} g_{T} ^! g_{T,+} \DD _{T'}  g_{T} ^! (\E)}
\ar[r] ^-{\mathrm{adj} _g} _-\sim
&
{g_{T,+} \DD _{T'} g_{T} ^!  (\E),}
}
\\
\label{def-dg!-g!d2}
 \xymatrix  @R=0,3cm    {
{g_{T,+} g_{T} ^! \DD _{T} (\E)}
\ar[r] ^-{\mathrm{adj} _g} _-\sim
\ar[d] ^-{\theta _g} _-\sim
&
{\DD _{T}  (\E)}
\ar[r] ^-{\mathrm{adj} _g} _-\sim
&
{\DD _{T} g_{T,+} g_{T} ^! (\E)}
\ar[d] ^-{\chi_g} _-\sim
\\
{g_{T,+} \DD _{T'} g_{T} ^!  (\E)}
\ar@{=}[rr]
&
&
{g_{T,+} \DD _{T'} g_{T} ^!  (\E),}
}
\end{gather}
sont égaux.
En effet, les composés des morphismes du haut sont égaux par fonctorialité (on établit la commutativité
du carré correspondant),
tandis que pour ceux du bas, cela découle fonctoriellement d'une propriété standard des foncteurs adjoints.
Or, le carré de gauche de \ref{def-dg!-g!d1} est l'image de \ref{def-dg!-g!d} par $g _{T,+}$.
La commutativité de celui de droite étant fonctorielle, il en résulte celle
de \ref{def-dg!-g!d1}. Ainsi, \ref{def-dg!-g!d2} est commutatif.
Comme $g$ est une immersion fermée, d'après l'analogue $p$-adique de Berthelot du théorème de Kashiwara,
la commutativité de \ref{def-dg!-g!d2} détermine $\theta _g$.
\end{vide}

\begin{prop}
  \label{tausigma-thg}
  Soient $g$, $g'$ : $\X '\rightarrow \X$ deux immersions fermées de $\V$-schémas formels lisses telles que
  $g _0 =g '_0$,
  $T$ un diviseur de $X$ tel que $ T' := T \cap X '$ soit un diviseur de $X'$.
  On dispose du diagramme commutatif :
\begin{equation}
\label{tausigma-thgdiag0}
  \xymatrix  @R=0,3cm  {
  {g ^! _T \DD _T (\E)}
    \ar[r] _-{\theta _{g,T,\E}} ^-\sim
    \ar[d] _-{\tau _{ g', g}} ^-\sim
    &
  {\DD _{T'} g^! _T (\E)}
  \\
  {g ^{\prime !} _T \DD _T (\E)}
  \ar[r] _-{\theta _{g',T,\E}} ^-\sim
    &
  {\DD _{T'} g ^{\prime !} _T (\E) .}
  \ar[u] _-{\tau _{ g', g}} ^-\sim
  }
\end{equation}
\end{prop}
\begin{proof}
Considérons le diagramme suivant :
\begin{equation}
\label{tausigma-thgdiag}
  \xymatrix  @R=0,3cm {
    &
    {g^! _T \DD _{T} g_{T,+} g^! _T (\E)}
    &&
    {g ^! _T g _{T,+} \DD _{T'} g ^! _T (\E)}
    \ar[ll] _-{g ^! _T \chi _g g ^! _T } ^-\sim
    \\
    {g ^! _T \DD _T (\E)}
    \ar[ur] ^-{g^! _T \DD _{T} \mathrm{adj} _{g,T,\E}} _-\sim
    \ar[rr] _-(0.3){\theta _{g,T,\E}} ^-(0.3)\sim
    \ar[dd] _-{\tau _{ g', g}} ^-\sim
    &&
    {\DD _{T'} g^! _T (\E)}
    \ar[ur] _-{\mathrm{adj} _{g,T,\DD _{T'} g^! _T (\E) }} ^-\sim
    \\
    &
    {g^{\prime !} _T \DD _{T} g '_{T,+} g ^{\prime !} _T (\E)}
        \ar'[u]_-{\tau _{ g', g}\sigma _{ g', g}\tau _{ g', g}} ^-\sim [uu]
    &&
    {g ^{\prime !} _T g '_{T,+} \DD _{T'} g ^{\prime !} _T (\E)}
    \ar'[l][ll] _-{g ^{\prime !} _T \chi _{g'} g ^{\prime !} _T } ^-\sim
    \ar[uu] _-(0.6){\tau _{ g', g}\sigma _{ g', g}\tau _{ g', g}} ^-(0.6)\sim
    \\
    {g ^{\prime !} _T \DD _T (\E)}
    \ar[ur] ^-{g ^{\prime !} _T \DD _{T} \mathrm{adj} _{g',T,\E}} _-\sim
    \ar[rr] _-{\theta _{g',T,\E}} ^-\sim
    &&
    {\DD _{T'} g ^{\prime !} _T (\E).}
    \ar[ur] _-{\mathrm{adj} _{g',T,\DD _{T'} g^{\prime !} _T (\E) }} ^-\sim
    \ar[uu] _-(0.3){\tau _{ g', g}} ^-(0.3)\sim
    }
\end{equation}
Il découle de \cite[2.4.4.2]{caro-frobdualrel} et
\cite[2.5.1]{caro-frobdualrel} (ou \cite[2.6.2]{caro-frobdualrel})
que les carrés de droite, de gauche et de derrière sont commutatifs.
Comme les deux horizontaux le sont par définition (\ref{def-dg!-g!d}),
le dernier carré de \ref{tausigma-thgdiag}, celui de devant,
est donc commutatif.
\end{proof}

\begin{prop}\label{trans-thg}
Soient $g'$ : $\X'' \rightarrow \X'$ et $g $ : $ \X '\rightarrow \X$
deux immersions fermées de $\V$-schémas formels lisses,
$T$ un diviseur de $X$ tel que $T':=g ^{-1} (T)$ (resp. $T'':=g ^{\prime -1} (T')$)
  soit un diviseur de $X'$ (resp. $X''$).
Pour tout
$\E \in D ^{\mathrm{b}}_{\mathrm{coh}} ( \overset{^\mathrm{g}}{} \smash{\D} ^{\dag}   _{\X,\Q}(\hdag T))$
à support dans $X''$,
le diagramme
$$\xymatrix  @R=0,3cm    {
{\DD _{T''} g ^{\prime !} _{T'} g^! _T (\E)}
&
{g ^{\prime !} _{T'} \DD _{T'}  g^! _T(\E)}
\ar[l] ^-{\theta _{g',T'}} _-\sim
&
{g ^{\prime !} _{T'} g^! _T \DD _T (\E)}
\ar[d] _-\sim
\ar[l] ^-{\theta _{g,T}} _-\sim
\\
{ \DD _{T''} (g \circ g ')^! _T (\E)}
\ar[u] _-\sim
&&
{(g \circ g ')^! _T  \DD _T(\E)}
\ar[ll] ^-{\theta _{g\circ g',T}} _-\sim
}
$$
est commutatif.
\end{prop}
\begin{proof}
Quitte à alourdir les notations, on suppose que le diviseur $T$ est vide et on n'indiquera pas {\og $(\E)$\fg}.
Considérons le cube ci-dessous :
  \begin{equation}
    \label{trans-thg-diag1}
  \xymatrix  @R=0,3cm   @C=-0,3cm {
  &&
  {g ^{\prime !} g ^! g _+ g' _+ \DD g ^{\prime !} g^!}
  \ar[rr] _-\sim ^-{\chi _{g'}}
  \ar'[d]'[dd][ddd] _-\sim
  &&
  {g ^{\prime !} g ^! g _+ \DD g' _+  g ^{\prime !} g^!}
  \ar[rr] _-\sim ^-{\chi _{g}}
  &&
  {g ^{\prime !} g ^! \DD g _+  g' _+  g ^{\prime !} g^!}
  \\
  &
  {g ^{\prime !}  g' _+ \DD g ^{\prime !} g^!}
  \ar[rr] _-\sim ^-{\chi _{g'}}
  \ar[ur] _-\sim ^-{\mathrm{adj} _g}
  &&
  {g ^{\prime !} \DD g' _+  g ^{\prime !} g^!}
  \ar[ur] _-\sim ^-{\mathrm{adj} _g}
  &&
  {g ^{\prime !} g ^! \DD g _+  g^!}
  \ar[ur] _-\sim ^-{\mathrm{adj} _{g'}}
  &
  \\
  {\DD g ^{\prime !} g^!}
  \ar[ur] _-\sim ^-{\mathrm{adj} _{g'}}
  &&
  {g ^{\prime !} \DD  g^!}
  \ar[ll] _-\sim ^-{\theta _{g'}}
  \ar[ur] _-\sim ^-{\mathrm{adj} _{g'}}
  &&
  {g ^{\prime !} g^! \DD }
  \ar[ur] _-\sim ^-{\mathrm{adj} _{g}}
  \ar[ll] _-\sim ^-{\theta _{g}}
  \ar[ddd] _-\sim
  &&
  \\
  &&
  {g ^{\prime !} g ^! g \circ g' _+ \DD g ^{\prime !} g^!}
  \ar'[rr][rrrr] _-\sim ^-{\chi _{g\circ g'}}
  \ar[ld] _-\sim
  &&&&
  {g ^{\prime !} g ^! \DD g \circ g' _+  g ^{\prime !} g^!}
  \ar[dl] _-\sim
  \ar[uuu] _-\sim
  \\
  &
  { g\circ g ^{\prime !} g \circ g' _+ \DD g\circ g ^{\prime !}}
  \ar'[rrr][rrrr] _-\sim ^-{\chi _{g\circ g'}}
  &&&&
  {g\circ g ^{\prime !} \DD g \circ g' _+  g\circ g ^{\prime !}}
  &
  \\
  {\DD g\circ g ^{\prime !}}
  \ar[ur] _-\sim ^-{\mathrm{adj} _{g \circ g'}}
  \ar[uuu] _-\sim
  &&&&
  {g\circ g ^{\prime !}  \DD }
  \ar[ur] _-\sim ^-{\mathrm{adj} _{g \circ g'}}
  \ar[llll] _-\sim ^-{\theta _{g\circ g'}}
  }
  \end{equation}
Le carré du fond de \ref{trans-thg-diag1} correspond, en omettant les termes $g ^!$ et $g ^{\prime !}$,
au diagramme signifiant
la transitivité des isomorphismes de la forme $\chi $. Celui-ci est donc commutatif.
Par fonctorialité ou définition (\ref{def-dg!-g!d}), on obtient la commutativité de la face du bas.
On vérifie de même celle de la partie gauche de la face supérieure de \ref{trans-thg-diag1}.
Sa partie droite s'identifie (i.e. les composés des contours respectifs sont égaux) au contour du diagramme
commutatif
$$\xymatrix  @R=0,3cm     {
{g ^{\prime !} \DD  g^!}
  \ar[r] _-\sim ^{\mathrm{adj} _{g}}
  &
  {g ^{\prime !} g ^! g _+ \DD  g^!}
  \ar[d] _-\sim ^{\chi _{g}}
  \ar[r] _-\sim ^{\mathrm{adj} _{g'}}
  &
  {g ^{\prime !} g ^! g _+ \DD g' _+  g ^{\prime !} g^!}
  \ar[d] _-\sim ^{\chi _{g}}
  \\
  {g ^{\prime !}   g^! \DD}
  \ar[r] _-\sim ^{\mathrm{adj} _{g}}
  \ar[u] _-\sim ^{\theta _{g}}
  &
  {g ^{\prime !} g ^! \DD g _+   g^!}
  \ar[r] _-\sim ^{\mathrm{adj} _{g'}}
  &
  {g ^{\prime !} g ^! \DD g _+  g' _+  g ^{\prime !} g^!.}
}$$
La face du haut de \ref{trans-thg-diag1} est donc commutative.

En appliquant $g ^{\prime !}   g^! \DD$ au carré de gauche de \ref{tr-transifdag}
(en remplaçant $f$ et $g$ par $g$ et $g'$), on obtient le rectangle supérieur de
\begin{equation}
  \label{trans-thg-diag2}
  \xymatrix  @R=0,3cm     {
  {g ^{\prime !}   g^! \DD}
  \ar[r] _-\sim ^{\mathrm{adj} _{g}}
  &
  {g ^{\prime !} g ^! \DD g _+   g^!}
  \ar[r] _-\sim ^{\mathrm{adj} _{g'}}
  &
  {g ^{\prime !} g ^! \DD g _+  g' _+  g ^{\prime !} g^!}
  \\
  {g ^{\prime !}   g^! \DD}
  \ar@{=}[u]
  \ar[rr] _-\sim ^{\mathrm{adj} _{g\circ g'}}
  \ar[d] _-\sim
  &&
  {g ^{\prime !} g ^! \DD g \circ g' _+  g \circ g ^{\prime !}  }
  \ar[u] _-\sim
  \\
  {g \circ g ^{\prime !} \DD  }
  \ar[r] _-\sim ^-{\mathrm{adj} _{g\circ g'}}
  &
  {g \circ g ^{\prime !} \DD g \circ g' _+  g \circ g ^{\prime !}  }
  &
  {g ^{\prime !} g ^! \DD g \circ g' _+  g ^{\prime !} g ^! .}
  \ar[u] _-\sim
  \ar[l] _-\sim
}
\end{equation}
Le diagramme \ref{trans-thg-diag2} est donc commutatif. Son contour, i.e.,
le carré de droite de \ref{trans-thg-diag1} l'est donc aussi.
De même, on obtient la commutativité du carré de gauche de \ref{trans-thg-diag1}.

Comme les flèches de \ref{trans-thg-diag1} sont des isomorphismes,
comme cinq des faces du cube de \ref{trans-thg-diag1} sont commutatifs, le
dernier, i.e., celui de devant, égal à celui de \ref{trans-thg}, est commutatif.
\end{proof}

\begin{lemm}\label{def-thetaf}
  Soient $f$ : $ \X ' \hookrightarrow \X$ une immersion de $\V$-schémas formels lisses,
  $T$ un diviseur de $X$ tel que $T' :=f ^{-1} (T)$ soit un diviseur de $X'$,
  $\E \in D ^{\mathrm{b}}_{\mathrm{coh}} ( \overset{^\mathrm{g}}{} \smash{\D} ^{\dag}   _{\X}(\hdag T) _\Q)$
 à support dans l'adhérence schématique $\overline{X'}$ de $X'$ dans $X$.

  Choisissons $\Y$ un ouvert de $\X$ tel que $f$ se factorise
  en une immersion fermée $g $ : $ \X ' \hookrightarrow \Y$.
En notant $j$ : $\Y \subset \X$ l'inclusion,
l'isomorphisme composé
\begin{equation}
\label{def-thetafeq1}
  f _T ^! \DD _T (\E) \riso
  g ^! _{T \cap Y} j _T ^! \DD (\E)
  =
  g ^! _{T \cap Y} \DD _{T \cap Y} j ^! _T (\E)
  \tilde{\underset{\theta _g}{\longrightarrow}}
  \DD _{T'} g ^! _{T \cap Y}  j _T ^! (\E)
  \riso
  \DD _{T'} f ^! _T  (\E)
\end{equation}
ne dépend pas du choix de la factorisation de $f$ par $g$ et $j$.
\end{lemm}
\begin{proof}
  Pour le vérifier, quitte à alourdir les notations, supposons le diviseur $T$ vide.
  Notons $\Y _0$ l'ouvert de $\X$ dont l'espace sous-jacent
  est $X \setminus (\overline{X'} \setminus X')$. Comme $\Y \subset \Y _0$
  et que le morphisme induit (par $g$) $\X '\hookrightarrow \Y _0$ est une immersion fermée,
  on remarque qu'il suffit de prouver cette indépendance lorsque
  $f$ est une immersion fermée. Supposons donc $f$ est une immersion fermée.
Considérons le diagramme ci-dessous :
\begin{equation}
  \label{def-thetafeq2}
  \xymatrix  @R=0,3cm   {
  &
  {j ^! f_{+} f  ^! \DD  (\E)}
\ar[rr] ^-{\mathrm{adj} _f}
\ar[dl] ^-{\theta _f}
\ar'[d][dd]
&&
{j ^!  \DD (\E)}
\ar[rr] ^-{\mathrm{adj} _f}
\ar'[d][dd]
&&
{j ^!  \DD  f _{+} f  ^! (\E )}
\ar[dl] ^-{\chi _f}
\ar[dd]
\\
{j ^! f _{+} \DD f ^!  (\E) }
\ar@{=}[rrrr]
&&&&
{j ^!  f _{+} \DD f ^!  (\E) }
\\
&
{g_{+} g  ^! \DD j ^! (\E)}
\ar[rr] ^-{\mathrm{adj} _g}
\ar[dl] ^-{\theta _g}
&&
{ \DD j ^! (\E)}
\ar'[r]^-{\mathrm{adj} _g} [rr]
&&
{ \DD  g _{+} g  ^! j ^! (\E )}
\ar[dl] ^-{\chi _g}
\\
{g _{+} \DD g ^! j ^!  (\E) }
\ar[uu]
\ar@{=}[rrrr]
&&&&
{g _{+} \DD g ^! j ^!  (\E) ,}
\ar[uu]
}
\end{equation}
dont les flèches verticales se déduisent des isomorphismes canoniques
$f ^! \riso g ^! j ^!$ et $j ^! f _+ \riso g _+$.
Il résulte de \ref{def-dg!-g!d2} que les faces horizontales sont commutatives.
La commutativité des faces de devant, derrière et de droite se vérifie par construction de $\chi$
et $\mathrm{adj}$.
Comme toutes les flèches de \ref{def-thetafeq2} sont des isomorphismes,
il en résulte celle du carré de gauche.
Via l'isomorphisme canonique $j ^! f _+ \riso g _+$ et grâce à l'analogue
$p$-adique du théorème de Kashiwara,
le carré de gauche de \ref{def-thetafeq2} dont on a enlevé les termes $j ^! f _+$
et $g _+$ est alors commutatif.
Cela se traduit par le fait que le composé \ref{def-thetafeq1} est $\theta _f$.
C.Q.F.D.

\end{proof}

\begin{nota}
  Avec les notations de \ref{def-thetaf},
  on désignera par $\theta _{f,T,\E}$ : $  f _T ^! \DD _T (\E) \riso
  \DD _{T'} f ^! _T  (\E)$
  l'isomorphisme \ref{def-thetafeq1}.
\end{nota}

\begin{prop}\label{taudual}
  Soient $f $, $f'$ : $\X ' \hookrightarrow \X$ deux immersions de $\V$-schémas formels lisses
  telles que $f _0 = f '_0$, $T$ un diviseur de $X$ tel que $T' :=f _0 ^{-1} (T)$ soit un diviseur de $X'$.
  On bénéficie du diagramme commutatif :
  \begin{equation}
\label{diag-taudual}
  \xymatrix  @R=0,3cm  {
  {f ^! _T \DD _T (\E)}
    \ar[r] _-{\theta _{f,T,\E}} ^-\sim
    \ar[d] _-{\tau _{ f', f}} ^-\sim
    &
  {\DD _{T'} f^! _T (\E)}
  \\
  {f ^{\prime !} _T \DD _T (\E)}
  \ar[r] _-{\theta _{f',T,\E}} ^-\sim
    &
  {\DD _{T'} f ^{\prime !} _T (\E). }
  \ar[u] _-{\tau _{ f', f}} ^-\sim
  }
\end{equation}
\end{prop}
\begin{proof}
  Notons $\Y _0$ l'ouvert de $\X$ dont l'espace sous-jacent
  est $X \setminus (\overline{X'} \setminus X')$,
  $j $ : $\Y _0 \subset \X$ l'inclusion canonique,
  $g$ et $g'$ : $\X ' \hookrightarrow \Y _0$ les immersions fermées
  factorisant respectivement $f$ et $f'$.
  Il suffit de prouver que le diagramme
\begin{equation}\label{diag1-taudual}
\xymatrix  @R=0,3cm  {
  {  f _T ^! \DD _T (\E) }
  \ar[r] _-\sim
  \ar[d] _-\sim ^{\tau _{f',f}}
  &
  {g ^! _{T \cap Y} j _T ^! \DD (\E)}
  \ar@{=}[r]
  \ar[d] _-\sim ^{\tau _{g',g} }
  &
  {g ^! _{T \cap Y} \DD _{T \cap Y} j ^! _T (\E)}
  \ar[r] ^-{\theta _g} _-\sim
  &
  {\DD _{T'} g ^! _{T \cap Y}  j _T ^! (\E)}
  \ar[r] _-\sim
  &
  {\DD _{T'} f ^! _T  (\E)}
  \\
  {  f _T ^{\prime !} \DD _T (\E)  }
  \ar[r] _-\sim
  &
  {g ^{\prime !} _{T \cap Y} j _T ^! \DD (\E)}
  \ar@{=}[r]
  &
  {g ^{\prime !} _{T \cap Y} \DD _{T \cap Y} j ^! _T (\E)}
  \ar[r] ^-{\theta _g} _-\sim
  \ar[u] _-\sim ^{\tau _{g',g} }
  &
  {\DD _{T'} g ^{\prime !} _{T \cap Y}  j _T ^! (\E)}
  \ar[r] _-\sim
  \ar[u] _-\sim ^{\tau _{g',g} }
  &
  {\DD _{T'} f ^{\prime !} _T  (\E)}
  \ar[u] _-\sim ^{\tau _{f',f} }
}
\end{equation}
  est commutatif.
  Or, grâce à \ref{tausigma-thg}, le deuxième carré de droite est commutatif.
  La commutativité du deuxième carré de gauche est évidente tandis que celle des deux autres
  se vérifie par transitivité des isomorphismes de la forme $\tau$.
\end{proof}

\begin{prop}\label{trans-theta-f2}
  Soient $f$ : $\X ' \hookrightarrow \X$ et $f '$ : $\X ''\hookrightarrow \X'$
  deux immersions de $\V$-schémas formels lisses,
  $T$ un diviseur de $X$ tel que $T' :=f ^{-1} (T)$ (resp. $T'':=f ^{\prime -1} (T')$
  soit un diviseur de $X'$ (resp. $X''$).
  Pour tout
  $\E \in D ^{\mathrm{b}}_{\mathrm{coh}} ( \overset{^\mathrm{g}}{} \smash{\D} ^{\dag}   _{\X}(\hdag T) _\Q)$
  à support dans l'adhérence schématique de $X''$ dans $X$,
  le diagramme
\begin{equation}
  \label{trans-theta-f2diag}
  \xymatrix  @R=0,3cm   @C=2cm  {
  {(f \circ f ')^{ !} _T \DD _T (\E) }
  \ar[d] _-\sim \ar[rr] _-\sim ^{\theta _{f\circ f',T}}
  &&
  {\DD _T (f \circ f') ^{ !} _T (\E)}
  \\
  { f ^{\prime !} _{T'}  f ^! _T \DD _T(\E)}
  \ar[r] _-\sim ^{\theta _{f,T}}
  &
  { f ^{\prime !}_{T'}  \DD _{T'} f ^! _T (\E)}
  \ar[r] _-\sim ^{\theta _{f',T'}}
  &
  {\DD _{T''} f ^{\prime !} _{T'} f ^! _T(\E)}
  \ar[u] _-\sim
}
\end{equation}
  est commutatif.
\end{prop}
\begin{proof}
Par \cite[4.3.12]{Be1}, on se ramène au cas où $T$ est vide.
  On désigne par $\Y$ (resp. $\Y'$) un ouvert de $\X$ (resp. $\X'$) tel que $f$ (resp. $f '$) se factorise
  en une immersion fermée $g $ : $ \X ' \hookrightarrow \Y$ (resp. $g '$ : $ \X '' \hookrightarrow \Y'$).
  Soit $\Y ''$ un ouvert de $\Y$ tel que $g ^{-1} (\Y'')=\Y'$.
On note $j $ : $\Y \subset \X$ et $j'$ : $\Y ' \subset \X'$ les inclusions canoniques. Par abus
de notations, on désignera encore par $j'$ : $\Y '' \subset \Y$ et $g$ : $\Y'\hookrightarrow \Y''$ les morphismes
canoniques.
Considérons le diagramme :
  \begin{equation}
  \label{trans-theta-f2-diag1}
    \xymatrix  @R=0,3cm     {
    { f ^{\prime !}  f ^! \DD(\E)}
    \ar[rr] _-\sim ^{\theta _{f}}
    &&
    { f ^{\prime !} \DD f ^! (\E)}
    \ar[rr] _-\sim ^{\theta _{f'}}
    &&
    {\DD f ^{\prime !} f ^! (\E)}
    \ar[d] _-\sim
    \\
    {g ^{\prime !}   j ^{\prime !} g ^! j ^!\DD (\E)}
    \ar[r] _-\sim
    \ar[u] _-\sim
    &
    {g ^{\prime !}   j ^{\prime !} g ^! \DD j ^! (\E)}
    \ar[r] _-\sim ^{\theta _g}
    &
    {g ^{\prime !}   j ^{\prime !} \DD g ^! j ^!(\E)}
    \ar[r] _-\sim
    \ar[u] _-\sim
    &
    {g ^{\prime !} \DD  j ^{\prime !} g ^! j ^!(\E)}
    \ar[r] _-\sim  ^{\theta _{g'}}
    \ar[d] _-\sim
    &
    {\DD g ^{\prime !} j ^{\prime !} g ^! j ^! (\E)}
    \ar[d] _-\sim
    \\
    {g ^{\prime !}   g ^!  j ^{\prime !} j ^! \DD (\E)}
    \ar[r] _-\sim
    \ar[u] _-\sim\ar[d] _-\sim
    &
    {g ^{\prime !}   g ^!  j ^{\prime !}\DD j ^! (\E)}
    \ar[r] _-\sim \ar[u] _-\sim \ar[d] _-\sim
    &
    {g ^{\prime !}   g ^! \DD  j ^{\prime !} j ^! (\E)}
    \ar[r] _-\sim ^{\theta _g}
    \ar[d] _-\sim
    &
    {g ^{\prime !} \DD  g ^! j ^{\prime !} j ^! (\E)}
    \ar[r] _-\sim ^{\theta _{g'}}
    &
    {\DD g ^{\prime !}  g ^! j ^{\prime !} j ^! (\E)}
    \\
    {  g \circ g ^{\prime !}   j ^{\prime !} j ^! \DD (\E)}
    \ar[r] _-\sim
    \ar[d] _-\sim
    &
    {  g \circ g ^{\prime !}   j ^{\prime !} \DD j ^! (\E)}
    \ar[r] _-\sim
    &
    {  g \circ g ^{\prime !} \DD  j ^{\prime !} j ^! (\E)}
    \ar[rr] _-\sim ^{\theta _{g \circ g'}}
    &
    &
    {\DD   g \circ g ^{\prime !}  j ^{\prime !} j ^! (\E)}
    \ar[u] _-\sim
    \\
    { g \circ g ^{\prime !}  j \circ j ^{\prime !} \DD (\E)}
    \ar[rr] _-\sim
    \ar[d] _-\sim
    &&
    { g \circ g ^{\prime !}  \DD  j \circ j ^{\prime !} (\E)}
    \ar[rr] _-\sim ^{\theta _{g \circ g'}}
    \ar[u] _-\sim
    &&
    {\DD   g \circ g ^{\prime !}   j \circ j ^{\prime !} (\E)}
    \ar[u] _-\sim
    \\
{f \circ f^{\prime !} \DD (\E) }
  \ar[rrrr] _-\sim ^{\theta _{f\circ f'}}
  &&&&
  {\DD f \circ f ^{\prime !} (\E),}
  \ar[u] _-\sim
    }
  \end{equation}
  dont le contour correspond à celui de \ref{trans-theta-f2diag} (on vérifie par transitivité que le composé de droite
  de \ref{trans-theta-f2-diag1} est le morphisme de droite de \ref{trans-theta-f2diag}, de même à gauche).
  Il dérive de \ref{def-thetaf} que les rectangles du haut et du bas
  sont commutatifs.
  La commutativité du rectangle de la deuxième ligne résulte du fait que les isomorphismes de la forme
  $\theta _g$ commute à la restriction à un ouvert. Celle du rectangle de la troisième ligne découle de
  \ref{trans-thg} tandis que celle du rectangle de gauche de l'avant dernière ligne est immédiate.
  Celle des carrés se vérifie par fonctorialité. Le diagramme \ref{trans-theta-f2-diag1} est donc
  commutatif.
\end{proof}

\subsection{Commutation au foncteur dual}
Nous reprenons les notations de \ref{notat-construc}.

\begin{prop}\label{propspetdualsansfrob}
Désignons par $E$ un isocristal sur $Y$ surconvergent le long de
$T _X$ et par $E ^\vee$ son dual. On a l'isomorphisme
canonique fonctoriel en $E$ : $ \sp _{X \hookrightarrow \PP,\,T  \,+} (E ^\vee ) \riso
\DD _{T} \circ \sp _{X \hookrightarrow \PP,\,T  \,+} (E )$.
\end{prop}
\begin{proof}
En remplaçant les foncteurs de la forme $j^\dag$ ou $(\hdag T)$ par
les foncteurs duaux respectifs,
il s'agit de reprendre la preuve de \ref{spjdag}.
Par exemple, avec les notations de \ref{prop1},
grâce à \ref{taudual} et \ref{trans-theta-f2}, on vérifie pour tout
$\E \in \mathrm{Coh} (X,\, \PP,\, T)$
la commutativité des diagrammes ci-dessous :
\begin{gather}
\label{propspetdual-lem1-diag01}
  \xymatrix  @R=0,3cm {
{u _{\alpha \beta} ^! \DD _{T \cap P_{\alpha \beta}}  ( \E | _{\PP _{\alpha\beta}})}
\ar[rr]
\ar[d] ^-{\DD _{T \cap P_{\alpha \beta}} \tau}
&&
{\DD _{T \cap X_{\alpha \beta}} u _{\alpha \beta} ^! ( \E | _{\PP _{\alpha\beta}})}
\ar[d] ^-{ \tau}
\\
{p _{1 }  ^{\alpha \beta !}   u _\alpha ^! \DD _{T \cap P _{\alpha }} ( \E | _{\PP _\alpha})}
\ar[r]
&
{p _{1 }  ^{\alpha \beta !}  \DD _{T \cap X_{\alpha }} u _\alpha ^! ( \E | _{\PP _\alpha})}
\ar[r]
&
{\DD _{T \cap X_{\alpha \beta}} p _{1 }  ^{\alpha \beta !} u _\alpha ^! ( \E | _{\PP _\alpha}),}
}
\\
\label{propspetdual-lem1-diag02}
  \xymatrix  @R=0,3cm {
{u _{\alpha \beta} ^! \DD _{T \cap P_{\alpha \beta}}  ( \E | _{\PP _{\alpha\beta}})}
\ar[rr]
\ar[d] ^-{\DD _{T \cap P_{\alpha \beta}} \tau}
&&
{\DD _{T \cap X_{\alpha \beta}} u _{\alpha \beta} ^! ( \E | _{\PP _{\alpha\beta}})}
\ar[d] ^-{ \tau}
\\
{p _{2 }  ^{\alpha \beta !}   u _\beta ^! \DD _{T \cap P _{\alpha }} ( \E | _{\PP _\beta})}
\ar[r]
&
{p _{2 }  ^{\alpha \beta !}  \DD _{T \cap X_{\beta }} u _\beta ^! ( \E | _{\PP _\beta})}
\ar[r]
&
{\DD _{T \cap X_{\alpha \beta}} p _{2 }  ^{\alpha \beta !} u _\beta ^! ( \E | _{\PP _\beta}),}
}
\end{gather}
qui correspondent aux analogues des carrés horizontaux de \ref{lem1spjdagdiag2}.
De plus, pour valider l'analogue du lemme \ref{lem3spjdag}, on utilise \ref{corospEdualf!}.

\end{proof}

\begin{nota}\label{notaotimesdag}
Soient $\E ^{(\bullet)}  , \, \FF ^{(\bullet)}
\in \smash[b]{\underset{^{\longrightarrow }}{LD }}  ^\mathrm{b} _{\Q, \mathrm{qc}}
( \smash{\widehat{\D}} _{\PP} ^{(\bullet)} (T ))$ (voir \cite[1.1]{caro_courbe-nouveau}.
On note
$$\E ^{(\bullet)}
 \smash{\overset{\L}{\otimes}}^{\dag} _{\O _{\PP } ( \hdag T ) _{\Q}}
  \FF ^{(\bullet)}
  :=
  (\E ^{(m)}  \smash{\widehat{\otimes}} ^\L _{\smash{\widehat{\B}} _{\PP} ^{(m)} (T )} \FF ^{(m)} ) _{m \in  \N }.$$
Pour une définition
des images directes et inverses extraordinaires à singularités surconvergentes le long d'un diviseur
de ces complexes quasi-cohérents,
on pourra consulter \cite[1.1.7]{caro_courbe-nouveau}.

  De manière analogue à \cite[4.2.2]{Beintro2},
on dispose du foncteur $\underset{\longrightarrow}{\lim}$ :
$\smash{\underset{^{\longrightarrow}}{LD}} ^{\mathrm{b}} _{\Q ,\mathrm{qc}}
(\smash{\widehat{\D}} _{\PP} ^{(\bullet)}(T))
\rightarrow
D  ( \smash{\D} ^\dag _{\PP} (\hdag T) _{\Q} )$. Celui-ci induit une équivalence de
catégorie entre $D ^\mathrm{b} _\mathrm{coh} ( \smash{\D} ^\dag _{\PP} (\hdag T) _{\Q} )$
et une sous-catégorie pleine de $\smash{\underset{^{\longrightarrow}}{LD}} ^{\mathrm{b}} _{\Q ,\mathrm{qc}}
(\smash{\widehat{\D}} _{\PP} ^{(\bullet)}(T))$, noté
$\smash{\underset{^{\longrightarrow}}{LD}} ^{\mathrm{b}} _{\Q ,\mathrm{coh}}
(\smash{\widehat{\D}} _{\PP} ^{(\bullet)}(T))$ (voir \cite[4.2.4]{Beintro2}).
Lorsque
$\E ^{(\bullet)}  , \, \FF ^{(\bullet)}
\in \smash[b]{\underset{^{\longrightarrow }}{LD }}  ^\mathrm{b} _{\Q, \mathrm{coh}}
( \smash{\widehat{\D}} _{\PP} ^{(\bullet)} (T ))$, en notant
$\E  := \underset{\longrightarrow}{\lim} \E ^{(\bullet)} $
et
$\FF  := \underset{\longrightarrow}{\lim} \FF ^{(\bullet)} $,
on pose
$$\E
 \smash{\overset{\L}{\otimes}}^{\dag} _{\O _{\PP } ( \hdag T ) _{\Q}}
  \FF
  :=
  \underset{\longrightarrow}{\lim}  \E ^{(\bullet)}
 \smash{\overset{\L}{\otimes}}^{\dag} _{\O _{\PP } ( \hdag T ) _{\Q}}
  \FF ^{(\bullet)} .$$
\end{nota}

\begin{vide}
  L'isomorphisme d'autodualité de \cite[2.1.1]{Becohdiff} s'écrit aussi
$\DD _{\PP,T} ( \R \underline{\Gamma} ^\dag _{X } \O _{\PP } ( \hdag T ) _{\Q} [d])
\riso
 \R \underline{\Gamma} ^\dag _{X } \O _{\PP } ( \hdag T ) _{\Q}  [d] $
 où $d :=-d _{X/P}$ la codimension de $X$ dans $P$.
 On en déduit :
$ \DD _{\PP,T} ( \R \underline{\Gamma} ^\dag _{X } \O _{\PP } ( \hdag T ) _{\Q} ) [-2d ]
\riso
\R \underline{\Gamma} ^\dag _{X } \O _{\PP } ( \hdag T ) _{\Q}$.
Cet isomorphisme ne semble pas compatible à Frobenius. En effet,
lorsque $X =P$, on retrouve l'isomorphisme canonique :
$ \DD _{\PP,T} ( \O _{\PP } ( \hdag T ) _{\Q} )
\riso
\O _{\PP } ( \hdag T ) _{\Q}$. Je n'ai pas de contre-exemple mais
la compatibilité à Frobenius de ce dernier isomorphisme me paraît inexacte.

  Pour tout $\E \in F\text{-}D ^{\mathrm{b}} _{\mathrm{coh}} (\smash{\D} ^\dag _{\PP } (\hdag T ) _\Q)$
  à support dans $X$, avec les notations \ref{notaotimesdag},
  on pose
  \begin{equation}\label{defiDD*}
    \DD ^* _{X,\PP,T} (\E) :=
\DD _{\PP,T} \left  ( (\DD _{\PP,T} ( \R \underline{\Gamma} ^\dag _{X } \O _{\PP } ( \hdag T ) _{\Q} ) [-2d ])
 \smash{\overset{\L}{\otimes}}^{\dag} _{\O _{\PP } ( \hdag T ) _{\Q}}
 \E \right ).
  \end{equation}

Or, comme $\E$ est aussi à support dans $X$, on obtient le deuxième isomorphisme :
$$\R \underline{\Gamma} ^\dag _{X } \O _{\PP } ( \hdag T ) _{\Q}
 \smash{\overset{\L}{\otimes}}^{\dag} _{\O _{\PP } ( \hdag T ) _{\Q}}
 \E
 \riso
 \R \underline{\Gamma} ^\dag _{X }  (\E) \riso \E.$$

Il en découle l'isomorphisme $\DD ^* _{X,\PP,T} (\E) \riso \DD _{\PP, T} (\E)$.
Ce dernier n'est pas à priori compatible à Frobenius.

\end{vide}

\begin{lemm}\label{proj+!}
  Soient $u$ : $\X \hookrightarrow \PP$ une immersion fermée de $\V$-schémas formels lisses,
  $T$ un diviseur de $P$ tel que $T _X:= u ^{-1} _0 (T)$ soit un diviseur de $X$,
$\E ^{(\bullet)} \in F\text{-}\smash[b]{\underset{^{\longrightarrow }}{LD }}  ^\mathrm{b} _{\Q, \mathrm{qc}}
( \smash{\widehat{\D}} _{\PP} ^{(\bullet)} (T ))$,
$\FF ^{(\bullet)} \in F\text{-}\smash[b]{\underset{^{\longrightarrow }}{LD }}  ^\mathrm{b} _{\Q, \mathrm{qc}}
( \smash{\widehat{\D}} _{\X} ^{(\bullet)} (T _X))$.
Avec les notations de \ref{notaotimesdag}, on dispose
  d'un isomorphisme compatible à Frobenius :
  $$u _{T,+} ( u ^! _T (\E ^{(\bullet)}) [-d _{X/P}]
  \smash{\overset{\L}{\otimes}}^{\dag} _{\O _{\X } ( \hdag T _X ) _{\Q}}
  \FF ^{(\bullet)} )
  \riso
  \E ^{(\bullet)}
 \smash{\overset{\L}{\otimes}}^{\dag} _{\O _{\PP } ( \hdag T ) _{\Q}}
u _{T,+} (\FF ^{(\bullet)}).$$
\end{lemm}
\begin{proof}
  Par les arguments habituels de complétion puis de passage à la limite sur le niveau,
  on se ramène au cas des schémas. Par passage de gauche à droite, cela correspond
  à \cite[1.4.1]{caro-frobdualrel} et \cite[1.4.10.2]{caro-frobdualrel}
  (on s'était contenté dans \cite[1.4]{caro-frobdualrel}
  de prouver le cas où $\E = \smash{\D} _X ^{(m)}$ et $\FF = \omega _X$, mais la preuve est identique
  pour le cas général).
\end{proof}

\begin{prop}\label{propspetdual}
Désignons par $E$ un $F$-isocristal sur $Y$ surconvergent le long de
$T _X$ et par $E ^\vee$ son dual.
On dispose de l'isomorphisme
canonique fonctoriel en $E$ et compatible à Frobenius :
$$ \sp _{X \hookrightarrow \PP,\,T  \,+} (E ^\vee ) \riso
\DD ^* _{X,\PP,T} \circ \sp _{X \hookrightarrow \PP,\,T  \,+} (E ).$$
\end{prop}
\begin{proof}
  Comme $X$ et $P$ sont lisses, par additivité, il ne coûte rien de supposer que $X$ et $P$ soient irréductibles.
Par application de $\DD _{\PP,T}$ (on dispose de l'isomorphisme de bidualité
qui est compatible à Frobenius \cite[II.3.5]{virrion})
et en échangeant $E$ et $E ^\vee$, il est équivalent d'obtenir un isomorphisme compatible à Frobenius de la forme :
\begin{equation}
  \label{propdiag1spetdual}
\DD _{\PP,T} ( \R \underline{\Gamma} ^\dag _{X } \O _{\PP } ( \hdag T ) _{\Q} ) [-2d ]
 \smash{\overset{\L}{\otimes}}^{\dag} _{\O _{\PP } ( \hdag T ) _{\Q}}
 \sp _{X \hookrightarrow \PP,\,T  \,+} (E ^\vee )
 \riso
\DD _{\PP,T} \sp _{X \hookrightarrow \PP,\,T  \,+} (E ).
\end{equation}
Supposons dans un premier temps que $X \hookrightarrow P$ se relève
en une immersion fermée $u$ : $\X \hookrightarrow \PP$ de $\V$-schémas formels lisses.
Via \cite[1.2.14.1]{caro_courbe-nouveau},
le terme de gauche de \ref{propdiag1spetdual} est isomorphe à
\begin{equation}
  \label{propdiag2spetdual}
\DD _{\PP,T} ( u _{T+} u _{T} ^!  \O _{\PP } ( \hdag T ) _{\Q} ) [-2d ]
 \smash{\overset{\L}{\otimes}}^{\dag} _{\O _{\PP } ( \hdag T ) _{\Q}}
u _{T+} \sp _{*} (E ^\vee ).
\end{equation}
Or, grâce au théorème de dualité relative compatible à Frobenius
(voir \cite{caro-frobdualrel}),
\begin{equation}
  \label{propdiag3spetdual}
\DD _{\PP,T} ( u _{T+} u _{T} ^!  \O _{\PP } ( \hdag T ) _{\Q} ) [-2d ]
\riso
u _{T+} \DD _{\X,T _X}  u _{T} ^! ( \O _{\PP } ( \hdag T ) _{\Q} ) [-2d ]
\riso
u _{T+} \DD _{\X,T _X} ( \O _{\X } ( \hdag T _X) _{\Q} ) [-d ]
\end{equation}
Avec \ref{proj+!}, comme $u _T ^! u _{T,+}\riso id$ (\cite[5.3.3]{Beintro2}),
le terme de gauche de \ref{propdiag1spetdual} devient isomorphe à
\begin{equation}
  \label{propdiag4spetdual}
   u _{T,+} ( \DD _{\X,T _X} ( \O _{\X } ( \hdag T _X) _{\Q} )
 \otimes _{\O _{\X } ( \hdag T _X ) _{\Q}}
\sp _{*} (E ^\vee )).
\end{equation}
D'un autre côté, le terme de droite de \ref{propdiag1spetdual} est isomorphe à
\begin{equation}
  \label{propdiag5spetdual}
   u _{T,+} \DD _{\X,T _X} \sp _{*} (E ).
\end{equation}

Dans le cas général où $X \hookrightarrow \PP$ ne se relève pas,
grâce à \ref{coro-frob-sp+}, il en résulte que les deux termes de
\ref{propdiag1spetdual}
sont dans l'image essentielle de $\sp _{X \hookrightarrow \PP,\,T  \,+}$.

Soit $\PP'\subset \PP \setminus T$ un ouvert affine de $\PP$ tel que $X \cap P'$ soit dense dans $X$.
Grâce à \cite[4.1.1]{tsumono} et à \cite{kedlaya_full_faithfull}, (2-ième ligne page 2),
\ref{EEhatsp+} et \ref{spjdag}, les deux termes de \ref{propdiag1spetdual}
sont isomorphes si et seulement s'ils le sont au dessus de $\PP'$.
On se ramène ainsi au cas où $T$ est vide et $X \hookrightarrow P$
se relève en une immersion fermée $u$ : $\X \hookrightarrow \PP$ de $\V$-schémas formels lisses.
Or, d'après \cite{caro_comparaison}, on dispose d'un isomorphisme compatible à Frobenius :
\begin{equation}
  \label{propdiag6spetdual}
 \DD _{\X} ( \O _{\X ,\Q} )
 \otimes _{\O _{\PP ,\Q}}
\sp _{*} (E ^\vee )
\riso
\DD _{\X} \sp _{*} (E ).
\end{equation}
Via la première partie de la preuve (\ref{propdiag4spetdual} et \ref{propdiag5spetdual}),
en appliquant $u _+$ à \ref{propdiag6spetdual}, on obtient à isomorphisme près
\ref{propdiag1spetdual} (avec $T$ vide).

\end{proof}

\begin{rema}
Si la conjecture \cite[1.2.1.(TC)]{tsumono} était validée,
  les preuves de \ref{spjdag} et de \ref{propspetdualsansfrob}
seraient moins techniques et délicates.
En effet,
en remplaçant \cite{kedlaya_full_faithfull} par \cite[1.2.1.(TC)]{tsumono},
de manière analogue à ce qui a été fait dans la preuve de \ref{propspetdual},
on se ramènerait au cas relevable.
\end{rema}

Le théorème qui suit fournit une application importante et prometteuse
de nos constructions.
Pour sa preuve, on pourra consulter \cite[2.2]{caro_unite}.
\begin{theo}
Soit $E$ un $F$-isocristal unité sur $Y$ surconvergent le long de
$T _X$. Le faisceau $\sp _{X \hookrightarrow \PP,\,T  \,+} (E )$ est $\smash{\D} ^\dag _{\PP, \Q}$-cohérent.
\end{theo}

\bibliographystyle{smfalpha}
\bibliography{bib1}

\end{document}